\documentclass[11pt]{article}

\usepackage{amsmath,amsthm,amsfonts}
\usepackage{graphicx}
\usepackage[usenames]{color}
\usepackage[T1,T2A]{fontenc}

\usepackage[english]{babel}
\usepackage{amssymb}

\usepackage{ifthen}

\newtheorem{theorem}{Theorem}[section]
\newtheorem{claim}{Proposition}[section]
\newtheorem{lemma}{Lemma}[section]
\newtheorem{corollary}{Corollary}[section]

\theoremstyle{definition}
\newtheorem{remark}{Remark}[section]
\newtheorem{exercise}{\quad}

\newtheorem{example}{Example}[section]

\newcommand{\supp}{{\rm supp}}
\newcommand{\rank}{{\rm rank}}
\newcommand{\tr}{{\rm tr}}
\newcommand{\wt}{{\rm wt}}
\newcommand{\cor}{{\rm cor}}
\newcommand{\per}{{\rm per}}

\begin{document}

\begin{center}

\Large

{\bf \huge Perfect Combinatorial Structures in Coding Theory and Cryptography}

\bigskip

Vladimir N. Potapov
\end{center}

\pagebreak

\tableofcontents

\pagebreak

\section*{Preface}

\addcontentsline{toc}{section}{Preface}

This book develops algebraic and combinatorial methods for studying
discrete structures. It brings together graph theory, Boolean functions,
Fourier analysis on finite groups, coding theory, perfect colorings and
perfect codes, association schemes, Latin squares, and related topics.
A central theme is the interaction between different representations of
the same object: combinatorial, algebraic, spectral, and
coding-theoretic. The main mathematical object studied in this book is
a perfect coloring of a graph, or, equivalently, an equitable partition
of a graph.

The book is intended for advanced undergraduate and graduate students
in mathematics and computer science, as well as for researchers in
discrete mathematics, combinatorics, coding theory, and related fields.

\smallskip

{\bf Prerequisites}

The reader is expected to have a basic knowledge of linear algebra.
Some familiarity with graph theory, coding theory, and algebra will be
useful for understanding the more advanced topics.

\smallskip

{\bf Methods and References}

The book is based primarily on methods from algebraic graph theory and
Fourier analysis on finite abelian groups. For an introduction to
algebraic graph theory, I recommend Godsil's monograph \cite{Godsil}.
For Fourier analysis and other methods in the theory of Boolean
functions, I recommend Carlet's monograph \cite{Carlet20}. I also
recommend the classical books on coding theory by van Lint \cite{Lint},
MacWilliams and Sloane \cite{MacW}, and Delsarte \cite{Delsarte73}.
A more detailed and rigorous treatment of a substantial part of the
material presented in this handbook can be found in \cite{CRCKP}.

\smallskip

{\bf Acknowledgments}

The author is grateful to Sergey Avgustinovich, who inspired me to
study perfect combinatorial structures; Denis Krotov, whose important
results form a substantial basis for this book; Anna Taranenko, who
carefully read previous versions of the text and provided many useful
comments; and Ev Sotnikova, who prepared the notes of my lectures.

I am also grateful to all participants of the weekly Friday meetings
at the Sobolev Institute of Mathematics: Evgeny Bespalov, Sergey
Malyugin, Ivan Mogilnykh, Alexey Perezhogin, Yuriy Tarannikov,  Alexandr Valyuzhenich, Anastasia
Vasil'eva, and Konstantin Vorob'ev. Their results
and comments have been valuable in preparing this book.

The English translation of the book was prepared using Google Translate, DeepSeek and ChatGPT.

\pagebreak

{\bf Notation}

Throughout the handbook, we use the following notation.

$\mathbb{R}$ and $\mathbb{C}$ denote the fields of real and complex
numbers, respectively.

$GF(q)$ denotes the Galois field of order $q$, where $q$ is a prime
power.

$\mathbb{Z}_q$ denotes the ring of integers modulo $q$ and, when
appropriate, its additive group.

$Q_q^n=\{0,1,\dots,q-1\}^n$ denotes the $q$-ary hypercube. This set may
be equipped with the Hamming metric, in which case $Q_q^n$ is the
Hamming graph $H(n,q)$, or with the structure of a vector space over
$GF(q)$, or with the structure of the abelian group $\mathbb{Z}_q^n$.

$G$ usually denotes a graph with adjacency matrix $M$.

$S$ denotes the quotient matrix of a perfect coloring or equitable
partition.

$\widehat{f}$ denotes the Fourier transform of $f$, normalized to be an
isometry.

\smallskip
 
{\bf Open problems}

Below  is a list of well-known open problems related to the subject of the book.
\begin{itemize}
\item Existence of $1$-perfect codes for non prime power alphabets  \hfill p.\,\pageref{PCprob}
\item Asymptotic number of $1$-perfect binary codes  \hfill p.\,\pageref{Vasprob}
\item Number of iterations of the Weisfeiler--Leman algorithm \hfill p.\,\pageref{WLprob}
\item Fourier Entropy -- Influence conjecture   \hfill p.\,\pageref{FEIprob}
\item Description of the vertices of the generalized Birkhoff polytope  \hfill p.\,\pageref{Berprob}
\item Existence of $3$ mutually orthogonal Latin squares of order $10$ \hfill p.\,\pageref{LSprob}
\item Existence of projective geometries  for non prime power orders   \hfill p.\,\pageref{PGprob}
\item Ryser -- Wanless conjecture on the existence of transversals   \hfill p.\,\pageref{Tprob}
\item Unrestricted main MDS conjecture \hfill p.\,\pageref{MDSprob}
\item Existence of Hadamard matrices for all orders divisible by $4$   \hfill p.\,\pageref{Hprob}
\end{itemize}

\pagebreak

\section{Graphs and Their Colorings}

\subsection{Graph Theory: Definitions and Notation}
Consider an undirected {\sl graph} $G$ without {multiple edges} or
loops  (a loop is an edge with the same vertex as both its beginning and end). Such
graphs are called {\sl simple}. Throughout this text, whenever we use the
term
"graph", we mean a simple graph, unless explicitly stated otherwise.
$V(G)$ denotes the set of {\sl vertices} of the graph $G$, and $E(G)\subset {V(G) \choose 2}$ denotes the set of {\sl edges} of the graph
$G=(V(G),E(G))$. Here and throughout, ${A \choose k}$
denotes the set of unordered $k$-element subsets of $A$.

A graph is called
{\sl connected} if any two of its vertices are connected by a {\sl path} of
edges. A connected graph corresponds to a discrete metric
space $X(V(G),d(G))$. The points of $X(V(G),d(G))$ are the vertices of the graph, and the {\sl distance} $d(G)$ between
two points
is the number of edges in the shortest path connecting them.
 The maximum distance between vertices of a graph is called the {\sl diameter} of the graph.

\begin{example}
Consider a set $Q_q=\{0,1,\dots,q-1\}$. Two $n$-tuples $x,y\in Q^n_q$ are connected by  an edge if $x$ and $y$  differ in only one position. This graph is called the {\sl Hamming graph} and denoted by $H(n,q)$. The distance  $d(x,y)$ in the corresponding metric  space $(Q^n_q, d)$  is called  the {\sl Hamming distance}.  $d(x,y)$ is equal to the number of positions at
which $n$-tuples $x,y\in Q^n_q$ differ. The metric
space $(Q^n_q, d)$ and the corresponding graph are also 
called an {\sl $n$-dimensional $q$-ary hypercube}; if $q=2$, then it is a {\sl Boolean $n$-cube}.
\end{example}

A {\sl spanning subgraph} of $G$  is a subgraph that contains every vertex of $G$ but only a subset of the edges.  An {\sl induced subgraph} of $G$ by $C\subset V(G)$ is another graph, formed from $C$ and all of the edges, from $G$, connecting pairs of vertices in $C$.

A {\sl clique} in a graph is an induced subgraph all of whose vertices
are pairwise connected by edges. A one-dimensional $q$-ary cube is a 
clique on $q$ vertices. In graph theory, a clique on $q$ vertices
has the standard notation $K_q$. However, for convenience, we will henceforth use the notations $Q_q$ and $Q^n_q$ in the multidimensional case, both for the metric space and for the corresponding graph, as well as for the set of its vertices.

The {\sl degree} of a graph vertex is the number of edges incident to it. In a simple graph, the degree of a vertex is equal to the number of adjacent vertices.

%\begin{definition}
A graph is called {\sl regular} if the degrees of all its vertices are equal.
%\end{definition}
If the degrees of all vertices in a graph are $r$, then the graph is called $r$-{\sl regular}. Let's  arbitrarily enumerate the vertices of the graph:
$V(G)=\{v_i\}$.

%\begin{definition}
{\sl The adjacency matrix} of a graph is defined by the equality

\[
M_{n\times n}(G)=(m_{ij})= \left\{
\begin{array}{l}
1,\ \mbox{if} (v_i,v_j) \in E(G);\\
0, \ \mbox{otherwise}.
\end{array}
\right.
\]
%\end{definition}

An undirected graph without multiple edges has a symmetric
$(0,1)$-adjacency matrix. If $G$ is a regular graph, then each
row and each column of the adjacency matrix contains the same
number of ones, equal to the degree of the graph.

%\begin{definition}
{\sl The complement} of a graph $G$ is the graph $\overline{G}$ with the same vertex set $V(\overline{G})=V(G)$ and an 
edge set $E(\overline{G})={V(G)\choose 2}\setminus E(G)$.
%\end{definition}
The adjacency matrix of   $\overline{G}$  satisfies the equality
$\overline{M}=J-I-M$.  Here and throughout, $I$ is the identity matrix and $J$ is a matrix all of whose elements are equal to $1$.   We will also use notations  $I_n$ and $J_n$ for  matrices of size $n\times n$.

\begin{example}
The simplest graph $Q_2$ consists of one edge. A cycle on four vertices,
also known as the Boolean $2$-cube, has the following adjacency matrix:

\begin{minipage}{0.5\textwidth}
%\vskip15mm
\center{\includegraphics[width=0.4\textwidth]{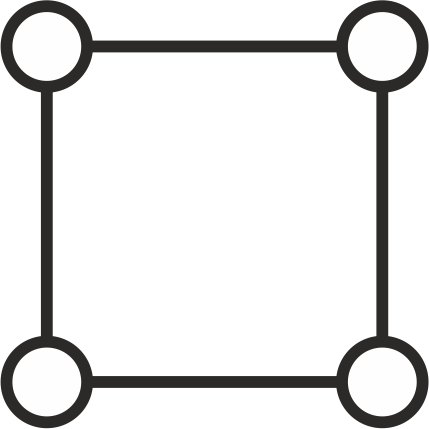}}
\end{minipage}
\begin{minipage}{0.5\textwidth}
$\begin{pmatrix}
0 & 1 & 1 & 0 \\
1 & 0 & 0 & 1 \\
1 & 0 & 0 & 1\\
0 & 1 & 1 & 0
\end{pmatrix}$
\end{minipage}
\end{example}

\begin{example}

The adjacency matrix of a Boolean cube $Q^3_2$ with eight vertices is given below. \\

\begin{minipage}{0.4\textwidth}
%\vskip15mm
\includegraphics[width=0.7\textwidth]{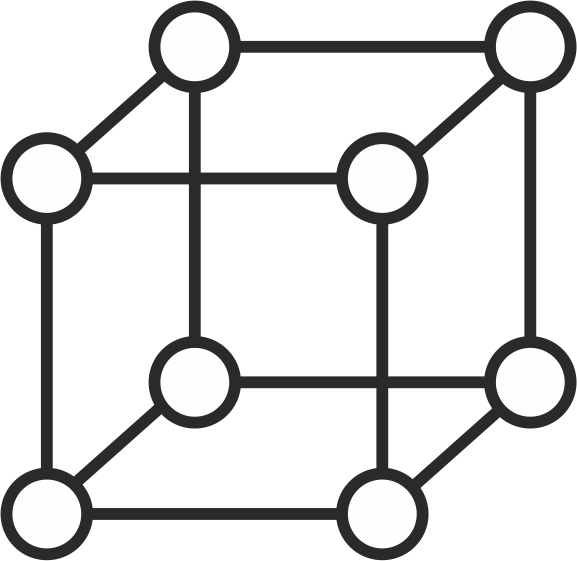}
% \center{\includegraphics[width=0.07\textwidth]{E3.png}}
\end{minipage}
\begin{minipage}{0.5\textwidth}
$\begin{pmatrix}
0 & 1 & 1 & 0 & 1 & 0 & 0 & 0 \\
1 & 0 & 0 & 1 & 0 & 1 & 0 & 0\\
1 & 0 & 0 & 1 & 0 & 0 & 1 & 0\\
0 & 1 & 1 & 0 & 0 & 0 & 0 & 1\\
1 & 0 & 0 & 0 & 0 & 1 & 1 & 0 \\
0 & 1 & 0 & 0 & 1 & 0 & 0 & 1 \\
0 & 0 & 1 & 0 & 1 & 0 & 0 & 1\\
0 & 0 & 0 & 1 & 0 & 1& 1 & 0
\end{pmatrix}$
\end{minipage}
\end{example}

The vertices of these Boolean cubes are numbered lexicographically.

{%\begin{definition}
A graph $G$ is called {\sl $k$-partite} if its vertices can
be divided into $k$ parts such that no two vertices in the same part are adjacent,
i.e., there are no edges connecting two vertices in the same part.
%\end{definition}
If we number the vertices of the graph so that each part is
numbered consecutively, then the adjacency matrix of a $k$-partite graph will
contain zero square submatrices on the diagonal. For example,
the adjacency matrix of a bipartite graph with this vertex numbering has the form
\begin{minipage}{0.25\textwidth}
$M=\begin{pmatrix}
O & B \\
B^{T} & O
\end{pmatrix}$,
\end{minipage}
where $O$ is a matrix of all zeros, and $B$ is the adjacency matrix
of the bipartite subgraph between two parts.
 For every $n$ the Boolean $n$-cube $Q^n_2$ is a bipartite graph.} In particular, after
renumbering the vertices, the adjacency matrix of a $2$-dimensional Boolean cube
takes the form $\begin{pmatrix}
0 & 0 & 1 & 1 \\
0 & 0 & 1 & 1 \\
1 & 1 & 0 & 0\\
1 & 1 & 0 & 0
\end{pmatrix}$.
A bipartite graph in which the degree of a vertex depends only on the part
of the graph is called {\sl biregular}.

%\begin{definition}
Two graphs are called {\sl isomorphic} if it is possible to establish
a one-to-one correspondence between their vertices such that
adjacency is preserved. 
%\end{definition}
Note that graphs are isomorphic if and only if their adjacency matrices map to each other by a synchronous permutation of rows and columns.

The set of vertices of a hypercube $Q_q^n$ that have a fixed value
of one coordinate, for example $x_n=0$, is called a {\sl hyperface} or axis-aligned {\sl hyperplane}.
The subgraph of the hypercube induced by the vertices of a hyperface is isomorphic
to a hypercube of one lower dimension, i.e., $Q_q^{n-1}$.
By fixing $m$ coordinates, we obtain an $(n-m)$-dimensional  {\sl face} that corresponds to the subgraph isomorphic to 
$Q^{n-m}_q$. Here and below,
we mean the induced subgraph, i.e., a subset of the graph's vertices
with edges  incident only to these vertices.

For brevity, we will call an arbitrary subset of the graph's vertices
a {\sl code}. This name stems from the fact that the main
application of the theory of perfect combinatorial structures presented below is
coding theory.

An {\sl automorphism} of a graph $G$ is a permutation of its
vertices $\pi:V(G)\rightarrow V(G)$ such that edges map to
edges, i.e., $$\{u,v\}\in E(G) \Leftrightarrow \{\pi(u),\pi(v)\}\in
E(G).$$%\end{definition}
The set of automorphisms of a graph is a group with respect to composition  and it is denoted
$Aut(G)$.
%If a graph is considered as a metric space, then
%its automorphisms are called isometries.
{Automorphisms of the code} $C\subset V(G)$
are automorphisms $\pi\in Aut(G)$ of the graph $G$ that preserve the code,
i.e., $\pi(C)=C$.

\begin{example}
Automorphisms of a Boolean $n$-cube are permutations of coordinates, additions  of arbitrary binary vectors  modulo $2$ and all their compositions.
\end{example}

By  definition, a permutation of  $V(G)$ is an automorphism  if  it  preserves the distance $1$, and therefore
all other distances in the discrete
metric space  $X(V(G),d(G))$  corresponding to $G$.
Thus, the set of automorphisms of a graph coincides with the set
of isometries of this space. So, graph automorphisms, especially those of a hypercube, are often called isometries, and isomorphic subsets
of graph vertices are called {\sl equivalent} (or {\sl isometric})  codes.

%\begin{definition}
A graph $G$ is called {\sl vertex-transitive} if, for any two
graph vertices $v,v'\in V(G)$, there exists $\pi\in Aut(G)$ such that
$v'=\pi(v)$. A code $C\subset V(G)$ is called vertex-transitive
if, for any two vertices of the code, there exists an automorphism of the code
mapping one vertex to the other.
%\end{definition}
%\begin{definition}
A vertex-transitive graph $G$ is called {\sl arc-transitive}  if for any two
edges of the graph $e,e'\in E(G)$, there exists $\pi\in Aut(G)$ such that
$e'=\pi(e)$, where $e=\{u,v\}$ and $\pi(e)=\{\pi(u),\pi(v)\}$.
%A subgraph $G'$ is called edge-
%transitive if, for any of its edges, there exists an automorphism of the code
%$V(G')$ that maps one edge to the other.
%\end{definition}
Each $n$-dimensional $q$-ary hypercube is an 
example of an edge-transitive graph. Their automorphism groups 
are discussed in the section \ref{Orbit}.

\subsection{Graph Colorings}

%\begin{definition}
Let $\{1,\ldots,k\}$ be a set of $k$ colors.  A {\sl coloring} or  $k$-{\sl coloring}  of vertices of a graph $G$ is a mapping $f: V(G) \rightarrow
\{1,\ldots,k\}$. We will use other color sets if it is convenient. 
%\end{definition}
In graph theory, one considers  not only vertices but also edges coloring. Further,   we always mean a vertex coloring, unless explicitly stated otherwise. $f^{-1}(i)$ is the set of vertices of the $i$th color. $P(f)=\{f^{-1}(i) : i=1,\dots,k\}$ is a {\sl partition} of $V(G)$.  Two colorings $f$ and $g$ of $G$ are {\sl equivalent} if there exists an automorphism $\pi\in Aut(G)$ such that $P(f)=P(g\circ \pi)$.

For each vertex $x\in V(G)$, the neighborhood of a vertex (the unit sphere centered at the vertex) is defined as follows: $$A_1(x)=\{y\in
V(G) \mid \{x,y\}\in E(G)\},$$ i.e.,  $A_1(x)$ consists of the vertices
adjacent to  $x$. A unit sphere together with its center is called  a unit ball: $B_1(x)=A_1(x)\cup \{x\}$.
Let $f(x)=i$, and denote $s_{ij}(x)=|A_1(x)\cap f^{-1}(j)|$.

%\begin{definition}
A coloring is called {\sl perfect} if the value $s_{ij}(x)=s_{ij}$
does not depend on the choice of vertex $x$ and color $i$.
%\end{definition}

Partitions corresponding to perfect colorings are called {\sl equitable}.

%\begin{definition}
The matrix 
$S_{k\times k}=(s_{ij})$ is called the {\sl quotient matrix}
of the perfect coloring.
%\end{definition}

\begin{example} A perfect coloring of the Petersen graph.  We  order the
 colors in the following way: white, green, and red.

\begin{minipage}{0.5\textwidth}
%\vskip20mm
\includegraphics[width=0.5\textwidth]{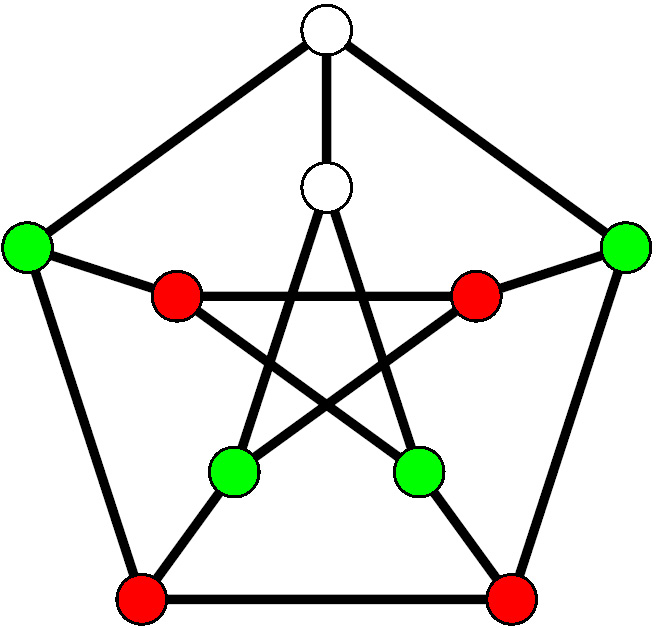}
\end{minipage}
\begin{minipage}{0.5\textwidth}
$S=\begin{pmatrix}
1 & 2 & 0 \\
1 & 0 & 2 \\
0 & 2 & 1
\end{pmatrix}$
\end{minipage}

\end{example}

\begin{example}\label{ex:bcube}
All nonequivalent perfect $2$-colorings   of the Boolean $3$-cube and their quotient matrices. The first color is black, the second is white.

\begin{minipage}{0.24\textwidth}
%\vskip12mm
\center{\includegraphics[width=0.5\textwidth]{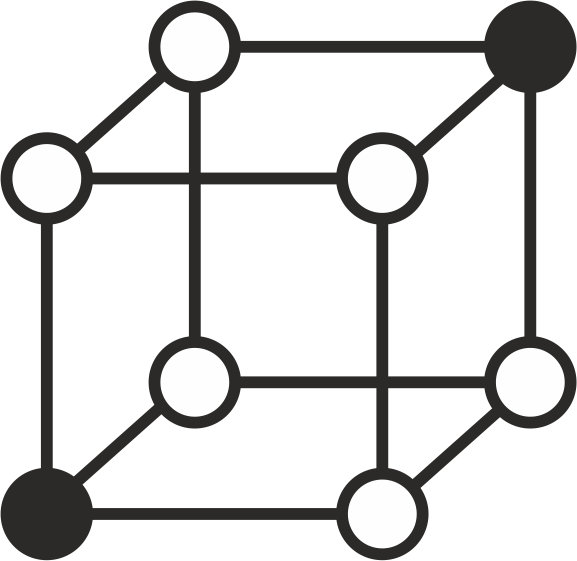}}
\end{minipage}
\begin{minipage}{0.24\textwidth}
$S=\begin{pmatrix}
0 & 3 \\
1 & 2
\end{pmatrix}$
\end{minipage}
\begin{minipage}{0.24\textwidth}
%\vskip12mm
\center{\includegraphics[width=0.5\textwidth]{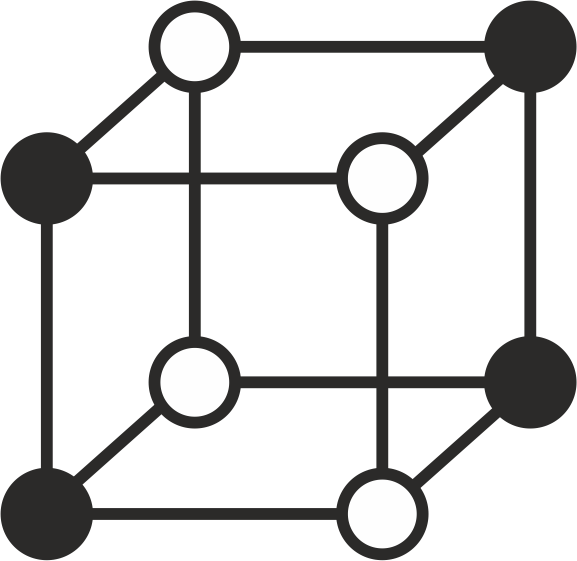}}
\end{minipage}
\begin{minipage}{0.24\textwidth}
$S=\begin{pmatrix}
1 & 2 \\
2 & 1
\end{pmatrix}$
\end{minipage}

\begin{minipage}{0.24\textwidth}
%\vskip12mm
\center{\includegraphics[width=0.5\textwidth]{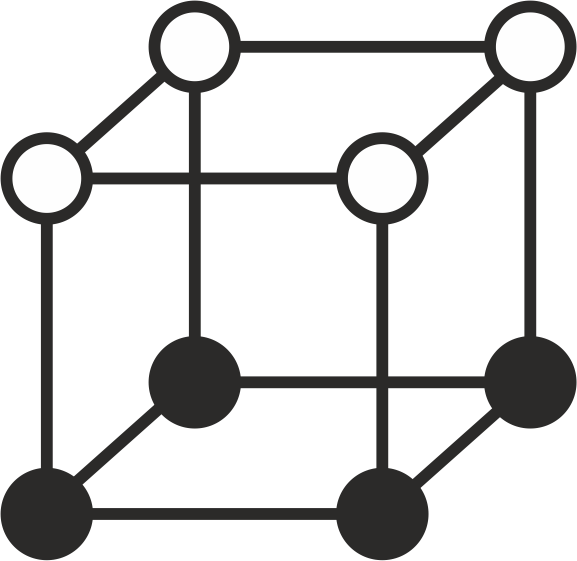}}
\end{minipage}
\begin{minipage}{0.24\textwidth}
$S=\begin{pmatrix}
2 & 1 \\
1 & 2
\end{pmatrix}$
\end{minipage}
\begin{minipage}{0.24\textwidth}
%\vskip12mm
\center{\includegraphics[width=0.5\textwidth]{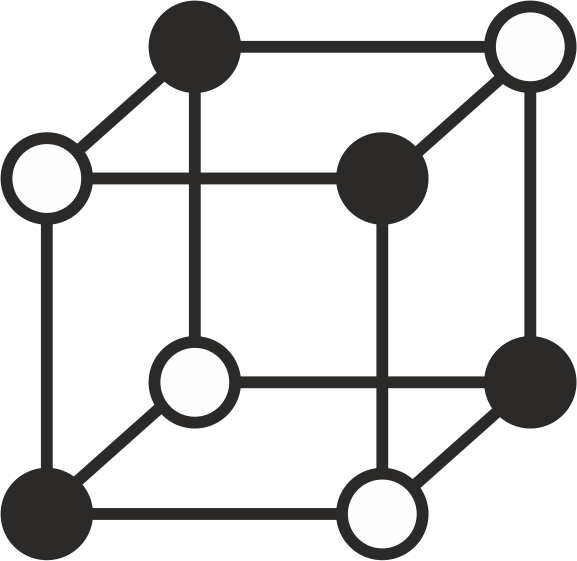}}
\end{minipage}
\begin{minipage}{0.24\textwidth}
$S=\begin{pmatrix}
0 & 3 \\
3 & 0
\end{pmatrix}$
\end{minipage}
\end{example}

\begin{claim}\label{cl:complete}
Every coloring of the complete graph is perfect.
\end{claim}
\begin{proof}
Consider two vertices $x$ and $y$ of the same color. Their neighborhoods  coincide, except for the vertices $y$  and $x$, respectively. But $x$ and $y$ have  the same color.
\end{proof}

As is often the case with useful and easily defined mathematical objects, concepts equivalent to a perfect graph coloring were independently introduced in various fields of discrete mathematics. In particular, the term "equitable partition" is used in the paper by Godsil, Koolen, and their students\,\cite{Godsil0}, \cite{Martin}. The equivalent term "partition design" was introduced in the paper by Delsarte et al.\,\cite{CCD} and has been used by specialists in the field of combinatorial design. This refers to the theory of group designs, where groups are sets of vertices of the same color, and blocks are unit spheres. 
 An equivalent term for a  perfect coloring,
"distributive coloring", was also used by  Vizing\,\cite{Vizing}. 

The term "intriguing set"\,\cite{Kelly} is used in the field of
finite geometry   to denote the set of vertices of the same color in a two-color perfect coloring.  The term "extremal graphical
design"\,\cite{Golubev} is used in the field of extremal graph problems.
The term "$(c_0,c_1)$-regular function"\  is used by specialists in Boolean functions to refer to
perfect $2$-colorings of the Boolean hypercube\,\cite{Taran02}. The term "perfect
coloring"\ was introduced by  Avgustinovich. This term
was first used  in a paper by his student Puzynina \cite{Puz}. The epithet "perfect"\ arose because the set of vertices
of the first color in a coloring with the quotient matrix $\begin{pmatrix}
0 & r \\
1 & r-1
\end{pmatrix}$
turns out to be a perfect code in the graph (see the next section).

Combinatorial configurations that can be described as
perfect colorings in some graph, multigraph, or
hypergraph are considered as {\sl perfect  combinatorial structures}.

Another, more widely studied, type of graph coloring is proper
coloring.
%\begin{definition}
A coloring is called {\sl proper} if no vertex has neighbors
of the same color as itself.
%\end{definition}
Obviously, each perfect coloring whose quotient matrix
has zeros on the main diagonal is proper. Most proper
colorings are quite chaotic. Unlike perfect colorings, proper colorings do not have to have a homogeneous local structure.

\subsection{Perfect Colorings and Information Transmission}\label{sec:it}

Perfect colorings often arise as a solution to some extremal problem, i.e., as a result of optimizing mappings with respect to some parameters. Several such problems, arising in coding theory and cryptography, will be considered below. Here, we will consider two simple examples related to data transmission over noisy communication channels.

Consider the following model of data transmission.
We encode all possible messages by the vertices of a graph $G$.
Let us connect two vertices $u, v \in V(G)$ by an edge in $G$ if the noisy communication channel can transform message $u$ into message
$v$. A communication channel is called symmetric if the probabilities of errors transforming $u$ into $v$ and $v$ into $u$ are equal. Such a communication channel
corresponds to a simple graph $G$.
%\begin{definition}
A subset of the graph's vertices is called a {\sl $1$-error-correcting code}
if balls of radius $1$ centered at the code vertices do not intersect.
%\end{definition}
If we  transmit only those messages that correspond to the vertices
of the code,  then  we can fix  any error in the received message. The greater the cardinality
 of the code, the more information can be transmitted over
the communication channel per unit time. Let all balls of radius $1$ in the graph
$G$ have the same cardinality $r+1$, i.e., the graph $G$ is $r$-regular. Then
the inequality $(r+1)|C|\leq |V(G)|$ holds,  where $|C|$ is the cardinality of a 1-error-correcting code $C$. This inequality is called the {\sl Hamming bound}. It is easy to see that a code
has maximum cardinality when the graph
is completely partitioned into balls of radius $1$ centered at the code
vertices. This means that a coloring in which the code vertices have one color and non-code vertices have another color is perfect with the quotient matrix
$\begin{pmatrix}
0 & r \\
1 & r-1
\end{pmatrix}$.
Here the set of vertices of the first color  is called a 
1-{\sl perfect}  code.   The theory of
error-correcting codes is described in \cite{Etzion}, \cite{MacW},
\cite{Lint},  and \cite{Sol}.

Next we suppose that it is sufficient  to know whether an error occurred in  the received
message. For example, we  assume that we can  request a retransmission of a corrupted message.
%\begin{definition}
A subset of graph vertices is called {\sl independent} if it does not contain adjacent vertices.
%\end{definition}
If only messages corresponding to
some independent set are transmitted over a communication channel, then all erroneous messages
are detected.

In an $r$-regular graph, an independent set cannot contain more than half of the vertices because any edge contains one vertex of the independent set at most. If a graph
 contains an independent set containing exactly half of the vertices, then the other half
of the vertices is also an independent set. So, the coloring of the two
parts of the graph in different colors is perfect with the quotient matrix
 $\begin{pmatrix}
0 & r \\
r & 0
\end{pmatrix}$. Moreover, the graph is bipartite by the definition.  In Section \ref{3.1}, we will obtain an upper bound for the cardinality
of an independent set for regular non-bipartite graphs and will show that when this bound is reached, the independent set also turns out to be a perfect coloring.

In Example \ref{ex:bcube}, the first coloring of the Boolean $3$-cube
corresponds to a perfect code, and the last one corresponds to a partition
of the $3$-cube into independent sets.

\subsection{Perfect and MDS Codes in Hypercubes}\label{s:MDS}

A subset of $Q^n_q$ is called a code with distance $d$ if
the Hamming distance between any two distinct vertices of the code is at least $d$. Let $q=p^k$ be a power of a prime $p$. Then, on
the set $Q_q$, we can define the structure of a Galois field $GF(q)$, and
the set $Q^n_q$ can be viewed as an $n$-dimensional vector
space over $GF(q)$. Linear subspaces of $Q^n_q$
are called $q$-ary {\sl linear codes}. 

Consider an arbitrary matrix $A$ of size $m\times n$ with elements
from $Q_q$ and the solution set of the system of equations $Ax=\bar{0}$ over
the field $GF(q)$.

\begin{claim}[\cite{Hamming50},\cite{Hamming}]\label{st_dop01}
If matrix $A$ does not contain collinear or zero columns, then
the set $C=\{x\in Q^n_q \ |\ Ax=\bar{0}\}$ does not contain distinct
vectors with Hamming distance less than $3$.
\end{claim}
\begin{proof}
Suppose that $Ax=Ay=0$. Then $A(x-y)=0$. Vector $(x-y)$ defines a nontrivial linear
combination  of the columns of $A$. Since the matrix has no
collinear columns, vector $(x-y)$ has at least three
nonzero coordinates. Therefore, $d(x,y)\geq 3$.
\end{proof}

\begin{claim}[\cite{Golay},\cite{Hamming50}]\label{c:codeHam}
If $n=\frac{q^t-1}{q-1}$ and $q$ is a prime power, then the
hypercube $Q^n_q$ has a $1$-perfect code $C$, i.e., there exists a
perfect coloring ${\bf 1}_{_C}$ with quotient matrix
$\begin{pmatrix}
0 & n(q-1) \\
1 & n(q-1)-1
\end{pmatrix}$.
\end{claim}
\begin{proof}
Consider a subset $\mathcal{A}\subset Q^t_q$ consisting of
all possible vectors whose first nonzero coordinate is equal to
$1$. It is clear that the set $\{ax \ |\ a\in GF(q)\setminus \{0\},
x\in \mathcal{A}\}$ consists of all nonzero vectors from $Q^t_q$,
and each vector is represented  uniquely. Then
$|\mathcal{A}|=\frac{|Q_q^t|-1}{q-1}=n$ and no two
vectors from $\mathcal{A}$ are collinear. Let $A$ be a matrix whose columns are
vectors from $\mathcal{A}$. By  Proposition
\ref{st_dop01},   $C=\{x\in Q^n_q \ |\
Ax=\bar{0}\}$ is a code with a distance of at least $3$. 

The matrix
$A$ contains an identity submatrix of order $t$, so it has rank exactly $t$. Therefore, the code $C$ is a linear
subspace of dimension $n-t$. Then $|C|=q^{n-t}$. Each
non-code vertex $x\not \in C$ has at most one code vertex at distance $1$,
since the distance between two codewords is at least $3$.
Since the number of vertices in the unit ball is $n(q-1)+1=q^t$, the number
of  vertices that are contained in $C$ or adjacent to vertices in
$C$ is $q^t|C|$. From the equality $|Q_q^n|=q^t|C|$ it follows that
every non-code vertex has exactly one neighbor in the code. Thus, the code
$C$ is 1-perfect. \end{proof}

The matrix $A$ is called the {\sl parity-check matrix} of the code $C$. A linear
$1$-perfect code $C$ defined above is called the {\sl Hamming code}. It can be shown that
all linear 1-perfect codes in a hypercube are equivalent and
can be obtained from $C$
by  linear transformations.

\begin{remark}
Vasil'ev \cite{Vas62} constructed the first nonlinear $1$-perfect codes. He proved that there are  at least  $2^{2^{n/2 -o(n)}}$ nonequivalent binary $1$-perfect codes with length $n=2^t-1$.  There exist many different constructions of  nonlinear $1$-perfect codes (see  \cite{Sol}, \cite{Etzion}).  It is a well-known problem to find an effective (i.e., less than  $2^{2^{n -o(n)}}$) asymptotic upper bound on the number of such codes. \label{Vasprob}
\end{remark}

\begin{example}
The parity-check matrix of a binary Hamming code of length $n=7$:
$A=\begin{pmatrix}
0 & 0 & 0 & 1& 1 & 1&1 \\
0 & 1 & 1 & 0 & 0 & 1 & 1\\
1 & 0 & 1 & 0 & 1 & 0 & 1
\end{pmatrix}$.
\end{example}

\begin{remark} From Proposition \ref{s:integer},  it is easy to see that the condition
$n=\frac{q^t-1}{q-1}$ is necessary for the existence of a
$1$-perfect code in the hypercube $Q^n_q$. For $q$ that is not  a prime power, the well-known  conjecture claims  that $1$-perfect codes in
$Q^n_q$, $n>1$, do not exist. \label{PCprob}
\end{remark} 

Consider a group structure with the operation $\circ$ on the set $Q_q$,
for example, an  abelian group $\mathbb{Z}_q=\mathbb{Z}/q\mathbb{Z}$. Define the set $$L=\{x\in Q_q^n \ |\ x_1\circ
x_2\circ\cdots\circ x_n=\varepsilon\}$$ where $\varepsilon$ is the identity element of the group.
\begin{claim}\label{s:integer11}
The function ${\bf 1}_{L}$ is a perfect coloring with the quotient matrix
 $\begin{pmatrix}
0 & n(q-1) \\
n & n(q-2)
\end{pmatrix}$. \end{claim}
\begin{proof}
Let $x\in L$ and $d(x,y)=1$ for some vertex $y\in Q^n_q$.
Suppose $y_i\neq x_i$. Then we have $x_1\circ\cdots \circ x_i=
(x_{i+1}\circ\cdots\circ x_n)^{-1}$, and $$y_1\circ\cdots \circ y_{i-1}\circ y_i=x_1\circ\cdots\circ x_{i-1} \circ y_i \neq
(x_{i+1}\circ\cdots\circ x_n)^{-1}. $$ Therefore, $y\not \in L$. From the
definition of a group, it follows that the equality $x_1\circ
x_2\circ\cdots\circ x_n =\varepsilon$ is equivalent to $(x_1\circ
x_2\circ\cdots\circ x_{n-1})^{-1}=x_n$. Consequently, all coordinates
of a vertex $x\in L$ except one can be chosen arbitrarily. Therefore,
$|L|=q^{n-1}$. Since $|L|q=|Q_q^n|$, each clique (one-dimensional
face)  of $Q^n_q$ should contain exactly one vertex from $L$. Each
vertex of  $Q^n_q$ is contained in $n$ cliques (one for each coordinate),
therefore, any non-code vertex has exactly $n$ neighbors in $L$.\end{proof}

As will be shown below, the first color of the  coloring from Proposition \ref{s:integer11}
is a maximal independent set in
$Q^n_q$. More precisely, it will be proved that the minimum eigenvalue
of the adjacency matrix of $Q_q^n$ is $-n$, so the cardinality
of the set of vertices of the first color reaches the
Delsarte--Hoffman bound (Proposition \ref{PCclaim11}).
Such subsets of the hypercube are called MDS codes (see Section  \ref{MDS})  with code
distance $2$.

\subsection{Problems}

\begin{exercise}\label{exe1} 
Prove that the sum of the eigenvalues of the adjacency matrix
of every simple graph is zero.
\end{exercise}

\begin{exercise} 
Prove that if the adjacency matrix of a bipartite graph
has an eigenvalue $\lambda$, then it also has an eigenvalue
$-\lambda$.
\end{exercise}

\begin{exercise}\label{exe3}
Prove that the adjacency matrices
of a connected regular graph $G$ and its complement have the same
eigenvectors. Calculate the eigenvalues of the adjacency matrix
of the complement if the eigenvalues of the adjacency matrix of
$G$ are known. 
\end{exercise}

\begin{exercise}
The Seidel matrix of a graph $G=(V,E)$ is defined by the equality
\[
Z_{n\times n}(G)= \left\{
\begin{array}{rr}
1 \mbox{, if}\ & (v_i,v_j) \in E(G);\\
-1\mbox{, if}\ & (v_i,v_j) \not \in E(G);\\
0\mbox{, if}\ & i=j.
\end{array}
\right.
\]
Prove that the Seidel matrix of a regular connected graph has the same
eigenvectors as the adjacency matrix of the graph. Calculate the
eigenvalues of the Seidel matrix if the eigenvalues
of the adjacency matrix of the graph are known. \end{exercise}

\begin{exercise}
Prove that a perfect coloring of a graph is a perfect
coloring of its complement.
\end{exercise}

\begin{exercise} Prove that a coloring of a regular graph $G$ is
perfect if the subgraphs induced by vertices of any two different colors are
biregular bipartite and  the  subgraphs induced by vertices of any
one color are regular.
\end{exercise}

\begin{exercise}
Prove that any linear set (code) in $Q^n_q$ (over the field
$GF(q)$) can be defined by a parity-check matrix $A$. Find a necessary
and sufficient condition on the columns of $A$ that ensures that
the code has distance $d$. \end{exercise}

\begin{exercise}
A subset  $S$   of a Boolean $n$-cube is called a Sidon set if
all pairwise sums of elements of $S$ are distinct. Prove that a linear
binary code has code distance of at least $5$ if and only if
the columns of its parity-check matrix, together with the all-zero vector,
form a Sidon set.
\end{exercise}

\section{Algebraic Properties of Perfect Colorings}

\subsection{Stochastic and  Symmetrizable Matrices}

A simple necessary condition for an integer matrix
to be the quotient  matrix of some perfect coloring
is the symmetry of the  set of positions of nonzero
elements (support) of the matrix. Indeed, if a vertex of the first color has neighbors
 of the second color, then a vertex of the second color also should also have neighbors of the first one. However, this condition is not sufficient. In this section, we  establish
necessary and sufficient conditions.

%\begin{definition}
A square matrix of nonnegative elements is called
{\sl stochastic} if the sums of its elements by rows (or by columns)
are equal to $1$. 
%\end{definition}

From the definitions, it is clear that if $S$ is the quotient matrix 
of a perfect coloring of an $r$-regular graph, then the matrix
$\frac{1}{r}S$ is stochastic.  Let $n_i$
be the number of vertices of the $i$th color, and
\begin{equation} \label{eq:colorededges0}
 \sum\limits_{i=1}^k n_i = n
\end{equation}
 be the number of vertices in the graph.
By calculating
the number of edges which are incident to vertices of colors $i$ and $j$, we obtain the following
relation:
\begin{equation} \label{eq:colorededges}
n_i s_{ij}=n_js_{ji}.
\end{equation}

We obtain a simple necessary condition for the existence of perfect colorings.
{
\begin{claim}\label{s:integer}
If there is a perfect coloring of a graph on $n$ vertices with a quotient matrix
$S$, then the system of equations (\ref{eq:colorededges0}), (\ref{eq:colorededges}) has
nonnegative integer solutions $n_i$.
\end{claim}
}

Let $N$ be a diagonal matrix of size $k\times k$, consisting of
elements $n_i$:
$$
N=
\begin{pmatrix}
n_1 & \ldots & 0\\
\vdots & \ddots & \vdots\\
0 & \ldots & n_k
\end{pmatrix}.
$$

From the equality (\ref{eq:colorededges}), we  conclude that $NS$
is a symmetric matrix, since we can rewrite the relation~(\ref{eq:colorededges}) in matrix form:
\[
\label{eq:matrixcolledges} S^{\mathsf T}N=NS.
\]
Matrices $S$ with this property are called {\sl symmetrizable}.
We  always can choose  a nonnegative integer   matrix $N$ if $S$ is nonnegative and
integer.

This implies that the following proposition is true.
\begin{claim}\label{s:S_sym}
The quotient matrix $S$  is symmetrizable
by a diagonal matrix, and there exists a basis consisting of  its eigenvectors.
 \end{claim}
\begin{proof}
Since $n_i>0$, let us consider the diagonal matrix
$N^{\frac12}$. We have the equality
$N^{-\frac12}S^{\mathsf T}N^{\frac12}=N^{\frac12}SN^{-\frac12}$.
Thus, the matrix $N^{\frac12}SN^{-\frac12}$ is symmetric, and,
consequently, there is a basis for the space $\mathbb{R}^k$ consisting
of its eigenvectors. The matrix $S$ is similar to the matrix
$N^{\frac12}SN^{-\frac12}$, and therefore there is also a basis of its
eigenvectors. \end{proof}

It turns out that the symmetrizability property of a stochastic matrix $S$ is not only necessary but also sufficient.  After multiplying by a common denominator
of its elements, $S$  becomes  a quotient matrix of some perfect coloring. More precisely, 

\begin{theorem}[\cite{Kir}, \cite{25}]\label{thAF}
Let the elements of an integer nonnegative symmetrizable matrix
$S$ satisfy the equalities $\sum_{j=1}^k s_{ij}=r$.
Then there exists an $r$-regular graph $G$ and its perfect coloring with quotient matrix $S$.
\end{theorem}
\begin{proof}
Take $i, j$, $1\leq i<j\leq k$.
We will show that if the equality $n_i
s_{ij}=n_js_{ji}$ holds, there exists a bipartite biregular graph with parts of cardinalities
 $n_i$ and $n_j$, and vertex degrees in the parts $s_{ij}$ and $s_{ji}$.
It is sufficient to construct the adjacency matrix for  such a biregular graph.  We need  a matrix  $B$  of size
$n_i\times n_j$ which  provides an adjacency of the parts of the  bipartite  graph.   Let the $k$th row
of matrix $B$ contain a block of ones at positions $(k-1)s_{ij}+1$ through
$ks_{ij}$, where both position numbers are considered modulo $n_j$.
We fill the remaining elements of matrix $B$ with zeros. Each row
of matrix $B$ contains exactly $s_{ij}$ ones by  construction. Furthermore, by the construction, the number of ones in the columns of matrix $B$ differs by no more than one. When the equality
(\ref{eq:colorededges}) holds, these numbers are equal to
$s_{ji}=n_i s_{ij}/n_j$.

By the condition, there exists a diagonal nonnegative matrix $N$ such that the matrix $NS$ is symmetric. We will assume that the elements of the matrix $N$ are even, since multiplication by $2$ preserves the equality
(\ref{eq:colorededges}). For each $i=1,\dots,k$, we take a set
$V_i$ of $n_i$ vertices, where $n_i$ is the $i$th diagonal element
of the matrix $N$. For each pair $i,j$, $i\neq j$, we draw edges so as
to obtain a bipartite biregular graph with parts $V_i$ and
$V_j$, and vertex degrees in parts $s_{ij}$ and $s_{ji}$. This is possible,
since $|V_i|s_{ij}=|V_j|s_{ji}$. For each $i=1,\dots,k$, we divide the vertices of the $i$th part into two equal parts, $V_i'$ and $V_i''$, and construct a bipartite $s_{ii}$-regular graph. This is possible by the  above construction,
since the required equality
$|V_i'|s_{ii}=|V_i''|s_{ii}$ holds. The desired graph has been constructed; it remains to
color the vertices of $V_i$ with the $i$th color for all $i=1,\dots,k$.
\end{proof}

Theorem \ref{thAF} can be generalized
to irregular graphs (see \cite{25} and \cite{Lei}).

We have established necessary and sufficient conditions to be a quotient matrix of a perfect
coloring of some graph. But can a given graph be colored
according to a given quotient matrix? This question remains unanswered in many specific cases.

\begin{remark}
Any stochastic matrix $T=(t_{ij})$ can be viewed as
the transition probability matrix of a Markov chain, where $t_{ij}$ are the conditional transition probabilities from state $i$ to state $j$.
Thus,  a perfectly colored graph
$G$ generates a Markov chain whose states
are the colors of the vertices, and the choice of the next  edge to move  along
is equiprobable. Note that the vector of  color frequencies 
$q=\left(\frac{n_1}{n}, \ldots, \frac{n_k}{n} \right)$ is
the vector of stationary probabilities of the states of the Markov chain,
since it satisfies the equation $q^{\mathsf T}T=q$.
However, not every discrete Markov chain can be interpreted as a random walk on an undirected colored graph. We can simulate only such chains  whose trajectory probabilities  are independent on the order in which the states are traversed: from right to left or from left to right. This property is equivalent to the symmetrizability of the matrix $T$.
\end{remark}

\subsection{Matrix Equation of a Perfect Coloring}

Let $f: V(G) \rightarrow \{1,\ldots,k\}$ be a graph coloring.
Consider a matrix $F=F[f]$ of size $n\times k$, in which column
$f_i$ is the characteristic function of color $i$:
\begin{equation}\label{eq:F}
f_i(x)=\left\{
\begin{array}{ll}
1 & \mbox{for}\ f(x)=i,\\
0 & \mbox{for}\ f(x)\ne i.
\end{array}
\right.
\end{equation}

The matrices $F[f]$ will also be called graph colorings, since
they are in  one-to-one correspondence with the functions $f: V(G) \rightarrow
\{1,\ldots,k\}$.  Note that each row of the matrix $F$ contains
exactly one $1$. We will assume throughout that
the function defined on the vertices of the graph is a column vector.

\begin{claim}[algebraic criterion of perfect coloring \cite{Martin},\cite{Haemers1},\cite{Godsil}]\label{s:perfcol_criteria}
Let $M$ be the adjacency matrix of $G$. $F=F[f]$ is a
perfect coloring
of a graph $G$ with quotient matrix $S$ if and
only if $MF=FS$.
\end{claim}
\begin{proof} To prove the equality, consider the $(v, j)$-th elements of the matrices on the left-
and right-hand sides of the equality $MF=FS$, $j=1,\dots,k$, $v\in V(G)$.
$(MF)_{v j}$ is the number of neighbors of color $j$ of vertex $v$.
Now consider the right-hand side of the equality. In the $v$th row of $F$, the entry $1$ is in the position
corresponding to the color $f(v)$.  We obtain that $(FS)_{v j}$ is the number of vertices of color $j$ adjacent to a vertex of color $f(v)$.
Thus, the right-hand side depends only on the color of vertex $v$.
\end{proof}

Note that if matrix $S$ is a solution to the equation $MF=FS$,
where $M$ is a nonnegative integer matrix and $F$
satisfies condition (\ref{eq:F}), then matrix $S$ is also nonnegative and integer.  Properties of arbitrary matrices
$F$ satisfying the equation $MF=FS$ were studied in \cite{Taranenko}.

\begin{corollary}[Lloyd's Theorem \cite{Lloyd}]\label{p:perf_color_properties}
Let  $F$  be  a perfect 
coloring with the  quotient matrix $S$ of a graph with adjacency matrix $M$.  Let $z$ be an eigenvector of $S$, i.e., $Sz=\lambda z$
and $z\ne \mathbf{0}$. Then $Fz$ is an eigenvector
of the adjacency matrix $M$ with the same eigenvalue $\lambda$.
\end{corollary}
\begin{proof}
By using Proposition~\ref{s:perfcol_criteria}
and the associativity of matrix multiplication, we obtain the following chain of equalities:
\[
MFz=FSz=F\lambda z=\lambda Fz.
\]
 If
$z\neq \mathbf{0}$ then $Fz\neq \mathbf{0}$ by the definition of $F$.
\end{proof}

Since $M$ is a symmetric matrix, all its eigenvalues
are real.
\begin{corollary}
Every eigenvalue of  a quotient matrix is 
real. \end{corollary}

Since $M$ is a symmetric matrix, there is a basis consisting of
its eigenvectors. That is, $\mathbb
R^n=X_{\lambda_1}\oplus\ldots\oplus X_{\lambda_m}$,
where $X_{\lambda_i}$ is the eigenspace of $M$,
corresponding to $\lambda_i$.

\begin{corollary}\label{p:perf_color_properties1}
Any  column vector of a perfect coloring $F$
 is contained in the subspace
$X_{\lambda_1}\oplus\ldots\oplus X_{\lambda_{m'}}$,
where $\lambda_1$, $\ldots$, $\lambda_{m'}$ are the eigenvalues of  $S$.
Moreover,  all columns of $F$  are not contained in $X_{\lambda_1}\oplus\ldots\oplus X_{\lambda_{m''}}$, where $m''<m'$.
\end{corollary}
\begin{proof}
By Proposition \ref{s:S_sym},  $S$ has a basis of eigenvectors. Therefore, the vector $e_i=(0\ldots \underset{i}{1} \ldots
0)$ can be expanded in this basis: $e_i=\sum \alpha_j z_j$, where
$z_j$
are the eigenvectors of the matrix $S$. We obtain:
\[Fe_i=F\sum \alpha_j z_j=\sum \alpha_j Fz_j\]
By Corollary \ref{p:perf_color_properties}, $Fz_j$ is an eigenvector of the matrix $M$,
thus, we have obtained the required expansion.
Since $\{e_i\}$ is a basis, all  eigenvectors of the matrix $S$ are used in the expansions. 
\end{proof}

Note that  column vectors of a perfect coloring $F$ are characteristic functions of the elements (cells)  of the equitable partition corresponding to the coloring.

A subspace $V$ is called {\sl invariant} with respect to
the matrix $M$ if $MV\subseteq V$. In particular, the eigensubspaces
of the matrix are invariant.   A {\sl projector} onto a linear subspace
is a linear operator acting  as the identity on this subspace and as zero 
 on its complement.  Consider the subspace
$V_F$ generated by the columns of the coloring matrix $F$. From the definition of a perfect coloring, it follows that

\begin{claim}\label{cl:invar}
$F$ is a perfect coloring if and only if $V_F$ is an invariant subspace.
\end{claim}
\begin{proof}
$F$ is a perfect coloring if and only if the rows
of the matrix $MF$ corresponding to vertices of the same color coincide. This means precisely that all columns of the matrix $MF$ are linear combinations
of columns of the matrix $F$.
\end{proof}
Let $F$ be a $k$-coloring of a graph with adjacency matrix $M$. Then $N=F^{\mathsf T}F$
is a diagonal matrix of size $k\times
k$, whose $(ii)$-th element is $n_i$, the number of vertices of color $i$.

\begin{claim}\label{claimprojector}
The linear operator defined by the matrix $P=FN^{-1}F^{\mathsf T}$ is a
projector to the subspace $V_F$ generated by the
characteristic vectors of the colors of the coloring $F$.
\end{claim}
\begin{proof}
Let $f_i$ be the $i$-th column of the matrix $F$, i.e.,
the characteristic vector of the $i$-th color. It is easy to see that
$Pf_i=f_i$ for any $i$, $i=1,\dots,k$. Since the rank of $F$
is $k$, the rank of  $P$ is no greater than $k$. Therefore,
the dimension of the image of the operator $P$ is $k$, and the vectors $f_i$
form a basis for the image. It is well-known that an operator is a projection  if and only if $P^2=P$.
We have $P^2=FN^{-1}F^{\mathsf T}FN^{-1}F^{\mathsf T}=FN^{-1}NN^{-1}F^{\mathsf T}=P$.
\end{proof}

It is easy to see that $P$ can be defined by the formula
\begin{equation}\label{eq:projector}
P(x,y)=\left\{
\begin{array}{ll}
\frac{1}{n_i} & \mbox{if vertices }\,x\  \mbox{and }\,y\  \mbox{have color}\,i,  \mbox{i.e.,}\ F(x,i)=F(y,i)=1,\\
0 & \mbox{otherwise}.
\end{array}
\right.
\end{equation}

\begin{claim}\label{claimprojector1} Let $G$ be a graph  with adjacency matrix $M$. A coloring $F$  of $G$ is perfect  if and only if $MP=PM$, where $P$ is the projector corresponding onto the subspace $V_F$.
\end{claim}
\begin{proof}
1. ($\Leftarrow$) Let us multiply the equality $MFN^{-1}F^{\mathsf T}=FN^{-1}F^{\mathsf T}M$ by
the matrix $F$ on the right. We have $MF=FN^{-1}F^{\mathsf T}MF$. Then, by
Proposition \ref{s:perfcol_criteria}, $F$ is a perfect coloring of $G$ with quotient matrix $S=N^{-1}F^{\mathsf T}MF$.

2. ($\Rightarrow$) Let $F$ be a perfect coloring with quotient matrix
$S$. Then, by Proposition \ref{s:perfcol_criteria}, we have the equality $MF=FS$. By the definition of $P$, we have $PF=F$. Then $MF=PFS=PMF$.
Multiplying the equality $MF=PMF$ on the left  by $N^{-1}F^{\mathsf T}$, we obtain
$MP=PMP$. Since  $P$ and $M$ are symmetric, we have
$PM=(MP)^{\mathsf T}=(PMP)^{\mathsf T}=PMP$.
\end{proof}

\subsection{Theorem on the Spectral Decomposition of  Adjacency Matrices}\label{2.3}

Let $X_{\lambda_i}$ be an eigenspace of a symmetric
matrix $M$ corresponding to an eigenvalue $\lambda_i$.
Consider some orthonormal basis
$\{u^i_1,\dots,u^i_{m_i}\}$ of $X_{\lambda_i}$ and a matrix
$U_i$ of size $n\times m_i$ whose columns are the vectors
$u^i_j$, $j=1,\dots,m_i$. We define a matrix $E_{\lambda_i}$ of size
$n\times n$ by the equality $E_{\lambda_i}=U_i U_i^{\mathsf T}$.

\begin{theorem}[on spectral decomposition \cite{Godsil}]
For any symmetric matrix $M$, the matrices $E_{\lambda_i}$
satisfy the conditions \\
{\rm (a)} $E^2_{\lambda_i}=E_{\lambda_i}$; $E_{\lambda_i}E_{\lambda_j}=0$
for $i\neq j$;\\
{\rm (b)} $\sum_iE_{\lambda_i}=I_n$; the matrix $E_{\lambda_i}$ does not depend on the choice of an orthonormal basis
of the space $X_{\lambda_i}$;\\
{\rm (c)} $ME_{\lambda_i}=\lambda_iE_{\lambda_i}$; For any polynomial
$p$, the following equality holds:\\
$p(M)=\sum_ip(\lambda_i)E_{\lambda_i}$.
\end{theorem}
\begin{proof}

(a) By the definition, $E_{\lambda_i}=U_i U_i^{\mathsf T}$. Since the basis is orthonormal, we have $U_i^{\mathsf T}U_i=I_{m_i}$. Then $E^2_{\lambda_i}=
U_i(U_i^{\mathsf T}U_i) U_i^{\mathsf T}= E_{\lambda_i}$. The eigenspaces
of a symmetric matrix corresponding to different eigenvalues
are orthogonal, therefore $U_i^{\mathsf T}U_j=0$ for $i\neq j$. Then
$E_{\lambda_i}E_{\lambda_j}=U_i(U_i^{\mathsf T}U_j) U_j^{\mathsf T}=0$.

(b) The squares of the eigenvalues of  $M$ are the eigenvalues
of  $M^2$. From the equality $E^2_{\lambda_i}=E_{\lambda_i}$
it follows that the eigenvalues of  $E_{\lambda_i}$ are either $0$ or
$1$. The rank of  $E_{\lambda_i}$ is equal to the dimension of the subspace
$X_{\lambda_i}$, since the rank of  $I_{m_i}=U_i^{\mathsf T}(U_iU_i^{\mathsf T})U_i$
is equal to $m_i$ and does not exceed the rank of  $E_{\lambda_i}=U_i
U_i^{\mathsf T}$. Then the eigenvalue $1$ has multiplicity $m_i$ and
$\tr(E_{\lambda_i})= m_i$. Let $E=\sum_iE_{\lambda_i}$. We have
$\tr(E)=\sum_i\tr(E_{\lambda_i})=\sum_im_i=n$. Furthermore, from item (a)
it follows that $E^2=E$ and  $E$ is symmetric, due to
the symmetry of the terms. This means that all  eigenvalues of $E$ is equal to $1$. Consequently, $E=I_n$. If we change
the orthonormal basis in one of the eigenspaces
$X_{\lambda_i}$, then in the equality $\sum_iE_{\lambda_i}=I_n$, only one term can
change. However, this is impossible in the equality.

(c) $ME_{\lambda_i}=(MU_i) U_i^{\mathsf T}=\lambda_i
U_iU_i^{\mathsf T}=\lambda_iE_{\lambda_i}$. Then
$M^kE_{\lambda_i}=\lambda^k_iE_{\lambda_i}$ for any natural
$k$. From item 2), we have $M^kI_n= \sum_i\lambda^k_iE_{\lambda_i}$. Then
for any polynomial $p$, we have the equality
$p(M)=\sum_ip(\lambda_i)E_{\lambda_i}$.
\end{proof}

\begin{corollary}\label{colE}
The matrices
$E_\lambda$ can be represented as polynomials of
$M$. \end{corollary}
\begin{proof}
From the spectral decomposition theorem, we have the equality $p(M)=\sum_i
p(\lambda_i)E_{\lambda_i}$, where $p$ is an arbitrary polynomial.
Consider the polynomial $q(t)=(\prod_i(t-\lambda_i))/(t-\lambda)$. We have
$q(M)=q(\lambda)E_\lambda$.
\end{proof}

{
\begin{corollary}\label{corEM}
The dimension of the linear space of matrices generated by powers
of a symmetric matrix $M$ is equal to the number of distinct eigenvalues
of $M$. \end{corollary}
\begin{proof}
The matrices $E_{\lambda_i}$ are linearly independent by item 1) and constitute a basis for this space by item 3) of the spectral decomposition theorem.
\end{proof}
}

Recall that  $J_n$ is the matrix of size $n\times n$, all elements of which
are equal to $1$. For the maximum eigenvalue of an $r$-regular connected graph, the eigenspace $X_r$
consists of vectors collinear with the vector ${\bf 1}$ (see Problem
\ref{exe8}). Therefore, $E_r=\frac{1}{n}J_n$. Then $J_n$ can be represented as polynomials in the adjacency matrix
$M$ by Corollary \ref{colE}. We will show that this property of a regular
graph is a criterion.

\begin{theorem}[Hoffman \cite{Godsil}]
A graph $G$ is regular and connected if and only if
there exists a polynomial $p$ such that $J_n=p(M)$, where $M$ is the adjacency matrix
of $G$.
\end{theorem}
\begin{proof}

$(\Leftarrow) $ If $J_n = p(M)$ for some polynomial $p$, then
$Mp(M) = p(M)M$ and $MJ_n = J_nM$, which means that the graph
$G$ is regular. If $G$ is disconnected, then there exist vertices  $i$
and $j$ such that the element $(i, j)$ in the matrix $M^d$ is zero for
any degree $d$.

$(\Rightarrow)$ Let $G$ be an $r$-regular connected graph and $p(x) =
(x-r)(x-\lambda_1)\cdots(x-\lambda_l) = (x-r)q(x)$, where
$r,\lambda_1,\dots,\lambda_l$ are all distinct eigenvalues
of the matrix $M$. Since $p(M) = 0$, we have $Mq(M) = rq(M)$. Therefore,
the columns of the matrix $q(M)$ are the eigenvectors with
the maximum eigenvalue $r$. Since $G$ is connected, all
these columns are collinear with the vector ${\bf 1}$. Since $q(M)$ is
symmetric, all its rows are equal to its columns; hence all entries of $q(M)$ are equal.
Therefore, $q(M) = cJ_n$
for some $c$. \end{proof}

\subsection{Perfect Colorings and Eigenfunctions}

Consider a perfect coloring of a connected $r$-regular graph in
two colors. Let $S=\begin{pmatrix}
a & b \\
c & d
\end{pmatrix}=\begin{pmatrix}
r-b & b \\
c & r-c
\end{pmatrix}$ be the matrix of parameters of this coloring.
The main parameters of the coloring are the numbers $b, c$, i.e., the numbers
of neighboring vertices of a different color. The eigenvalues
of this matrix are  $\lambda_1=r-(b+c)$ and $\lambda_2=r$.
Note that the eigenvalues of the quotient  matrix of size $2\times 2$ are integer and 
 the adjacency matrix of  the $r$-regular graph has the same  eigenvalues.

We calculate the number of vertices of each color
($n_1$ and $n_2$, respectively) by the system of equations:\\
$n_1 b=n_2 c$
is the number of two-color edges and  \\
$n_1+n_2=n=|V(G)|$.\\
From this, we obtain the equalities $n_1 =
\frac{nc}{b+c}=\frac{n}{1+\frac{b}{c}}$ and $n_2 =
\frac{nb}{b+c}=\frac{n}{1+\frac{c}{b}}$.

{ In this case, a necessary condition for the existence of a perfect
graph coloring (Proposition \ref{s:integer}) means the divisibility
of the numbers $nc$ and $nb$ by $b+c$.}

Let $A$ be the set of vertices of the first color and let $\mathbf{1}$
be a vector of ones, which is an eigenvector of the matrix $S$ with eigenvalue $r$. From Corollary \ref{p:perf_color_properties1}, we obtain
that ${\bf 1}_{A} = f + \alpha\mathbf{1}$, where $f$ is the eigenfunction corresponding to the number $\lambda_1$. We have
$n_1=({\bf 1}_{A},\mathbf{1})= \alpha(\mathbf{1},\mathbf{1})=\alpha
n$, from which we obtain:
\[
f(x)=\left\{
\begin{array}{ll}
1-\frac{c}{b+c} & x\in A;\\
-\frac{c}{b+c} & x\notin A.
\end{array}
\right.
\]

Therefore, a perfect 2-color coloring can be associated
with an eigenfunction that takes only two different values. We will now show that the converse is also true.
For colorings with a larger number of colors, such a correspondence
arises under certain additional conditions.  By Corollary \ref{p:perf_color_properties} we have 

\begin{corollary}\label{z:colorlin0}
Consider a perfect $k$-coloring of a regular graph $G$ with quotient matrix $S$. Let $v=(v_1,\dots,v_k)$ be an eigenvector of $S$, where all components of the vector $v$ are distinct.
Then the function $Fv$
is an eigenfunction of graph $G$, which takes values $v_i$ at points of color $i$.
\end{corollary}

This means that $Fv$ determines the same coloring as the initial perfect $k$-coloring.

{
\begin{claim}\label{c:colorlin}
Let $f$ be a $k$-coloring of a regular graph $G$  and
let there be $k$ linearly independent eigenfunctions of the adjacency matrix
of the graph, which can be obtained as linear combinations
of the characteristic functions of the colors. Then $f$ is a
perfect coloring.
\end{claim}
\begin{proof}
Let the eigenfunctions $(h_1,\dots, h_{k})$ be linear
combinations of the characteristic functions of the colors $(f_1,\dots, f_{k})$, i.e., $(h_1,\dots, h_{k}) = (f_1,\dots, f_{k}) A$ or
$H=FA$, where the eigenfunctions are the columns of the matrices $H$. $A$
is a nonsingular matrix because the eigenfunctions $h_1,\dots, h_{k}$
are linearly independent by assumption. Since $h_1,\dots, h_{k}$ are eigenvectors of the adjacency matrix $M$ of  $G$, we have the equality
$MH=H\Lambda$, where $\Lambda$ is a diagonal matrix and
$\lambda_{ii}$ is an eigenvalue corresponding to the vector $h_i$.
We obtain the equality $MFA=FA\Lambda$.
Then, by
Proposition \ref{s:perfcol_criteria}, the coloring $f$ is
perfect with quotient matrix $A\Lambda A^{-1}$. \end{proof}
}

\begin{corollary}[\cite{FdF2}]\label{corol5}
Let $f$ take only two values on the vertices of an $r$-regular
graph $G$ and $f=g+c\mathbf{1}$, where $g$ is an eigenfunction of the adjacency matrix
of $G$ with an eigenvalue distinct from $r$. Then $f$ is a
perfect $2$-coloring.
\end{corollary}
\begin{proof}
Let $v_0$ be a vertex. We define the characteristic
functions of the colors $f_1(v)=1$ if $f(v)=f(v_0)$ and $f_1(v)=0$ if
$f(v)\neq f(v_0)$; $f_2=\mathbf{1}-f_1$. Then the sum $f_1+f_2=\mathbf{1}$
is an eigenfunction with eigenvalue $r$ and $g=f-c\mathbf{1}$ is a
linear combination of $f_1$ and $f_2$. Then $f$
is a perfect $2$-coloring.
\end{proof}

\subsection{Orthogonality of Perfect Colorings}

\begin{claim}\label{c:color_distr_proportion}
Let $G$ be a connected $r$-regular graph, and let $g$ and $h$ be perfect
colorings of $G$, such that the sets of eigenvalues of the quotient matrices
corresponding to $g$ and $h$ intersect only at the maximum
eigenvalue $r$. Then, for any color $i$ of coloring $g$ and
any color $j$ of coloring $h$, the proportion of vertices of color $j$ among the vertices
of color $i$ is the same as the proportion of vertices of color $j$ in the entire graph.
\end{claim}
\begin{proof}
Denote by $B$ the set of vertices of color $i$ and by $A$ the set
of vertices of color $j$. We define the function
$f=\mathbf{1}_{B} - \frac{|B|}{|V(G)|}\mathbf{1}$. By the construction,
$(f,\mathbf{1})=0$. By Corollary \ref{p:perf_color_properties1}
the function $f$ is contained in the union of the subspaces $X_{\lambda_1}\oplus\cdots \oplus X_{\lambda_k}$, where $\lambda_1,\dots,\lambda_k$  are eigenvalues
of the coloring $g$ without the eigenvalue $r$. Then
$(f,{\bf 1}_{A})=0$, since by the hypothesis ${\bf 1}_{A}$ lies in
a subspace that is orthogonal to the function $f$ (see Corollary
\ref{p:perf_color_properties1}). Consequently,
$(f,{\bf 1}_{A})=|B\cap A|-\frac{|B||A|}{|V(G)|}=0$. Then
$\frac{|B\cap A|}{|B|}=\frac{|A|}{|V(G)|}$.
\end{proof}

In probability theory, two partitions $\{A_i\}$ and $\{B_j\}$
of the event space are called independent if the events $A_i$ and $B_j$ are pairwise independent for any $i$ and $j$, i.e., the equalities
$Pr(A_iB_j)=Pr(A_i)Pr(B_j)$ hold. Thus,
equitable partitions whose spectra intersect only at
the maximum eigenvalue can be considered as
independent partitions of the graph vertices (here $Pr(A)=|A|/|V(G)|$).

By the proof  of Proposition \ref{c:color_distr_proportion} we see that the following statement is also true. 

\begin{claim}\label{c:color_distr_proportion0}
Let two sets of eigenspaces  containing the characteristic functions of  vertex sets
$A$ and $B$, respectively, intersect only in a subspace
consisting of constant functions. Then
$\frac{|A|}{|V(G)|}=\frac{|A\cap B|}{|B|}$.
\end{claim}

For hypercubes, Proposition \ref{c:color_distr_proportion0} is
a special case of the Proposition  \ref{c:plansh}, which will be
formulated and proved in Section \ref{Four}.

%\textcolor{green}{
\begin{claim}\label{PCclaim101}
Let $f,g:V(G)\rightarrow \{0,1\}$ be perfect $2$-colorings
of a connected regular graph $G$. Suppose that the vertex sets
of colors $f^{-1}(1)$ and  $g^{-1}(1)$ in the coloring
$g$ are disjoint ($f^{-1}(1)\cap g^{-1}(1) =\varnothing$). Then the coloring $f+g$
is perfect.
\end{claim}
\begin{proof}
%From Proposition \ref{c:color_distr_proportion} it follows that
%the eigenvalues of the quotient matrices of the colorings $f$ and $g$ coincide.
Since
\[
\frac{|f^{-1}(1)\cap g^{-1}(1)|}{|f^{-1}(1)|}=0
\]
while
\[
\frac{|g^{-1}(1)|}{|V(G)|}>0,
\]
Proposition \ref{c:color_distr_proportion} implies that the spectra of the two quotient matrices
must have a nontrivial common eigenvalue (i.e., they do not intersect only at 
 $r$). Since each quotient matrix of a perfect $2$-coloring has exactly one nontrivial eigenvalue, these nontrivial eigenvalues must coincide.
Therefore, the colorings $f$ and $g$ correspond to eigenfunctions
$f_1=f-\frac{|f^{-1}(1)|}{|V(G)|}{\mathbf 1}$ and
$g_1=g-\frac{|g^{-1}(1)|}{|V(G)|}{\mathbf 1}$ with the same
eigenvalue. We denote
$c=\frac{|f^{-1}(1)|+|g^{-1}(1)|}{|V(G)|}$. The sum $h=f_1+g_1$ is also
an eigenfunction and takes two values $1-c$ and $-c$.
Therefore, by Corollary \ref{corol5}, the function $h$ is a perfect
coloring. The function $f+g=h+c{\mathbf 1}$ differs from $h$ by a
constant, and therefore it is a perfect coloring with
values $0$ and $1$.
\end{proof}

\subsection{Problems}

\begin{exercise}
Let $P$ be a projector corresponding to a coloring $F$  (see Proposition
\ref{claimprojector}). Prove that a coloring $F$ of a graph with adjacency matrix $M$
is perfect if and only if
the equality $MF=PMF$ holds. In other words, the columns of $MF$ are
eigenvectors of the projector $P$.
\end{exercise}

\begin{exercise}
Prove that for any invariant subspace of a symmetric
matrix $A$, there exist matrices $F$ and $S$ such that the image
$F(\mathbb{R}^k)$ is this invariant subspace and
the equality $AF=FS$ holds.
\end{exercise}

\begin{exercise}
 Let $M$ be an adjacency matrix of a regular graph and let $S$ be a quotient matrix of a perfect coloring of this graph.
Prove that the polynomial $\det(S-xI)$ divides the polynomial $\det(M-xI)$. \end{exercise}

\begin{exercise}\label{exe8}
 Let $M$ be an adjacency matrix of an $r$-regular graph and let $S$ be a quotient matrix of a perfect coloring of this graph.
Prove that the maximal eigenvalues of the matrices $M$ and $S$
 are equal to $r$  and that the eigenspace $X_{r}$ of the {connected graph} is one-dimensional.
\end{exercise}

\begin{exercise}
We define the support of the matrix $A$ as the matrix $\supp(A)$ by the equality
$(\supp(A))_{ij}=\left\{
\begin{array}{ll}
1 & \mbox{for}\ (A)_{ij}\neq 0,\\
0 & \mbox{for}\ (A)_{ij}=0.
\end{array}
\right.$ Prove that any symmetric $(0,1)$-matrix
is the support of some symmetrizable stochastic matrix.
\end{exercise}

\begin{exercise}\label{exe13}
Let $d$ be the diameter of a graph. Prove that its adjacency matrix
has at least $d+1$ distinct eigenvalues.
\end{exercise}

\begin{exercise}
Prove that for any finite set of eigenvalues, there are only
finitely many (possibly $0$) connected regular graphs with such
set of eigenvalues.
\end{exercise}

\begin{exercise}
Prove that half of the vertices of any $1$-perfect code in a Boolean
hypercube have even weight (number of ones), and the other half have odd weight. \end{exercise}

\begin{exercise}
Let $v=(v_1,\dots,v_k)$ be an eigenvector of quotient matrix $S$ of a $k$-perfect coloring $F$ of a connected regular graph $G$, where all components of the vector $v$ are distinct. Prove
that the linear span of the columns of $F$ coincides with the linear
span of the set of vectors $f^0,f^1,\dots,f^{k-1}$, where $f=Fv$ and  components  of $f^j$ are $j$th powers of the corresponding  components  of $f$.
\end{exercise}

\begin{exercise}\label{1p:Pot04b}
 Let $G$ be a regular graph with adjacency matrix~$M$.
An $M$-invariant subspace~$V$ is closed with respect to the Hadamard (pointwise) product if and only if it is generated by
the characteristic functions of the colors of a perfect coloring of~$G$.
\end{exercise}

\section{Extremal Properties of Perfect Colorings}

\subsection{Extremal Property of Perfect $2$-Colorings}

Many years ago, Avgustinovich  postulated the following
principle in his report. Each perfect coloring is a solution of some optimization
problem and on the other hand, a tight (reached a theoretical bound)
solution of any optimization problem in combinatorics is equivalent
to a perfect coloring with some parameters that depend on the
problem. We illustrate this principle below.

Let $f:V(G)\rightarrow \{1,2\}$ be a coloring
of graph $G$. Denote by $a(x)$ the number of vertices
of the first color adjacent to the vertex $x$ of the first color,
and by $c(y)$ the number of vertices of the first color adjacent to  a vertex $y$ of the second color.

\begin{claim}\label{PCclaim10}
Suppose that for some coloring $f$ of graph $G$ and for all  vertices $x$
and $y$, we have $a(x)\leq a$ and $c(y)\leq c$. Then the proportion of vertices
of the first color in an $r$-regular graph $G$ is at most
$\frac{c}{r-a+c}$. If this proportion is attained, then $f$
is a perfect coloring of the graph $G$ with quotient matrix $S=\begin{pmatrix}
a & r-a \\
c & r-c
\end{pmatrix}$.
\end{claim}
\begin{proof}
Let $n=|V(G)|$, $n_i=|f^{-1}(i)|$, and $f_i$ be the characteristic
function of the $i$-th color. From the assumptions on the quantities $a(x)$ and $c(y)$, we obtain the inequalities
\begin{equation}\label{PCe10}
(Mf_1,f_1)\leq an_1,\quad (Mf_1, \mathbf {1} - f_1)=(f_1, Mf_2)
\leq cn_2= c(n-n_1)
\end{equation}
and $(Mf_1,\mathbf {1})=rn_1$.
Since $$
(Mf_1,\mathbf{1})=(Mf_1, \mathbf{1} - f_1)+(Mf_1,f_1), $$ we  obtain 
$rn_1\leq 
c(n-n_1)+an_1$.

Suppose that $a(x)< a$ for at least one vertex of the first color  or that $c(y)<c$
for at least one vertex of the second color. Then at least  one of the
inequalities (\ref{PCe10}) is strict.
In this case, we obtain $rn_1<
c(n-n_1)+an_1$ that is, $\frac{n_1}{n}< \frac{c}{r+c-a}$. Consequently, if the required proportion is attained, then $a(x)= a$ for
any vertices of the first color and $c(y)=c$ for any vertices of the second
color.
\end{proof}

\begin{corollary}
Let $G$ be an $r$-regular graph and let  $A$ be some subset of its
vertices. Suppose that any ball of radius $1$ contains at most one
vertex from $A$. If $\frac{|A|}{|V(G)|}\geq \frac{1}{r+1}$, then $A$
is a perfect code.
\end{corollary}
\begin{proof}
We assume that the vertices of $A$ are colored the first color, and
the rest are colored the second. Then, by assumption, $a=0$ and $c=1$.  Therefore, by Proposition \ref{PCclaim10}, the coloring
is perfect with quotient matrix $\begin{pmatrix}
0 & r \\
1 & r-1
\end{pmatrix}$.
\end{proof}

\subsection{Cuts and Independent Sets in Graphs}\label{3.1}

Let $A$ and $B$ be subsets of the  vertex set of a simple graph $G$. We denote by
$e(A,B)=|\{\{a,b\}\in E(G) : a\in A, b\in B\}|$  the number
of edges between $A$ and
$B$.  Note that if $a,b\in A\cap B$, then the edge $\{a,b\}$ is counted
twice in $e(A,B)$. So, $e(A,B)=(\mathbf{1}_{A}, M\mathbf{1}_{B})$, where $M$
is the adjacency matrix of $G$.

\begin{lemma}[Expander Mixing Lemma \cite{AlonCh}, \cite{Dev}]\label{EMLclaim}
Let $G$ be an $r$-regular connected graph and let $\lambda$ be an eigenvalue of its adjacency matrix $M$ of maximum absolute value among all eigenvalues different from $r$, i.e., $|\lambda|>|\lambda'|$ for all eigenvalues
$\lambda'$, $\lambda'\neq \lambda$, $\lambda'\neq r$.
Then, for any subsets of vertices $A$ and $B$,  we have
$$|e(A,B)-\frac{r|A||B|}{|V(G)|}|\leq
|\lambda|\sqrt{|A||B|(1-|A|/|V(G)|)(1-|B|/|V(G)|)}.$$ Moreover,
if equality holds, then $\mathbf{1}_{A}$ is a perfect
$2$-coloring with eigenvalue $\lambda$ and
$B=A$ or $B=V(G)\setminus A$.
\end{lemma}
\begin{proof}
The characteristic functions $\mathbf{1}_{A}$ and $\mathbf{1}_{B}$ can be represented as
$\mathbf{1}_{A}=\sum_i\alpha_i\phi_i$ and
$\mathbf{1}_{B}=\sum_i\beta_i\phi'_i$, where $\phi_i$ and $\phi'_i$ are the normalized eigenfunctions of the matrix $M$ corresponding
to the eigenvalue $\lambda_i$. The normalized eigenfunction
corresponding to the eigenvalue $r$ is the constant
$\phi_0=\phi'_0=\mathbf{1}/\sqrt{n}$, $n=|V(G)|$. We have
$\alpha_0=(\mathbf{1}_{A},\phi_0)=|A|/\sqrt{n}$,
$\beta_0=(\mathbf{1}_{B},\phi_0)=|B|/\sqrt{n}$. Then
\begin{equation}\label{eqEx1}
(\mathbf{1}_{A}-\alpha_0\phi_0, M(\mathbf{1}_{B}-\beta_0\phi_0))=
\sum\limits_{i\neq 0}\lambda_i\alpha_i\beta_i(\phi_i,\phi'_i).
\end{equation}

From the Cauchy--Schwarz  inequality we have
\begin{equation}\label{eqEx2}
|\sum\limits_{i\neq 0}\lambda_i\alpha_i\beta_i(\phi_i,\phi'_i)|\leq
|\lambda|\sum\limits_{i\neq 0}|\alpha_i\beta_i|\leq
|\lambda|(\sum\limits_{i\neq 0}\alpha^2_i)^{1/2}(\sum\limits_{i\neq
0}\beta^2_i)^{1/2}.
\end{equation}

On the other hand
$$(\mathbf{1}_{A}-\alpha_0\phi_0, M(\mathbf{1}_{B}-\beta_0\phi_0))=(\mathbf{1}_{A},
M\mathbf{1}_{B})-\beta_0(\mathbf{1}_{A}, M\phi_0)-(\alpha_0\phi_0, M(\mathbf{1}_{B}-\beta_0\phi_0))
$$
\begin{equation}\label{eqEx33}
 =e(A,B)-\frac{r|A||B|}{n},
\end{equation}
since $(\phi_0, M(\mathbf{1}_{B}-\beta_0\phi_0))=(M\phi_0, \mathbf{1}_{B}-\beta_0\phi_0)=
(r\phi_0, \mathbf{1}_{B}-\beta_0\phi_0)=0$.

From the equality $(\mathbf{1}_{A},\mathbf{1}_{A})= \sum_i\alpha^2_i$ it follows that
$\sum\limits_{i\neq 0}\alpha^2_i=
(\mathbf{1}_{A},\mathbf{1}_{A})-\alpha^2_0=|A|- |A|^2/n$. A similar
inequality holds for $B$. Then from 
(\ref{eqEx1})--(\ref{eqEx33}) we have

$$\left|e(A,B)-\frac{r|A||B|}{n}\right|=
|(\mathbf{1}_{A}-\alpha_0\phi_0, M(\mathbf{1}_{B}-\beta_0\phi_0))|\leq$$
$$\leq |\lambda|\sqrt{(|A|- |A|^2/n)(|B|- |B|^2/n)}.$$

Equality in the Cauchy--Schwarz inequality holds only
if the vector $\mathbf{1}_{A}-\alpha_0\phi_0$ is collinear
with the vector $\mathbf{1}_{B}-\beta_0\phi_0$, i.e. $ \alpha_i=\beta_i$ or $
\alpha_i=-\beta_i$ for all $i\neq 0$. Then
equality in the first inequality of (\ref{eqEx2}) holds only if
 $\mathbf{1}_{B}-\beta_0\phi_0$ is an eigenfunction, with
eigenvalue $\lambda$. By Corollary \ref{corol5}, we obtain that ${\mathbf 1}_{B}$
is a perfect $2$-coloring. If the vector
$\mathbf{1}_{A}-\alpha_0\phi_0$ is codirectional with the vector
$\mathbf{1}_{B}-\beta_0\phi_0$, then $B=A$; if it is 
directed oppositely, then $B=V(G)\setminus A$. \end{proof}

Let's consider in detail the cases where $B=V(G)\setminus A$ and $B=A$.

%\begin{definition}
A {\sl cut} of a graph $G$ is a pair of vertex sets $A\subset
V(G)$ and $B=V(G)\setminus A$, or the set of edges between them.
The {\sl cut size} is the number $e(A,B)$. 
%\end{definition}
In the graph theory it is 
considered problems of finding both minimum and maximum
cuts of a graph.

This following inequalities were proved in 
\cite{AlonM} and  \cite{Tanner} but Golubev \cite{Golubev} proved
that equality is attained if and only if the corresponding partition is equitable. 

\begin{corollary}[\cite{Golubev}]\label{cut}
Let $G$ be an $r$-regular connected graph and let
$\lambda_k<\lambda_{k-1}<\dots<\lambda_{1}<\lambda_0=r$
be the eigenvalues of its adjacency matrix $M$. Let $A\subset
V(G)$, $B=V(G)\setminus A$. Then the following inequalities hold:
$$ \frac{(r-\lambda_1)|A||B|}{|V(G)|}\leq e(A,B)\leq \frac{(r-\lambda_k)|A||B|}{|V(G)|}.$$
Moreover, if  equality holds in either  inequality, then $\{A, B\}$ is an equitable partition.
\end{corollary}
\begin{proof}
The proof is similar to the proof of Lemma \ref{EMLclaim}. In this
case we have ${\mathbf 1}_{B}=\sqrt{n}\phi_0-{\mathbf 1}_{A}$, therefore the equality
(\ref{eqEx1}) takes the form
$$({\mathbf 1}_{A}-\alpha_0\phi_0, M({\mathbf 1}_{B}-\beta_0\phi_0))=
-\sum\limits_{i\neq 0}\lambda_i\alpha_i^2$$ and
$$-\lambda_1\sum\limits_{i\neq 0}\alpha_i^2\leq -\sum\limits_{i\neq
0}\lambda_i\alpha_i^2\leq -\lambda_k\sum\limits_{i\neq
0}\alpha_i^2.$$ Using  equality (\ref{eqEx33}), we obtain
required inequality. Obviously, equality holds 
only if ${\mathbf 1}_{A}-\alpha_0\phi_0$ is an
eigenfunction of  $M$ with the eigenvalue $\lambda_1$
or $\lambda_k$.
\end{proof}

Corollary \ref{cut} implies

\begin{corollary}
The size of a maximum cut of  $G$ is at most
$(r-\lambda_k)|V(G)|/4$, and if equality holds, then the cut corresponds to a
perfect $2$-coloring associated with the minimal eigenvalue.
\end{corollary}

\begin{corollary}[Haemers \cite{Haemers}]\label{haemers}
Let $G$ be an $r$-regular connected graph and let
$\lambda_k<\lambda_{k-1}<\dots<\lambda_{1}<\lambda_0=r$
be the eigenvalues of its adjacency matrix $M$. Let $C\subset
V(G)$, $C\neq V(G)$, and let $e(C,C)$ denote twice the number of edges induced by
$C$. Then the following inequalities hold:
$$ \lambda_k|C|+\frac{(r-\lambda_k)|C|^2}{|V(G)|}\leq e(C,C)\leq
\lambda_1|C|+\frac{(r-\lambda_1)|C|^2}{|V(G)|}.$$
Moreover, if equality holds  either inequalities, then $\{C, V(G)\setminus C\}$ is an equitable partition. 
\end{corollary}
\begin{proof}
We proceed similarly to the proof of Lemma \ref{EMLclaim}. In this
case, we have ${\mathbf 1}_{C}={\mathbf 1}_{B}={\mathbf 1}_{A}$, so the equality
(\ref{eqEx1}) takes the form $$({\mathbf 1}_{C}-\alpha_0\phi_0,
M({\mathbf 1}_{C}-\beta_0\phi_0))= \sum\limits_{i\neq
0}\lambda_i\alpha_i^2$$ and $$\lambda_k\sum\limits_{i\neq
0}\alpha_i^2\leq \sum\limits_{i\neq 0}\lambda_i\alpha_i^2\leq
\lambda_1\sum\limits_{i\neq 0}\alpha_i^2.$$ Using the equality
(\ref{eqEx33}), we obtain the required inequality. Obviously, the equality holds only if
${\mathbf 1}_{A}-\alpha_0\phi_0$ is an eigenfunction of  $M$
with eigenvalue $\lambda_1$ or $\lambda_k$.
\end{proof}

Recall that
an independent subset of graph vertices is a subset that does not
contain adjacent vertices.
%\begin{definition}
A subset  $C\subseteq V(G)$  is called
{\sl dominating} if every vertex of $G$  either belongs to $C$ or
is adjacent to a vertex of $C$.
%\end{definition}
From the definition, it is clear that every inclusion-maximal independent set is dominating.
It is easy to see that every color of a perfect $2$-coloring of a connected graph generates a dominating set.
Moreover, if a regular graph has a
$1$-perfect code, then it is a minimum dominating
set.

\begin{claim}[Delsarte--Hoffman bound \cite{Hoff}]\label{PCclaim11}
Let $G$ be an $r$-regular graph and let $\lambda_{min}$ be the minimum eigenvalue of its
adjacency matrix $M$. Then the cardinality of any independent
set in $G$ does not exceed
$\frac{\lambda_{min}|V(G)|}{\lambda_{min}-r}$. 
An independent set attains this bound if and only if its characteristic function is a perfect coloring with quotient matrix
 $S=\begin{pmatrix}
0 & r \\
-\lambda_{min} & r+\lambda_{min}
\end{pmatrix}$.
\end{claim}
\begin{proof}
Without loss of generality, we can assume that the graph $G$ is connected,
since otherwise the connected components could
be considered separately. The definition of an independent set
is equivalent to the equality $e(C,C)=0$. Then, from Corollary \ref{haemers}, we have the desired inequality $
\lambda_{min}|C|+\frac{(r-\lambda_{min})|C|^2}{|V(G)|}\leq 0.$

Furthermore, Corollary \ref{haemers} asserts that the equality
$\frac{|C|}{|V(G)|}=\frac{\lambda_{min}}{\lambda_{min}-r}$ holds
only if ${\mathbf 1}_{C}$ is a perfect coloring. The coefficient $a=s_{11}$ of  $S$ is equal to $0$ by the definition of independent set.
The entry $c=s_{21}$ of  $S$ can be determined by counting the edges joining vertices of different colors:
 $|C|r=c(|V(G)|-|C|)$.

Consider a perfect $2$-coloring with quotient matrix $S=\begin{pmatrix}
0 & r \\
-\lambda_{min} & r+\lambda_{min}
\end{pmatrix}$.
Since $s_{11}=0$, the first color of a perfect coloring with parameter matrix $S$ is an independent set. The number of vertices in the first
color is easily calculated from the relation $|C|r=-\lambda_{min}
(|V(G)|-|C|)$.
Thus, the cardinality of the first color class attains the Delsarte--Hoffman bound.
 \end{proof}

\begin{remark} Instead of the  adjacency matrix $M$, we can consider any other
nonnegative symmetric matrix $M'$ whose row sums are all equal to  $r'$ 
and whose support is contained in the support of $M$. Under this condition, for any indicator function $f$, $(Mf,f)=0$ implies that $(M'f,f)=0$.
The proof of the Delsarte--Hoffman theorem does not use the fact that $M$ is a  $(0,1)$-matrix, it uses  only the nonnegativity
and symmetry of  $M$. Therefore, every independent set in $G$ satisfies 
the  bound, where $\lambda'_{min}$
is the minimum eigenvalue of $M'$.
\end{remark}

\subsection{Delsarte Cliques in Transitive Graphs}

Recall that a clique in a graph
is a set of vertices that are pairwise adjacent.
It is easy to see that a clique in a graph corresponds to an independent set in its complement, and vice versa.
 For an $r$-regular connected graph $G$ on $n$ vertices,
it is easy to show that the minimum eigenvalue of its complement
is $-\lambda_2-1$, where $\lambda_2<r$ is the maximum
nontrivial eigenvalue of $G$, and $n-r-1$
is the degree of the complement (see Problem \ref{exe3}). Therefore, by the Delsarte--Hoffman bound,
 every clique in $G$
 has cardinality at most
$\frac{(\lambda_2+1)n}{\lambda_2+n-r}$. Another, usually stronger, bound can be obtained for the size of a clique in an arc-transitive graph.

\begin{theorem}[on the Delsarte clique \cite{Godsil}]\label{PCclaim13}
Every clique in an arc-transitive $r$-regular graph $G$ has cardinality at most
 $1-\frac{r}{\lambda_{min}}$, where $\lambda_{min}$ is the minimum
eigenvalue of $G$.
\end{theorem}
\begin{proof}
Let $u_1,\dots,u_m$ be an orthonormal basis of the eigensubspace $X_{\lambda_{min}}$. Consider some automorphism $\pi$
of $G$. Obviously, each automorphism of a graph maps one orthonormal basis of the eigenspace $X_{\lambda_{min}}$ to another. Let $Q_{\pi}$ denote the corresponding orthogonal matrix of size $m\times m$.

Consider the vectors
$w(x)=(u_1(x),\dots,u_m(x))$, where $x\in V(G)$. We will claim that the
scalar product $(w(x),w(y))$ is the same for every pair of adjacent vertices $x$ and
$y$. Indeed, consider the automorphism $\pi$ of the graph $G$,
mapping the edge $\{x,y\}$ to the edge $\{\pi(x),\pi(y)\}$. Then
$(u_1(x),\dots,u_m(x))Q_\pi=(u_1(\pi(x)),\dots,u_m(\pi(x)))$. We have
$(w(x),w(y))=(w(x)Q_\pi,w(y)Q_\pi)=(w(\pi(x)),w(\pi(y)))$. Let
$(w(x),w(y))=\beta$ for any adjacent vertices $x,y\in V(G)$.
Similarly, since $G$ is vertex-transitive, we obtain $(w(x),w(x))=\alpha$ for every
vertex $x\in V(G)$.

From the definition of  eigenfunctions, we have the equality $
\sum\limits_{y\in A_1(x)}u_i(y)=\lambda_{min} u_i(x),$ i.e., $
\sum\limits_{y\in A_1(x)}w(y)=\lambda_{min} w(x)$. Taking the inner product of the previous equality with $w(x)$, we obtain the equation
$r\beta=\lambda_{min}\alpha$.

Let $C\subset V(G)$ be a clique in graph $G$. Consider the Gram matrix
$W$ of size $|C|\times |C|$, whose elements are
$(w(x),w(y))$, $x,y\in C$.
%(since the corresponding quadratic form defines the scalar
%product of linear combinations of vectors $w(x)$).
As proven
above, $W=\begin{pmatrix}
\alpha & \beta & \beta & \dots & \beta\\
\beta & \alpha & \beta& \dots & \beta\\
\dots & \dots & \dots & \dots & \dots \\
\beta & \dots & \dots & \beta & \alpha
\end{pmatrix}$.

The Gram matrix is positive semidefinite, so
$(W\overline{1},\overline{1})\geq 0$, and therefore $\alpha+
(|C|-1)\beta\geq 0$. The minimum eigenvalue
$\lambda_{min}$ is negative (see Problem \ref{exe1}). Combining this with
 $r\beta=\lambda_{min}\alpha$, we obtain $|C|\leq
1-\frac{r}{\lambda_{min}}$.
\end{proof}

%\begin{definition} 
A clique of cardinality $1-\frac{r}{\lambda_{min}}$ in an
$r$-regular graph $G$ is called a {\sl  Delsarte clique}.
%\end{definition}

\begin{corollary}
Let $G$ be an arc-transitive graph which contains Delsarte cliques and
 an independent set $C$. If  $C$ attains the Delsarte--Hoffman bound,  then each Delsarte clique in $G$ contains exactly one vertex
from $C$, and each vertex from $C$ belongs to the same number
of Delsarte cliques.
\end{corollary}
\begin{proof}
Since the graph is transitive, each vertex belongs to the same
number of Delsarte cliques.
 We denote this number by $s$.  Then
the number of Delsarte cliques is equal to
$s|V(G)|/(1-\frac{r}{\lambda_{min}})$.
 Since $C$ is an independent set, each Delsarte clique contains at most one vertex of $C$.
 Then $|C|s\leq s|V(G)|/(1-\frac{r}{\lambda_{min}})$. Since $C$ attains the Delsarte--Hoffman bound, each clique contains a vertex from $C$.
\end{proof}

\subsection{Problems}

\begin{exercise}
Let $G$ be an $r$-regular graph and $f:V(G)\rightarrow\mathbb{C}$.
Prove the equality $((rI-M)f,f)=\sum\limits_{\{x,y\}\in
E(G)}|f(x)-f(y)|^2$.
\end{exercise}

\begin{exercise}
Let $f:V(G)\rightarrow\{0,1\}$ be a coloring
of an $r$-regular graph $G$. An edge of $G$ is called bichromatic if
its endpoints are colored differently. Prove the inequality 
$(r-\lambda)|f^{-1}(1)||f^{-1}(0)|\geq m|V(G)|$, where $m$ is the number
of bichromatic edges and $\lambda$
is the minimum eigenvalue
of the adjacency matrix of $G$. Prove that if equality is reached, then the coloring $f$ is perfect.
\end{exercise}

\begin{exercise}
{ Let $G$ be an $r$-regular graph and $\lambda$ be the minimum
eigenvalue of its
adjacency matrix $M$.  Suppose that every vertex of $A$ has at most $t<r$ neighbors in $A$. Prove that
$|A|\leq\frac{(t-\lambda)|V(G)|}{r-\lambda}$.  Moreover, if equality holds,  then ${\bf 1}_A$ is a perfect $2$-coloring of the graph $G$. }
\end{exercise}

\begin{exercise}
The chromatic number $\chi(G)$ of a graph $G$ is the minimum
number of colors in its proper coloring. Prove that $\chi(G)\geq
\frac{\lambda-r}{\lambda}$ for an $r$-regular graph $G$ with
minimum eigenvalue $\lambda$. Prove that if
$\chi(G)= \frac{\lambda-r}{\lambda}$, then every proper coloring
of the graph $G$ with $\chi(G)$ colors is perfect. \end{exercise}

\begin{exercise}
{ Let $G$ be an $r$-regular graph and $\lambda$ be the second largest eigenvalue  of its
adjacency matrix $M$. Prove that the cardinality of any clique in $G$ does not exceed
$\frac{\lambda+1}{1-(r-\lambda)/|V(G)|}$.}
\end{exercise}

\section{Lattice of Perfect Colorings}

\subsection{Merging of Colors}\label{SklCol}

For any set $V$,
 the collection of all partitions of 
$V$, ordered by refinement, is a {\sl lattice}. The maximum element of this lattice is one-element partition $\{V\}$, and the minimum element  is the partition consisting of all one-element subsets of $\{V\}$. 
For two partitions $B$ and $C$ of the same set~$V$,
the finest partition of~$V$ that is coarse  than both $B$ and $C$ are its refinements is denoted by $B \lor C$; the coarsest partition of~$V$ that   refines  both $B$ and $C$  is denoted by $B \wedge C$.  So,  $\lor$ is
the least upper bound  and  $\wedge$ is the greatest lower bound in the natural lattice of partitions.
A partially ordered set containing the least upper bound (or greatest lower bound) for any pair of its elements is called a {\sl semilattice}.

Let $G$ be a graph.  Partitions of $V(G)$ correspond to colorings of $G$. We will regard the elements of partitions as colors. Recall  that a partition of  $V(G)$ arising from a perfect coloring is called equitable. 
Consider a sublattice consisting of equitable partitions of $V(G)$. 

\begin{claim}\label{cl:mincolor}
Let $C^1$ and $C^2$ be equitable partitions of $V(G)$. Then $C^1 \lor C^2$ is an equitable partition of~$V(G)$.
\end{claim}
\begin{proof}
Denote by $F^i$ a matrix whose columns are indicator functions of elements of $C^i$, where $i=1, 2$.
Let $L_i$ be the subspace generated by the columns of the matrix
$F^i$. By Proposition \ref{cl:invar}, 
a partition $C$ is equitable if and only if the subspace generated by the indicator functions of its cells is invariant under the adjacency matrix of $G$.
Clearly, the intersection of invariant
subspaces $L_0= L_1\cap L_2$ is an invariant subspace.  It is easy to see that a basis of $L_0$ consists of indicator functions of  the cells of   $C^1 \lor C^2$.   Therefore,  $C^1 \lor C^2$ is an equitable partition of~$V(G)$ by 
Proposition \ref{cl:invar}.
\end{proof}

\begin{corollary} For every graph $G$ the sublattice of  equitable partitions is a semilattice.
\end{corollary}

In the lattice of partitions, the operation of combining several cells of a partition into one is called a  {\sl merging}. By merging, we obtain a coarse partition  from the original one.

Let $B$ be a partition, $|B|=k$, and let $F_B$ be a matrix whose columns are indicator functions of elements (colors) of $B$. 
 Consider a $(0,1)$-matrix $H=(h_{ij})$ of size $k\times m$, $m<k$, which
describe the merging the original $k$ colors into
$m$ new colors: $h_{ij}=1$ if and only if the $i$-th color of the original
coloring is merged into the $j$-th color in the new coloring. The matrix $H$ is called  a {\sl merging matrix}.

\begin{claim}\label{merge} Applying the merging matrix  $H$  of a perfect coloring with
quotient matrix $S$, we obtain a perfect coloring with quotient matrix
$S'$ if and only if $SH=HS'$.
\end{claim}
\begin{proof}
Let $F$ be the matrix of a perfect coloring with quotient matrix
$S$ of a graph with adjacency matrix $M$. Then $FH$ is the matrix
of the coloring obtained by merging colors. Then, from the equality $MF=FS$
and the condition $SH=HS'$, we have
the equalities $MFH=FSH=FHS'$. By the algebraic criterion of  perfect coloring (Proposition \ref{s:perfcol_criteria})
it follows that the merged coloring  is  perfect.

Conversely, from the equalities $MFH=FHS'$ and $MFH=FSH$, as well as the non-singularity of the matrix $F$, we obtain $SH=HS'$.
\end{proof}

The Cartesian product of perfectly colored graphs generates
 perfect colorings with a larger number of colors (see Proposition \ref{cartprod:1}). 
 Below, we consider sufficient conditions under which the number of colors of a perfect coloring can be reduced, that is, conditions under which certain colors can be merged.

%\begin{definition}
Two colors $i$ and $j$ of a perfect coloring are said to be
{\sl equivalent}  if the $i$-th and $j$-th rows in the quotient matrix
of the coloring coincide, with the possible exception of the $i$-th and $j$-th
coordinates in each row.
%\end{definition}

%\textcolor{green}{
\begin{claim} The equivalence relation of colors of a perfect coloring satisfies the
transitivity property. \end{claim}
\begin{proof}
Suppose
that the first color is equivalent to the second and third colors in a perfect
coloring with quotient matrix $S$. Then, from the definition
of equivalent colors, it follows that $s_{13}=s_{23}$ and
$s_{12}=s_{32}$. Let $n_i$
be the number of vertices of color $i$. Then we have the relations $n_1 s_{13}=n_3s_{31}$,
$n_2s_{21}=n_1 s_{12}$, and $n_3 s_{32}=n_2s_{23}$  (see equality
(\ref{eq:colorededges})). Multiplying the right and left sides of the equalities and applying the equalities $s_{13}=s_{23}$ and
$s_{12}=s_{32}$, we obtain $s_{21}= s_{31}$. This means the second and third
colors are equivalent.
\end{proof}
%}

\begin{claim}\label{cor:scleiv}
Consider an $r$-regular graph $G$.
Let the coloring $f'$ be obtained from the perfect coloring $f$ of $G$
by merging of a pair of equivalent colors. Then $f'$ is a perfect
coloring. \end{claim}
\begin{proof}
By the definition of equivalent colors, vertices of $i$-th and $j$-th colors have the same number of neighbors of other fixed color. 
Therefore, it remains only to show that every vertex in the merged color class has the same number of neighbors within this class.
 Let $S$ be a quotient matrix of the initial coloring. The equality  $s_{ii}+s_{ij}=s_{jj}+s_{ji}$ follows from the equation $\sum\limits_{m}s_{jm}=r$ for every $j$.   
\end{proof}

It is easy to see that the merger of any subset of an equivalence class is again equivalent to the remaining colors of that class. Consequently, any collection of equivalent colors in a perfect coloring can be merged, and the resulting coloring is again perfect.

\subsection{Weisfeiler--Leman Algorithm}\label{Vizing}

There is a natural algorithm that constructs the coarsest equitable refinement 
of a given partition of the vertex set of a graph. 
Thus, for any partition of the vertex set, this algorithm produces the coarsest equitable refinement of the given partition (see Theorem~\ref{th:Vizing}).
In the modern literature, this algorithm is often 
called the Weisfeiler--Leman (WL) algorithm (to be exact, the one-dimensional  Weisfeiler--Leman algorithm) \cite{WL:68}, although it was independently discovered many times, including earlier works, see, e.g.~\cite{Morgan:65}.
The algorithm is used in the software that 
recognizes the isomorphism of graphs, see~\cite{McKay:80}.

We describe the algorithm in terms of colorings.
 Suppose we have a graph colored with several colors.
Consider the vertices of a fixed color.
Partition them into groups consisting of vertices with the same neighborhood color composition.
Perform this operation simultaneously for all color classes. Color each group with its own (new)
color.  The
recoloring process stops if a perfect coloring is obtained,
i.e., no color is divided into groups. Otherwise, we repeat
the recoloring procedure.
 Note that no two vertices that already have different colors can be assigned the same color during the algorithm.
 Since at each step of the algorithm the number
of colors increases, but cannot exceed the number of vertices in the graph,
we will always end up with a perfect coloring in the final
graph.
It is possible that, in the end, all vertices receive distinct colors.
Obtaining a sharp upper bound on the number of iterations of the WL algorithm remains an open problem.
There are examples of graphs for which the number of steps of the WL algorithm is approximately equal to the number of vertices and many times greater than the graph diameter, even for a regular graph (if the initial coloring is not monochromatic). \label{WLprob}

\begin{theorem}[Vizing \cite{Vizing}]\label{th:Vizing}
Let $g:V(G)\rightarrow\{1,\ldots,k\}$ be a coloring of some graph $G$.
% $A_i=f^{-1}(i)$, $B=V(G)\setminus (A_0\cup\dots\cup A_s)$, $s<k-2$.
Let $f$ be a perfect coloring such that every color class of $g$ is obtained by merging color classes of $f$.
Then the perfect coloring obtained by
applying the WL algorithm to $g$ does not split the colors of $f$. In other words, the partition produced  by WL algorithm coarser than the partition corresponding to $f$. 
\end{theorem}
\begin{proof} We will prove the theorem by induction.
Assume that after the $i$-th step of the algorithm, every color class of the current coloring is obtained by merging color classes of $f$.
We will show that this holds true at the $(i+1)$-th step.
The induction base is valid by assumption.
Consider two vertices $x_1,x_2\in f^{-1}(j)$ of the same color in
coloring $f$.
Then they have the same neighborhood color composition  in coloring $f$, and therefore, by the induction hypothesis,
they will have the same neighborhood color composition in the current coloring as well.
Then the next step of the WL algorithm colors the vertices $x_1$ and
$x_2$ with the same color. \end{proof}

\begin{remark}
Essentially, the WL algorithm
finds an $M$-invariant subspace
closed with respect to the Hadamard (pointwise) product, where $M$ is the adjacency matrix.
Such subspaces are equivalent to equitable partitions (see Problem~\ref{1p:Pot04b}).
Initially, we have the subspace spanned by the characteristic functions of the colors of the starting coloring.
Each step acts on the subspace by~$A$
and then finds % the basis of 
the Hadamard-product closure of the resulting subspace. 
The algorithm stops when we get a subspace closed both under the action of~$M$ and under the Hadamard product.
\end{remark}

The WL algorithm can be used to check whether a set $A$ of vertices in a graph
is a color of some perfect coloring or not. It suffices to apply
the algorithm to a $2$-coloring, where the first  color is $A$ and the second color is  the remaining vertices. If the algorithm at some step
splits the vertices of the first color, then Vizing's theorem implies
that the set $A$ is not a color of any perfect
coloring.
Denote the perfect coloring obtained by the WL algorithm
by $f_A$. We have

\begin{corollary}
If $A$ is a color class of a perfect coloring $f$, then $f$ can be obtained from $f_A$ by splitting the remaining color classes.
\end{corollary}

\begin{remark}
Consider random $r$-regular graph with $N$ vertices, where $N-r_0\geq r\geq r_0$ for sufficiently large   $r_0$.
Then the following holds with high probability ($N\rightarrow \infty$): for every non-trivial initial partition  of the vertex set of $G$,
 WL algorithm runs at most $2{\rm diam}(G)+3$ steps on $G$ and outputs an $N$-colouring (see \cite{Isaev}).
 This means that a random regular graph with high probability does not
have an equitable partition, except those with $1$ part or $N$ parts.
\end{remark}

\subsection{Induced Colorings}\label{Vizing2}

Consider a bipartite graph. Suppose that one part of the graph is colored and the other part is left uncolored. The coloring obtained after one step of the WL algorithm is called an {\sl induced coloring}.

\begin{claim}[\cite{Lis}]
If a coloring of one part of a bipartite graph can be extended to a perfect coloring of the entire graph, then the induced coloring is also perfect.
\end{claim}

\begin{proof}
Suppose that a coloring of one part of a bipartite graph can be extended to a perfect coloring $g$ of the entire graph. If some color $i$ of $g$ appears in both parts of the graph, then the vertices of color $i$ in the two parts can be separated into distinct colors, and the resulting coloring remains perfect. Therefore, without loss of generality, we may assume that the two parts of the graph use disjoint sets of colors in the coloring $g$.

Consider the coloring $f$ that coincides with $g$ on the first part of the graph and is monochromatic on the second part. By Theorem~\ref{th:Vizing}, the colors of $f$ in the first part are not split at any step of the WL algorithm. Hence, the color compositions of the neighborhoods of vertices in the second part do not change during the algorithm. Consequently, the WL algorithm terminates after the first step and produces the induced coloring.
\end{proof}

Note that an induced coloring uses disjoint sets of colors on the two parts of the graph. However, some colors belonging to different parts may still be merged. Induced colorings can be used to reduce the search space when searching for perfect colorings of bipartite graphs.

\subsection{Coverings of Graphs}\label{sec_cov}

The most efficient method for transferring perfect
colorings from one graph to another is the covering.

%\begin{definition}
Let $G_1$ and $G_2$ be regular graphs of the same degree.
A {\sl covering} $\phi: V(G_1)\rightarrow V(G_2)$ is a mapping
 that induces a one-to-one correspondence  in each sphere.  In other words,  for all $ x\in V(G_1)$ 
 the restriction of $\phi$ to the neighborhood $A_1^1(x)$ is a bijection onto $A_1^2(\phi(x))$.
%\end{definition}

It follows immediately from the definition that the composition of two coverings is again a covering.
 Moreover, a composition of a covering with a perfect coloring is again a perfect coloring.

\begin{claim}\label{cover1}
Let $\phi: V(G_1)\rightarrow V(G_2)$ be a covering and
$f: V(G_2)\rightarrow \{1,\ldots,k\}$ be a perfect coloring of $G_2$.
Then $f\circ\phi: V(G_1)\rightarrow\{1,\ldots,k\}$ is a perfect coloring of $G_1$ with the same quotient matrix.
\end{claim}
\begin{proof}
Let the colors of the vertices $x_1$, $x_2\in V(G_1)$ be the same, i.e., $f\circ\phi(x_1)=f\circ\phi(x_2)$.
Then the vertices $\phi(x_1)$ and $\phi(x_2)$ also have the same color.
The  neighborhood color compositions $f(A_1^2(\phi(x_1)))$ and  $f(A_1^2(\phi(x_2)))$ of the vertices $\phi(x_1)$ and $\phi(x_2)$ 
are the same, since $f$ is a perfect coloring.
Because $\phi$ induces a bijection between the neighborhoods of a vertex and its image,
the  neighborhood color compositions of $x_1$ and $x_2$ in the coloring $f\circ\phi$ also coincide. 
\end{proof}

\begin{example}\label{exphi}
Define a covering $\phi:Q_2^3\rightarrow Q_4$ according to the rule: each pair
of opposite vertices of the Boolean cube corresponds to one
vertex of the complete graph $Q_4$.
\begin{equation*}
\begin{array}{c}
\phi^{-1}(0)=\{(000),\,(111)\},\\
\phi^{-1}(1)=\{(100),\,(011)\},\\
\phi^{-1}(2)=\{(010),\,(101)\},\\
\phi^{-1}(3)=\{(001),\,(110)\}.
\end{array}
\end{equation*}
As mentioned above (Proposition \ref{cl:complete}), every coloring of the complete graph
is perfect. The perfect $2$-colorings of the Boolean cube obtained by this
covering and $2$-colorings of the complete graph are in the figure below.
\end{example}

\begin{minipage}{0.24\textwidth}
%\vskip10mm
\center{\includegraphics[width=0.5\textwidth]{E3_color_0312.png}}
\end{minipage}
\begin{minipage}{0.24\textwidth}
$S=\begin{pmatrix}
0 & 3\\
1& 2
\end{pmatrix}$
\end{minipage}
\begin{minipage}{0.24\textwidth}
%\vskip10mm
\center{\includegraphics[width=0.5\textwidth]{E3_color_1221.png}}
\end{minipage}
\begin{minipage}{0.24\textwidth}
$S=\begin{pmatrix}
1&2\\
2&1
\end{pmatrix}$
\end{minipage}

From the definitions it follows that  a covering $\phi: V(G_1)\rightarrow
V(G_2)$ can be viewed as a perfect coloring of the vertices of $G_1$ by vertices of $G_2$ considered as colors. The converse is also true.
Let $f$ be a perfect coloring  of $G_1$ whose quotient matrix of the coloring is the adjacency matrix of  $G_2$.
Then the mapping $f:V(G_1)\rightarrow V(G_2)$
is a covering.

Since a covering can be viewed as a perfect coloring,
the covering $\phi$ corresponds to a rectangular matrix $\Phi$, satisfying $M_1\Phi=\Phi M_2$, where $M_1$ and $M_2$
are the adjacency matrices of graphs $G_1$ and $G_2$, respectively. Therefore, Proposition
\ref{cover1} is a special case of Proposition
\ref{merge} on color merging. Indeed, from the equalities $M_2F=FS$ and
$M_1\Phi=\Phi M_2$ it follows that $M_1\Phi F=\Phi FS$.  Moreover,
we obtain that the set of eigenvalues of the covered graph is a subset of the set of eigenvalues  of the covering graph (see Corollary  \ref{p:perf_color_properties}).
An isomorphism between two graphs is equivalent to a one-to-one covering. In
this case, the matrix $\Phi$ is a square permutation matrix, and we have the equalities $M_1\Phi=\Phi M_2$,
$M_2\Phi^{-1}=\Phi^{-1} M_1$. 
An automorphism is an isomorphism from a graph to itself. 

The examples above (perfect colorings, mergings, isomorphism and coverings) motivate us to consider the algebraic
equation $M_1F=FM_2$ in its own itself.
A matrix $M_2$  of size $k\times k$ is called a {\sl divisor} of matrix  $M_1$   of size $n\times n$, where $k\leq n$ if there exist   a $(0,1)$-matrix  $F$ containing exactly one $1$ in every row such that   $M_1F=FM_2$.
We define a natural partial order on the set of  nonnegative  integer symmetrizable matrices by setting $M_2\leq M_1$ whenever $M_2$ is a divisor of $M_1$. By the next theorem and its generalization,  this  partial order forms a semilattice.

\begin{theorem}[\cite{Angluin}]\label{Angl}
For any two $r$-regular graphs $G_1$ and $G_2$, there exists a graph
$G_0$ that can cover each of them.
\end{theorem}
\begin{proof}
First, we construct a covering $\phi$ of $G_1$ by the bipartite graph
$G'_1$. The vertex set of the bipartite graph $G'_1$ with parts $V$ and $V'$ consists of
twice the vertex set of $G_1$, i.e., $V(G'_1) =V \cup V'$
and there is a mapping $\phi_1:V \cup V'\rightarrow V(G_1)$,
one-to-one on $V$ and on $V'$.
Vertices $u\in V$ and $u'\in V'$ of $G'_1$ are connected by an edge
if and only
if the vertices $\phi_1(u)$ and $\phi_1(u')$ are adjacent in $G_1$.
It is easy to see that the mapping $\phi_1$ is a covering by
construction. Similarly, we define a bipartite graph $G'_2$ that
covers $G_2$ and a covering $\phi_2$.

Now we define a graph $G_0$ that covers $G'_1$ and
$G'_2$. By a corollary of König's theorem, which will be proved in Section \ref{7.2}, the edges of a bipartite $r$-regular
graph can be colored with $r$ colors such that each vertex is incident to exactly one edge of each  color. We color
the edges of $G'_1$ and $G'_2$ in this manner. Let
$V(G_0)=V(G'_1)\times V(G'_2)$. Vertices $(u_1,u_2)$ and $(v_1,v_2)$
are connected by an edge in graph $G_0$ if and only if vertices
$u_1$ and $v_1$ in graph $G'_1$ are connected by an edge of the same color
as vertices $u_2$ and $v_2$ in graph $G'_2$. By the construction,
the mapping $\psi_1:V(G_0)\rightarrow V(G'_1)$, defined
by the equality $\psi_1(u_1,u_2)=u_1$, is a covering of graph $G'_1$. Indeed, each edge color corresponds to unique vertex
$v_1$ adjacent to $u_1$ in the graph $G'_1$, and one vertex
$(v_1,v_2)$ adjacent to $(u_1,u_2)$ in the graph $G_0$. Similarly,
the mapping $\psi_2:V(G_0)\rightarrow V(G'_2)$, defined
by the equality $\psi_2(u_1,u_2)=u_2$, is a covering of the graph $G'_2$.
The compositions $\psi_1\circ \phi_1$ and $\psi_2\circ \phi_2$ are
coverings of the graphs $G_1$ and $G_2$, respectively, by the graph $G_0$.
\end{proof}

Theorem \ref{Angl} can be generalized to irregular graphs. Namely, two graphs have a common covering if and only if they admit a perfect coloring with the same quotient matrix
(see \cite{Lei}).

Note that the vertices of the graphs $G_1$ and $G_2$ have the same degrees,
but the number of vertices in the graphs can differ greatly. The graphs (or one
of them) can have an infinite number of vertices. Thus,
a covering of a "small" graph  $G_2$ by a "large" graph $G_1$ allows us to lift a perfect
coloring of   $G_2$  to a perfect coloring of $G_1$ with the same quotient matrix. 

From equality (\ref{eq:colorededges}) for the parameters of a perfect
coloring, it follows that the ratio of the number of vertices of  two colors
$\frac{n_i}{n_j}$ is equal to the ratio of the elements of the quotient matrix
$\frac{s_{ji}}{s_{ij}}$ if $s_{ij}\neq 0$. In the case of a covering,
the quotient  is the adjacency matrix
of the graph being covered, whose elements are either $0$ or $1$.
Therefore, under the covering map $\phi: V(G_1)\rightarrow V(G_2)$, the preimages of all vertices of $G_2$ have the same cardinality
$k=|V(G_1)|/|V(G_2)|$. Such a covering is called a {\sl $k$-covering}.
Note that a covering of the complete graph $Q_q$ by  a graph $G$ is both a proper and a perfect coloring of $G$ with $q$
colors.

\begin{claim}
If there exists a covering of the complete graph
$Q_{q+1}$ by a $q$-regular graph $G$, then $G$ contains a $1$-perfect code.
\end{claim}
\begin{proof}
Since every coloring of a complete graph is a perfect coloring,
there exists a perfect coloring of the graph $Q_{q+1}$ with parameter matrix
$S=\begin{pmatrix}
0 & q \\
1 & q-1
\end{pmatrix}$. By Proposition \ref{cover1}, there exists a perfect
coloring of the graph $G$ with the same parameter matrix. Then the first color
of this coloring is a $1$-perfect code by definition.
\end{proof}

\subsection{Problems}

\begin{exercise}
Let $P$ be a equitable partition of a bipartite connected
graph.
Prove that if  at least one cell of the partition is contained entirely in one part,
then this is true for vertices of every color. 
\end{exercise}

\begin{exercise}
Let $P$ be a equitable partition of a bipartite 
graph. Prove that  if  one cell $P_i\in P$  is contained  in both parts of the graph then this cell can be split into two cells such that the partition
remains equitable. 
 \end{exercise}

\begin{exercise}
Suppose there is a coloring of the graph $G$ with colors from the sets $A$ and
$B$ such that merging all the colors from $A$ into a single color while preserving
all the colors from $B$, and conversely, merging all the colors from
$B$ into a single color while preserving all the colors from $A$, both result in perfect
colorings. Prove that the original coloring is perfect.
\end{exercise}

\begin{exercise}
Let $\phi: V(G_1)\rightarrow V(G_2)$ be a covering. Prove that
the number of common neighbors of vertices $u,v\in V(G_1)$ is no greater than the number
of common neighbors of vertices $\phi(u),\phi(v)\in V(G_2)$. Deduce that
the size of a maximum clique in the covering graph is no greater than in
the covered graph.
\end{exercise}

\begin{exercise}
Suppose that $G_1$ covers $G_2$. Prove that if
$G_2$ has a perfect matching (see Section \ref{7.3}),
then the graph $G_1$ also has a perfect matching.
\end{exercise}

\begin{exercise}
Prove that an infinite $r$-regular tree  covers every
$r$-regular graph. \end{exercise}

\begin{exercise}
Give an example of a graph and an initial partition whose color-indicator functions  
lie in an   invariant subspace of the graph, but for which the equitable partition obtained by the WL algorithm lie in an invariant subspace with a new eigenvalue.
\end{exercise}

\section{Cartesian Products of Graphs}\label{Cart}

\subsection{Eigenfunctions of the Cartesian Product of Graphs}

The {\sl Kronecker product} of matrices $A$ of size
$n_1\times m_1$ and  $B$ of size $n_2\times m_2$ is a matrix
$A\otimes B$ of size $n_1n_2\times m_1m_2$, consisting of the elements\\
$(A\otimes B)_{(i_1, i_2),(j_1, j_2)}=A_{i_1,j_1}B_{i_2,j_2}$.
For  matrices of compatible  sizes, we have  $(A\otimes
B)\cdot (C\otimes D)=(AC)\otimes(BD)$ and $(A\otimes B)+(A\otimes
D)=A\otimes (B+D)$.

%\begin{definition}
The {\sl Cartesian product} of two graphs $G_1=(V_1,E_1)$ and
$G_2=(V_2,E_2)$ is a graph $G=G_1\square G_2$ such that
$V(G)=V_1\times V_2$, 
$E(G)=\Bigl\{\{(v_1,v_2),(u_1,u_2)\} :
(v_1=u_1 \text{ and } \{v_2,u_2\}\in E_2)
\text{ or }
(v_2=u_2 \text{ and } \{v_1,u_1\}\in E_1)
\Bigr\}.$
%\end{definition}

Let $|V_1|=n$ and $|V_2|=m$. Then the adjacency matrix $M$ of graph $G$
can be expressed in terms of the Kronecker products of the adjacency matrices $M_1$
and $M_2$ of graphs $G_1$ and $G_2$ and identity matrices. With an appropriate
numbering of the vertices of  $G$, the equality $M=M_1 \otimes
I_{m} + I_{n} \otimes M_2$ holds. 
Thus, the matrix $M$ is obtained from $M_1$ by replacing each
entry $1$ with the block $I_m$, each entry $0$ with the zero
matrix of order $m$, and then adding the block-diagonal matrix
with copies of $M_2$ on the main diagonal.

In particular, the $n$-dimensional $q$-ary cube $Q^n_q$ is defined as the
Cartesian product of cubes of lower dimension, namely
$Q^n_q=Q^{n-1}_q\square Q_q$. We obtain the following
recurrence representation of its adjacency matrix
\[
M_{n}=
\begin{pmatrix}
M_{n-1} & I_{q^{n-1}}&\cdots & I_{q^{n-1}} \\
I_{q^{n-1}} & M_{n-1}&\cdots & I_{q^{n-1}} \\
\cdots & \cdots&\cdots &\cdots\\
I_{q^{n-1}} &  I_{q^{n-1}} &\cdots & M_{n-1}
\end{pmatrix}
\]

The definition of the Cartesian product of graphs implies that the subgraph
of  $G=G_1\square G_2$ induced by the set of vertices
$V(G_1)\times \{y\}$ is isomorphic to  $G_1$ for any vertex $y\in
V(G_2)$. Similarly, the subgraph induced by $\{x\}\times V(G_2)$ is isomorphic to $G_2$.

From the definition of the Cartesian product of graphs, for any function
$f(x,y): V(G) \rightarrow \mathbb C$, we have the equality $$M f(x,y) =
\sum\limits_{x' \in A_1^1(x)} f(x',y) +
\sum\limits_{y' \in
A_1^2(y)}f(x,y'),$$
where $A_1^1(x)$ and $A_1^2(y)$ are the unit spheres in  $G_1$ and
$G_2$, respectively.

\begin{claim}\label{c:direct_product_eigen}
Let $f_1,\ldots,f_n$ be an orthogonal basis of eigenfunctions of graph $G_1$ with eigenvalues
$\alpha_1,\ldots,\alpha_n$; and $g_1,\ldots,g_m$ be an orthogonal basis of eigenfunctions of graph $G_2$ with eigenvalues
$\beta_1,\ldots,\beta_m$. Then $\{f_i\otimes g_j\}$ is an orthogonal basis of eigenfunctions of graph $G=G_1\square G_2$
with eigenvalues $\alpha_i+\beta_j$. \end{claim}
\begin{proof}

$M(f_i\otimes g_j) = (M_1 \otimes I_{m})\cdot(f_i\otimes g_j) +
(I_{n} \otimes M_2)\cdot(f_i\otimes g_j)$ \\ $=M_1f_i\otimes
g_j+f_i\otimes M_2g_j=\alpha_i f_i\otimes g_j+\beta_j f_i\otimes
g_j=(\alpha_i+\beta_j)f_i\otimes g_j.$

Now let's check the orthogonality of the basis. The well-known equality
\[
(f_{i_1}\otimes g_{j_1},
 f_{i_2}\otimes g_{j_2})
=
(f_{i_1},f_{i_2})
(g_{j_1},g_{j_2}),
\]
 immediately implies orthogonality.
The products $f_{i}\otimes g_{j}$ form a basis, since
the number of these functions is $|V(G)|=nm$.
\end{proof}

Consider the graph $Q_q$, isomorphic to a clique on $q$ vertices. A basis of  eigenfunctions of $Q_q$ is given by
$
\Phi_0=(1,1,1,\ldots,1)\mbox{  with eigenvalue } q-1,
$
and for any $i,j\in\{1,\dots,q\}$, $i\neq j$,
$
\Phi_{ij}=(0,\dots,\underset i{1},0,\dots,0,\underset
j{-1},0,\ldots,0)
$
 with eigenvalue $-1$.
Among the listed functions, we can choose a basis for the entire space,
therefore, the graph $Q_q$ has no other eigenvalues.

\begin{corollary}\label{cl:spec}
The eigenvalues of $Q^n_q$ are $(q-1)n-qk$, where $k=0,\dots,n$.
\end{corollary} 
\begin{proof}
By the definition  $Q^n_q=Q^{n-1}_q\square Q_q$.  The eigenvalues of  $Q_q$ are $q-1$ and $-1$.
Suppose that we know the eigenvalues of $Q^{n-1}_q$. By Proposition \ref{c:direct_product_eigen},   every eigenvalue of $Q^n_q$ is equal to $(q-1)(n-1)-qk+q-1=(q-1)n-qk$ or  $(q-1)(n-1)-qk-1=(q-1)n-q(k+1)$, where $k=0,\dots,n-1$.
\end{proof}

Next we will show how, from the eigenfunctions of the Cartesian product of graphs and
one of the factors, we can obtain the eigenfunction of the second
factor.

%\subsection*{Lemma on $\phi$-Projection}
\begin{lemma}[contraction theorem  \cite{Taranenko}]\label{l:phi_projection}
Let $G=G_1\square G_2$; $f(x,y)$ is an eigenfunction of $G$ with
eigenvalue $\lambda$; $\phi_0$ is an eigenfunction of 
$G_2$ with eigenvalue $\alpha_0$. Then
$c_0(x)=\sum\limits_{y\in V(G_2)} f(x,y)\cdot \phi_0(y)$ is an eigenfunction of $G_1$ with eigenvalue $(\lambda -
\alpha_0)$.
%Here $\alpha_0$ is an arbitrary eigenvalue.
\end{lemma}
\begin{proof}
It is known that any vector can be included in some orthogonal
basis containing this vector. Therefore, without loss of generality,
we can assume that $\phi_0,\ldots,\phi_{m-1}$
form an orthonormal basis of eigenfunctions of  $G_2$ with
eigenvalues 
$\alpha_0,\ldots,\alpha_{m-1}$.  For
any fixed $x\in V(G_1)$,   the function $f$ can be expanded in the basis:
$f(x,y)=\sum\limits_{i=0}^{m-1} c_i(x) \phi_i(y)$, since by
 fixing $x$ we obtain a graph
isomorphic to $G_2$. Thus we have
$f=\sum\limits_{i=0}^{m-1}c_i\otimes\phi_i$ 
and
$\sum\limits\limits_{y\in V(G_2)}f(x,y) \phi_i(y)=c_i(x)$ due to
orthonormality of the basis.

Since $f$ is an eigenfunction of $G$, it holds the equality
$$(M_1 \otimes
I_{m} + I_{n} \otimes
M_2)\left(\sum\limits_{i=0}^{m-1}c_i\otimes\phi_i\right)=\lambda\left(\sum\limits_{i=0}^{m-1}c_i\otimes\phi_i\right).$$

Since $\phi_i$, $i=0,\dots, m-1$, are eigenfunctions of $G_2$, we obtain a new equality

\begin{equation}\label{eq:phi_proj}
\sum\limits_{i=0}^{m-1}(M_1c_i+\alpha_ic_i)\otimes\phi_i=\sum\limits_{i=0}^{m-1}\lambda
c_i\otimes\phi_i.
\end{equation}

It is easy to see that scalar product with one of the factors
in the Cartesian product yields the following equalities:
$$g(x)(h_1,h_2)=\sum\limits_{{y\in V(G_2)}}g(x)h_1(y)h_2(y)=
\sum\limits_{{y\in V(G_2)}}(g\otimes h_1)(x,y)\cdot h_2(y).$$
Therefore, scalar product of both sides of the equality (\ref{eq:phi_proj})
by the function $\phi_0$ yields the equality $M_1c_0+\alpha_0c_0=\lambda
c_0$.
\end{proof}

Lemma \ref{l:phi_projection} is used for construction of eigenfunctions of hypercubes in \cite{ValVor}.
Consider
some eigenfunction $f$ with eigenvalue $\lambda$ in
 $Q^{n+1}_q=Q^n_q\square Q_q$. By   Lemma  \ref{l:phi_projection} the sum of the  $f$
 over parallel  hyperfaces
$g(x_1,\dots,x_n)=\sum\limits_{z\in Q_q}f(x_1,\dots,x_n,z)$ 
is an eigenfunction of $Q^n_q$ with eigenvalue
$\lambda - q+1$.  Similarly, for two hyperfaces
of the same direction, the difference
$g_{ij}(x_1,\dots,x_n)=f(x_1,\dots,x_n,i) - f(x_1,\dots,x_n,j)$ 
is an eigenfunction of  $Q^n_q$ with eigenvalue
$\lambda + 1$. 

\subsection{Perfect Colorings of Cartesian Products of Graphs}

Let us consider methods for constructing perfect colorings of Cartesian
products of graphs by perfect colorings of factors.

\begin{claim}[\cite{Taranenko}]\label{cartprod:1}
Let $f_i$ be a perfect $k_i$-coloring of graph $G_i$ with quotient matrix $S_i$, $i=1,2$. Then $f_1\times f_2$ is a $(k_1k_2)$-perfect coloring of  $G=G_1\square G_2$ 
with quotient matrix $S=S_1 \otimes I_{k_2} + I_{k_1} \otimes
S_2$.
\end{claim}
\begin{proof}
Consider the neighborhood of the vertex $(x,y)$. It consists of vertices $(x,y')$,
$y'\in S^2_1(y)$ of colors $(f_1(x),f_2(y'))$ and vertices $(x',y)$,
$x'\in S^1_1(x)$ of colors $(f_1(x'),f_2(y))$. Moreover, the multisets of colors
$f_2(y')$ and $f_1(x')$ depend only on the colors $f_2(y)$ and $f_1(x)$, respectively,
and therefore only on the color of  $(x,y)$. Thus, $f_1\times f_2$ is a perfect coloring.

Below we calculate the quotient  matrix of the resulting coloring and independently
verify that this coloring is perfect. As noted earlier, the adjacency matrix of graph $G_1\square G_2$ has the form $M=M_1 \otimes
I_{m_2} + I_{m_1} \otimes M_2$, where $M_i$
is the adjacency matrix of graph $G_i$, $m_i=|V(G_i)|$, $i=1,2$. Let
$F_1$
be the matrix of size $m_1\times k_1$ corresponding to the perfect
coloring $f_1$ of graph $G_1$, and $F_2$ be the matrix of size $m_2\times k_2$ corresponding to the perfect coloring $f_2$ of graph $G_2$. It is easy to see that
the matrix $F_1\otimes F_2$  corresponds to the coloring $f_1\times
f_2$. From Proposition
\ref{s:perfcol_criteria}, we have $M_iF_i=F_iS_i$, where $i=1,2$. Using
the following equality $(A_1\otimes B_1)(A_2\otimes B_2)=
(A_1A_2)\otimes (B_1B_2)$ for Kronecker products of matrices
of suitable size, we obtain the equalities
$$(M_1 \otimes I_{m_2} + I_{m_1} \otimes M_2)(F_1\otimes F_2)=
(F_1S_1)\otimes F_2+ F_1\otimes (F_2S_2)=(F_1\otimes F_2)(S_1
\otimes I_{k_2} + I_{k_1} \otimes S_2).$$ Thus, the coloring
$f_1\times f_2$ is perfect by the algebraic criterion
for a perfect coloring (Proposition \ref{s:perfcol_criteria}).
\end{proof}

Consider the special case when $f_2$ is a trivial one-coloring. Let $f_1$
be a perfect coloring of the graph $G_1$ with quotient matrix $S$.
Then $f(x,y)=f_1(x)$ is a perfect coloring of the graph $G_1\square
G_2$ with quotient matrix $S+rI$, where $r$ is the degree of $G_2$.

\begin{claim}\label{cartprod:2}
Let $g:V(G)\rightarrow\{0,1\}$ be a perfect coloring of a graph $G$
with quotient matrix $\begin{pmatrix}
a & b \\
b & a
\end{pmatrix}$ and $h:V(H)\rightarrow\{0,1\}$ be a perfect coloring of graph $H$ with quotient matrix
$\begin{pmatrix}
c & d \\
d & c
\end{pmatrix}$. Then the coloring $f(x,y)=g(x)\oplus h(y)$,
where $\oplus$ is addition modulo $2$, is perfect with
the quotient matrix $\begin{pmatrix}
a+c & b+d \\
b+d & a+c
\end{pmatrix}$ in the graph $G\square H$.
\end{claim}
\begin{proof}
Consider the number of neighbors of color $1$ of a vertex $(x,y)$. If
$g(x)=h(y)=0$, then there are $b$ neighbors of $(g(x'),h(y))$ and $d$ neighbors
$(g(x),h(y'))$ of color $1$. If $g(x)=h(y)=1$, then there are $b$
neighbors of $(g(x'),h(y))$ and $d$ neighbors of $(g(x),h(y'))$ of color $1$.
The remaining cases are analogous.
\end{proof}

\begin{claim}\label{claim:cartprodcover}
If the mappings $\phi_i: V(G_i)\rightarrow V(H_i)$, $i=1,2$, are coverings,
then the mapping $\bar\phi:V(G_1\square G_2)\rightarrow V(H_1\square H_2)$,
where $\bar\phi(x,y)=(\phi_1(x),\phi_2(y))$ is also a covering.
\end{claim}
\begin{proof}
We have $M_{i}\Phi_i=\Phi_i N_{i}$, where the matrix $\Phi_i$ corresponds
to the covering $\phi_i$, $M_{i}$ is the adjacency matrix of $G_i$, and
$N_{i}$ is the adjacency matrix of  $H_i$, $i=1,2$. By definitions,
the mapping $\bar\phi$ has matrix $\Phi_1\otimes\Phi_2$. The matrix
$M_{1} \otimes I_{m_{2}} + I_{m_{1}} \otimes M_{2}$ is the adjacency matrix of the Cartesian product of graphs $G_1$ and $G_2$ with
adjacency matrices $M_{1}$ and $M_{2}$, and the matrix $N_{1} \otimes
I_{n_{2}} + I_{n_{1}} \otimes N_{2}$ is the adjacency matrix of the Cartesian product of graphs $H_1$ and $H_2$ with the adjacency matrices
$N_{1}$ and $N_{2}$ . Therefore we have $$(M_{1} \otimes I_{m_{2}} +
I_{m_{1}} \otimes M_{2})(\Phi_1\otimes\Phi_2)=M_{1}\Phi_1 \otimes
\Phi_2 + $$ $$+ \Phi_1 \otimes M_{2}\Phi_2=
(\Phi_1\otimes\Phi_2)(N_{1} \otimes I_{n_{2}} + I_{m_{1}} \otimes
N_{2}), $$ which,  by the algebraic criterion of perfect coloring (Proposition \ref{s:perfcol_criteria}), proves that $\bar\phi$ is a covering.
\end{proof}

\begin{example}\label{exBesp} %(Bespalov)
The mapping $\bar\phi_n :({Q}^{3}_2)^n
\rightarrow {Q}^n_4$, $\bar\phi_n(x_1,\dots,x_n)=(\phi(x_1),\dots,\phi(x_n))$, where
$\phi$ is defined in Example \ref{exphi}, is a covering.
\end{example}

\subsection{$k$-Fold Coverings and Examples}

%\begin{definition}
Let $G_1$ be an $mr$-regular graph and $G_2$ an $r$-regular graph.
An {\sl $m$-fold covering} $\phi: V(G_1)\rightarrow V(G_2)$ is a mapping
such that for all $ x\in V(G_1)$ we have $\phi(S^1_1(x))=S^2_1(\phi(x))$
and  for all $ y\in S^2_1(\phi(x))$ we have  $|\phi^{-1}(y)\cap S^1_1(x)|=m$.
%\end{definition}

\begin{claim}
Let $\phi: V(G_1)\rightarrow V(G_2)$ be an $m$-fold covering and let
$f: V(G_2)\rightarrow \{1,\ldots,k\}$ be a perfect coloring of $G_2$ with quotient matrix
$S$.
Then $f\circ\phi: V(G_1)\rightarrow\{1,\ldots,k\}$ is a perfect coloring of $G_1$ with quotient matrix
$mS$.
\end{claim}
\begin{proof}
From the definition of an $m$-fold covering,  $\phi$ correspond to  
 $(0,1)$-matrix $\Phi$ such that  $M_1\Phi=\Phi(mM_2)$, where $M_1$ and $M_2$
are the adjacency matrices of the graphs $G_1$ and $G_2$, respectively. By the algebraic criterion of perfect coloring (Proposition \ref{s:perfcol_criteria}) we have $M_2F=FS$.
Therefore, $M_1\Phi F=\Phi(mM_2)F=\Phi F(mS)$.
\end{proof}

\begin{example}\label{ex:pch}
Let $pch:Q^k_2\rightarrow Q_2$ be the
 perfect coloring with parameter matrix $
\begin{pmatrix}
0 & k \\
k & 0
\end{pmatrix}$, i.e., $pch(x)=1$,
if $x$ contains an odd number of ones, and $pch(x)=0$ if $x$
contains an even number of ones. Then the mapping
$p_m^k:Q^{km}_2\rightarrow Q_2^m$, defined by the equality
$p_m^k(x_1,\dots,x_m)=(pch(x_1),\dots,pch(x_m))$, $x_i\in Q^k_2$,
is a $k$-fold covering.
\end{example}

\begin{example}\label{ex:covercube}
Let  $*$  be a group  operation  on the set $Q_q$. Consider the mapping $\ell:Q^n_q\rightarrow Q_q$ defined by the equality $\ell(x_1,\dots,x_n)=x_1* x_2*\cdots* x_n$. It is easy to show that $\ell$
is an $n$-fold covering of  $Q_q$. Let $id$ be a coloring of the graph $Q_q$ with $q$ different colors (with quotient  matrix 
$J_q-I_q$), then quotient matrix  of $id\circ \ell$ is
$n(J_q-I_q)=
\begin{pmatrix}
0 & n & \ldots & n\\
n & 0 & \ldots & n\\
\vdots & \vdots & \ddots & \vdots\\
n & \ldots & \ldots & 0
\end{pmatrix}$. 
Such a coloring is called {\sl latin coloring}, and a hypercube with such a coloring
is called {\sl latin hypercube}. \end{example}

\begin{example}\label{MDSexam}
Consider a $2$-coloring $h$ of $Q_q$: one vertex receives the first color, the rest are colored by the second one. The coloring $h\circ \ell$, where
the mapping $\ell$ is defined in Example \ref{ex:covercube}, is a
perfect coloring of the graph $Q^n_q$ with quotient matrix
$\begin{pmatrix}
0 & n(q-1) \\
n & n(q-2)
\end{pmatrix}$. 
Thus, the vertices of the first color form an
MDS code (see Section \ref{s:MDS}, and the non-code vertices can be viewed as the union
of the remaining colors of the latin coloring $h\circ \ell= h\circ (id \circ
\ell)$.
\end{example}

\begin{example}
Consider the covering $\bar\phi_n :{Q}^{3n}_2
\rightarrow {Q}^n_4$ from Example \ref{exBesp} and an MDS code in the hypercube ${Q}^n_4$, i.e., the coloring $f$
with parameter matrix $S=\begin{pmatrix}
0 & 3n \\
n & 2n
\end{pmatrix}$. By Proposition \ref{cover1}, it follows that $f\circ \bar\phi_n$ is a perfect coloring
 of  ${Q}^{3n}_2$ with quotient matrix $S$.
A perfect $2$-coloring with the same quotient matrix can also be defined as
$\phi\circ p_3^n$, where the mapping $p_3^n$ is defined in Example
\ref{ex:pch}, and the mapping $\phi$ is defined in Example \ref{exphi}.
Perfect colorings with this quotient matrix attains the Fon-Der-Flaass bound
on correlation immunity (see Section \ref{11.6}).
\end{example}

Further constructions of perfect colorings of hypercubes can be found in
\cite{Besp}.

\subsection{Relevant Variables and Wegener's Theorem}

%\begin{definition}
Let $T_0,T_1, \dots, T_n$ be arbitrary nonempty sets. Let $f:T_1\times\cdots\times T_n\rightarrow T_0$ a function. A variable $x_i$, $1\leq i\leq
n$, is called {\sl relevant} (essential, or effective) if there
exist $a_j \in T_j$, $j\neq i$ and $b, c\in T_i$
such that
$$f(a_1,\dots,a_{i-1},b, a_{i+1},\dots, a_n)\neq f(a_1,\dots,a_{i-1},c, a_{i+1},\dots,
a_n).$$
%\end{definition}

 Otherwise, the variable is said to be irrelevant. The coloring of a Cartesian
product of graphs $G=G_1\square\cdots\square G_n$ can be considered as a function of $n$ variables. 
Thus, it is natural to consider the relevant variables of such colorings.

A coloring of a Boolean $n$-cube $f:Q_2^n\rightarrow
\{1,\dots,k\}$ can be viewed as a function in $n$
Boolean variables. We denote by $\eta(f)$ the number of relevant
variables of $f$. For each vertex $x\in Q_2^n$, let $c(x)$ denote the number of neighbors
 whose color differs from the color of $x$. Denote $c(f)=\max_x
c(x)$.

\begin{theorem}[Wegener \cite{Wegener}]
For any coloring $f$ of a Boolean $n$-cube, the inequality
$\eta(f)\leq c(f)2^{2c(f)-2}$ holds.
\end{theorem}
\begin{proof}
First, we prove the following statement by induction. If every vertex of a nonempty induced subgraph of the Boolean $n$-cube has degree at least $r\leq n$, then the subgraph contains at least  $2^r$ vertices. 
   For $n=1$, the statement is obvious. Let us prove the induction step. Assume the statement is true for the $n$-dimensional cube. Consider a vertex set $T$ in $Q^{n+1}_2$ such that
every vertex of the induced subgraph on $T$ has  degree  at least $r$. Choose an arbitrary edge of the subgraph induced by  $T$. Without
loss of generality, we assume that this is an edge in the $n$-th direction.
Consequently, the sets $T_0=Q_2^n\square \{0\}\cap T$ and
$T_1=Q_2^n\square \{1\}\cap T$ are nonempty. Subgraphs induced by
 $T_0$ and $T_1$ have vertex degrees at least $r-1$,
since any vertex in $T_0$ is adjacent to at most one
vertex in $T_1$ and vice versa. Then, by the induction hypothesis, we have
$|T|=|T_0|+|T_1|\geq 2\cdot 2^{r-1}$.

We proceed to the proof of the theorem. Without loss of generality, we can
assume that the $n$-th variable of $f$ is relevant.
 Then there exists a nonmonochromatic  edge  $\{x,y\}$
of the $n$-th direction. The vertices incident with this edge have different colors: $f(x)=\alpha$ in 
$Q_2^{n-1}\square \{0\}$ and $f(y)=\beta\neq \alpha$ in 
$Q_2^{n-1}\square \{1\}$. Let us estimate the number of nonmonochromatic edges in the $n$-th direction adjacent to
$\{x,y\}$. At most
$c(f)-1$ vertices in $Q_2^n\square \{0\}$ adjacent to vertex $x$
have a color different from $\alpha$. Similarly,
at most $c(f)-1$ vertices in $Q_2^n\square
\{1\}$ adjacent to vertex $y$ have a color different from $\beta$. Vertices adjacent to vertices $x$ and $y$ are pairwise connected by edges in the $n$-th direction. 
Therefore, among the $n-1$ edges parallel to $\{x,y\}$, at least
 $n-1-2(c(f)-1)$
 are nonmonochromatic. Consider the subgraph of the graph
$Q_2^{n-1}\square \{0\}$ consisting of the vertices that are incident to
nonmonochromatic edges in the $n$-th direction. As we showed above, the degree
of the vertices of this subgraph is at least $n-2c(f)+1$. Consequently, the number of
 nonmonochromatic edges in the $n$-th direction and any other essential
direction is at least $2^{n-2c(f)+1}$.

By the definition of $c(f)$, there are at most
$c(f)2^{n-1}$ nonmonochromatic edges in the Boolean $n$-cube. Therefore
$\eta(f)2^{n-2c(f)+1}\leq c(f)2^{n-1}$.
\end{proof}

\begin{corollary}
For any $c,b\in \mathbb{N}$, there exists $n_0\in \mathbb{N}$ such that for $n> n_0$, any perfect coloring of the Boolean $n$-cube with
quotient matrix $\begin{pmatrix}
n-b & b \\
c & n-c
\end{pmatrix}$ has irrelevant variables. \end{corollary}
\begin{proof}
By the definition of the perfect coloring, $c(f)=\max\{c,b\}$.
Thus, by Wegener's theorem, the number of relevant variables
in the coloring does not exceed the constant
$n_0=\max\{c,b\}2^{2\max\{c,b\}-2}$.
\end{proof}

Any function with an irrelevant variable is a trivial extension of a function of fewer variables.
 It is easy to see that if a perfect coloring of a Boolean
hypercube with quotient matrix $\begin{pmatrix}
n-b & b \\
c & n-c
\end{pmatrix}$ has an irrelevant variable, then it is
 the Cartesian product of a monochromatic edge and a perfect
coloring with quotient matrix $\begin{pmatrix}
n-b-1 & b \\
c & n-c-1
\end{pmatrix}$. 
Therefore, a constructive description of perfect $2$-colorings of Boolean $n$-cubes requires determining the maximum dimension $n$ for which there exists a perfect coloring with quotient matrix
$\begin{pmatrix}
n-b & b \\
c & n-c
\end{pmatrix}$  without irrelevant variables.

\subsection{Problems}

\begin{exercise}
Suppose there exists a perfect coloring of $Q_2^k$ with
quotient matrix $S$. Prove that there exists a perfect coloring of $Q_2^{mk}$ with quotient matrix $mS$.
\end{exercise}

\begin{exercise}
Suppose there exists a perfect coloring of  $Q_q^k$ with quotient matrix
$S$. Prove that there exists a coloring of
$Q_{pq}^{k}$ with quotient matrix $pS+k(p-1)I$.
\end{exercise}

\begin{exercise}
Find the eigenvalues of the graph $Q^m_2\square Q^k_3$.
\end{exercise}

\begin{exercise}
Formulate and prove an analogue of Wegener's Theorem for the ternary hypercube.
\end{exercise}

\section{Amply Regular Graphs}\label{3.3}

\subsection{Average Degree of Vertices}

Let $M$ be the adjacency matrix of some graph
$G$. Then $M{\bf1}$
is the vector of vertex degrees of $G$, $n=({\bf1},{\bf1})$ is the number of vertices of $G$
and $(M{\bf1},{\bf1})/n$ is the average  degree of vertices of $G$.  The $i$-th coordinate of 
$M(M{\bf1})=M^2{\bf1}$ is equal to the sum  of the  neighbors
of the $i$-th vertex (if a vertex is adjacent  to $k$ neighbors, then it
is counted $k$ times). 
Equivalently, it is the number of walks of length $2$ ending at the $i$-th vertex.
 Thus, $(M^2{\bf1},{\bf1})/n$ is
the average sum of the degrees of the neighbors of the vertices in  $G$.

\begin{claim}\label{KBclaim}
For any graph, the square of the average degree of  vertices does not exceed
the average of the sums of the degrees of the vertex's neighbors. Moreover, if these are equal, then
the graph is regular. \end{claim}
\begin{proof}
From the Cauchy--Schwarz inequality we have
$$(M{\bf1},{\bf1})^2\leq (M{\bf1},M{\bf1})({\bf1},{\bf1}),$$ hence
$$\left(\frac{(M{\bf1},{\bf1})}{({\bf1},{\bf1})}\right)^2\leq
\frac{(M^2{\bf1},{\bf1})}{({\bf1},{\bf1})}.$$ Moreover, equality
is achieved only if the vectors $M{\bf1}$ and ${\bf1}$ are collinear.
It is equivalent to the regularity of the graph.
\end{proof}

Consider a $2$-coloring $f:V(G)\rightarrow \{0,1\}$ of a regular graph
$G$. Let $G_i$ be the subgraph induced by the vertices
of color $i\in \{0,1\}$. For each of these subgraphs, Proposition \ref{KBclaim} holds. Define $g=f\cdot Mf$ and
$g'=f'\cdot Mf'$, where $f'={\bf 1}-f$ and  $\cdot$ denotes
 entrywise (Hadamard) product. From the Cauchy--Schwarz inequality, we have $(g,f)^2\leq
(g,g)(f,f)$ and $(g',f')^2\leq (g',g')(f',f')$. Let us show that equality in both inequalities is simultaneously reached if and only if the coloring of $f$ is perfect.

Remove all vertices of color $0$ from the graph. Then, for the induced
subgraph $G_1$, the functions $f$ and $g$ play the roles of ${\bf1}$ and
$M{\bf1}$. In the case of equality, all vertices of color $1$ have the same
number of neighbors of color $1$, and hence, by the regularity of  $G$,
the same number of neighbors of color $0$. The same  arguments apply to the
characteristic function $f'$ of color $0$. Conversely, if
the coloring of $f$ is perfect, then the vectors $f$ and $g$, as well as $f'$ and $g'$,
are collinear. Then Cauchy--Schwarz inequalities become  equalities.

Below we give another example of the extremal property of a perfect coloring.
Denote by $M_2=M_2(G)$ the  adjacency matrix of  distances  $2$  of the graph $G$,
i.e.,
\[
M_2(x,y)=\left\{
\begin{array}{ll}
1\mbox{,} & d(x,y)=2;\\
0\mbox{,} & \mbox{otherwise.}
\end{array}
\right.
\]
As noted above, the  $(x,y)$-entry of $M^2$ is the 
 number of walks of length $2$ from $x$ to $y$. In particular, traversing
an edge twice gives a walk of length $2$ that starts and ends at the same vertex.

%\begin{definition}
A regular graph is called {\sl amply regular} if every pair of  adjacent
vertices has the same number of common neighbors and every pair of vertices at
distance $2$  has the same number of common neighbors.
%\end{definition}

In other words, every pair of adjacent vertices in such a graph has $\lambda$
common neighbors, and every pair of nonadjacent vertices has $0$ or $\mu>0$ common
neighbors. Consider a graph that is amply regular. It is clear that
between any two vertices at  distance of $1$ there is the same number
$a_1=\lambda$ of walks of length $2$, and between any two vertices at  distance of $2$ there is the same number $c_2=\mu$ of walks of length $2$. Therefore,
$M^2=c_{2}M_{2}+ a_1M+ b_0I, $ where $b_0$
is the degree of the graph. The notations $c_2,a_1,b_0$ are consistent with
the standard notation for distance-regular graphs (see
Section \ref{4.1}). Thus, 
\begin{equation}\label{e1.2}
M_2=p_2M^2+p_1M+p_0I, 
\end{equation}
for
some $p_0,p_1,p_2$.  Note that for an amply $r$-regular graph $r=-p_0/p_2$. Let $f$ be an arbitrary $k$-coloring of the graph
$G$, let $f_i$ be the characteristic function of color $i$, and let $n_i$ be the number of vertices of color $i$. As above,  we denote by $s_{ij}(x)$  the number of neighbors of $x$  colored by $j$. Here we imply  that  $x$ has 
$i$-th color. We define the {\sl average quotient matrix} $S_f=({\bar s}_{ij})$, where $ {\bar s}_{ij}=
\frac{1}{n_i}\sum\limits_{f(x)=i} s_{ij}(x),$
 i.e.,  $\bar s_{ij}$ is the average number of neighbors  of color $j$ of  a vertex of color $i$.

\begin{theorem}[Krotov \cite{Krotov12}]\label{th:Krot}
Let $G$ be an amply regular graph, i.e.,
its adjacency matrix $M_2$ with respect to distance $2$ satisfies the equality
$M_2=P(M)=p_2M^2+p_1M+p_0I$, and let
$f$ be a $k$-coloring of $G$  with
average quotient matrix
$S_f$.
Then, for any color $i$, the number of pairs of vertices of color $i$ at distance
$2$ from each other does not exceed $n_i(p(S_f))_{ii}/2$, and if equality is reached for
each color $i$, then $f$ is a perfect
$k$-coloring with quotient matrix $S_f$. \end{theorem}
\begin{proof}
Let $f_i$ be the characteristic function of $i$-th
color. Every function on  $V(G)$ can be decomposed into a
sum of functions supported on individual color classes. Then we can define the functions $g^i_j$ such that
$Mf_i=\sum_jg^i_j$, where $\supp (g^i_j)\subseteq \supp (f_j)$.

By definition
$(g^i_j,f_j)= (Mf_i,f_j)=\bar s_{ij}n_i=\bar s_{ji}n_j$. By
the Cauchy--Schwarz  inequality, we have $(g^i_j,f_j)^2\leq
(g^i_j,g^i_j)(f_j,f_j)$, and equality is reached only if
the vectors $g^i_j$ and $f_j$ are collinear. Since
$g^i_j(x)=s_{ji}(x)f_j(x)$ for any vertex of color $j$, we have
 $g^i_j=\bar s_{ji}f_j$ if the vectors $g^i_j$ and $f_j$
are collinear.

Twice the number of pairs of vertices of color $i$ at distance 2 from each other
is equal to $(M_2f_i,f_i)$. Therefore,
$$(M_2f_i,f_i)= (P(M)f_i,f_i)=(p_0f_i+p_1Mf_i+p_2M^2f_i,f_i)$$ 

$$=p_0(f_i,f_i)+p_1(Mf_i,f_i)+p_2(Mf_i,Mf_i)$$
$$=p_0n_i+p_1\bar s_{ii}n_i+p_2(\sum_jg^i_j,\sum_jg^i_j)=p_0n_i+p_1\bar s_{ii}n_i+p_2\sum_j(g^i_j,g^i_j)$$
$$\geq p_0n_i+p_1\bar s_{ii}n_i+p_2\sum_j(g^i_j,f_j)^2/n_j$$
$$=p_0n_i+p_1\bar s_{ii}n_i+p_2\sum_j\bar s_{ij}\bar s_{ji}n_in_j/n_j=n_i(P(S_f))_{ii}.$$

Moreover, equality holds  only if $g^i_j=\bar
s_{ij}f_j$ for all $i,j$, i.e., if $f$ is a perfect
coloring with quotient matrix $S_f$. \end{proof}

Denote by $\sigma(C)$ the average number of  neighbors lying
in  $C\subset V$, i.e., $\sigma(C)=e(C,C)/|C|=(M{\mathbf{1}}_C,{\mathbf{1}}_C)/|C|$;
and denote by $\sigma_2(C)$ the average number of vertices of $C$ at
distance~$2$, i.e.,
$\sigma_2(C)=(M_2{\mathbf{1}}_C,{\mathbf{1}}_C)/|C|$. Below we formulate Theorem \ref{th:Krot} for $2$-colorings.

\begin{corollary}\label{krotov} Let $G$ be an amply $r$-regular graph with
polynomial~$P$, and  let $C\subset V(G)$. If $\sigma(C)= a$ and
$\sigma(V\backslash C)= d $, then $\sigma_2(C)\leq (P(S))_{11}$ and
$\sigma_2(V\backslash C)\leq (P(S))_{22} $, where $S=\begin{pmatrix}
a & r-a \\
r-d & d
\end{pmatrix}$.
Moreover,  both inequalities are equalities if and only if
${\mathbf{1}}_C$ is a perfect $2$-coloring with quotient matrix~$S$.
\end{corollary}

By Corollary  \ref{krotov}, we  obtain the following criterion
for perfect $2$-colorings with the minimum eigenvalue in amply
regular graphs.

\begin{corollary}\label{corExt}
Let $G$ be an amply $r$-regular graph with polynomial~$P$,  and let
$C\subset V$ be an independent set.   Then $\sigma_2(C)\leq
-p_2r(\lambda_{min}+1)$, where~$\lambda_{min}$ is
the minimum eigenvalue of~$G$. Moreover, equality holds whenever ${\mathbf{1}}_C$ is
a perfect $2$-coloring with the eigenvalue~$\lambda_{min}$.
\end{corollary}
\begin{proof}
Let  $S$ be the  average  quotient matrix of  ${\mathbf{1}}_C$.
 It is clear that  $\sigma(C)=0$ and $\sigma(V\backslash C)=\frac{(|V(G)|-|C|)r-|C|r}{|V(G)|-|C|}=
r(1-\frac{|C|}{|V(G)|-|C|})$ for each independent set~$C$. A straightforward computation shows 
 that $P(S)_{11}=p_2(r(r-\sigma(V\backslash C))-r)$. By Theorem \ref{th:Krot}, it
holds that $\sigma_2(C)\leq p_2(\frac{r^2|C|}{|V(G)|-|C|}-r)$. By the Delsarte--Hoffman
bound, we get $\frac{r^2|C|}{|V(G)|-|C|}\leq -\lambda_{min}r$.
Thus $\sigma_2(C)\leq -p_2r(\lambda_{min}+1)$.  For any perfect
$2$-coloring with the eigenvalue~$\lambda_{min}$, the
equality $\sigma_2(C)= -p_2r(\lambda_{min}+1)$ follows from Corollary \ref{krotov} and the Delsarte--Hoffman
bound
(Proposition \ref{PCclaim11}). \end{proof}

\subsection{Generalizations of Delsarte--Hoffman Bound}

By the definition,  a  quotient matrix
determines how many vertices of one color are  adjacent to a vertex of another color.
In an amply regular graph, the quotient matrix also determines the number of vertices of one color that lie at distance $2$ from a vertex of the other color.
Let $\begin{pmatrix}
a & b\\
c & d
\end{pmatrix}$  be a quotient matrix and 
let $C$ be  the set of first
color vertices. Then we have $\sigma(C)=a$,
$\sigma_2(C)=(p_2(a^2+bc)+p_1a+p_0)$, and
$|C|=\frac{c|V(G)|}{b+c}$. Define
$\beta=p_2(a^2+bc)+p_1a+p_0$, then
$\frac{c}{b+c}=\frac{bc}{b^2+bc}=\frac{\beta-P(a)}{p_2b^2+\beta-P(a)}$.

Next we will show
that, in an amply regular graph an upper bound on  $\sigma_2(C)$ yields  an upper bound on the cardinality of $|C|$. Moreover, when this
upper bound  is attained,  the characteristic
function of  $C$ is a perfect $2$-coloring.

\begin{theorem}\label{thEPPC}
Let $G$ be an amply $r$-regular graph with polynomial
$P$, and let $C\subset V$. If  $a=\sigma(C)$ and $\beta=\sigma_2(C)$
then $|C|\leq \frac{(\beta-P(a))n}{p_2(r-a)^2+\beta-P(a)} $.
Moreover, if $|C|= \frac{(\beta-P(a))n}{p_2(r-a)^2+\beta-P(a)} $
then
 ${\mathbf{1}}_C$ is a perfect $2$-coloring of~$G$.\end{theorem}

\begin{proof} Without loss of generality, we suppose that  $G$ is
connected.
Otherwise, the theorem can be proved separately for each connected component.
 Expand ${ {\mathbf{1}}}_{C}$
in an orthonormal basis of eigenfunctions of $G$:
${\mathbf{1}}_{C}=\sum_i\alpha_i\phi_i$, where $\phi_i$ is an
eigenfunction of~$M$ with   eigenvalue~$\lambda_i$. Without loss of
generality, we assume  $\|\phi_i\|_2=1$ for all~$i$. The
eigenfunction with   eigenvalue~$r$ is equal to
$\phi_0=\mathbf{1}_V/\sqrt{n}$.  From $({
{\mathbf{1}}}_{C},{\mathbf{1}}_{C})= \sum_i\alpha^2_i$ we obtain
$$
{\rm (I)} \qquad \sum\limits_{i\neq 0}\alpha^2_i=
({\mathbf{1}}_{C},{ {\mathbf{1}}}_{C})-\alpha^2_0=|C|- \varrho|C|,
$$
 where $\varrho=|C|/n$. From $(M{\mathbf{1}}_{C},{\mathbf{1}}_{C})=
 a|C|$ and $
(M{\mathbf{1}}_{C},{\mathbf{1}}_{C})=\sum_i\alpha^2_i\lambda_i $, it
follows
$$
{\rm (II)} \quad r\varrho|C|+\sum\limits_{i\neq
0}\alpha^2_i\lambda_i= a|C|.$$ From (\ref{e1.2}) and  the hypothesis
of the theorem we obtain
$$(M^2{\mathbf{1}}_{C},{\mathbf{1}}_{C})=
\frac{1}{p_2}((M_2-p_1M-p_0I){\mathbf{1}}_{C},{\mathbf{1}}_{C})=
\frac{|C|}{p_2}(\beta-p_1a-p_0).  $$ Therefore,
$$
{\rm (III)} \quad r^2\varrho|C|+\sum\limits_{i\neq
0}\alpha^2_i\lambda^2_i= \frac{|C|}{p_2}(\beta-p_1a-p_0).
$$

Consider the linear combination
${\rm (III)}-2\theta{\rm
(II)}+\theta^2{\rm (I)}$.  Then for any $\theta \in \mathbb{R}$ we
obtain the  inequalities
$$r^2\varrho|C|-2 r\varrho\theta|C|+\sum\limits_{i\neq
0}\alpha^2_i(\lambda_i-\theta)^2=
\frac{|C|}{p_2}(\beta-p_1a-p_0)-2a\theta|C|+|C|(1-\varrho)\theta^2.$$
Since $\sum\limits_{i\neq
0}\alpha^2_i(\lambda_i-\theta)^2\geq 0$, we obtain
\begin{equation}\label{eqEx8}
r^2\varrho-2 r\varrho\theta\leq
\frac{1}{p_2}(\beta-p_1a-p_0)-2a\theta+(1-\varrho)\theta^2,
\end{equation}
$$\varrho\leq\frac{\frac{1}{p_2}(\beta-P(a))+(a-\theta)^2}{(r-\theta)^2}.$$
Let $\theta=a-\frac{\beta-P(a)}{p_2(r-a)}$. Then we conclude that
$$\varrho\leq\frac{(a-\theta)(r-a)+(a-\theta)^2}{(r-a+a-\theta)^2}=\frac{a-\theta}{r-\theta}=   \frac{\beta-P(a)}{p_2(r-a)^2+\beta-P(a)}.$$

It is clear that (\ref{eqEx8}) holds with equality if and only if
${\mathbf{1}}_{C}=\phi+\alpha_0\varphi_0$, where $\phi$ is an eigenfunction with
eigenvalue $\theta$. By Proposition \ref{corol5}, this implies  that
${\mathbf{1}}_{C}$ is a perfect $2$-coloring. 
\end{proof}

 For $a=0$ and $\beta=0$, the new bound coincides
 with the
Hamming bound $|C|\leq \frac{|V(G)|}{r+1} $. If $C$ is an independent
set, then $a=0$, $P(0)=p_0=-rp_2$ and $\varrho\leq
1/(1+\frac{p_2r^2}{\beta+p_2r})$. This bound and the Delsarte--Hoffman bound
are reached simultaneously on perfect $2$-colorings with the minimum
eigenvalue.  In this case, it holds that $\sigma_2(C)=
-p_2r(\lambda_{min}+1)$. By Corollary~\ref{corExt},  we
have $\sigma_2(C)\leq -p_2r(\lambda_{min}+1)$ for any
independent set~$C$. Consequently, it is meaningful to consider only the
case   $\beta< -p_2r(\lambda_{min}+1)$.

\begin{corollary}\label{corExt2}
Let $C$ be an independent set  in an amply $r$-regular graph $G$ with
polynomial~$P$. If $\beta=\sigma_2(C)<
-p_2r(\lambda_{min}+1)$, then the bound $\frac{|C|}{|V(G)|}\leq
1/(1+\frac{p_2r^2}{\beta+p_2r})$ is  better than the Delsarte--Hoffman
bound.
\end{corollary}
\begin{proof}
$1/(1+\frac{p_2r^2}{\beta+p_2r})<1/(1+\frac{p_2r^2}{-p_2r\lambda_{min}})=\frac
{-\lambda_{min}}{r-\lambda_{min}}$. 
\end{proof}

\section{Distance-Regular Graphs}

\subsection{Intersection Array and Krawtchouk Polynomial}\label{4.1}

Let $G$ be a simple graph. For a
fixed vertex $x_0$ we consider the {\sl distance
partition} of the graph $G$ with respect to $x_0$: $$A_t(x_0)=\{x\in V(G) \mid
d(x,x_0)=t\}.$$

%\begin{definition}
A simple connected graph $G$ is called {\sl distance-regular} if, for every vertex $x_0$, the distance partition of $G$  is equitable and  its quotient matrix is independent of the choice of $x_0$. 
%\end{definition}

 Consider layer $A_t(x_0)$ in detail.  A vertex $x\in A_t(x_0)$ is adjacent to
$a_t$ vertices in its layer, to $b_t$ vertices in the next layer, and to
$c_t$ vertices in the previous layer. In a distance-regular graph,
these  parameters $a_t, b_t, c_t$ are independent of the choice of the initial vertex $x_0$
and the vertex $x\in A_t$ from a layer.

%\begin{definition}
The set of values $(b_0,\dots,b_{d-1};c_1,\dots,c_d)$ is called
the {\sl intersection array} of a distance-regular graph of diameter $d$.
%\end{definition}

Any distance-regular graph is regular; the degree
of each vertex in it is $b_0=r$. Strictly speaking, one should say $a_t +b_t +c_t=r$, $0\leq t\leq d$, with conventions $c_0=0$ and $b_d=0$.

From the definition, the distance coloring  of a distance-regular graph with respect to any vertex
$x_0$ is perfect with quotient matrix 
 $$S[G]=\begin{pmatrix}
0 & b_0 & 0 & 0& \dots & 0\\
c_1 & a_1 & b_1 & 0 & \dots & 0\\
0 & c_2 & a_2 & b_2 & \dots & 0 \\
\dots & \dots & \dots & \dots & \dots & \dots \\
0 & \dots & 0 & 0 & c_d & a_d
\end{pmatrix}.$$

The quotient matrix of the distance coloring is called the intersection matrix of $G$.
Therefore, every eigenvalue of $S[G]$ is an  eigenvalue of $G$.

\begin{example}\label{ex:hamm_cube_bool_parameters}
The Boolean $n$-cube $Q_2^n=\{0,1\}^n$ is a distance-regular
graph of diameter $n$ with parameters $a_t=0$, $b_t=n-t$, $c_t=t$.
\end{example}

{
\begin{example}
The hypercube $Q_q^n$ is a distance-regular graph of diameter
$n$ with parameters $a_t=(q-2)t$, $b_t=(n-t)(q-1)$, $c_t=t$.
\end{example}
}

{
\begin{example}\label{ex:John}
{The Johnson graph} $J(n,k)$ is defined as follows.  We take all binary $n$-tuples of 
 weight $k$ as a
set of vertices. Vertices are connected by an edge if the Hamming distance between
two $n$-tuples is equal to $2$. The Johnson graph $J(n,k)$ is
distance-regular with parameters 
$b_t=(k-t)(n-k-t)$, $c_t=t^2$, $a_t=k(n-k)-(k-t)(n-k-t)-t^2=t(n-2t)$.
\end{example}
}

Constructions and properties of special
distance-regular graphs is available in \cite{BCN}.

\begin{remark} A natural question arises: is a distance-regular graph uniquely determined
by its intersection array? It turns out
not. For example, $Q_4^2$ and the Shrikhande graph (see
Example \ref{exSh} below) have the same intersection array
$(b_0,b_{1};c_1,c_2)=(6,3;1,2)$. Moreover, Cartesian products (see Section \ref{Cart})
of these graphs with the same number of factors are nonequivalent distance-regular graphs with matching intersection arrays.
\end{remark}

\begin{remark} 
Listed above distance-regular graphs are transitive.
However, these two classes are incomparable: there exist distance-regular graphs that are not transitive and transitive graphs that are not distance-regular (see \cite{BCN}).
\end{remark}

%\begin{remark}
A  {\sl distance-$t$  graph} is a graph with adjacency matrix $M_t$.
A  {distance-$t$  graph}, obtained from a distance-regular graph, is usually no longer distance-regular. But there are
exceptions. For example, the distance-$2$ graph obtained from  the Boolean $n$-cube for any $n$ is
distance-regular.
%\end{remark}

Let $M_t$ be the adjacency matrix of the   distance-$t$ graph, that is
\[
M_t(x,y)=\left\{
\begin{array}{ll}
1\mbox{,} & d(x,y)=t;\\
0\mbox{,} & \mbox{otherwise.}
\end{array}
\right.
\]

The main property of a distance-regular graph is the following. Matrices  $M_t$ for any $t$ can be represented as  a polynomial in the adjacency matrix $M$ of the graph.
Note that the $(ij)$th element of $M^k$ is equal to  the number of walks  of length $k$ between
  $i$th and $j$th vertices.  These walks are not required to be simple and may 
traverse the same edge several times.  For example, $(M^2)_{ii}$ is the twice number of the degree of $i$th vertex.

\begin{example}\label{ex:hamm_cube_bool_distmatrix}
Let $M$ be the adjacency matrix of the Boolean $n$-cube $Q^n_2=\{0,1\}^n$.
Any vertex at distance $2$ from the original can be reached
in two ways: by moving
first along the edge of direction $i$, then $j$,  $j\neq i$, and vice versa.

$M_0=I_{2^n}$, $M_1=M$,

$M^2 = b_0I_{2^n}+c_2M_2= n I_{2^n} +2M_2$. \end{example}

\begin{claim}\label{remDR}
A graph of diameter $d$ is distance-regular if and only if
each of its distance-$t$ adjacency matrices,
$t=2,\dots,d$, is a polynomial of degree $t$ in its
adjacency matrix $M$.
\end{claim}
\begin{proof}
Fix a vertex $y$ of a distance-regular graph,
take a vertex $x_i$ at distance $t$ from it, and consider its neighbors.
For the adjacency matrices 
of a distance-regular graph, the following equality holds:
\begin{equation}\label{eqST0}
M_t M= a_tM_t +b_{t-1}M_{t-1}+c_{t+1}M_{t+1}.
\end{equation}

Indeed, in the matrix $M_t M$, position $(i,j)$ indicates the number
of ways to get from vertex $x_i$ to $x_j$ by first moving to an intermediate vertex $z$ adjacent to $x_j$, located at a distance
$t$ from $x_i$, and then from there to vertex $x_j$. Clearly, the element
of the matrix $(a_tM_t)_{ij}$ is equal to the number of ways to complete such a path
if $d(x_i,x_j)=t$; $(b_{t-1}M_{t-1})_{ij}$, if $d(x_i,x_j)=t-1$;
$(c_{t+1}M_{t+1})_{ij}$, if $d(x_i,x_j)=t+1$. Thus, the recurrence relations (\ref{eqST0}) hold.
By
induction, we obtain $M_t=P_t(M)$, where $P_t$ is a polynomial of degree
$t$.  This shows that the multiplication of adjacency matrices by different
distances is commutative, i.e., $M_tM_s=M_sM_t$. Moreover, from (\ref{eqST0}) follows the recurrence relation
\begin{equation}\label{eqST11}
xP_t(x)=a_tP_t(x)+b_{t-1}P_{t-1}(x)+c_{t+1}P_{t+1}(x).
\end{equation}

Let's prove that
the converse is also true.  If $M_t=P_t(M)$ then by induction we see that $M^{t+1}$ and $M_t M$ are  linear combinations of $I, M,\dots, M_{t+1}$. 
Since an edge changes the distance to a fixed vertex by at most one, every walk contributing to $M_tM$ ends at distance $t-1$, $t$, or $t+1$. So,  $M_tM$ is a nonnegative linear combination of $M_{t-1}, M_t, M_{t+1}$. We obtain equation (\ref{eqST0}) for some coefficients $a_t, b_{t-1}, c_{t+1}$.
 Let $e_0$  be the indicator function of the initial vertex.  Then $M_te_0$ is the characteristic function of $t$th color and  $(M_tMe_0)_j=(MM_te_0)_j$ indicates how many neighbors  of  the $t$th color  the $j$th vertex has.    From  (\ref{eqST0})  we obtain that if  $j$th vertex has $t$th color then the number of such vertex is $a_t$; if  $j$th vertex has $(t-1)$th color then the number of such vertex is $b_{t-1}$ and 
 if  $j$th vertex has $(t+1)$th color then the number of such vertex is $c_{t+1}$.  Therefore, we prove that the distance coloring is perfect and its quotient does not depend on the initial vertex.
\end{proof}

In the case of a Boolean $n$-cube, the polynomial $P_{t}$ (\ref{eqST11}) is called 
the  Krawtchouk polynomial.  A strict form of   Krawtchouk polynomial will be calculated in Section  \ref{10.2}. 
 If we consider a $q$-ary
cube, then the polynomial $P_{t}$ is called the $q$-Krawtchouk polynomial
$P_{t}[n,q](y)$. For Johnson graphs,
the polynomial $P_{t}$ is called the Eberlein polynomial.  For other distance-regular graphs we treat $P_t$ as  a
Krawtchouk-type polynomial.

\begin{theorem}[Lloyd, Shapiro and Zlotnik, Martin]\label{th:shapiro_zlotnik_martin}
Let $G$ be a distance-regular graph, $G_t$ be the $t$-distance graph obtained from $G$, i.e., $V(G_t)=V(G)$,
$M(G_t)=M_t$. If $f$ is a perfect coloring of $G$ with
quotient matrix $S$, then $f$ is a perfect coloring of $G_t$ with quotient matrix $P_t(S)$.
\end{theorem}
\begin{proof}
We have $MF=FS$, where $F$ is the matrix corresponding to the coloring $f$.
From the equality $M_t=P_t(M)$, we have $M_t F=P_t(M) F$.
Furthermore, we have the equations $M^t F=M^{t-1}MF=M^{t-1}(FS)=FS^t$. Then
$P_t(M) F=FP_t(S)$. By Proposition \ref{s:perfcol_criteria},  we see that $f$ is a perfect coloring of $G_t$ with quotient matrix $P_t(S)$.
\end{proof}

\begin{corollary}\label{c:szm_spherecoloring}
For any perfect coloring of a distance-regular graph,
the color composition of a sphere of radius $t$ depends only on the color of the center
of the sphere. \end{corollary}

%\begin{remark}
This fact was originally established by Lloyd \cite{Lloyd} for
perfect codes in the Boolean $n$-cube. Shapiro and Zlotnik
\cite{Shapiro} found the weight distribution (see below) of a perfect
code. Martin established  this claim for any  distance-regular graph in his  thesis \cite{Martin}.
%\end{remark}

%\begin{remark}
%How many vertices of what color is important for coding theory.
%\end{remark}

Theorem \ref{th:shapiro_zlotnik_martin} provides a way to prove
the nonexistence of a perfect coloring. Substituting the quotient matrix $S$
of a perfect coloring into the polynomial $P_t$, we should obtain a matrix
whose entries are  nonnegative integer. Furthermore, the vector
of color proportions calculated from the matrix $P_t(S)$ and multiplied by the size
of the  graph must also be integral.
To the best of our knowledge, no examples are known in which a matrix $S$ satisfies the conditions
of Theorem \ref{thAF}, but  some matrix $P_t(S)$ fails to have nonnegative integer entries .

\subsection{Strongly Regular Graphs}

%\begin{definition}
A distance-regular graph of diameter $2$ is called {\sl strongly
regular}.
%\end{definition}

It is easy to see that intersection array $(b_0,b_1; c_1,c_2)$ of  a strongly regular graph calculated from only three parameters: degree $r$ of the graph, the number $\lambda$ of common neighbors of two adjacent vertices and the number $\mu$ of common neighbors of two non-adjacent vertices. Thus, $b_0=r$, $b_1=r-\lambda-1$, $c_1=1$, $c_2=\mu$.

{
\begin{claim}[\cite{GodsilR}] A connected regular graph $G$ with exactly three distinct
eigenvalues is strongly regular.
\end{claim}
\begin{proof}
Let $M$ be the adjacency matrix of an $r$-regular graph $G$. The matrix $M$
is the root of its characteristic polynomial $P$, therefore
$(M-rI)P(M)=0$, where $P(x)=(x-\lambda_1)(x-\lambda_2)$ and
$\lambda_1<\lambda_2<r$ are the eigenvalues of $M$. From
the equality $(M-rI)P(M)=0$ it follows that  $M^3$ is equal to a linear
combination of the matrices $I$, $M$, and $M^2$. Then any walk of length $3$ in
graph $G$ connects vertices that are at distance $0$, $1$, or
$2$. Consequently, the graph $G$ has diameter $2$, i.e.,  $M_2=J-M-I$. Since the graph
$G$ is connected, it follows that the matrix $M$ has a unique eigenvector
${\bf 1}$ with eigenvalue $r$. Then, from the equation
$(M-rI)P(M)=0$, it follows that all columns of the matrix $P(M)$ are collinear
with the vector ${\bf 1}$. Similarly, all rows of the matrix $P(M)$ are collinear
with the vector ${\bf 1}$. Then $P(M)=\alpha J$, where all components of 
$J$ are equal to $1$. We have  $P(M)=\alpha(M_2+M+I)$.
\end{proof}
}

If the complement of a strongly regular graph $G$ with parameters $(v=|V(G)|, r, \lambda, \mu)$ is connected then it is
also a strongly regular graph with parameters $(v,v-r-1,v-2-2r+\mu, v-2r+\lambda)$.
The eigenvalues of a strongly regular graph with parameters $(v, r, \lambda, \mu)$ are 
(see, e.g.~{\cite[Theorem\,9.1.3]{BHaemers}})
$$
 r, \quad \frac{\lambda-\mu \pm \sqrt{(\lambda-\mu)^2+4(r-\mu)}}2,
$$
with multiplicities, respectively,
$$
 1, \quad \frac12
 \left(
 {v-1 \mp \frac{2r+(v-1)(\lambda-\mu)}{\sqrt{(\lambda-\mu)^2+4(r-\mu)}} }
 \right).
$$

%\subsection{Krawtchouk Polynomials}\label{4.2}

\subsection{Minimum Support of Eigenfunctions}\label{4.4}

In this section, we will explore how large the minimum difference between two perfect colorings can be. At first we need to consider the cardinality  of an eigenfunction's support.
The {\sl support}
of a function $f:V(G)\rightarrow \mathbb{C}$ is the set $\supp(f)=\{x\in V(G) \mid f(x)\ne
0\}$.
{
%\begin{definition}
The {\sl weight distribution} of a function $f:V(G)\rightarrow \mathbb{C}$
with respect to a vertex $a\in V(G)$ of a graph $G$ is the sequence
$\Big(f(a),\dots,\sum\limits_{d(y,a)=t}f(y),\dots,\sum\limits_{d(y,a)=d}f(y)\Big)$,
where $d$ is the diameter of $G$.
%\end{definition}}

\begin{claim}\label{c:WDB_bound0}
Let $f$
be an eigenfunction of a distance-regular graph $G$ with eigenvalue
$\lambda$, i.e., $Mf=\lambda f$.
Then $\sum\limits_{d(y,a)=t}f(y)=P_t(\lambda) f(a)$.
\end{claim}
\begin{proof}
Similarly to the proof of  Theorem \ref{th:shapiro_zlotnik_martin}, we have
the equalities $M^i f=\lambda^i f$ and $M_t f = P_t(M) f = P_t(\lambda) f$.
Consequently,
$\sum\limits_{d(y,a)=t}f(y)=(M_tf)(a)=P_t(\lambda) f(a)$.
\end{proof}

{ Thus, in a distance-regular graph, the sum of the values
of an eigenfunction over the sphere depends only on the eigenvalue and
the value of the function at the center of the sphere. For definiteness, we can assume
that $f(a)=1$. Then the set $[1, P_1(\lambda),\ldots,P_d(\lambda)]$
is the weight distribution (spectrum) of an eigenfunction with
eigenvalue $\lambda$.}

Any function $g:V(G)\rightarrow \mathbb{C}$ can be represented as
a linear combination of the graph eigenfunctions:
$g=\sum_i\alpha_if_i$, where $f_i$ is an eigenfunction
of  $G$ with eigenvalue $\lambda_i$.
Then we have

\begin{corollary} Let $g$ be defined on
the vertices of a distance-regular graph $G$. Then
$\sum\limits_{d(y,a)=t}g(y)=\sum_i\alpha_iP_t(\lambda_i) f_i(a)$.
\end{corollary}

\begin{claim}[\cite{KMP}]\label{c:WDB_bound}
Let $f$ be an eigenfunction with eigenvalue $\lambda$ in a
distance-regular graph. Then $|\supp(f)| \ge
\sum\limits^d_{t=0} |P_t(\lambda)|$, where $d$
is the diameter of the graph.
\end{claim}
\begin{proof}
Choose a vertex $a$ such that $f$ has its maximum absolute value, that is, $|f(a)|
= \underset{x\in V(G)}{\mathrm{max}} |f(x)|$. From Proposition
\ref{c:WDB_bound0} we have $\sum\limits_{d(y,a)=t}f(y)=P_t(\lambda)
f(a)$, which yields: $$ |\{y \mid d(y,a)=t, f(y)\ne 0\}|\geq  \sum\limits_{d(y,a)=t}\frac{|f(y)|}{|f(a)|}\geq
|P_t(\lambda)|.$$
\end{proof}

Consider perfect colorings $f_1$ and $f_2$
of a distance-regular graph $G$ with two colors $0$ and $1$ with
the same quotient matrix.
Since $f_1$ and $f_2$ have the same quotient matrix,
their difference is orthogonal to the constant function and is an
eigenfunction corresponding to the nontrivial eigenvalue $\lambda$ of the
quotient matrix.
The number of vertices, where  $f_1$ and $f_2$ differ can be estimated by Proposition \ref{c:WDB_bound}:
$|\supp(f_1-f_2)| \ge \sum\limits_{t=0}^d|P_t(\lambda)|$.

A detailed review of minimum supports of eigenfunctions of graphs is available in \cite{Sot}.

%\begin{definition} 
The function $u:V(G)\rightarrow \{0,\pm1\}$ is called a {\sl bitrade (ball bitrade)} in graph $G$ if for any ball $B(x)$, $x\in V(G)$,
of radius $1$
the following holds\\
$
\left\{
\begin{array}{ll}
|B(x)\cap u^{-1}(1)|=1,\\
|B(x)\cap u^{-1}(-1)|=1,
\end{array}
\right.
$
or $u(B(x))=\{0\}$.
%\end{definition}
In other words, the function takes the values $0$ and $\pm1$ and
either takes exactly one value each of $1$ and $-1$ in a ball of radius $1$,
or only  zeros in this ball. The difference
of the characteristic functions of two $1$-perfect codes is a bitrade.
It is easy to see that a bitrade is an eigenfunction with eigenvalue $-1$. Indeed, if the ball $B(x)$  contains only zeros, then $Mu(x)=0=-u(x)$, otherwise the ball contains one $1$ and one $-1$, so
$u(x)+Mu(x)=0$.
By Proposition \ref{c:WDB_bound}, we have $|\supp(u)| \ge \sum\limits_{t=0}^d|P_t(-1)|$,
where $d$ is the graph diameter.

{ Denote by $A_t$  the set of vertices at distance
$t$ from some fixed vertex $x_0$
in  $G$.

\begin{claim}\label{c:dopeigenfun}
Let $f_1, f_2$ be eigenfunctions of a distance-regular
graph $G$ of diameter $d$ with the same eigenvalue $\lambda$,
$P_d(\lambda)\neq 0$. If the functions $f_1$ and $f_2$ coincide on
the set ${A_{t-1}\cup A_{t}}$, where $t<d/2$, then
$f_1|_{A_i}=f_2|_{A_i}$ for $i\leq t$.
\end{claim}
\begin{proof}
Let $x'\in A_i$ and $f_1(x')\neq f_2(x')$. Denote by $A'_d$  the set
of vertices that are at distance $d$ from some fixed
vertex $x'$. The triangle inequality implies that $A'_d\cap
A_i=\varnothing$ for $i\leq t$. Indeed, if $y\in A'_d\cap
A_i$ then $d(y,x')\leq i+i<d$. From the condition $P_d(\lambda)\neq 0$ and
Proposition \ref{c:WDB_bound0} it follows that $\sum\limits_{x\in
A'_d}f_1(x)\neq \sum\limits_{x\in A'_d}f_2(x)$. Consider the function
$f_3$ defined by the equality
$$f_3(x)=\left\{
\begin{array}{ll}
f_1(x), & x\in A_i, i\leq t\\
f_2(x), & x\in A_i, i>t.\\
\end{array}
\right.
$$
Since  $f_1$ and $f_2$ coincide on the set $A_{t-1}\cup A_{t}$,
every vertex outside $A_{t-1}\cup A_{t}$,
 has all its neighbors in layers on which $f_3$
agrees with either $f_1$ or $f_2$. Then $f_3$ is an eigenfunction with eigenvalue
$\lambda$. By Proposition \ref{c:WDB_bound0} it follows that
$$\sum\limits_{x\in
A'_d}f_3(x)=P_d(\lambda)f_3(x')=P_d(\lambda)f_1(x')=\sum\limits_{x\in
A'_d}f_1(x).$$ But by construction, $f_3|_{A'_d}=f_2|_{A'_d}$. We obtain a
contradiction.
\end{proof}

\begin{corollary}
If $P_d(-1)\neq 0$, then the support of a bitrade has 
diameter at least $\lceil d/2\rceil-2$.
\end{corollary}

In particular,  the diameter of  $Q^n_2$ is $n$,
the eigenvalues are $\lambda_k=n-2k$, $k=0,\dots,n$, and
$P_d(\lambda_k)=(-1)^{k+1}$ (see Section \ref{10.2}).  Then any bitrade has diameter at least $(n-3)/2$ if $n$ is odd and $(n-4)/2$ if $n$ is even.

\subsection{Orthogonality Relations for  Krawtchouk-Type Polynomials}

 Let $s_t=|A_t|$. Since a constant
is an eigenfunction with eigenvalue $r$ of an $r$-regular graph, by Proposition \ref{c:WDB_bound0} we have
$s_t=P_t(r)$.

{
\begin{claim}\label{c:polyKrv}
Let $\lambda$ be an eigenvalue of a distance-regular graph
$G$ of diameter $d$ and  $$z(\lambda)=(z_0(\lambda), \dots, z_d(\lambda) ), \qquad z_t(\lambda)=P_t(\lambda)/s_t, \qquad t=0,\dots,d.$$
Then the vector $z(\lambda)$ is an eigenvector of the intersection matrix
 $S[G]$ of $G$ with eigenvalue $\lambda$.
\end{claim}
\begin{proof}
Consider some eigenfunction $f$ of graph $G$ with
eigenvalue $\lambda$. By the definition of an eigenfunction for
an arbitrary vertex $x\in A_t$, we have the equality
\begin{equation}\label{dop44}\lambda f(x)= \sum\limits_{y: y\in
A_{t-1},d(x,y)=1}f(y)+ \sum\limits_{y: y\in
A_{t},d(x,y)=1}f(y)+\sum\limits_{y: y\in
A_{t+1},d(x,y)=1}f(y).\end{equation}
Let us sum both sides
of the equality over all vertices $x\in A_t$. 
Since every vertex of $A_t$ has exactly $c_t$ neighbors in $A_{t-1}$, the number of edges between $A_{t-1}$ and $A_t$ is $c_t|A_t|$.
By  distance regularity, 
 each vertex  $y\in A_{t-1}$ is adjacent to the same number $\frac{c_t|A_t|}{|A_{t-1}|}$ of vertices in $A_t$. We have
the equality $\sum\limits_{x\in A_t}\sum\limits_{y: y\in
A_{t-1},d(x,y)=1}f(y)=\frac{c_t|A_t|}{|A_{t-1}|}\sum\limits_{y\in
A_{t-1}}f(y)$. Similarly, transforming the other terms and
applying Proposition \ref{c:WDB_bound0}, we obtain the equality
$$\lambda P_{t}(\lambda)=\frac{c_t|A_t|P_{t-1}(\lambda)}{|A_{t-1}|}
+\frac{a_t|A_t|P_{t}(\lambda)}{|A_{t}|}+\frac{b_t|A_{t}|P_{t+1}(\lambda)}{|A_{t+1}|},$$
which is equivalent to the required equality $\lambda
z_t(\lambda)=c_tz_{t-1}(\lambda)+a_tz_t(\lambda)+b_tz_{t+1}(\lambda)$.
\end{proof}

\begin{remark}\label{remSG}
Since a graph $G$ of diameter $d$ has at least $d+1$ distinct
eigenvalues (see Problem \ref{exe13}),  $S[G]$ has
$d+1$ distinct eigenvalues, and the sets of eigenvalues
of $G$ and $S[G]$ coincide.
\end{remark}

\begin{corollary}
Let $\bar{z}^\lambda(x)=z_t(\lambda)$ if $x\in A_t$. Then
$\bar{z}^\lambda$
is an eigenfunction of $G$ with eigenvalue $\lambda$.
\end{corollary}
\begin{proof}
The equation $\lambda
z_t(\lambda)=c_tz_{t-1}(\lambda)+a_tz_t(\lambda)+b_tz_{t+1}(\lambda)$
is equivalent to equality (\ref{dop44}), which coincides with
the definition of an eigenfunction.
\end{proof}

\begin{corollary}[Orthogonality of Krawtchouk-type Polynomials]
We introduce the notation
$\widetilde{z}_t(\lambda)=\sqrt{s_t}z_t(\lambda)=P_t(\lambda)/\sqrt{s_t}$.
If $\lambda_1\neq \lambda_2$
are eigenvalues of a distance-regular graph $G$, then the vectors
$\widetilde{z}(\lambda_1)$ and $\widetilde{z}(\lambda_2)$
are orthogonal.
\end{corollary}
\begin{proof}
From the orthogonality of the eigenfunctions $\bar{z}^{\lambda_1}$ and
$\bar{z}^{\lambda_2}$ of a graph $G$ with different eigenvalues, we have the equalities

\begin{align*}
0
&=
\sum_{x\in V(G)}
\bar z^{\lambda_1}(x)\bar z^{\lambda_2}(x)
=
\sum_{t=0}^d
s_t\,z_t(\lambda_1)z_t(\lambda_2)
\\
&=
\sum_{t=0}^d
s_t
\frac{P_t(\lambda_1)}{s_t}
\frac{P_t(\lambda_2)}{s_t}
=
\sum_{t=0}^d
\widetilde z_t(\lambda_1)
\widetilde z_t(\lambda_2).
\end{align*}
\end{proof}

Equivalently,
$
\sum_{t=0}^d
\frac{P_t(\lambda_i)P_t(\lambda_j)}{s_t}
=0,
\qquad i\neq j.
$

\subsection{Completely Regular Codes}\label{4.5}

%\begin{definition}\label{d:completely_regular_code}
A subset of the vertices of a graph is called a {\sl completely regular code}
if the distance coloring with respect to it
is a perfect coloring.  
%\end{definition}

In this case, we assume that the code is the zero
color, the set of vertices at distance $1$ from the code is the first color,
and so on. It is clear that every one-element subset of a distance-regular graph is a completely regular code.
Every color of  $2$-perfect coloring in any regular graph is completely
regular code by definition. It is easy to show that in a
distance-regular graph, any perfect code with arbitrary
code distance (see Section \ref{PC}) is completely regular.
The theory of completely regular codes is presented in detail in \cite{Koolen} and \cite{CRC}.

{
%\begin{definition}\label{d:code_spectra}
Consider a graph $G$ and a vertex $x_0\in V(G)$. Recall that by $A_i$
we denote the set of vertices at distance $i$ from the vertex
$x_0$. Let $C$
be a code and $|A_i\cap C|$
be the weight of its intersection with the $i$-th layer $A_i$. Then the {\sl weight
distribution (spectrum)} of the code (relative to the vertex $x_0$,
$A_0=\{x_0\}$) is the weight distribution of its
characteristic function, i.e., the set $[|A_0 \cap C|, |A_1 \cap
C|, \ldots, |A_d \cap C|]$.}
%\end{definition}

By the Theorem \ref{th:shapiro_zlotnik_martin}, in distance-regular graphs  the weight distribution of a completely regular
code depends only on the parameters of the distance coloring  and a color of the initial vertex relative to which the distance refinement  is considered. In particular, all $1$-perfect codes in a
distance-regular graph that
contain the initial vertex have the same weight distribution.

Let $G$ be a distance-regular graph, and let $C$ be a code in $G$.
For $t\geq 0$, let  $C_t=\{x\in V(G) : d(x,C)=t\}$, where $d(x,C)=\min_{y\in C}d(x,y)$. Let $\rho$ be a largest integer such that  $C_\rho \neq \varnothing$. The number $\rho$ is called the {\sl covering radius} of the code $C$.

The following theorem is actually the standard characterizations of completely regular codes (see \cite{BCN}).

\begin{theorem}\label{thcritDRG}
A code $C$ in a distance-regular graph $G$ is completely
regular if and only if, for every $t$ and every $i\geq t$,
each vertex $x\in C_t$ has the same number
of vertices from $C$  at distance $i$ from $x$.
\end{theorem}
\begin{proof}
$(\Rightarrow)$ Follows from Theorem \ref{th:shapiro_zlotnik_martin}.

$(\Leftarrow)$ From the matrix equation of a perfect coloring
(Proposition \ref{s:perfcol_criteria} or Proposition \ref{cl:invar}) it follows that for any $j$
it suffices to prove the equality
\begin{equation}\label{eqC_t}
M\mathbf{1}_{C_j}=\sum\limits_{i=0}^ds_{ji}\mathbf{1}_{C_i}
\end{equation}
with some coefficients $s_{ji}$. By the hypothesis, for every $i$ the value 
$M_i\mathbf{1}_{C}(x)$ depends only on  the distance $t=d(x,C)$. Thus, for any $i$, $i=0,\dots,d$, we have the equality
$$M_i\mathbf{1}_{C}=\sum\limits_{j=0}^dp_{ji}\mathbf{1}_{C_j}.$$
By the definition of $C_t$, we have $p_{ti}=0$ for $t>i$ and $p_{tt}>0$.
For $t=0$ (\ref{eqC_t}) is true  by the hypothesis.
Let's assume it is proven for $j<t$. Then
\begin{equation}\label{eqC_t1}
M_{t}M{\mathbf 1}_{C}=MM_{t}\mathbf{1}_{C}=\sum\limits_{j=0}^{t}p_{jt}M\mathbf{1}_{C_j}=
\sum\limits_{j=0}^{t-1}p_{jt}\left(\sum\limits_{i=0}^ds_{ji}{\mathbf 1}_{C_i}\right)+
p_{tt}M\mathbf{1}_{C_t}.
\end{equation}
From equality (\ref{eqST0}) for a
distance-regular graph, we have
$$M_t M\mathbf{1}_{C}= a_tM_t\mathbf{1}_{C}
+b_{t-1}M_{t-1}\mathbf{1}_{C}+c_{t+1}M_{t+1}\mathbf{1}_{C}.$$
Therefore,
$$M_{t}M\mathbf{1}_{C}=\sum\limits_{i=0}^dq_{ti}\mathbf{1}_{C_i},$$
for some coefficients $q_{ti}$. 
 Since $p_{tt}>0$, the coefficient of
$M{\mathbf 1}_{C_t}$ in (\ref{eqC_t1}) is nonzero.
Therefore, $M\mathbf{1}_{C_t}$ is a linear combination of the vectors
${\mathbf 1}_{C_0},\dots,\mathbf{1}_{C_d}$, proving (\ref{eqC_t}) for $j=t$.
%From (\ref{eqC_t1}) and the last equality,  (\ref{eqC_t}) follows for $j=t$. 
Thus, the induction step
is proved. \end{proof}
{
\begin{example}\label{ex:face0}
Consider the hypercube $Q^n_q$.
%By definition, $Q^n_q=Q^{k}_q\square
%Q^{n-k}_q$.
The set $C=V(Q^{k}_q)\times \{a\}$ is a $k$-dimensional
face of  $Q^n_q$. Every face of $Q^n_q$
of dimension $k$ is a completely regular code with the quotient matrix
 $S=S[Q^{n-k}_q]$.
\end{example}}

{
\begin{example}\label{ex:face}
A famous example of completely regular code is Preparata code (see \cite{Preparata} and
\cite{Lint}), which consist of the vertices of the first color in a perfect
coloring of the Boolean $n$-cube with quotient matrix $\begin{pmatrix}
0 & n & 0 & 0 \\
1 & 0 & n-1 & 0\\
0 & 2 & n-3 & 1\\
0 & 0 & n & 0
\end{pmatrix}$. Preparata codes exist in Boolean $n$-cubes  if
 $n=2^{2t}-1$ and $t$ is an integer.
\end{example}}

{
\begin{example}\label{ex:face1}
Another well-known example of completely regular code is constructed in \cite{CCZ}. It consists of the vertices of the first color in a perfect
coloring of the Boolean $n$-cube with quotient matrix $\begin{pmatrix}
0 & n & 0 & 0 \\
1 & 0 & n-1 & 0\\
0 & 2 &\frac{n-7}{2} & \frac{n+3}{2}\\
0 & 0 & \frac{n-1}{2} & \frac{n+1}{2}
\end{pmatrix}$. Preparata codes exist in Boolean $n$-cubes  if
 $n=2^{2t-1}-1$ and $t$ is an integer.
\end{example}}

Examples of completely regular in Johnson graphs are available in \cite{AM10} and \cite{AM11}.

\begin{claim}\label{c:compl_reg_code_color}
Let $G$ be a distance-regular graph and let $C$ be  a completely
regular code in $G$ with covering radius $d$. Suppose that   ${\mathbf 1}_{C}=\sum\limits_{i=1}^n a_i f_i$, where  $f_i$ are  eigenfunctions of $G$.  Then, for any $t$,
$t\leq d$, ${\mathbf 1}_{C_t}$ is a linear combination of the same eigenfunctions, i.e.,
${\mathbf 1}_{C_t}=\sum\limits_{i=1}^n b_i f_i$.
\end{claim}
\begin{proof}
The proposition is true if $t=0$.  Suppose it is true for $j<t$. By the definition of  completely
regular codes, we have equality  $$M_t{\mathbf 1}_{C}=\sum\limits_{j=0}^{t-1}p_{jt}{\mathbf 1}_{C_j}+p_{tt}{\mathbf 1}_{C_t},$$
where $p_{tt}>0$.
By Proposition \ref{c:WDB_bound0} we have equalities
$M_tf_i=P_t(M)f_i=P_t(\alpha_i)f_i$. Then by the induction hypothesis  we obtain the required result. \end{proof}

Let $G$ be a regular graph and $C\subset V(G)$
a completely regular code. Consider the distance refinement of vertices with respect to
$C$: $C_j=\{x\in V(G) \ |\ \min_{y\in C}
d(x,y)=j\}$. Clearly, $C_0=C$. By the definition of a completely regular
code, $C_j$ is the set of vertices of the $j$th color in some perfect
coloring $h$ with quotient matrix $S_{_C}$. Let
$f:V(G)\rightarrow \{1,\dots,k\}$
be another perfect coloring of $G$. We will call the weight
distribution of a coloring $f$ with respect to $C$ the sets of vectors
$w_0=(|C_0\cup f^{-1}(1)|,\dots, |C_0\cup f^{-1}(k)|),\dots,
w_j=(|C_j\cup f^{-1}(1)|,\dots, |C_j\cup f^{-1}(k)|),\dots$.}

\begin{theorem}[Krotov \cite{Krotov11}]
The weight distribution $w_0,\dots,w_j,\dots$ of a perfect coloring
$f$ of  $G$ with respect to a completely regular code $C\subset
V(G)$ can be recovered from $w_0$.
\end{theorem}
\begin{proof}
Let $M$ be the adjacency matrix of graph $G$ and $S_f$ the quotient matrix of
 $f$. Then, by Proposition \ref{s:perfcol_criteria} and
the definition of a completely regular code, we have the equalities $MF=FS_f$ and
$MH=HS_{_C}$, where $H_{ij}=1$ only if the $i$-th vertex of $G$
is at distance $j$ from  $C$.

Consider the matrix $W=F^{\mathsf T}H$. By definition, its columns
are  $w_j$. We have the equalities $S_f^{\mathsf T}F^{\mathsf T}H=F^{\mathsf T}M^{\mathsf T}H=F^{\mathsf T}MH=F^{\mathsf T}HS_{_C}$. By the condition, we know the first (with
distance $0$ from the code $C$) column of the matrix $W$.

Since $S_{_C}$ is a distance coloring matrix, it is
tridiagonal: $$S_{_C}=\begin{pmatrix}
a_0 & b_0 & 0 & 0& \dots & 0\\
c_1 & a_1 & b_1 & 0 & \dots & 0\\
0 & c_2 & a_2 & b_2 & \dots & 0 \\
\dots & \dots & \dots & \dots & \dots & \dots \\
0 & \dots & 0 & 0 & c_d & a_d
\end{pmatrix},$$ where $d$ is the maximum possible distance from the code
$C$ to an arbitrary vertex in $G$. Moreover, $c_i\neq 0$ and
$b_i\neq 0$ for all $i$, since some vertices of consistent layers in
the distance-from-code refinement should be adjacent. We reconstruct the  columns $w_1,w_2,\dots,w_d$ of  $W$ recursively. From
the equation $WS_{_C}=S_f^{\mathsf T}W$, we have
$a_0w_0+c_1w_1=S_f^{\mathsf T}w_0$, i.e., $w_1=(S_f^{\mathsf T}w_0 - a_0w_0)/c_1$. Similarly, for column $w_i$,
$i=1,\dots,d-1$, we have the equality
$b_{i-1}w_{i-1}+a_iw_i+c_{i+1}w_{i+1}=S_f^{\mathsf T}w_i$. From
these equalities, we can sequentially calculate all columns of the matrix
$W$, and hence reconstruct the weight distribution of the coloring $f$
with respect to the code $C$.
\end{proof}

Consider the Cartesian product $G_1\square G_2$, in which
$G_1$ is an $r$-regular graph, and
$G_2$ is a distance-regular graph with quotient matrix
$S[G_2]$. Consider the vertex set $C=V(G_1)\times \{a\}$, $a\in
G_2$.

\begin{claim}
The distance coloring $f$ of the graph $G_1\square G_2$ with respect to 
$C$ ($f(x,y)=d(y,a)$) is perfect with quotient matrix
$S[G_2]+rI$. Thus, the code $C$ is completely regular in
 $G_1\square G_2$.
\end{claim}
\begin{proof}
Let a vertex $y\in V(G_2)$ be at distance $i$ from a vertex
$a$. By the definition of a distance-regular graph,
it is adjacent to vertices $c_i,a_i,b_i$ at distances $i-1,i,i+1$ from
$a$, respectively, where  $c_i,a_i,b_i$ are from the intersection array of $G_2$. By the definition of the Cartesian product,
a vertex $(x,y)$ of graph $G_1\square G_2$ has $c_i$ neighbors of color $i-1$, $b_i$ neighbors of color
$i+1$, and $a_i+r$ neighbors of color $i$.
\end{proof}

For the code $C\subset V(G)$, we define $s^*(C)$ such that $s^*(C)+1$
is the minimum number of eigenfunctions of the graph $G$, a linear
combination of which is the characteristic function ${\mathbf 1}_C$.

\begin{theorem}[\cite{BGKM}]
Let $G$ be a distance-regular graph of diameter $d$. Then\\
{\rm (a)} for any code $C$, the equality ${\rm
diam}(C)+s^*(C)\geq d$ holds;\\
{\rm (b)} if ${\rm diam}(C)+s^*(C)= d$, then $C$ is completely
regular.
\end{theorem}
\begin{proof}
Let $M_t$ be the adjacency matrix of $G$ over distance $t$.
Consider a matrix $B$ of size $(d+1)\times V(G)$, whose rows
are the vectors $M_t{\mathbf 1}_{C}$, $t=0,\dots,d$. Let
$C_{0},\dots,C_{\rho}$, $d-{\rm diam}(C)\leq \rho\leq d$, be the distance refinement of vertices of graph $G$
with respect to code $C$.

(a) From the Proposition \ref{c:WDB_bound0}, the row rank
of matrix $B$ is $s^*(C)+1$. From the definition of the diameter of a code, it follows
that in the columns of matrix $B$ corresponding to vertices from $C_i$, only positions from $i$ to ${\rm diam}(C)+i$ are nonzero. After ordering the columns according to the partition $C_{0},\dots,C_{\rho}$ the matrix $B$ contains upper-triangle submatrix of order $d+1-{\rm diam}(C)$.
Therefore, the column rank of matrix $B$ is at least $d+1-{\rm diam}(C)$.

(b) If ${\rm diam}(C)+s^*(C)= d$, then the rank of matrix $B$ is equal to
the number of distinct supports of the columns. Hence, we conclude that columns with
the same supports are collinear. Since the sum of the elements in
each column is $|C|$, columns with the same supports
coincide. By Theorem \ref{thcritDRG} this condition is equivalent to the definition
of a completely regular code.
\end{proof}

\subsection{Delsarte Cliques in Distance-Regular Graphs}

Let $\lambda$ be an eigenvalue of a distance-regular graph
$G$ of diameter $d$. Let $u_1,\dots,u_m$ be an orthonormal basis
of the eigenspace $X_{\lambda}$ of $G$.
Consider the vectors
$w^\lambda(x)=(u_1(x),\dots,u_m(x))$, where $x\in V(G)$. Their scalar
products are the elements of the matrix $E_{\lambda}$ defined in Section \ref{2.3}.

From Corollary \ref{colE} and Proposition \ref{remDR} we have

\begin{claim}
If $G$ is a distance-regular graph of diameter $d$, then the matrices
$E_\lambda$ can be represented as linear combinations of the distance-$t$ adjacency matrices of $G$,
where $t=0,\dots,d$.
\end{claim}

\begin{corollary}\label{corElambda}
If $G$ is a distance-regular graph of diameter $d$, then the element
$E_\lambda(x,y)$ depends only on the distance between vertices $x$ and
$y$ in $G$.
\end{corollary}

We define the vector $\gamma(\lambda)=(\gamma_0,\dots,\gamma_d)$, where
$\gamma_t=E_\lambda(x,y)$ if $d(x,y)=t$. Fix an arbitrary
vertex $x_0\in V(G)$ as the initial vertex and the distance refinement
$[A_0,\dots,A_d]$ of $G$ with respect to it. Recall the notation
from Section \ref{4.4}: $z_t(\lambda)=P_t(\lambda)/s_t$,
$t=0,\dots,d$, where $s_t=|A_t|$.

\begin{claim}\label{c:ElemE}
Let $G$ be a distance-regular graph of diameter $d$. Then
the vectors $z(\lambda)$ and $\gamma(\lambda)$ are collinear.
\end{claim}
\begin{proof}
Consider an orthonormal basis
of the eigenspace $X_{\lambda}$ of $G$ corresponding
an eigenvalue $\lambda$. By the definition of an eigenfunction, for each $x\in A_t$, we have the equality $$\lambda u_i(x)=
\sum\limits_{y: y\in A_{t-1},d(x,y)=1}u_i(y)+ \sum\limits_{y: y\in
A_{t},d(x,y)=1}u_i(y)+\sum\limits_{y: y\in
A_{t+1},d(x,y)=1}u_i(y).$$ Then

\begin{equation}\label{e:ElemE1}\lambda w^\lambda(x)= \sum\limits_{y: y\in
A_{t-1},d(x,y)=1}w^\lambda(y)+ \sum\limits_{y: y\in
A_{t},d(x,y)=1}w^\lambda(y)+\sum\limits_{y: y\in
A_{t+1},d(x,y)=1}w^\lambda(y).\end{equation}

Consider the scalar product of both sides of (\ref{e:ElemE1}) by the vector
$w^\lambda(x_0)$.  We obtain the equality
\begin{equation}\label{eqdopDC}
\lambda
\gamma_t(\lambda)=c_t\gamma_{t-1}(\lambda)+a_t\gamma_t(\lambda)+b_t\gamma_{t+1}(\lambda),
\end{equation}
where $a_t,c_t,b_t$ are elements of  the intersection array of $G$.

Then the vector $\gamma(\lambda)$ is an eigenvector of the intersection matrix $S[G]$ of  $G$ with eigenvalue $\lambda$.  By
Remark \ref{remSG}, every eigenspace of  $S[G]$ is one-dimensional.  By Proposition \ref{c:polyKrv}, the vector $z(\lambda)$
is contained in the same eigenspace of $S[G]$.  \end{proof}

Next we prove an analog of Theorem \ref{PCclaim13} for a distance-regular
graph.

\begin{theorem}[on the Delsarte clique \cite{Godsil}]\label{PCclaim130}
Let $G$ be an $r$-regular graph. If $G$ is a distance-regular  graph  then the cardinality of any clique in $G$
does not exceed $1+\frac{r}{-\lambda_{min}}$, where $\lambda_{min}$ is the minimum
eigenvalue of $G$. \end{theorem}
\begin{proof}
Let $u_1,\dots,u_m$ be an orthonormal basis of the eigensubspace $X_{\lambda_{min}}$.
Consider the vectors
$w(x)=(u_1(x),\dots,u_m(x))$, where $x\in V(G)$. Their scalar
products are elements of the matrix $E_{\lambda_{min}}$. By
Corollary \ref{corElambda}, their scalar products are the same for
all pairs of adjacent vertices $x$ and $y$, as well as for pairs of coinciding
vertices. Let $(w(x),w(y))=\beta$ for any pair of adjacent vertices $x,y\in
V(G)$ and $(w(x),w(x))=\alpha$ for any vertex $x\in V(G)$.

From the definition of an eigenfunction, we have the equality\\ $
\sum\limits_{y\in A_1(x)}u_i(y)=\lambda_{min} u_i(x),$ i.e., $
\sum\limits_{y\in A_1(x)}w(y)=\lambda_{min} w(x)$. The scalar product of both sides of
the previous equality by $w(x)$ gives that
$r\beta=\lambda_{min}\alpha$.

Let $C\subset V(G)$ be a clique in  $G$. Consider the Gram matrix
$W$ of size $|C|\times |C|$, whose elements are
$(w(x),w(y))$, $x,y\in C$. We have \\
$W=\begin{pmatrix}
\alpha & \beta & \beta & \dots & \beta\\
\beta & \alpha & \beta& \dots & \beta\\
\dots & \dots & \dots & \dots & \dots \\
\beta & \dots & \dots & \beta & \alpha
\end{pmatrix}$.

The Gram matrix is positive semidefinite, so
$(W\overline{1},\overline{1})\geq 0$, therefore, $\alpha+
(|C|-1)\beta\geq 0$. Since the minimum eigenvalue
$\lambda_{min}$ is negative (see Problem \ref{exe1}), this and the equality $r\beta=\lambda_{min}\alpha$ yield $|C|\leq
1+\frac{r}{-\lambda_{min}}$.
\end{proof}

Recall that a clique of cardinality $1-\frac{r}{\lambda_{\min}}$ in an
$r$-regular graph $G$ is called a Delsarte clique.

\begin{theorem}[\cite{Godsil}]\label{PCclaim131}
Every Delsarte clique in a distance-regular graph of diameter $d$
is a completely regular code with covering radius $d-1$.
\end{theorem}
\begin{proof}
Let $C$ be a Delsarte clique and $C_t$ be the set of vertices at
distance $t$ from $C=C_0$. By Theorem \ref{thcritDRG}, it suffices to verify that for any $t$ and $i\geq t$, each vertex $x\in C_t$
has the same number of vertices from $C$ at distance $i$. We need to check the number of vertices only at distance $t$,
since the definition of a clique implies that the distance can only be $t$ or
$t+1$.

Below, we will use the notation  of Theorem
\ref{PCclaim130}. From the proof of Theorem \ref{PCclaim130}, it is clear
that for every Delsarte clique, the equality $\alpha+ (|C|-1)\beta=
0$ holds, i.e., $(\sum\limits_{x\in C} w(x), \sum\limits_{x\in C}
w(x))= 0$. Consequently, $\sum\limits_{x\in C} w(x)=\bar 0$.

By definition, $E_{\lambda_{min}}(x,y)=(w(x),w(y))$. From Corollary
\ref{corElambda}, we see that $(w(x),w(y))$ depends only on
the distance between vertices $x$ and $y$. Let $(w(x),w(y))=\gamma_t$
if $d(x,y)=t$. Consider $y\in C_t$. We have equality
$$0=(w(y),\sum\limits_{x\in C} w(x))= (w(y),\sum\limits_{x\in C,
d(x,y)=t} w(x))+ (w(y),\sum\limits_{x\in C, d(x,y)=t+1} w(x))=$$
$$\hfill = \gamma_tk_t(y)+ \gamma_{t+1}(|C|-k_t(y))=(\gamma_t-\gamma_{t+1})k_t(y)+ \gamma_{t+1}|C|,$$ where $k_t(y)=|\{x\in
C, d(x,y)=t\}|$. If $\gamma_t-\gamma_{t+1}=0$, then $\gamma_{t+1}=0$ and $\gamma_{t}=0$.
By using (\ref{eqdopDC}) and induction, we obtain  $\gamma_{i}=0$ for all $i<t$. It contradicts  Proposition \ref{c:ElemE} ($P_0(\lambda)=1$). Therefore, 
$\gamma_t-\gamma_{t+1}\neq 0$ for $t=0,\dots,d-1$. Then $k_t(y)=
\frac{\gamma_{t+1}|C|}{\gamma_{t+1}-\gamma_{t}}$ is independent of $y$.
Furthermore, the graph has no vertices at distance $d+1$. Therefore,
$\gamma_dk_d(y)=0$.
\end{proof}

\begin{corollary}
Let $C$ be a Delsarte clique in a distance-regular graph $G$
with minimal eigenvalue $\lambda_{min}$. Then the function
${\mathbf 1}_{C}$ is orthogonal to the eigenspace
$X_{\lambda_{min}}$ of $G$.
\end{corollary}
\begin{proof}
Since $u_1,\dots,u_m$ is an orthonormal basis of the eigenspace
$X_{\lambda_{min}}$, we have
$({\mathbf 1}_{C},u_i)=\sum\limits_{x\in C} (w(x))_i$. Then the statement 
is equivalent to the equality $\sum\limits_{x\in C} w(x)=\bar 0$.
\end{proof}

\begin{example}
A maximal clique in $Q^n_q$ is isomorphic to the complete graph $K_q=Q_q$. As mentioned above the minimum eigenvalue of $Q^n_q$  is $-n$. Since $|Q_q|=q=1+\frac{(q-1)n}{n}$, every maximal cliques in  $Q^n_q$ is a Delsarte clique. 
\end{example}

\subsection{Middle Layer of Boolean Hypercube as a Testing Set}

%\begin{definition}
A graph $G$ of diameter $d$ is called {\sl antipodal} if for any
vertex $x\in V(G)$, there is a unique
vertex $y\in V(G)$ at maximum distance $d$ from $x$, called the antipode  (or opposite vertex) of $x$.
%\end{definition}
For example, a Boolean $n$-cube is antipodal since at
distance $n$ from a vertex $x$ there is a unique vertex
$x\oplus \bar{1}$. Note that $q$-ary hypercubes for $q>2$ are not
antipodal graphs.

\begin{claim}\label{c:27}
For a perfect coloring in an antipodal distance-regular
graph, the color of any vertex is uniquely determined by the color
of its opposite vertex. \end{claim}
\begin{proof}
By Theorem \ref{th:shapiro_zlotnik_martin}, a perfect vertex coloring
of a graph remains perfect in a graph with a distance equal to the diameter
$d$ of the graph. In an antipodal graph, there is exactly one vertex at distance $d$ from each vertex, i.e., the  distance-$d$ graph is a perfect matching, and the color of a vertex is determined by the color of its neighbor.
\end{proof}

Recall that the number of units in $n$-tuple $x\in Q^n_2$ is called the {\sl weigh} of $x$ and denoted by $\wt(x)$. $\{x\in  Q^n_2 : \wt(x)=t\}$ is a $t$th layer of  $Q^n_2$.  The  layer of weight $\frac{n}{2}$ (if $n$ is even) and  
layers of weight $\frac{n-1}{2}$ and  $\frac{n+1}{2}$ (if $n$ is odd) are called the {\sl middle} layers.

\begin{theorem}[\cite{Avgust}]\label{thAvg}
Let $n$ be odd, $f_1$ and $f_2$ be perfect colorings
of $Q^n_2$ with the same quotient matrix, and for any
vertex $x\in Q^n_2$ of weight $\wt(x)=\frac{n-1}{2}$ (or
$\wt(x)=\frac{n+1}{2}$), we have $f_1(x)=f_2(x)$. Then $f_1=
f_2$.
\end{theorem}
\begin{proof}
Vertices opposite to vertices of weight $\wt(x)=\frac{n-1}{2}$ have
weight $\frac{n+1}{2}$. Therefore, by Proposition \ref{c:27} it follows that
there is a correspondence $\varphi$ between the colors such that $f_1(x\oplus
\bar{1})=\varphi(f_1(x))=\varphi(f_2(x))=f_2(x\oplus \bar{1})$, where
$i=1,2$, $\wt(x)=\frac{n+1}{2}$ (or $\wt(x)=\frac{n-1}{2}$).

Consider the coloring $f(x)=\left\{
\begin{array}{ll}
f_1(x)\mbox{,} & \wt(x)\le\frac{n+1}{2}\\
f_2(x)\mbox{,} & \wt(x)\geq \frac{n-1}{2}
\end{array}
\right. $. This definition of $f$ is correct by the previous arguments. The coloring $f$ is perfect with the same quotient matrix
as $f_1$ and $f_2$, since in every ball of radius $1$ around a vertex of weight $t$ involves only vertices of weights $t-1$ and $t+1$. So, on every  radius $1$ ball $f$  coincides with
$f_1$ or $f_2$. From Proposition \ref{c:27} we have
$f_1(x)=f(x)=\varphi(f(x\oplus \bar{1}))=\varphi(f_2(x\oplus
\bar{1}))=f_2(x)$ and for every $x\in Q^n_2$ with weight $\wt(x)\le\frac{n-1}{2}$. Therefore,
$f_2(x)=f_1(x)$ for every $x\in Q^n_2$.
\end{proof}

%\begin{definition}
A {\sl testing set} for a function from a certain class is a
subset of its domain such that the function can be completely reconstructed if one know values of the function on this subset.
%\end{definition}

Here we are discussing only "unconditional" testing sets, when
the  set should not depend on the function. Sometimes, for example if we need better upper bound on the number of some combinatorial configurations, it is useful
to consider "conditional" testing sets, when the choice of the next
testing argument depends on the function value at previous testing
arguments.

By Theorem  \ref{thAvg} the middle layer  is a testing set for   perfect colorings of the Boolean $n$-cube when $n$ is odd. Moreover, Theorem
\ref{thAvg} provides best known upper bound on the number of $1$-perfect
codes in the Boolean $n$-cube. Since a $1$-perfect code is uniquely determined by its intersection with the middle layer, the number of $1$-perfect codes is at most $2^{n \choose (n-1)/2}$. Testing sets for  perfect colorings of nonbinary hypercubes will be considered in Section \ref{secTest}.

\subsection{Problems}

{\begin{exercise} Prove that for the Krawtchouk-type polynomials $P_t$ of a distance-regular graph $G$ the equality
$(x-b_0)\sum\limits_{t=0}^sP_t(x)=c_{s+1}P_{s+1}(x)-b_sP_s(x)$ holds, if
$s+1$ does not exceed the diameter of $G$.
\end{exercise}}

{\begin{exercise}\label{exer} Prove the equality\\
{ $2n\sum\limits_{t=0}^s(-1)^{t+1}{n \choose
t}=(-1)^{s+1}\left((s+1){n \choose {s+1}}+(n-s){n \choose
{s}}\right)$.}
\end{exercise}}

\begin{exercise} Prove that if a graph is the complement of a strongly regular
graph and it is connected, then it is strongly regular.
\end{exercise}

{\begin{exercise} Describe all graphs with exactly two eigenvalues
of its adjacency matrix.
\end{exercise}}

\begin{exercise}
Prove that for the quotient matrix $S[G]$ of each
distance-regular graph $G$, the sequence $b_i$ is decreasing
(non-strictly), while the sequence $c_i$ is increasing (non-strictly).
\end{exercise}

{
\begin{exercise}
Prove that for every distance-regular graph, the equality $|A_2|P_1(-1)=|A_1|P_2(-1)$ holds, where $A_i$ is a sphere of radius $i$.
\end{exercise}}
\begin{exercise}
Prove that for every distance-regular graph, the number
of positive values of an eigenfunction with eigenvalue
$\alpha$ is not less than $\sum\limits_{t=0}^d \max\{P_t(\alpha),0\}$.
\end{exercise}

\begin{exercise}
Let $C$ be a completely regular code in a distance-regular
graph. Prove that the quotient matrix of distance coloring 
with respect to the code $C$ has eigenvalues that coincide with
the eigenvalues of the minimal set of eigensubspaces
containing the function ${\mathbf 1}_C$.
\end{exercise}

{
\begin{exercise}
Let $f$ be a perfect coloring of the Boolean cube. Faces of the Boolean $n$-cube
that have the same set of fixed coordinates but completely different
fixed coordinate values are called opposite.
Prove that if two faces of a Boolean $n$-cube have the same color composition
(i.e.,  equal numbers of vertices of each color in both faces), then
the color composition of the opposite  faces  them also match.
\end{exercise} }

{
\begin{exercise}
Let $C\subset Q^n_2$ ($n$ is odd) be a completely regular
code. Prove that  $C$ is completely recoverable from the set
of its vertices of weight $\frac{n+1}{2}, \frac{n+1}{2}+1,\dots,
\frac{n+1}{2}+r$, where $r$ is the covering radius of $C$.
\end{exercise} }

\begin{exercise}
Let $G$ be a distance-regular graph and the characteristic
function ${\mathbf 1}_{C}$ of the set $C\subset V(G)$ is
the sum of $k$ eigenfunctions of $G$ with $k$ distinct
eigenvalues. Prove that 1) if the set $C$
is a color of some perfect coloring, then it has at least
$k$ colors; 2) if after $k-1$ steps of the WL algorithm
the function ${\mathbf 1}_{C}$ yields a coloring with $k$ colors, then
the resulting coloring is perfect. \end{exercise}

\begin{exercise}
Let $G$ be a distance-regular graph and let $C$ is a color of some equitable partition of $G$.
Prove that for two vertices of any other color $C'$, the number of paths
of arbitrary fixed length to vertices from $C'$
is the same. Give an example showing that the coincidence of the path lengths
to the vertices of $C$ is not a sufficient condition for the WL algorithm, starting from a $2$-color partition 
$\{C,V(G)\setminus C\}$, to color these vertices by the same color.
\end{exercise}

\begin{exercise}
Let $P$ be  an equitable partition of an antipodal distance-regular graph. Prove that for every cell  $P_i\in P$ there exists $P_j\in P$ such that  if $x\in P_i\in P$ and vertex $y$ is antipodal to $x$ then $y\in P_j$. It is possible that $i$ is equal to $j$.
\end{exercise}

\section{Automorphisms of Graphs and Perfect Colorings}\label{Orbit}

\subsection{Orbit Colorings}

Let $G$ be a graph and let $H$ be a subgroup of its 
  automorphism group, i.e., $H\leq Aut(G)$.
%\begin{definition}
An {\sl orbit} with respect to $H$ is a minimal set with respect to
inclusion $O\subset V(G)$ such that $\pi (O)=O$ for
any automorphism $\pi \in H$.
%\end{definition}

Since $H$ is a subgroup, each orbit $O$ coincides with the image
 $$H(v)=\{\pi(v)\
|\ \pi \in H\}$$ of any element $v\in O$. Indeed, if $v',v\in O$ then $v'=\pi_1\circ \dots\circ \pi_k(v)$,  where $\pi_1,\dots,\pi_k\in H$ and, consequently, $\pi_1\circ \dots\circ \pi_k\in H$. So, the vertex set
of a graph is partitioned into disjoint orbits.

\begin{claim}\label{claim:orbit}
Let $P=\{O_1,\dots,O_k\}$ be a partition of  $V(G)$ into orbits
with respect to the subgroup $H$. Then partition $P$ is equitable.
\end{claim}
\begin{proof}
Let elements $v$ and $v'$ be contained in the same orbit. Then $v'=\pi(v)$ for some automorphism $\pi\in
H$.  Consider an arbitrary vertex $u$ adjacent to $v$.
Then, by the definition of an automorphism, the vertex $u'=\pi(u)$ is adjacent to
the vertex $\pi(v)$. Therefore, vertices $u$ and $u'$ belong
to the same orbit. The mapping $\pi$ is bijective,
therefore, the color composition of the neighborhoods of vertices $v$ and $v'$ coincide.
\end{proof}

{ The perfect coloring defined in Proposition \ref{claim:orbit}
will be called an {\sl orbit coloring}.}

\begin{corollary}
Let $C\subset V(G)$ be a transitive code in graph $G$. Then $C$
is a color of some orbit coloring.
\end{corollary}

{
\begin{example}\label{ex:circ}
Consider a graph whose vertices are the numbers $0,\dots,8$, and
vertices $x$ and $y$ are connected by an edge if $(x-y)\mod 9\in
\{\pm1,\pm2\}$. It is easy to see that the mapping
$\varphi(x)=(x+3)\mod 9$ is an automorphism of the graph, and generates the
subgroup
$H=\{\varepsilon,\varphi,\varphi^2\}$.  Here and throughout, $\varepsilon$ is the identity
element of the group (the identity permutation). With respect to $H$, the graph
is partitioned into three orbits, according to the residue modulo
$3$. The coloring $h(x)=x\mod 3$ is perfect with quotient matrix
 $\begin{pmatrix}
0 & 2 & 2 \\
2 & 0 & 2 \\
2 & 2 &0
\end{pmatrix}$.
\end{example}}

Describing all orbital colorings reduces to finding the group
of graph automorphisms. Isometries (automorphisms) of hypercubes are described by Markov's theorem. To state it, we need some
definitions.

%\begin{definition} An {\sl isotopy} is an isometry of  $Q^n_q$,
defined by the equality
$\overline{\sigma}(x_1,x_2,\dots,x_n)=(\sigma_1x_1,\sigma_2x_2,\dots,\sigma_nx_n)$,
where each $\sigma_i$ is a permutation of the elements of the set $Q_q$.
%\end{definition}

%\begin{definition}
A {\sl parastrophy} is an isometry of  $Q^n_q$ defined by
the equality
${\tau}(x_1,x_2,\dots,x_n)=(x_{\tau1},x_{\tau2},\dots,x_{\tau n})$,
where $\tau$ is a permutation of the set $\{1,\dots,n\}$.
%\end{definition}

\begin{theorem}[Markov \cite{Markov}, \cite{Miraf}]\label{Markov}
Any isometry of a hypercube is the composition of an isotopy and a
parastrophy. More precisely, the full automorphism group of the Hamming graph $H(n, q)$  is the wreath product of the symmetric group $S_{q}$ by the symmetric group $S_{n}$, which is isomorphic to the semidirect product ${\mathbb Z}_2^n\rtimes S_n$ for $q=2$.
\end{theorem}

Some isometries of $Q^n_q$ can be viewed as affine
transformations in a vector space over a finite field $GF(q)$, if
$q$
is a prime power. In this case we can identify $Q_q$ and $GF(q)$ and regard $Q^n_q$ as an $n$-dimensional  
vector space over $GF(q)$. Consider a linear transformation
of the hypercube whose matrix has exactly one $1$ in each
row and column, and the remaining elements are zero. Clearly,
such  a linear transformation is a parastrophy. Consider a
translation map that adds a fixed vector to each vertex of the hypercube. Clearly, such  affine transformation
is an isotopy. It can be shown (see Problem \ref{exer6}) that for
$q=2,3$, all isometries of the hypercubes $Q^n_q$ are affine
transformations.

A group action on a graph is called {\sl regular}
if for every ordered pair of vertices of the graph there is a unique element of the group that maps the first vertex to the second.
{A graph on $q$ vertices such that a cyclic subgroup
$\mathbb{Z}_q=\mathbb{Z}/q\mathbb{Z}$ of its automorphism group acting
regularly
on its vertices, is called  {\sl circulant}. The graph considered in
Example \ref{ex:circ} is  circulant.}

A description of perfect colorings for some circulant graphs with two colors
can be found in the papers \cite{Horosh} and \cite{ParLis}.

A vertex-transitive
graph has exactly one orbit with respect to its
automorphism group. Transitive graphs with a subgroup-rich
automorphism group admit many orbit colorings.

\subsection{Perfect Colorings of Cayley Graphs}\label{Orbit}

The main examples
of transitive graphs are Cayley graphs.

%\begin{definition}
 Let $Gr$ be a group and let $A\subset Gr$ be a {\sl set of generators}
such that $A^{-1}=A$ (otherwise we obtain a directed graph),
$\varepsilon\not\in A$ (otherwise we get loops). A {\sl Cayley graph} $Cay(Gr,A)$
is a graph whose vertices are elements of the group, and two
vertices $x,y\in Gr$ are adjacent by an edge when $y=xa$ for some
$a\in A$. 
%\end{definition}

%\begin{remark}
The elements  $g\in Gr$ define automorphisms $\varphi_g$ of the graph
$Cay(Gr,A)$ by the rule $\varphi_g(x)=gx$. The group $Gr$ is a
subgroup of the automorphism group of  $Cay(Gr,A)$. Thus,
every Cayley graph is vertex-transitive. For each generator set $A$, the Cayley graph
$Cay(\mathbb{Z}_q,A)$ is circulant.
%\end{remark}

\begin{minipage}{0.4\textwidth}
%\vskip15mm
\includegraphics[width=0.7\textwidth]{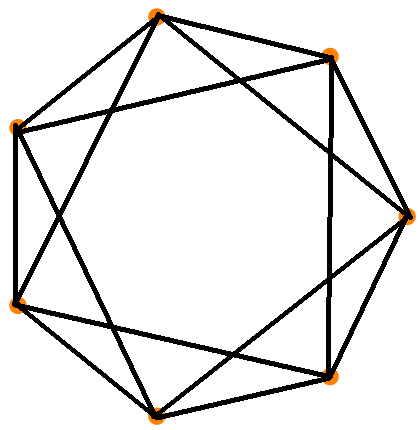}
% \center{\includegraphics[width=0.07\textwidth]{E3.png}}
\end{minipage}

The figure shows the graph $Cay(\mathbb{Z}_7,A)$, where $A=\{\pm 1, \pm
2\}$.

It is easy to show that the Cayley graph of the Cartesian square of $Gr\times
Gr$ with generator set $A\times \{\varepsilon\}\cup
\{\varepsilon\}\times A$ is the Cartesian product of the Cayley graphs
$Cay(Gr,A)\square Cay(Gr,A)$.
%\begin{remark}
The complete graph $Q_q$ can be viewed as the Cayley graph of the group
$\mathbb{Z}_q$ with generator set consisting of all nonzero
elements of the group. The hypercube $Q^n_q$ can be viewed as the Cayley graph
$Cay(\mathbb{Z}_q^n,A)$, where $A$
is the set of $n$-tuples containing exactly one nonzero element. 
%\end{remark}

{ A special case of the Cayley graph arises when 
the generator set  $A=A^{-1}$ as the set of representatives of cosets
of some subgroup $H$, i.e., $A\simeq(Gr/H)\backslash
\{\varepsilon\}$. Then, any element $g\in Gr$ of the group is uniquely
representable as $g=ha$ or $g=h\in H$, i.e., each
vertex of the graph $Cay(Gr,A)$ is  adjacent to exactly one element of  $H$,
or itself belongs to the subgroup $H$. Thus, the set $H$
forms a $1$-perfect code in $Cay(Gr,A)$. As proved in the statement below, $H$ can be any
normal subgroup. A {\sl normal} subgroup is a subgroup of $H$ for which $gH=Hg$ for any $g\in Gr$.

\begin{claim}
Let $H$ be a normal subgroup of $Gr$. Then, from the cosets
of $H$, we can choose a set of representatives of $A$ satisfying
the property $A=A^{-1}$.
\end{claim}
\begin{proof}
Since $H$
is a normal subgroup, the equality $Ha=Hb$ implies the equality
$Ha^{-1}=Hb^{-1}$ ($b=ha\Rightarrow b^{-1}=a^{-1}h^{-1}=h'a^{-1}$),
which means that elements of $A$ can be chosen in pairs
$\{a,a^{-1}\}$ when $a$ and $a^{-1}$ belong to different cosets
of $H$. If a pair of mutually inverse elements $a$ and
$a^{-1}$ belongs to the same coset, then it contains an element
of order $b^2= \varepsilon$, which can be chosen in the set $A$. Indeed, let $a^{-1}=ha$
for some $h\in H$, i.e., $\varepsilon=ha^2$. Consequently,
$a^2$ and all even powers of $a$ are contained in $H$, and all
odd powers are contained in the coset $Ha$. Since the group
$H$ is finite, there exists a minimal integer $s>0$ such that
$a^{2s}=\varepsilon$. Moreover, the parity of $s$ contradicts
the condition that $a$ and $a^{-1}$ belong to different cosets. Therefore, $a^s\in Ha$. \end{proof}

\begin{example}
Consider the symmetric group $S_n$ and a subgroup $H$ consisting of
permutations fixing some element. Without loss
of generality, we can assume that $\pi(1)=1$ for all $\pi\in H$. Let
$A=(S_n/H)\backslash \{\varepsilon\}$. It is easy to see that
$A=\{(1i) \ |\ i=2,\dots, n\}$ and $H\simeq S_{n-1}$. The graph
$Cay(S_n,A)$ is called a {\sl star graph}. The subgroup $H$ forms a $1$-perfect code in the
star graph. 
\end{example}

Other sufficient conditions for a subgroup of a given group to be a  perfect code in a Cayley graph are studied in \cite{Huang}. Examples of non-group
$1$-perfect codes in star graphs can be found in \cite{Mog20}.

A {\sl homomorphism}  is a mapping $\varphi:  Gr_1 \rightarrow Gr_2$,
satisfying the equality $\varphi(x*y)=\varphi(x)\circ\varphi(y)$
for any elements $x,y\in Gr_1$. Here $*$ and $\circ$ are group
operations in $Gr_1$ and $Gr_2$, respectively.

\begin{claim}
Let $Gr_1$ and $Gr_2$ be some groups, and there is a
homomorphism $\varphi:Gr_1\rightarrow
Gr_2$, which is bijective on the generator set
$A\subset Gr_1$. Then $\varphi$ is a covering
of the graph $Cay(Gr_2,\varphi(A))$ by the graph $Cay(Gr_1,A)$.
\end{claim}
\begin{proof}
Let $g\in Gr_1$. Then $\{ga : a\in A\}=S^1_1(g)$. From
the definition of homomorphism, we have $\varphi(ga)=\varphi(g)\varphi(a)$.
Thus, $\varphi(S^1_1(g))=S^2_1(\varphi(g))$, and the mapping
$\varphi$ is bijective on $S^1_1(g)$ for any $g\in Gr_1$ by
the assumption.
\end{proof}

The mapping $f:Gr_1\times Gr_2\rightarrow Gr_1$, defined
by $f((a,b))=a$
is a natural
homomorphism from the direct product $Gr_1\times Gr_2$ to
a factor. Let us consider other examples.

\begin{example}\label{ex:coverZ}
The mapping $\varphi(n)= n\,(\mbox{mod }4)$ is a homomorphism
of the group $\mathbb{Z}$ to the group $\mathbb{Z}_4$. Then the mapping
$\varphi: Cay(\mathbb{Z},\pm1)\rightarrow Cay(\mathbb{Z}_4,\pm 1)$
is a covering.
\end{example}

%\begin{example}
The groups $\mathbb{Z}_2\times \mathbb{Z}_2$ and $\mathbb{Z}_4$ are not
isomorphic. However, their Cayley graphs (with any two generators)
are both  isomorphic to the square $Q_2\square Q_2$. The isomorphism
$\iota:Cay(\mathbb{Z}_2\times \mathbb{Z}_2, \{(0,1),(1,0)\})
\leftrightarrow Cay(\mathbb{Z}_4,\pm1)$ such that
$\iota((0,0))=0$ is called the {\sl Gray map}.
%\end{example}

Consider the group $\mathbb{Z}^n_4$. Let $e_i=(0,\dots,\underset
i{1},\dots,0)$. We have graph isomorphisms $Cay(\mathbb{Z}^n_4,\{\pm
e_i\})\leftrightarrow Cay(\mathbb{Z}^{2n}_2,\{e_i\})\leftrightarrow
Q_2^{2n}$. From Proposition \ref{claim:cartprodcover} and Example
\ref{ex:coverZ} it follows that there exists a covering of the graph $Q_2^{2n}$
by the graph $Cay(\mathbb{Z}^n,\{\pm
e_i\})$ of the $n$-dimensional rectangular grid. Consequently, from any perfect coloring of the Boolean  $2n$-cube, one can obtain a perfect coloring of the $n$-dimensional rectangular grid with the same quotient matrix.

\begin{example}\label{exSh}
Consider the covering $\psi$ of the graph $Q_5$ by a two-dimensional rectangular
grid $(\mathbb{Z}^2, \{(\pm1,0), (0,\pm1)\})$,
defined by the formula $\psi(x,y)=x+2y(\mod 5)$. Every coloring of the complete
graph $Q_5$ is perfect. Hence, the two-dimensional rectangular
grid can be colored perfectly with the following quotient matrices $\begin{pmatrix}
0 & 4 \\
1 & 3
\end{pmatrix}$,
$\begin{pmatrix}
1 & 3 \\
2 & 2
\end{pmatrix}$,
$\begin{pmatrix}
0 & 2 & 2 \\
1 & 1 & 2 \\
1 & 2 &1
\end{pmatrix}$.
\end{example}

Consider the  {\sl  Shrikhande graph} 
$$Sh=Cay(\mathbb{Z}^2_4,\{(1,0),(-1,0),(0,1),(0,-1),(1,1),(-1,-1)\}).$$
We define the mapping \[\theta: V(Cay(\mathbb{Z}^3_4,\{\pm e_i\}))\rightarrow V(Sh).\] by the equality $\theta(x_1 e_1+x_2 e_2+x_3
e_3)=x_1(1,0)+x_2(0,1)+x_3(1,1)$.
To  verify that $\theta$ is a
covering, it  suffices to consider the image of a neighborhood of one
point, since both graphs are transitive. As noted above, $Cay(\mathbb{Z}^3_4,\{\pm
e_i\})$ is isomorphic to $Q^6_2$. Therefore, there is a
covering of the Shrikhande graph by the Boolean $6$-cube. It is not difficult to find a perfect $2$-coloring of the Shrikhande graph with
quotient matrix
$\begin{pmatrix}
1 & 5 \\
3 & 3
\end{pmatrix}$.

Then  $\mathbb{Z}^6_2$ admits a perfect $2$-coloring with the same quotient matrix.
This coloring was first found by Tarannikov. {It is an example of a
coloring that attains the Fon-Der-Flaass bound
of correlation immunity (see Section \ref{11.6}).}

\begin{example}\label{ex:trans_graph}
Consider the symmetric group $S_n$ and the generator set $T_n$ consisting of all transpositions.  $Cay(S_n,T_n)$ is called the {\sl transposition graph}.
Choose an arbitrary subset $\{a_1,\dots,a_k\}$ of cardinality $k,
1<k<n,$ in the domain of  permutations. We define
a mapping $\nu: V(Cay(S_n,T_n))\rightarrow V(J(n,k))$ by the following rule.
An $n$-tuple $\nu(\sigma)$ contains $1$  exactly at the positions $\sigma(a_i)$,
$i=1,\dots,k$, $\sigma\in S_n$. It is easy to see that 
$\nu$ covers the  Johnson graph (see Example
\ref{ex:John}) with additional ${k \choose 2} + {n-k \choose 2}$
loops at each
vertex of the Johnson graph. %$\overline{J(n,k)}$
\end{example}

\subsection{Problems}

\begin{exercise}
Prove that the $k$-dimensional rectangular grid covers 
$Q^k_3$.
\end{exercise}

\begin{exercise}
Prove that the $k$-dimensional rectangular grid has a $1$-perfect
code.
\end{exercise}

\begin{exercise}\label{exer6}
Using Markov's theorem prove  that for $q=2$ and $q=3$, all
isometries of the hypercubes $Q^n_q$ are affine transformations.
\end{exercise}

\begin{exercise}
Let $q$ be a prime power and let $A$ be a matrix of size
$k\times n$ over the field $GF(q)$. Define the color class of color
$b\in Q^n_q$ as $\{x\in Q^n_q : Ax=b\}$. Prove that this
coloring of  $Q^n_q$ is perfect.
\end{exercise}

\begin{exercise}
Let $q$ be a prime number. Prove that there exists a covering $\varphi: Cay(\mathbb{Z}^{\frac{n(q-1)}{2}}_q,\{\pm e_i : i=1,\dots, (q-1)/2\})\rightarrow Q_q^n.$
\end{exercise}

\begin{exercise}\label{exerCayCol}
Consider the symmetric group $S_n$ and the generator set
$A\subset S_n$, $A^{-1}=A$, such that $|A|=n$. For each permutation $\pi\in
A$, define the permutation matrix $P_\pi$, where $P_\pi(i,j)=1$ if
$\pi(i)=j$, and $P_\pi(i,j)=0$ otherwise. Consider a
coloring of $Cay(S_n,A)$ with $n$ colors according to the rule: a permutation
$\sigma\in S_n$ has color $i$ if $\sigma(0)=i$. Prove that
the resulting coloring is perfect with quotient matrix
$S=\sum\limits_{\pi\in A}P_\pi$.
\end{exercise}

\begin{exercise}
Prove that for $n\equiv 1,3\pmod 6$ in the permutation graph
$Cay(S_n,T_n)$ (see Example  \ref{ex:trans_graph}) there exist perfect $2$-colorings 
with quotient matrix\\
$S=\left(\begin{array}{cc}
3+{n-3 \choose 2}& 3(n-3) \\
3 &  3(n-3)+{n-3 \choose 2}\\
\end{array} \right)$.
\end{exercise}

\section{Fourier Transform}\label{Four}

\subsection{Characters of Abelian Groups}

We denote by $\mathbb{V}(K)\stackrel{\text{def}}{=}\{f:K\rightarrow\mathbb{C}\}$ the set of functions acting from the group $(K,\circ)$ to the field of
complex numbers.

%\begin{definition}
A mapping $\phi\in\mathbb{V}(K)$ is called  a {\sl   character} of a group if it  
  satisfies the following conditions:
    \begin{equation}\label{eq:char}
    \forall x\in K\quad |\phi(x)|=1;\qquad
    \forall x,y\in K\quad \phi(x\circ y)=\phi(x)\cdot{\phi(y)}.
    \end{equation}
%\end{definition}

Thus, any character is a homomorphism from the
group $K$ to the group of points on the unit circle
$T=\{x\in\mathbb{C}:\quad|x|=1\}$ with the operation of multiplication. The group $T$
is isomorphic to the unitary group $\mathrm{U}(1)$ and the special
orthogonal group $\mathrm{SO}(2)$.

Let $e$ be an identity element of the group, then $\phi(x)=\phi(ex)=\phi(e)\phi(x)$. So, $\phi(e)=1$.
For a finite group the condition $|\phi(x)|=1$ is actually redundant, since every element of the group has finite order. Then from $mx=e$ it follows $\phi(x)^m=\phi(mx)=1$.

%\textbf{Properties of characters:}

We show that the characters form an abelian group, which we will
denote by $X(K)$. Define the product of characters by the equality
$$(\phi_1\cdot\phi_2)(x)=\phi_1(x)\phi_2(x).$$ Let us verify that the
product of characters $\phi_1,\phi_2\in X(K)$ is a character:
$$|(\phi_1\cdot\phi_2)(x)|=|\phi_1(x)||\phi_2(x)|=1;$$
$$(\phi_1\cdot\phi_2)(x\circ y)=\phi_1(x\circ y)\phi_2(x\circ y)
\hfill =\phi_1(x)\phi_1(y)\phi_2(x)\phi_2(y)
=(\phi_1\cdot\phi_2)(x)(\phi_1\cdot\phi_2)(y).$$

Associativity and commutativity of the multiplication of characters
follow from the associativity and commutativity of the multiplication of
complex numbers. The identity of the group $X(K)$ is the identically
one function  $\mathbf{1}$, where $\mathbf{1}(x)=1$ for every $x\in K$. The inverse element is the
character $(\phi)^{-1}=\bar\phi$, which is the complex conjugate of the original one.
The character group $X(K)$ is called the dual group of $K$.

%It is easy to see that if $A<K$, then $X(A)<K$.
% and is denoted by $\widehat{K}=X(K)$.

Henceforth we will consider only abelian groups $(K, +)$, i.e.,
we assume the group operation $+$ to be commutative. Consider the characters
of the character group $X(K)$. Define the character $\Theta_x\in X(X(K))$
by the equality $\Theta_x(\phi)=\phi(x)$ for all $\phi\in X(K)$. 
Note that characters do not distinguish between 
conjugate elements of the group $\phi(a^{-1}ba)=\phi(b)$. 
However, for a finite abelian group (more generally, a locally compact abelian group) it is possible to prove that for every pair of elements $x,y\in K$ there exists $ \phi\in X(K)$  such that $\phi(x)\ne\phi(y)$. Therefore, $x\ne y \Leftrightarrow
\Theta_x\ne\Theta_y$ and $K$
is a subgroup of the group $X(X(K))$.

\begin{theorem}[Pontryagin--van Kampen duality theorem \cite{Morris}, \cite{RS}]\label{thPK}

    If $K$ is a locally compact commutative group, then $K\simeq X(X(K))$.
\end{theorem}

{\bf Explanation of the term "locally compact".}  A topological group is a group equipped with a topology for which multiplication and inversion are continuous. A group equipped with such a topology
is called a topological group. The characters of a topological group
are only the continuous (with respect to this topology) functions
satisfying conditions (\ref{eq:char}). A topological group
is locally compact if every element of the group has a compact
neighborhood. In the case of any finite group and the group $\mathbb{Z}$,
the discrete topology possesses the required properties. For
$\mathbb{R}$ and $\mathbb{C}$, the standard topology given by the
absolute value of the number is suitable.

In the case when $K\simeq X(X(K))$, one can introduce an operator
$\Phi:K\times X(K)\rightarrow \mathbb{C}$ such that for
a fixed $y\in X(K)$ the function $\Phi(x,y)$ is a character of the
group $K$, and for a fixed $x\in K$ it is a character of the group
$X(K)$.

\begin{example}
 Consider the group $(\mathbb{R},+)$. As mentioned above,  it is locally compact with respect to the
 natural topology defined by the absolute value of a number. It is easy to check that
functions of the form $\phi_y(x)=e^{ixy}$, $y\in \mathbb{R}$, satisfy
conditions $(\ref{eq:char})$. Indeed,
$\phi_y(x_1+x_2)=e^{i(x_1+x_2)y}=e^{ix_1y}\cdot
e^{ix_2y}=\phi_y(x_1)\cdot\phi_y(x_2)$, $|\phi_y(x)|=|e^{ixy}|=1$.
It can be shown that there are no other continuous functions on $\mathbb{R}$
satisfying conditions $(\ref{eq:char})$.

From the equality $\phi_{y_1}(x)\cdot\phi_{y_2}(x)=e^{ixy_1}\cdot
e^{ixy_2}=e^{ix(y_1+y_2)}=\phi_{y_1+y_2}(x)$ it follows that
$X(\mathbb{R},+) \simeq (\mathbb{R},+)$.\end{example}

\begin{example}\label{ex:charZ}
 Consider the group $(\mathbb{Z},+)$, as a subgroup of
$(\mathbb{R},+)$. We have $\phi_y (n)=e^{iyn}$, $n\in \mathbb{Z}$.
Since $\phi_{y+2\pi} (n)=\phi_y (n)$, to each character $\phi_y$
one can assign a point on the unit circle
having angle $y$. Using the natural isomorphism
$\phi_y\leftrightarrow e^{iy}$, $y\in \mathbb{R}$,
we obtain that
$X({\mathbb{Z},+})\simeq U(1)\simeq SO(2)$. 
Strictly speaking, we have only shown that $U(1)\leq
X({\mathbb{Z},+})$. Since $X(U(1))\simeq  \mathbb{Z}$ by
 the duality theorem, in both cases we have an isomorphism.
It is clear that the
character group will not change (up to equivalence) upon changing
the argument $y$; therefore, one can define $\phi_z (n)=e^{2\pi izn}$,
$z\in [0,1)$. By the duality theorem we have
$X(U(1))\simeq(\mathbb{Z},+)$. Indeed, if we consider
the operators $\phi_z (n)$ as functions of the variable $z$, we have
the equality $\phi_z (n_1)\phi_z (n_2)=\phi_z(n_1+n_2)$.\end{example}

\begin{example}\label{ex:charZ_q}
 Consider the cyclic group
$\mathbb{Z}_q=\mathbb{Z}/q\mathbb{Z}$. Now we need
characters of the group $(\mathbb{Z},+)$ that do not distinguish between elements that are comparable
modulo $q$, i.e., such that $\phi_y (n+q)=\phi_y(n)$ for any
$n\in \mathbb{Z}$. Again we seek characters in the form $\phi_y (n)=e^{iyn}$.
Then $y(n+q)=yn+2\pi k$, i.e., $y=\frac{2\pi k}{q}$. Put
$\xi=e^{2\pi i/q}$; then the characters of the group $\mathbb{Z}_q$ can be
represented in the form $\phi'_k(n)=\xi^{n\cdot k}=e^{\frac{2\pi
ink}{q}}$. Since $\phi'_{k_1}(n)\cdot
\phi'_{k_2}(n)=e^{\frac{2\pi in}{q}(k_1+k_2) } =e^{\frac{2\pi
in}{q}((k_1+k_2)\mod q)}=\phi'_{(k_1+k_2)\mod q}(n)$, we have
$X({\mathbb{Z}}_q)\simeq \mathbb{Z}_q$.
\end{example}

\begin{claim}\label{c:char_sum}
        If $\phi$ is a character of a finite group $K$, then the
        equalities
        \[
        \sum\limits_{x\in K}\phi(x)=0,\mbox{ if }\phi\neq {\bf 1};
        \]
        \[
        \sum\limits_{x\in K}\phi(x)=|K|,\mbox{ if }\phi= {\bf 1}
        \]
        hold.
\end{claim}
\begin{proof}
   For any $a\in K$ the equalities
    \[
    \sum\limits_{x\in K}\phi(x)=\sum\limits_{x\in K}\phi(x+a)=\phi(a)\sum\limits_{x\in K}\phi(x)
    \]
    hold.
    If $\phi\not\equiv 1$, then there exists $a\in K$ for which $\phi(a)\neq 1$.
    Hence it follows that
     $\sum\limits_{x\in K}\phi(x)=0$. Otherwise $\sum\limits_{x\in
     K}\phi(x)=|K|$.
\end{proof}

From Theorem \ref{thPK} and Proposition \ref{c:char_sum} we have

\begin{corollary}\label{c:char_sum1}
     For a finite abelian group $K$, the
        equalities
        \[
\begin{array}{ll}
        \sum\limits_{\phi\in X(K)}\phi(x)=|X(K)|, & \mbox{ if } x \mbox{ is the identity element of the group } K;\\
        \sum\limits_{\phi\in X(K)}\phi(x)=0\, & \mbox{otherwise}\\
\end{array}
                        \]
        hold.
\end{corollary}

\begin{claim}\label{c:char_sum11}
Let $\xi=e^{2\pi i/q}$, $m=1,\dots,q-1$. Then $\sum\limits_{k=0}^{q-1}\xi^{mk}=0$. If $q$ is a prime number and $\sum\limits_{k=0}^{q-1}\alpha_k\xi^{mk}=0$  then $\alpha_k=\alpha_0$ for all $k$.
\end{claim}
\begin{proof}
The first statement is a special case of Proposition \ref{c:char_sum1} for $\mathbb{Z}_q$. It also follows from the equations $$0=\xi^{mq}-1=(\xi^m-1)\sum\limits_{k=0}^{q-1}\xi^{mk}.$$
The second  one follows from the well-known property that $\xi^m$ is not a root of any equation of degree less than $q-1$. 
\end{proof}

In  $\mathbb{V}({K})$ one can define a scalar
product. In the case of a finite group $K$, 
$\mathbb{V}({K})$ is a finite-dimensional space and the scalar product in
$\mathbb{V}({K})$ can be defined by the equality
$(f,g)=\sum\limits_{x\in K}f(x)\overline{g(x)}$.

\begin{claim}\label{c:charact_orthonorm_basis}
    The set of characters of a finite abelian  group forms an orthogonal basis in the space $\mathbb{V}({K})$.
\end{claim}
\begin{proof}
    We begin by proving the orthogonality of the characters. Let $a\in K$. From the definition of characters and the invariance of the inner product under shifts of the argument, we have the equalities\\
    $(\phi_1,\phi_2)=\sum\limits_{x\in K}\phi_1(x)\overline{\phi_2(x)}=
    \sum\limits_{x\in K}\phi_1(x+a)\overline{\phi_2(x+a)}=\phi_1(a)\cdot\overline{\phi_2(a)}\cdot(\phi_1,\phi_2)$.
   % $\sum\limits_{x\in K} \phi_1(x)\cdot \overline{\phi_2(x)}=\sum\limits_{x\in K}\phi_1(a+x)\cdot\overline{\phi_2(a+x)}=\sum\limits_{x\in K}\phi_1(a)\cdot\overline{\phi_2(a)}\cdot\phi_1(x)\overline{\phi_2(x)}=\phi_1(a)\overline{\phi_2(a)}\sum\limits_{x\in K}\phi_1(x)\cdot\overline{\phi_2(x)}$
    If $\phi_1\neq \phi_2$, then there exists $a\in K$ such that  $\phi_1(a)\overline{\phi_2(a)}\ne 1$. Consequently
     $(\phi_1,\phi_2)=\sum\limits_{x\in K}\phi_1(x)\cdot\overline{\phi_2(x)}=0$.

     Now we prove that the system of characters is complete. We have
     the inequality
    $\mathrm{dim}\,\mathbb{V}({K}) = |{K}|\ge |X(K)|$.
    The reverse inequality follows from the Pontryagin--van Kampen duality theorem. Indeed,
    consider the characters on the group $X(K)$:
    $\mathrm{dim}\,\mathbb{V}(X({K})) = |X({K})|\ge |X(X(K))|=|K|$, since $X(X(K))\simeq K$.
    Thus, the number of characters coincides with the dimension of the space $\mathbb{V}({K})$.
  %  If $\mathbb{K}$ is not finite, consider the space of square-integrable functions $L_2(K)$. That is, instead of the sum, we consider the integral. In the case of $\mathbb{R}$: $\int\limits_{\mathbb{R}} e^{ixy}dx$ is not integrable.
\end{proof}

\begin{corollary} For a finite abelian group $K$ we have $|X(K)|=|K|$.
\end{corollary}

As can be seen from the proof, the characters of any, in particular non-abelian, group are orthogonal. But in the non-abelian case they do not form a basis.
Note that if the constructed set of characters constitutes a basis
of the space $\mathbb{V}({K})$, then there are no other characters of the group $K$.
It follows that in Examples \ref{ex:charZ} and
\ref{ex:charZ_q} all characters of the groups $U(1)$ and
$\mathbb{Z}_q$ are defined.

\begin{claim}\label{c:charact_cartesian_product}
    Let  $K_1$ and $K_2$ be finite abelian groups and $K=K_1\times K_2$ be the direct product of the groups.
    Then $X(K)=X({K}_1)\times X({K}_2)$.
\end{claim}

\begin{proof}
        Consider the characters $\phi_{z_1}\in X({K}_1)$, $\phi_{z_2}\in X({K}_2)$. Their
     product
    $(\phi_{z_1}\cdot\phi_{z_2})(x_1,x_2)\stackrel{\text{def}}{=}\phi_{z_1}(x_1)\cdot\phi_{z_2}(x_2)$
     is a character of the group $K$.
   Conditions $(\ref{eq:char})$ are easily verified directly.
   The number of constructed characters is equal to $|X(K_1)||X(K_2)|$. It coincides
   with the cardinality of the group $X(K)$ equal to
   $|K_1||K_2|$.
\end{proof}

It is known that any finite abelian group $K$ is isomorphic to the Cartesian
product (direct sum) of cyclic groups $\mathbb{Z}_q$, $q\in
\mathbb{N}$. Then from Proposition \ref{c:charact_cartesian_product}
and Example \ref{ex:charZ_q} it follows that $K\simeq X(K)$ for any
finite abelian group $K$.

\begin{claim}\label{cl:subchar}
Let $K_1$ be a subgroup of a finite abelian group $K$. Then for any
character $\phi\in X(K_1)$ there exists a character $\widetilde{\phi}\in
X(K)$ such that $\phi=\widetilde{\phi}|_{K_1}$.
\end{claim}
\begin{proof}
Consider the factor group $K/K_1$. Each element $x\in K$
can be uniquely represented as a sum $x=x_1+x_2$, where
$x_1\in K_1$, $x_2\in K_2$, and the group $K_2$ is isomorphic to the factor group
$K/K_1$. Define $\widetilde{\phi}(x)=\phi(x_1)$. It is easy
to verify that $\widetilde{\phi}$ satisfies condition
(\ref{eq:char}).
\end{proof}

Let $K_1\subset K$. Define the subgroup $K_1^\bot$ in $X(K)$
by the equality
$$K_1^\bot=\{z\in X(K) \ |\  \phi_z(x)=1,\ \mbox{if}\ x\in K_1\}.$$
If a finite abelian group $K$ is representable as a direct
product of groups $K=K_1\times K_2$, then one can consider
the subgroups $K_1\simeq K_1\times \{e_2\}$ and $K_2\simeq \{e_1\}\times
K_2$ as subsets of $K$, where $e_1$ and $e_2$ are the identity
elements of the groups $K_1$ and $K_2$. From Proposition
\ref{c:charact_cartesian_product} we have $X(K)=X(K_1)\times
X(K_2)$. It is easy to see that in this case $K_1^\bot\simeq X(K_2)$.
Let us prove a similar fact in a slightly more general case.

\begin{claim}\label{claim:factor}
 Let $K_1$ be a subgroup of a finite abelian group $K$. Then $K_1^\bot\simeq X(K/K_1)$.
\end{claim}
\begin{proof}
Let $K_2$ be a subgroup of $K$ such that $K_2\simeq K/K_1$ and each element $x\in K$
can be uniquely represented as a sum $x=x_1+x_2$, where
$x_1\in K_1$, $x_2\in K_2$. 
For every  $\phi\in K_1^\bot$ we have $\phi(x_1+x_2)=\phi(x_2)$, i.e., there exists a homomorphism  $K_1^\bot\rightarrow X(K_2)$.

 For every $\phi\in X(K_2)$ define $\widetilde{\phi}(x)=\phi(x_2)$. By Proposition \ref{cl:subchar} we have  an injective map $X(K_2)\rightarrow K_1^\bot$ acting as  $\phi \rightarrow \widetilde{\phi}$. Thus, it is an isomorphism.
\end{proof}

\begin{corollary}\label{claim:factor1}
Let $K_1$ be a subgroup of a finite abelian group $K$.  Then
$|K_1^\bot|\cdot|K_1|=|K|$.
\end{corollary}

\begin{example}
As shown in Example $\ref{ex:charZ_q}$, the characters of the group
$\mathbb{Z}_q$ are the functions $\phi_y(x)=\xi^{x\cdot y}$, where
$\xi=e^{2\pi i/q}$ and $x,y\in \mathbb{Z}_q$. Then the characters of the group
$\mathbb{Z}^n_q$ are the functions $\xi^{\langle x, y\rangle}$, where
$x,y\in \mathbb{Z}^n_q$ and $\langle x,
y\rangle=x_1y_1+\dots+x_ny_n\mod q$.
\end{example}

Note that in the character formula $\xi^{\langle x\cdot y\rangle}$,
the inner product $\langle x\cdot
y\rangle=x_1y_1+\dots+x_ny_n$ can be considered both modulo $q$
and in integers, since $\xi^a=\xi^{a\mod q}$.

Let $K$ be a finite abelian group. $A\subset K$, $0 \notin A$,
$-A=A$. Recall that $\mathrm{Cay}(K,A)$
is the Cayley graph of the group $K$ with generating set $A$.

\begin{claim}[\cite{Bab}]\label{c:char_egenfunc}
    A character $\phi$ of the group $K$ is an eigenvector of the adjacency matrix $M$ of the graph $\mathrm{Cay}(K,A)$
    with eigenvalue
    $\sum\limits_{a\in A}\phi(a)$.
\end{claim}
\begin{proof}

    $(M\phi)(x)=\sum\limits_{a\in A}\phi(x+a)=
    \sum\limits_{a\in A}\phi(x)\phi(a)=\phi(x)\left(\sum\limits_{a\in A}\phi(a)\right)$.
\end{proof}

\begin{corollary}\label{c:baschar}
Every eigenspace of the adjacency matrix of the graph
$\mathrm{Cay}(K,A)$ has a basis consisting of characters of the group $K$.
\end{corollary}

\begin{claim}
    A character $\phi$ of the group $K$ is a perfect coloring (here we use complex numbers as the colors) of the graph $\mathrm{Cay}(K,A)$.
\end{claim}
\begin{proof}
A vertex of color $\phi(x)$ is adjacent to vertices of colors $\{\phi(x+a)\ |\
a\in A\}=\{\phi(x)\phi(a)\ |\ a\in A\}$. Therefore, if the colors of two
vertices coincide, i.e., $\phi(x_1)=\phi(x_2)$, then the color composition
of their neighborhoods also coincides.
\end{proof}

\begin{corollary}
The characters $\phi_z$ of the group $\mathbb{Z}^n_q$ are perfect
colorings of the hypercube $Q_q^n$ with quotient matrix depending
only on $\mathrm{wt}(z)$. \end{corollary}

\subsection{Krawtchouk Polynomials }\label{10.2}

The hypercube $Q_q^n$ is a Cayley graph of the group $\mathbb{Z}^n_q$ with
the generating set $A$ consisting of tuples of weight $1$. We find the
eigenvalues corresponding to the characters
$\phi_z(x)=\xi^{\langle x\cdot z\rangle}$, where $\xi=e^{2\pi i/q}$. Recall that
$\mathrm{wt}(z)=d(z,\bar{0})=|\operatorname{supp}(z)|$ is the weight of the tuple $z$. We have
the equalities

\[
\sum\limits_{a\in A}\phi_z(a)=\sum\limits_{a\in
A}\xi^{<a,z>}=\sum\limits_{\substack{a\in A\\ \supp(z)\cap
\supp(a)=\varnothing}}\xi^{<a,z>}+\sum\limits_{\substack{a\in A\\
\operatorname{supp}(z)\cap \operatorname{supp}(a)\ne\varnothing}}\xi^{<a,z>}.
\]

We compute the terms separately. Any element $a\in A$ has the form
$a=(0,\dots,0,a_i,0,\dots,0)$ for some $a_i=b\in\mathbb{Z}_q$,
$b\neq0$. For any fixed $z\in \mathbb{Z}_q$, we represent
both sums as
\[
\sum\limits_{\substack{a\in A\\  \supp(z)\cap
\supp(a)=\varnothing}}\xi^{<a,z>}=\sum\limits_{i:z_i=0}\sum\limits_{a\in
A, a_i\neq 0}\xi^{<a,z>}; \]
\[
\sum\limits_{\substack{a\in A\\
\operatorname{supp}(z)\cap
\operatorname{supp}(a)\ne\varnothing}}\xi^{<a,z>}=\sum\limits_{i:z_i\neq
0}\sum\limits_{a\in A, a_i\neq 0} \xi^{<a,z>}.\]
For a fixed $i$ we have
\begin{equation}\label{eqsum} \sum\limits_{a\in A, a_i\neq 0}\xi^{<a,z>}=\sum\limits_{b\in
\mathbb{Z}_q\setminus \{0\}}1=q-1,\ \mbox{if}\ z_i=0;
\end{equation}
$$\sum\limits_{a\in A, a_i\neq 0}\xi^{<a,z>}=\sum\limits_{b\in
\mathbb{Z}_q\setminus \{0\}} \xi^{z_i b}=-1,\ \mbox{if}\   z_i\neq
0;$$ where the last equality follows from Proposition
\ref{c:char_sum}. Then
\[
\sum\limits_{a\in
A}\phi_z(a)=(n-\operatorname{wt}(z))(q-1)+(-1)\operatorname{wt}(z)=n(q-1)-\operatorname{wt}(z)q.
\]
Thus, the eigenvalues of the hypercube $Q_q^n$ are the numbers $-n,
-n+q, \dots, n(q-1)$. This agrees with Corollary \ref{cl:spec}.

Now we compute the values of the Krawtchouk polynomials $P_t[n,q]$ for hypercubes
$Q_q^n$. From Proposition \ref{c:WDB_bound0} we have
$$P_t[n,q](n(q-1)-\mathrm{wt}(z)q)=\phi_z(\bar{0})
\left(\sum\limits_{x:\mathrm{wt}(x)=t}\phi_z(x)\right).$$ Without
loss of generality, we assume that $z=(1,\dots,1,0,\dots,0)$,
$\mathrm{wt}(z)=k$. Divide the set of vectors of weight $t$ into groups,
each containing vectors with the same support. Suppose that in some
group, exactly $m$ nonzero coordinates of the vectors lie in the
$1$-st through $k$-th positions. From formulas
(\ref{eqsum}), the sum $\sum_x\phi_z(x)$ over all vectors $x$ from such
a group is equal to $$(q-1)^{t-m}\sum\limits_{x}\xi^{x_{i_1}+\cdots+
x_{i_m}}=$$ $$=(q-1)^{t-m}\left(\sum\limits_{x_{i_1}\neq
0}\xi^{x_{i_1}}\right)\cdots \left(\sum\limits_{x_{i_m}\neq
0}\xi^{x_{i_m}}\right)=(-1)^m(q-1)^{t-m}.$$ For each $m$, the number of such
groups is equal to ${{n-k}\choose{t-m}}{{k}\choose{m}}$.
%For each such vector $\phi_z(x)=
% From the equality $\sum\limits_{x:x_i\neq0, i=1\dots
%m}\xi^{x_1+\cdots+ x_m}=(\sum\limits_{x_1\neq 0}\xi^{x_1})\cdots
%(\sum\limits_{x_m\neq 0}\xi^{x_m})=(-1)^m$ we have
Then the equality
$$P_t[n,q](n(q-1)-kq)=
%=\sum\limits_{m=0}^t\sum\limits_{x_i\neq 0}\xi^{x_1+\cdots+
%x_m}(q-1)^{t-m}{{n-k}\choose{t-m}}{{k}\choose{m}}=$$ $$=
\sum\limits_{m=0}^t
(-1)^m(q-1)^{t-m}{{n-k}\choose{t-m}}{{k}\choose{m}},
$$ holds, where for $k< m$ or
$ n-k<t-m$ the corresponding terms are equal to zero.

Further, we find the generating function for the Krawtchouk polynomials. Denote
by $y$ an independent variable. For any $x\in Q^n_q$, the equality
$y^{\mathrm{wt}(x)}=y^{\mathrm{wt}(x_1)}\cdots y^{\mathrm{wt}(x_n)}$ holds. As before,
we assume that $z=(1,\dots,1,0,\dots,0)$, $\mathrm{wt}(z)=k$.
Using Proposition \ref{c:WDB_bound0} and formulas (\ref{eqsum}),
we obtain the equalities
$$\sum\limits_{x\in Q^n_q}\phi_z(x)y^{\mathrm{wt}(x)}=
\sum\limits_{t=0}^ny^t\left(\sum\limits_{x:\mathrm{wt}(x)=t}\phi_z(x)\right)=
\sum\limits_{t=0}^ny^tP_t[n,q](n(q-1)-\mathrm{wt}(z)q);$$
$$\sum\limits_{x\in Q^n_q}\phi_z(x)y^{\mathrm{wt}(x)}=\sum\limits_{x\in Q^n_q}
\xi^{x_1}\cdots\xi^{x_k}y^{\mathrm{wt}(x_1)}\cdots y^{\mathrm{wt}(x_n)}=$$ $$=
\left(\prod\limits_{i=1}^{k}\sum\limits_{x_i\in
Q_q}\xi^{x_i}y^{\mathrm{wt}(x_i)}\right)
\left(\prod\limits_{i=k+1}^{n}\sum\limits_{x_i\in
Q_q}y^{\mathrm{wt}(x_i)}\right)=(1-y)^k(1+(q-1)y)^{n-k}.$$ Thus,

\begin{equation}\label{eqKraw}
(1+(q-1)y)^{n-k}(1-y)^k=
\sum\limits_{t=0}^n P_t[n,q](\lambda_k)y^t,
\end{equation} where
$\lambda_k=n(q-1)-kq$.

Let us express the powers of the adjacency matrix in terms of the Krawtchouk polynomials. Let $M$ be the adjacency matrix of a distance-regular graph. As noted above, $(M^r)_{ij}$
is the number of walks of length $r$ between vertices $x_i$ and $x_j$ in the graph. The coefficients of the Krawtchouk-type   polynomials $P_t$ define a
triangular matrix with a nonzero diagonal, so there are inverse
relations:
$$M^r = \alpha^r_0 I_{N} + \alpha^r_1 M_1+ \cdots+ \alpha^r_r M_r$$
with some coefficients $\alpha^r_i$, $i=1,\dots r$.

In the case of the Boolean $n$-cube from the previous equality and equalities
$M_iM=(n-i+1)M_{i-1}+(i+1)M_{i+1}$ for $i=1,\dots r$ we obtain
equality
$$M^{r+1} = \alpha^r_0 M + \alpha^r_1 (nI_{2^n}+2M_2)+ \cdots+
\alpha^r_r((n-r+1)M_{r-1}+(r+1)M_{r+1}).$$ 
We have recurrent
relations $\alpha^{r+1}_i=i\alpha^r_{i-1}+(n-i)\alpha^r_{i+1}$ for
$i=1,\dots r$ and $\alpha^r_{r+1}=0$. It is not difficult to verify that the values $\alpha^r_i
= ({\cosh}^{n-i}(x)\,{\sinh}^{i}(x))^{(r)} \mid_{x=0}$ satisfy the same recurrence relations. 
Indeed $$\frac{d}{dx} ({\cosh}^{n-i}(x)\,{\sinh}^{i}(x))=i\cosh^{n-i+1}(x)\sinh^{i-1}(x)+(n-i)\cosh^{n-i-1}(x)\sinh^{i+1}(x)$$
which yields exactly the same recurrence after taking the $r$-th derivative at $x=0$.
 Thus, to find the coefficients $\alpha^r_i$ for the Boolean $n$-cube, it suffices to calculate the $r$-th derivative at $x=0$ of the product
of the powers of the hyperbolic cosines and sines.

\subsection{Definition and Properties of the Fourier Transform}\label{10.4}

%\begin{definition}
Let $\phi_z$ be the characters of a group $K$. The {\sl Fourier transform} of a function $f\in V({K})$ is the function
$\widehat{f}(z)=\mu\cdot(f,\phi_z)$ mapping the character group $X(K)$ into $\mathbb{C}$. The factor $\mu\in \mathbb{R}$ is chosen so that the Fourier transform preserves the norm of the function:
$\|\widehat{f}\|_2=\|f\|_2=\sqrt{(f,f)}$.
%\end{definition}

\begin{example}
If $K=(\mathbb{R},+)$, then the scalar product is defined by
\[(a,b)=\int\limits_{\mathbb{R}}a(x)\overline{b(x)}\,dx\]
and
\[\widehat{f}(z)=\frac{1}{\sqrt{2\pi}}\int\limits_{\mathbb{R}}f(x)e^{-ixz}dx.\]
\end{example}

\begin{example}
If $K=U(1)$, then the scalar product is defined by
\[(a,b)=\int\limits_0^{2\pi}a(x)\overline{b(x)}\,dx\]
and
\[\widehat{f}(n)=\frac{1}{\sqrt{2\pi}}\int\limits_0^{2\pi}f(x)e^{-inx}dx.\]
In this case the result of the Fourier transform can be regarded as the set of coefficients of the Fourier series expansion.
\end{example}

\begin{example}
If $K=\mathbb{Z}^n_q$, then the scalar product is defined by
\[(a,b)=\sum\limits_{x\in \mathbb{Z}^n_q}a(x) \overline{b(x)}\]
and
\[\widehat{f}(z)=\frac{1}{q^{n/2}}\sum\limits_{x\in \mathbb{Z}^n_q}f(x)\overline{\xi}^{\langle x,z\rangle}. \]

In the case $K=\mathbb{Z}_q$, i.e.\ when $n=1$, this Fourier transform is called the {\sl discrete Fourier transform}.

Consider the special case $K=\mathbb{Z}_2^n$. We have $\xi=-1$ and
\[\widehat{f}(z)=\frac{1}{2^{n/2}}\sum\limits_{(x_1,\ldots,x_n)\in \mathbb{Z}_2^n}f(x_1,\ldots,x_n)\cdot (-1)^{x_1 z_1}\cdots (-1)^{x_n z_n}
= \frac{1}{2^{n/2}}\sum\limits_{x\in \mathbb{Z}_2^n}
f(x)(-1)^{<x,z>}. \]
\end{example}

Note that in all the listed cases, except $K=U(1)$, the group $X(K)$ is isomorphic to the group $K$. As noted above, the group $X(K)$ is isomorphic to $K$ for any finite abelian group. In this case one can regard the Fourier transform as acting from $\mathbb{V}(K)$ to $\mathbb{V}(K)$. In what follows we always assume that $K$ is a finite abelian group and the Fourier transform acts from $\mathbb{V}(K)$ to $\mathbb{V}(K)$ according to the formula
\[\widehat{f}(z)=\frac{1}{|K|^{1/2}}(f,\phi_z)= \frac{1}{|K|^{1/2}}\sum\limits_{x\in K}f(x)\overline{\phi_z(x)}. \]

\begin{remark}
The rows of the generalized Hadamard matrix $H_{q,n}$ are the vectors of values of the characters on the group $\mathbb{Z}_q^n$. Therefore the Fourier transform of a function on the group $\mathbb{Z}_q^n$ coincides with the multiplication of the matrix $H_{q,n}$ by the vector of function values, i.e.\ $\widehat{f}=\frac{1}{q^{n/2}}H_{q,n}f$.
\end{remark}

\begin{claim}\label{c:fourier_basis}  The following equality holds:
\[f(x)=\frac{1}{|K|^{1/2}}\sum\limits_{z\in K} \widehat{f}(z)\phi_z(x). \]
\end{claim}
\begin{proof}
From the definition of a character we have $\phi_z(x)\phi_z(-x)=\phi_z(0)=1$, hence $\phi_z(-x)=\overline{\phi_z(x)}$.
Then we have the equality $\|\phi_z\|_2^2=\sum_{x\in K}\phi_z(x)\overline{\phi_z(x)} =\sum_{x\in K}\phi_z(x){\phi_z(-x)}=|K|$. From Proposition \ref{c:charact_orthonorm_basis} it follows that the functions $\phi_z/|K|^{1/2}$ form an orthonormal basis. From the definition of the Fourier transform it follows that the numbers $\widehat{f}(z)$ are the coefficients of the expansion of $f$ with respect to the basis of characters.
\end{proof}

\begin{claim}[Plancherel theorem]\label{c:plansh}
For functions on finite abelian groups the equality $(f,g)=(\widehat{f},\widehat{g})$ holds.
\end{claim}
\begin{proof}
Apply complex conjugation to the equality from Proposition \ref{c:fourier_basis}:
\[\overline{g(x)}=\frac{1}{|K|^{1/2}}\sum\limits_{z\in K}\overline{ \widehat{g}(z)\phi_z(x)}. \]

Using the orthogonality of the characters, we obtain
$$(f,g)=\frac{1}{|K|}\sum\limits_{z\in K}
\widehat{f}(z)\overline{\widehat{g}(z)}\|\phi_z\|^2=\sum\limits_{z\in
K}
\widehat{f}(z)\overline{\widehat{g}(z)}=(\widehat{f},\widehat{g}).$$
\end{proof}

From Proposition \ref{c:plansh} with $f=g$ follows Parseval's identity:
    \[\|f\|_2^2=(f,f)=\frac{1}{|K|}\sum\limits_{z\in
 K}|\widehat{f}(z)|^2 \|\phi_z\|_2^2=\sum\limits_{z\in
 K}|\widehat{f}(z)|^2=\|\widehat{f}\|_2^2. \]

\begin{claim} \label{c:double fourier}
    $\widehat{\widehat{f}\,}(x)=f(-x)$.
\end{claim}
\begin{proof}
From the definition of the Fourier transform we have
\[\widehat{\widehat{f}\,}(x)=\frac{1}{|K|^{1/2}}(\widehat{f},\phi_x)=
 \frac{1}{|K|^{1/2}}\sum\limits_{z\in
 K}\widehat{f}(z)\overline{\phi_x(z)}.\]
From Proposition \ref{c:fourier_basis} we have
\[f(-x)=\frac{1}{|K|^{1/2}}\sum\limits_{z\in
 K} \widehat{f}(z)\phi_z(-x). \]
From the Pontryagin–van Kampen duality theorem it follows   (it  also follows from Proposition  \ref{c:charact_cartesian_product} and Example \ref{ex:charZ_q})  that characters satisfy
$\Phi(x,z)=\phi_z(x)=\phi_x(z)$. And from the definition of a character we have $\phi_z(-x)=\overline{\phi_z(x)}$.
\end{proof}

\begin{claim}\label{claim:ef}
A nonzero function $f:Q^n_q\rightarrow \mathbb{C}$ is an eigenfunction with eigenvalue $\lambda=n(q-1)-kq$ if and only if $\widehat{f}(z)=0$ for all $z\in Q^n_q$ with $\wt(z)\neq k$.
\end{claim}
\begin{proof}
As shown above, the characters form a basis of the space $\mathbb{V}(K)$ and are eigenfunctions of the hypercube $Q^n_q$. Hence one can form bases of the eigenspaces of the hypercube from them.
Since the eigenspace corresponding to $\lambda=n(q-1)-kq$ is spanned by precisely those characters with
  $\wt(z)= k$, the proposition follows.
\end{proof}

\begin{claim}[Properties of the Fourier transform]\label{cl:Four1-3} \quad
\\
 $(a)$ $ \widehat{f+g}=\widehat{f}+\widehat{g};\quad
 \widehat{c\cdot f}=c\cdot\widehat{f}$, where $c$ is a constant \quad
 (linearity);\\
$(b)$ $\widehat{f(-x)}(z)=\widehat{f}(-z)$, and if $f$ is a real-valued function, then
$\widehat{f}(-z)=\overline{\widehat{f}(z)}$;\\
 $(c)$
$\widehat{f(x+a)}(z)=\phi_z(a)\widehat{f}(z)$ and
$\widehat{f}(z+a)=\widehat{f\cdot\overline{\phi_a}(x)}(z)$ \quad
(shift of the argument).
\end{claim}
\begin{proof}
The linearity of the Fourier transform (a) is obvious from the definition.
Proof of (b):
\[\widehat{f(-x)}(z)=\frac{1}{|K|^{1/2}}\sum\limits_{x\in
K}f(-x)\overline{\phi_z(x)}=\frac{1}{|K|^{1/2}}\sum\limits_{y\in
K}f(y)\overline{\phi_z(-y)}.\]

From the duality theorem and the definition of a character we have the equalities
$\phi_z(-y)=\phi_{-z}(y)=\overline{\phi_z(y)}$. Then

\[\widehat{f(-x)}(z)=\frac{1}{|K|^{1/2}}\sum\limits_{y\in
K}f(y)\overline{\phi_{-z}(y)}=\widehat{f}(-z).\] If
$f(y)=\overline{f(y)}$, then
\[\widehat{f(-x)}(z)=\frac{1}{|K|^{1/2}}\sum\limits_{y\in
K}\overline{f(y)\phi_{z}(-y)}=\frac{1}{|K|^{1/2}}\sum\limits_{y\in
K}\overline{f(y)}\phi_{z}(y)=\overline{\widehat{f}(z)}.\]

 Proof of the shift property (c).
Changing the variable $x$ to $y=x+a$, we obtain the equalities

\[ \widehat{f(x+a)}(z)=\frac{1}{|K|^{1/2}}\sum\limits_{x\in K}f(x+a)\overline{\phi_z(x)}=\frac{1}{|K|^{1/2}}\sum\limits_{y\in K}f(y)\overline{\phi_z(y-a)}= \]

here we use the equalities
\[ \overline{\phi_z(y-a)}=\overline{\phi_z(y)}\cdot\overline{\phi_z(-a)}=\overline{\phi_z(y)}\cdot\phi_z(a)
\]

therefore

\[
=\frac{\phi_z(a)}{|K|^{1/2}}\sum\limits_{y\in K}
f(y)\overline{\phi_z(y)}=\phi_z(a)\widehat{f}(z).
\]

Proof of the second identity in (c):
\[
\widehat{f\cdot\overline{\phi_a}(x)}(z)=\frac{1}{|K|^{1/2}}\sum\limits_{x\in
K}
f(x)\overline{\phi_a(x)}\cdot\overline{\phi_z(x)}=\frac{1}{|K|^{1/2}}\sum\limits_{x\in
K}f(x)\overline{\phi_{a+z}(x)}=\widehat{f}(z+a).
\]
\end{proof}

Define the Dirac delta function on the group $K$:

\[
\delta(x)=\left\{
\begin{array}{ll}
|K|^{1/2}, & x={0}\\
0, & x\neq 0.
\end{array}
\right.
\]

Apply the Fourier transform to a character:

\[
\widehat{\phi_z(x)}(y)=\frac{1}{|K|^{1/2}}(\phi_z(x),\phi_y(x))=\delta(y-z).\]
 We have the equality
$\widehat{\phi_0}=\widehat{\bf 1}=\delta$. By Proposition \ref{c:double fourier} it follows the converse equality $\widehat{\delta}= {\bf 1}$.

\begin{claim}[Uncertainty principle \cite{Tao}]\label{c:uncertainty_heisenberg}
For any function $f\in \mathbb{V}(K)$ the inequality
    \[|K|\le |\supp(f)| \cdot |\supp(\widehat{f})|\]
holds.
\end{claim}
\begin{proof}
From Parseval's identity it follows that
    $$\|f\|_2^2=\sum\limits_{z\in K} |\widehat{f}(z)|^2\le
    (\max\limits_{z}|\widehat{f}(z)|^2)
    \cdot |\supp(\widehat{f})|.$$

    From the equality $|\phi_z(x)|^2=1$ for every $x\in K$ and the Cauchy--Schwarz inequality
    $
    (u,v)^2\le \|u\|^2_2\cdot\|v\|^2_2
    $ we obtain
    \[|\widehat{f}(z)|^2=\frac{1}{|K|}\Bigl|\sum\limits_{x\in\,\supp(f)}f(x)
    \overline{\phi_z(x)} \Bigr|^2\le \frac{1}{|K|}\|f\|_2^2\cdot|\supp(f)|. \]

    If the inequality holds for every $z\in K$, then it is also true for the maximal $|\widehat{f}(z)|^2$.
    Hence
    \[
    \|f\|_2^2\le \frac{1}{|K|} \|f\|_2^2 \cdot |\supp(f)| \cdot
    |\supp(\widehat{f})|.
    \]
\end{proof}

Equality in the uncertainty principle inequality is attained, for example, when $f$ is a character or when $f$ is the characteristic function of a subgroup of $K$ multiplied by a character.

\begin{claim}[Poisson summation formula]\label{fsP}
Let $K_1$ be a subgroup of a finite abelian group $K$. For any function $f\in \mathbb{V}(K)$ the equality
$$|K|^{1/2}\sum\limits_{x\in K_1} f(x)= |K_1|\sum\limits_{z\in
K_1^\bot}\widehat{f}(z)$$
holds.
\end{claim}
\begin{proof}

$$\sum\limits_{z\in K_1^\bot}\widehat{f}(z)=\frac{1}{|K|^{1/2}}
\sum\limits_{z\in K_1^\bot}\sum\limits_{x\in\,K}f(x)
    \overline{\phi_z(x)}=\frac{1}{|K|^{1/2}}\sum\limits_{x\in\,K}f(x)\sum\limits_{z\in K_1^\bot}\overline{\phi_z(x)}.$$

As noted above, an element $x\in K$ can be regarded as an element of $X(X(K))$, i.e.,  every $x\in K$ induces a character of the subgroup $K_1^\bot$, $K_1^\bot\subset X(K)$ by Proposition \ref{cl:subchar}. Moreover, every $x\in K_1$ induces the identity element of the subgroup $K_1^\bot$ by the definition of the group $K_1^\bot$. By using the isomorphism $K\simeq X(X(K))$, from Proposition \ref{claim:factor} we obtain
$|(K_1^\bot)^\bot|=|X(K)|/|K_1^\bot|=|K|/|K_1^\bot|=|K_1|$. Then $x\in K$ induces the identity element   of  $K_1^\bot$ only if $x\in K_1$. By Corollary \ref{c:char_sum1} we have
$$\sum\limits_{z\in
    K_1^\bot}\overline{\phi_z(x)}=0,\ \text{if}\ x\not \in K_1,$$
    $$\sum\limits_{z\in
    K_1^\bot}\overline{\phi_z(x)}=|K_1^\bot|,\ \text{if}\ x \in K_1.$$
     Then
$$\frac{1}{|K|^{1/2}}\sum\limits_{x\in\,K}f(x)\sum\limits_{z\in
    K_1^\bot}\overline{\phi_z(x)}= \frac{|K_1^\bot|}{|K|^{1/2}}\sum\limits_{x\in K_1}
    f(x).$$
By Proposition \ref{claim:factor} we have
$\frac{|K_1^\bot|}{|K|^{1/2}}=\frac{|K|^{1/2}}{|K_1|}$.
\end{proof}

\subsection{Fourier Transform and Convolution}

%\begin{definition}
Let $K$ be a finite group and $f,g: K\rightarrow \mathbb{C}$. We define the {\sl convolution} of two functions as follows:
\[
(f*g)(y)=\sum\limits_{x\in K} f(x)\cdot g(y-x).
\]
%\end{definition}

Note that for real functions convolution is defined similarly, only the sum is replaced by an integral.  We continue to assume that $K$ is a finite abelian group. The following properties of convolution are well known.

\begin{claim}[Properties of convolution]\quad \\
\begin{enumerate}
    \item $f*g=g*f$;
    \item $f*(g*h)=(f*g)*h$;
    \item $f*(g_1+g_2)=f*g_1+f*g_2$;
    \item $f*(\alpha g)=\alpha(f*g)$;
    \item $f*\delta=|K|^{1/2} f$.
\end{enumerate}
\end{claim}

\begin{claim}\label{c:convolution}
    $\widehat{f*g}=|K|^{\frac{1}{2}}\widehat{f}\cdot\widehat{g}$.
\end{claim}
\begin{proof}
    By the definition of characters we have $\phi_z(y)=\phi_z(x)\phi_z(y-x)$.
    Then from the definition of convolution we obtain the equalities
    \[
    \widehat{f*g}(z)=
    \frac{1}{|K|^{1/2}}\sum\limits_{y\in K}\overline{\phi_z(y)}
    \sum\limits_{x\in K}f(x)g(y-x)=\]
    \[=\sum\limits_{x\in K}f(x)\overline{\phi_z(x)}\frac{1}{|K|^{1/2}}
    \sum\limits_{y\in K} g(y-x)\overline{\phi_z(y-x)}= \sum\limits_{u\in K} g(u)\overline{\phi_z(u)}=
    |K|^{1/2}\cdot \widehat{f}\cdot\widehat{g}.
    \]
\end{proof}

\begin{claim}\label{c:convolution2}
     $\widehat{f}*\widehat{g}=|K|^{\frac12}\widehat{f\cdot g}$.
\end{claim}
\begin{proof}
Substituting the functions $\widehat{f}$ and $\widehat{g}$ for $f$ and $g$ in the equality
$\widehat{f*g}=|K|^{\frac{1}{2}}\widehat{f}\cdot\widehat{g}$ and applying the Fourier transform to both sides and using the fact that $\widehat{\widehat{f}\,}(z)=f(-z)$, we obtain
    \[
    (\widehat{f}*\widehat{g})(-z)=|K|^{1/2}\widehat{[f(-x)\cdot
    g(-x)]}(z).
    \]
    Using the equality $\widehat{f(-x)}(z)=\widehat{f}(-z)$ from Proposition \ref{cl:Four1-3}(b), we have $\widehat{[f(-x)\cdot g(-x)]}(z)=\widehat{[(f\cdot
    g)(-x)]}(z)=\widehat{[(f\cdot
    g)(x)]}(-z)$. Replacing $z$ to $-z$, we obtain the required equality.
\end{proof}

%\begin{definition}
 The transform $\widehat{(-1)^f}$ of a Boolean-valued function $f:K\rightarrow\{0,1\}$ is called the {\sl Walsh--Hadamard transform}.
%\end{definition}

\begin{theorem}[Titsworth \cite{Tits}]\label{th:Tits} For any function $f:K\rightarrow\{0,1\}$ we have
    \[
    \widehat{(-1)^f}*\widehat{(-1)^f}=|K|^{1/2}\delta
    \]
   \end{theorem}
\begin{proof}
By Proposition \ref{c:convolution2} we obtain
    \[
    \widehat{(-1)^f}*\widehat{(-1)^f}=|K|^{1/2}\widehat{(-1)^{f}
    \cdot(-1)^{f}}=|K|^{1/2} \widehat{\mathbf{1}}=|K|^{1/2} \delta.
    \]
\end{proof}

Thus, for the function $h= \widehat{(-1)^f}$, where $f$ is an arbitrary Boolean-valued function, we have
\[
h*h(y)=\sum\limits_{z\in K} h(z)\cdot h(y-z)=0\quad\text{for }y\ne 0.
\]
 From Proposition \ref{cl:Four1-3}(b) we have $h(-z)=\overline{h(z)}$, since $(-1)^f$ is a real-valued function. Then
\[
(h(z),h(z-y)) =\sum\limits_{z\in K} h(z)\cdot \overline{h(z-y)}=0\quad\text{for }y\ne 0.
\]
Since $h(\bar{0})=|K|^{1/2}$, the vectors are nonzero.
It follows that the family of functions $\{h(z-y) \mid y\in K\}$ is an orthogonal basis of the space $\mathbb{V}(K)$.

\subsection{Problems}

\begin{exercise}
Prove that all eigenvalues of any Cayley graph on the groups $\mathbb{Z}_2^n$ and $\mathbb{Z}_4^n$ are integers.
\end{exercise}

\begin{exercise}
Prove the equality $P_k[n,q](\lambda_m) |A_m|=P_m[n,q](\lambda_k) |A_k|$, where $\lambda_t=n(q-1)-tq$ and $A_t$ is the sphere of radius $t$ in $Q^n_q$.
\end{exercise}

\begin{exercise}
Prove the equality
$\sum\limits_{k=0}^tP_k[n,q](\lambda_m)=P_t[n-1,q](\lambda_{m-1})$.
\end{exercise}

\begin{exercise}
Prove the equality\\
$\sum\limits_{m=0}^nP_k[n,q](\lambda_m)=\frac{1}{q}{n+1 \choose
k+1}((q-1)^{k+1}-(-1)^{k+1})$.
\end{exercise}

\begin{exercise}\label{ex:Pnn}
Prove that $|P_n[n,q](\lambda_m)|=1$ for every $m=0,\dots,n$.
\end{exercise}

\begin{exercise}
Prove that for any $x\in Q_2^n$ the number of Boolean functions $f$ for which the Walsh--Hadamard coefficient $\widehat{(-1)^f}(x)$ equals zero is ${2^n \choose 2^{n-1}}$.
\end{exercise}

\begin{exercise}
Let $f:Q_2^n\rightarrow\mathbb{C}$ be a function that equals zero at vertices of odd weight. Prove that its Fourier coefficients coincide at antipodal points of the Boolean hypercube, i.e. $\widehat{f}(z)=\widehat{f}(z\oplus\bar{1})$.
\end{exercise}

\begin{exercise}
Prove that for any Boolean function $f:Q_2^n\rightarrow Q_2$ the $r$-th moments $\sum\limits_{z\in Q_2^n}\widehat{f}^r(z)\geq 0$ are nonnegative. 
Prove that, among all Boolean functions with fixed weight $\wt(f)$, the minimum of the $r$-th moment is attained by the characteristic function of the set of columns of a parity-check matrix of a code with minimum distance $d+1$.
\end{exercise}

\begin{exercise}
Prove that multiplication by the parity function $(-1)^{\wt (x)}$ transforms an eigenfunction of the Boolean hypercube with eigenvalue $\lambda$ into an eigenfunction with eigenvalue $-\lambda$.
Describe a transformation that turns perfect colorings of the Boolean hypercube with quotient matrix $\begin{pmatrix}
a & b\\
b & a
\end{pmatrix}$
into perfect colorings of the Boolean hypercube with quotient matrix
$\begin{pmatrix}
b & a\\
a & b
\end{pmatrix}$.
\end{exercise}

\begin{exercise}
Let $C\subset Q^n_2$ and suppose $\widehat{{\mathbf 1}_{C}}(z)\neq 0$ only if the weight $\wt(z)$ is odd or $z=\bar{0}$. Prove that $|C|=2^{n-1}$, except for the cases $C=Q^n_2$ and $C=\varnothing$.
\end{exercise}

\begin{exercise}\label{exer1}
Let $A_t$ be the set of vertices of weight $t$ in $Q^n_q$. Prove that
$\widehat{{\mathbf 1}_{A_t}}(z)=\frac{1}{q^{n/2}}P_t[n,q](n(q-1)-\wt(z)q)$.
\end{exercise}

\begin{exercise}\label{exer76}
Let $f:Q^n_2\rightarrow \mathbb{R}$ be an eigenfunction of the Boolean hypercube with eigenvalue $\lambda_k=n-2k$. Prove that
$|{\rm supp}({f})|\geq 2^{(n+|\lambda_k|)/2}$. Give examples of functions for which equality holds.
\end{exercise}

\begin{exercise}
Let $g,f:Q_q^n\rightarrow \mathbb{C}$. Define the matrix $A_f$ by $A_f(x,y)=f(x-y)$. Prove that $A_fA_g=A_{f*g}$. Prove that $\det A_f= q^{n/2}\prod\limits_{y\in Q_q^n}\widehat{f}(-y)$.
\end{exercise}

\begin{exercise}
Consider the Boolean $n$-cube 
whose edges are directed from vertices of smaller weight to adjacent vertices of larger weight.
The adjacency matrix of the directed graph is skew-symmetric; it contains $1$ in position $i,j$ if the edge is directed from vertex $i$ to vertex $j$, and $-1$ if the edge is directed from vertex $j$ to vertex $i$. Find the eigenvalues and eigenvectors of the adjacency matrix of the directed $n$-cube.
\end{exercise}

\begin{exercise}
Prove that if $f\in \mathbb{V}(K)$ is an eigenfunction of the graph $Cay(K,A)$, then $f*g$
is either an eigenfunction or the identically zero function of $Cay(K,A)$ for any function $g\in \mathbb{V}(K)$.
\end{exercise}

\begin{exercise}
Prove that if $f,g\in \mathbb{V}(K)$ are eigenfunctions of the graph $Cay(K,A)$ with different eigenvalues, then their convolution $f*g$ is the identically zero function.
\end{exercise}

\begin{exercise}\label{exer:Ess}
Let $f: Q^n_q\rightarrow \mathbb{C}$. Prove that $f$ depends essentially on the variable $x_i$ if and only if there exists a vertex $z\in Q^n_q$, $z_i\neq 0$, for which $\widehat{f}(z)\neq 0$.
\end{exercise}

\begin{exercise}
Let $f: Q^n_q\rightarrow \mathbb{C}$ take only two distinct values. Prove that there exist constants $c_1, c_2 \in \mathbb{C}$ such that $\widehat{f}*\widehat{f}=c_1 \widehat{f}+ c_2\delta$.
\end{exercise}

\begin{exercise}
Let $C\subset Q^n_p$, $p$ a prime. Prove that if $|C|$ is not divisible by $p$, then $\supp(\widehat{{\mathbf 1}_{C}})=Q^n_p$.
\end{exercise}

\section{Correlation Immune Functions and Orthogonal Arrays}

\subsection{Functions Balanced on Faces}

In this section, we will consider only the $q$-ary hypercubes $Q_q^n$. Accordingly, the characters have the form
$\phi_z(x)=\xi^{\langle x,z\rangle}$, where $\xi=e^\frac{2\pi i}{q}$.
However, many results in this section can be generalized to arbitrary Cayley graphs of finite abelian groups.

Recall that a face $\Gamma$ of dimension $\dim\,\Gamma = k$ in
the $n$-dimensional $q$-ary cube $Q_q^n$
is a subset of its vertices
which consists of vectors in which $n-k$
coordinates fixed while  the remaining $k$
coordinates are arbitrary.

Two faces of the same dimension are said to have the same
direction if they have the same sets of fixed and free coordinates.
 A face with exactly one fixed coordinate
is called a hyperface. The subgraph of the hypercube $Q_q^n$ induced by a face of dimension $k$ is isomorphic to the hypercube $Q_q^k$.

%\begin{definition}
    A function $f: Q_q^n\rightarrow \mathbb{C}$ is called
    {\sl balanced on faces of dimension} $k$,
    if for every face $\Gamma$ of dimension $k$,
    the equality $\sum\limits_{x\in\Gamma}f(x)=0$ holds.
%\end{definition}
Note that if a function is balanced on faces of dimension $k$, then
it is also balanced on faces of every dimension greater than $k$.

\begin{claim}\label{c:weighted0}
Let a function $f:Q_2^n\rightarrow \mathbb{C}$ be balanced on faces
of dimension $k$. Then $f$ is uniquely determined from its
values in the ball of radius $k-1$, i.e., from the values at vertices of weight
at most $k-1$.
\end{claim}
\begin{proof}
We prove the proposition by induction. Consider an arbitrary vertex $z$
of weight $k$. There is a face of dimension $k$ that contains $z$ and the all-zero vector.
Clearly, $z$ is the unique vertex
of weight at least $k$ in this face.  Since the function $f$ is balanced, value
$f(z)$ is uniquely determined. If $f$ is balanced on faces of
dimension $k$, then it is balanced on faces of dimension $k+i$, $i\geq 0$.
Therefore, after finding the values $f(z)$ for all vertices $z$ of weight
$k$, we can find the values $f(z)$ for all vertices of weight $k+1$, and
so on.
\end{proof}

The functions balanced on faces of a fixed dimension form a
linear subspace of $\mathbb{V}(Q^n_q)$. We next show that a basis
of this space can be formed from the characters of the group
$\mathbb{Z}_q^n$.

\begin{claim}\label{c:weighted}
    If $\wt(z)=k$, then the character $\phi_z$ is balanced on faces of dimension $n-k+1$
    and there exists a family of faces of dimension
     $(n-k)$ having the same direction such that $\phi_z$ is constant on all faces of that direction.
\end{claim}
\begin{proof}
Without loss of generality, assume that
$z=(z_1,\ldots,z_k,0,\ldots,0)$, $z_i\ne0$, $i=1,\dots,k$.
Consider an arbitrary face $\Gamma$ of dimension $n-k+1$.  Such a face necessarily
has at least one free coordinate $x_i$ among the first $k$. Without loss of generality,
assume that $i=1$. By Proposition \ref{c:char_sum11} we have
$\sum\limits_{x\in Q_q}\xi^{x z_i}=0$. Then

\[
\sum\limits_{x\in\Gamma}\phi_z(x)=
\sum\limits_{x_{i_2},\ldots,x_{i_{n-k}}}\sum\limits_{x_1}\phi_z(x)
=\sum\limits_{x_{i_2},\ldots,x_{i_{n-k}}}\xi^{\langle(z_2,  \ldots,
z_n),(x_2,\ldots, x_n)\rangle}\sum\limits_{x_1}\xi^{z_1 x_1} = 0.
\]

Next consider a face $\Gamma$ of dimension $n-k$ whose first $k$ coordinates are fixed: $x_1=a_1,\dots, x_k=a_k$. For
any vertex $x\in \Gamma$ the equality
\[
\phi_z(x)= \xi^{\langle(z_1, \ldots, z_k,0,\ldots,0),(x_1, \ldots, x_n)\rangle}
=\xi^{z_1a_1+
 \cdots +z_k a_k}
\]
holds.
\end{proof}

\begin{claim}\label{c:weighted1}
  A function $f: Q_q^n\rightarrow \mathbb{C}$ is balanced on faces of dimension
    $k$ if and only if
      $\widehat{f}(z)=0$ for all $z\in Q_q^n$ of weight $\wt(z)\le n-k$.
\end{claim}

\begin{proof}
 ($\Leftarrow$)     Assume that $\widehat{f}(z)=0$ for $\mathrm{wt}(z)\le n-k$.
    From Proposition \ref{c:fourier_basis} we have
    \[
    f(x)=\frac{1}{q^{n/2}}\sum\limits_{\mathrm{wt}(z)>n-k}\widehat{f}(z)\phi_z(x).
    \]

  By Proposition \ref{c:weighted}, all characters $\phi_z(x)$
  with $\wt(z)\ge n-k+1$ are balanced on faces of dimension $k$.

 ($\Rightarrow$)   Consider
    \[
    \widehat{f}(z)=\frac{1}{q^{n/2}}(f,\phi_z),\, \mbox{ where}\,
    \mathrm{wt}(z)\le n-k.
    \]
    By Proposition \ref{c:weighted}, $\phi_z$ is constant
    on faces $\Gamma$ of dimension $k$ (or greater than $k$ if $\wt(z)<n-k$) of some direction.
    Since $f$ is balanced on every face $\Gamma$, we have $\sum\limits_{x\in \Gamma} f(x)\overline{\phi_z(x)}=
    \overline{\phi_z(a)}\sum\limits_{x\in \Gamma} f(x)=0$. Since the vertex set of the hypercube is partitioned into the faces of any fixed direction, summing over all such faces yields
     $\widehat{f}(z)=0$.
\end{proof}

\subsection{Correlation Immune Functions: Definition and Applications}

Let $Im(f)$ be the set of values taken by the function $f$.

%\begin{definition}
    A function $f$ is called {\sl correlation immune
    of order} $k$ if, for every $b\in Im(f)$,  the value $b$ occurs equally often 
     in every face of dimension $n-k$.
%\end{definition}

If a function is uniformly distributed in faces of small
dimension, then this same property will also hold for faces of
larger dimension. For this reason, it is convenient to introduce the notation
$\cor(f)$ for the maximum order of correlation immunity
of the function $f$.

It follows from the definition of an MDS code (see Section \ref{MDS}) that for the characteristic function of an MDS code $C$ with code
distance $d$ we have that $\cor({\mathbf 1}_{C})=n-d+1$.  Clearly, if an MDS code exists, then this is the
minimum cardinality nonempty set with such correlation immunity.

\begin{example}
For a constant function in an $n$-dimensional hypercube we have  $\cor(const)=n$.
The Boolean parity function $f(x_1,\dots,x_n)=x_1+\cdots+x_n\mod 2$ it has
$\cor(f)=n-1$.
\end{example}

\begin{example}
Latin colorings (see Example \ref{ex:covercube}) in $Q^n_q$ are
correlation-immune functions of order $n-1$.
\end{example}

Correlation immune functions are used in cryptography to
simulate a random distribution of values at the vertices of the $q$-ary
hypercube.
Consider the problem of determining the value of a function when only part of its argument is known.
Suppose that an adversary knows the function
$f:Q^n_q\rightarrow Im(f)$ and the values of $k$ out of $n$ arguments. Their
aim is to determine the value of the function. If
$\cor(f)\geq k$, then in each face of dimension $n-k$ the function
takes any value $b$ with the same frequency as in the whole
hypercube. Therefore, knowing the values of $k$ variables provides no additional information about the value of the function. One can also consider the problem
where the adversary knows the value $f(x)=b$ and
the values of $k$ out of $n$ variables, and their aim is to determine
the remaining coordinates of the argument. If $\cor(f)\geq k$, then the information available to the
adversary does not help at all to find the values of the
remaining variables, since in each face of dimension $n-k$
the function takes any value $b$ with the same frequency.
The very name "correlation immune" comes
from the fact that the set of $k+1$ random variables (the value of the function and $k$
arbitrary variables) is independent (uncorrelated).

In cryptographic algorithms, Boolean-valued and even Boolean functions are most often considered.
The use of the Boolean parity check function, which possesses the maximum correlation immunity for non-constant functions, has limitations, since in this
case the frequencies of $0$ and $1$ are equal. Furthermore, an a priori unknown function
creates additional difficulties for the adversary. The possibility of
constructing sufficiently diverse correlation immune Boolean-valued functions with different proportions of ones will be considered
below.

\subsection{Correlation Immunity and the Fourier Transform}

\begin{claim}\label{c:weighted2}
    Let $f: Q_q^n\rightarrow \mathbb{C}$ and $\cor(f)=k$.
    Then the function $g=f-\frac{\mathrm{wt}(f)}{q^n}\mathbf{1}$, where $\mathrm{wt}(f)=\sum_x
    f(x)$,
    is balanced on faces of dimension $n-k$.
\end{claim}
\begin{proof}
    By assumption, the sum of the values of the function $f$ in each face of dimension
    $n-k$ is the same, and adding a constant function does not change this
    property. By definition, the function $g$ is balanced in the hypercube;
    consequently, it is balanced on faces of dimension $n-k$.
\end{proof}

\begin{theorem}[\cite{Carlet20}]\label{t:cor-im}
    For a function $f:Q_q^n\rightarrow \{0,1\}$, the equality $\cor(f)=k$ holds
    if and only if
    $\widehat{f}(z)=0$ for all $ z\in Q_q^n$ of weight $0<\mathrm{wt}(z)\le k$ and
    $\widehat{f}(v)\neq 0$ for some $v\in Q_q^n$ of weight $\mathrm{wt}(v)=k+1$.
\end{theorem}
\begin{proof}
    ($\Rightarrow$) Let $g=f-\frac{\mathrm{wt}(f)}{q^n}\mathbf{1}$.
 By Proposition \ref{c:weighted2}, the function $g$ is balanced
 on faces of dimension $n-k$. Then $\widehat{g}(z)=0$
  for all $ z\in Q_q^n$ of weight $\mathrm{wt}(z)\le k$ by Proposition \ref{c:weighted1}. The Fourier transform of $\mathbf{1}$ is a
  $\delta$-function; therefore, the functions $\widehat{f}$
  and $\widehat{g}$ differ only at $\bar{0}$.

  If $\widehat{f}(v)= 0$ for all $v\in Q_q^n$ of weight
  $\mathrm{wt}(v)=k+1$, then the function $g$ is balanced also on faces
  of dimension $n-k-1$ by Proposition \ref{c:weighted1}. Then the sums
  $\sum_{x\in \Gamma}f(x)$ are the same for all faces $\Gamma$
  of dimension $n-k-1$. Since $f$ takes only two values $0$ and $1$,
  it must take the value $1$ the same number of times in each face of dimension
  $n-k-1$. Then $\cor(f)\ge k+1$, a contradiction.

    ($\Leftarrow$) It is easy to see that $\widehat{g}(\bar{0})=0$. From the condition we have
    $\widehat{g}(z)=0$ for $0<\mathrm{wt}(z)\le k$. Then from Proposition
    \ref{c:weighted1}
    we obtain that $g$ is balanced on faces of dimension $n-k$.
    By Proposition \ref{c:weighted1}, the existence of a nonzero Fourier coefficient
    $\widehat{g}(v)\neq 0$ for some $v\in Q_q^n$ of weight $\mathrm{wt}(v)=k+1$ implies
    that $g$ is not balanced on some face
    of dimension $n-k-1$.
    Since the function $g$ takes only two
    values, we have
    $\cor(g)=k$. The functions $f$ and $g$ differ by a constant,
    hence $\cor(f)=\cor(g)$.
\end{proof}

In this theorem, only the cardinality of the image of the function matters.
 We will also apply the theorem   when $Im(f)=\{-1,1\}$.
If the function $f$ takes more than two distinct values (i.e., colors),
then $\cor(f)$ equals the minimum of the correlation immunities
of the characteristic functions of the colors. Indeed, if all colors
are uniformly distributed over faces of some dimension, then so is
each individual color. If at least one of the colors is not
uniformly distributed over faces of some dimension $k$, then the
function as a whole, by definition, is not correlation immune
of order $n-k$.

 Theorem \ref{t:cor-im} can be viewed as an
analogue of Proposition \ref{c:color_distr_proportion} and it is a corollary
of the orthogonality of the eigenspaces of the adjacency matrix of the graph with
distinct eigenvalues. Indeed, the characteristic
function of a face of dimension $n-k$ is a linear combination of
the characters $\phi_z$, where $\wt(z)\leq k$. Conversely, the
characters are linear combinations of the indicators of faces of the corresponding
dimension. Hence every linear combination of the characters
$\phi_z$ with $\wt(z)> k$ is orthogonal to a face of dimension $n-k$.
Conversely, every function orthogonal to all faces of dimension
$n-k$ can be represented as such a linear combination.

\subsection{Correlation Immunity of Perfect Colorings}

\begin{claim}\label{c:cor_immun_for_perfect_coloring0}
    Let $f$ be a perfect coloring of $Q^n_q$ with quotient
    matrix $S$. Then $\widehat{f}(z)\neq 0$ only if
$(n-\mathrm{wt}(z))(q-1)-\mathrm{wt}(z)$ is an eigenvalue of $S$.
\end{claim}
\begin{proof}
 The function $f$ can be
represented as a linear combination of eigenfunctions of 
$Q^n_q$ whose eigenvalues coincide with the eigenvalues of  $S$ (see
Corollary \ref{p:perf_color_properties1}). By
Corollary \ref{c:baschar}, the characters form a basis of every
eigenspace of $Q^n_q$. Therefore, every such eigenfunction is a linear combination of characters.  The character $\phi_z$ has the eigenvalue
$(n-\mathrm{wt}(z))(q-1)-\mathrm{wt}(z)$ (see Section \ref{10.2}).
\end{proof}

If $f$ is the characterictic function of a cell of an   equitable partition corresponding to a perfect coloring, then Proposition \ref{c:cor_immun_for_perfect_coloring0} remains true.

\begin{claim}\label{c:cor_immun_for_perfect_coloring01}
Let $f$ be a perfect coloring of  $Q^n_q$ and
$\lambda_1<n(q-1)$ the largest nontrivial eigenvalue
of its quotient matrix. Then $\cor(f)=
\frac{n(q-1)-\lambda_1}{q}-1$.
\end{claim}
\begin{proof}
Let $\{C_1, \dots, C_k\}$ be an equitable partition corresponding to $f$.  By Proposition \ref{p:perf_color_properties1}, for every $i$, $i=1,\dots,k$, 
${\bf 1}_{C_i} $ is a linear combination of eigenfunctions of $Q^n_q$ with eigenvalues from the spectrum of the quotient matrix.  Moreover, for every eigenvalue $\lambda$ of the quotient matrix there exists a cell $C_i$ such that  its linear expansion  has a nontrivial summand with eigenvalue $\lambda$ in the linear expansion. Thus, the proposition follows from Proposition \ref{c:cor_immun_for_perfect_coloring0}.
\end{proof}

\begin{claim}\label{c:cor_immun_for_perfect_coloring}
    Let $f$ be a perfect $2$-coloring of $Q^n_q$ with quotient matrix $S=\begin{pmatrix}
    n(q-1)-b & b\\
    c & n(q-1)-c
    \end{pmatrix}$. Then $\cor(f)=\frac{b+c}{q}-1$.
\end{claim}
\begin{proof}
The eigenvalues of $S$ are
$\lambda_0=n(q-1)$ and $\lambda_1=n(q-1)-(b+c)$.  Solving  $(n-\mathrm{wt}(z))(q-1)-\mathrm{wt}(z)=\lambda_0=n(q-1)$ gives 
$\mathrm{wt}(z)=0$, while solving
$(n-\mathrm{wt}(z))(q-1)-\mathrm{wt}(z)=\lambda_1=n(q-1)-(b+c)$,
gives $\mathrm{wt}(z)=\frac{b+c}{q}$. Then by Theorem
\ref{t:cor-im} and Proposition \ref{c:cor_immun_for_perfect_coloring0}
we obtain $\cor(f)=\frac{b+c}{q}-1$.
\end{proof}

\begin{corollary}
For the characteristic function $f$ of a 1-perfect code in 
$Q_q^n$, the equality $\cor(f)= \frac{(n-1)(q-1)}{q}$ holds.
\end{corollary}

Constructing perfect $2$-colorings is one method of obtaining correlation immune functions with prescribed parameters.
 It will be shown below that correlation immune functions
of the maximum possible order for some parameters are perfect
colorings.

The following theorem is a generalization of Theorem \ref{thAvg}
on recovering a perfect coloring from the average layer. Here, however,
two layers are required to recover the function. A generalization
of Theorem \ref{thAvg} on recovering a perfect coloring from one
layer under other additional conditions is given in \cite{Vas12}.

\begin{theorem}\label{th:2level}
Let $f$ be a perfect coloring of  $Q^n_2$ with
quotient matrix $S$. Let $v$ be an eigenvector of $S$
corresponding to the eigenvalue $\lambda_k$ of $S$,
where $-n<\lambda_k=n-2k<0$. If all entries of  $v$
are distinct, then $f$ is uniquely determined by its values on
vectors of weight $n-k$ and $n-k-1$.
\end{theorem}
\begin{proof}
By Remark \ref{z:colorlin0}, there exists an eigenfunction $f'$
of the hypercube with eigenvalue $\lambda_k$, whose distinct values
correspond to distinct colors of $f$. Since $f'$ is an
eigenfunction, $\widehat{f'}(z)=0$ for $\wt(z)\neq k$ by
Proposition \ref{claim:ef}. Then by Propositions \ref{c:weighted0}
and \ref{c:weighted1}, $f'$ is uniquely determined by its
values at vertices of weight at most $n-k$. By Proposition
\ref{c:dopeigenfun}, recovering $f'$ inside the ball of radius
$n-k$ requires only the values of $f'$ on the spheres of radii $n-k$ and
$n-k-1$, provided that $n-k<n/2$. The latter inequality holds by the
hypothesis.
\end{proof}

Note that in Theorem \ref{th:2level}, the condition $-n<\lambda_k=n-2k<0$
can be replaced by $0<\lambda_k=n-2k<n$, and the spheres of radius
$n-k$ and $n-k-1$ by the spheres of radii $k$ and $k-1$, since
the Hadamard (pointwise) multiplication of an eigenfunction with eigenvalue
$\lambda$  and the function $g(x)=(-1)^{\wt(x)}$ yields an eigenfunction with eigenvalue
$-\lambda$. This follows from the equality
$\widehat{((-1)^{\wt(x)}f)}(y)=\widehat{f}(y\oplus \bar 1)$.

\subsection{Orthogonal Arrays and Rao's Inequality}

%\begin{definition}
An {\sl orthogonal array} of type $OA(N, k, q, t)$  is an $N\times k$ array whose entries come from a fixed finite set of $q$ symbols (without loss of generality,  $Q_q$),   such that in every subset of $t$ columns of the array, every ordered $t$-tuple of symbols  occurs the same number of times. The number $t$ is called the {\sl strength} of the orthogonal array. The number of repeats is usually denoted by $\lambda$. Obviously, $\lambda=N/q^t$.
%\end{definition}

An orthogonal array of general form may contain repeated
rows and corresponds to the value table of a function $f:Q_q^n\rightarrow
\{0,1,\dots,m\}$: if $f(x)=s$, then the row $x$ appears in the array
$s$ times. By the definition of  orthogonal arrays, the function $f$ differs from a function balanced on
$(n-k)$-dimensional faces by a constant.

A simple orthogonal array   (i.e., without repeated rows) of strength $k$ is a table
whose set of rows is precisely the support of $f$,
where $f:Q_q^n\rightarrow \{0,1\}$
is a correlation immune function of order $k$.

An { orthogonal array}  of type $OA(N, k, q, t)$   of minimum cardinality $N=q^t$  is precisely the list of codewords of an MDS code with length $k$ and minimum code distance $k-t+1$ (see Section \ref{MDS}).

\begin{example}
A binary orthogonal array of strength $2$:

$\begin{array}{lllllll}
1& 1 &1 &1 &1 &1 & 1\\
1 &1 &1 &0 &0 &0 &0 \\
1 & 0 &0 &1 &1 & 0& 0 \\
1 &0 &0 &0 &0 &1 &1\\
0 &1& 0 &1 &0& 1 &0\\
0 & 1 &0 &0 &1 &0 &1\\
0 &0 &1 &1 &0 &0 &1\\
0 &0 &1 &0 &1 &1 &0
\end{array}.$
\end{example}

Correlation immune functions in the form of orthogonal arrays
are used for experiment design. Suppose we have a
program with $n$ parameters, each of which can take
$q$ values. To test the program for all
possible parameter sets, one would need to run $q^n$ computational
experiments, which may require unreasonably many resources.
Suppose it is sufficient  to verify the program's operation for any
values of any $k$ parameters. Then we need to set
the parameters as rows of an orthogonal array of strength $k$ of minimal
possible size. A similar solution can be applied not only
for testing the correctness of the program, but also when collecting
data to evaluate the performance of the program. Of course, to solve such a problem, one needs orthogonal arrays of as small a
size as possible or, equivalently, correlation immune functions with as
small a support size as possible.

Several theorems proved below for Boolean-valued correlation immune functions can be generalized to functions taking nonnegative integer values and to the corresponding orthogonal arrays.

\begin{theorem}[Rao's inequality \cite{Rao}]\quad

Let $f:Q_q^n\rightarrow \{0,1\}$, $f\neq \mathbf{0}$ and $k=\cor(f)$. Then,\\
 if $k$ is even, then $|\supp(f)|\geq
 \sum\limits_{i=0}^{k/2}{ n \choose i}(q-1)^i$;\\
 if $k$ is odd, then
$|\supp(f)|\geq
 { n-1 \choose (k-1)/{2}}(q-1)^{\frac{k+1}{2}}+
 \sum\limits_{i=0}^{(k-1)/2}{ n \choose i}(q-1)^i$.
\end{theorem}
\begin{proof}
Since $f$ takes only values $0$ and $1$, and $|\phi_u(x)|=1$ for every $x\in Q_q^n$, 
%and the characters take values of the form $\xi^a$, for any $u,v\in Q_q^n$
we have
\[
(f\phi_v,f\phi_u)=(f\phi_v,\phi_u)=(f,\phi_{-v}\phi_u)=
(f,\phi_{u-v})=q^{n/2}\widehat{f}(u-v).
\]

From Theorem \ref{t:cor-im} it follows that $\widehat{f}(z)=0$ for all
$z\in Q_q^n$ of weight $0<\mathrm{wt}(z)\le k$. Consequently,   all functions $f\phi_v$ for $0\le \mathrm{wt}(v)\le k/2$
are pairwise orthogonal.

 There can be at most $N$ pairwise orthogonal vectors
of length $N$. We have
$\supp(f\phi_v)=\supp(f)$. Therefore, the number of distinct vectors $f\phi_v$, where
$0\le \mathrm{wt}(v)\le k/2$, is at most $|\supp(f)|$. Hence
for even $k$ we obtain the inequality $|\supp(f)|\geq
 \sum\limits_{i=0}^{k/2}{ n \choose i}(q-1)^i$.

 In the case when $k$ is odd, we can consider all functions
$f\phi_v$ for $0\le \mathrm{wt}(v)\le (k-1)/2$ and all functions
$f\phi_u$, $u_1\neq 0$, $\mathrm{wt}(u)= (k+1)/2$. Then the weight of the differences of any two such vectors  does not exceed $k$.
\end{proof}

 Note that for any Boolean-valued function $f$, the number
$q^{\cor(f)}$ divides $|\supp(f)|$, since in each face of
dimension $n-\cor(f)$ contains the same number of ones. 

\subsection{Total Influence  of Boolean Functions}

Let $f:Q_q^n\rightarrow \{0,1\}$ be a characteristic function of a set $C\subset Q_q^n$. Denote by $I[f]$ the number
of edges whose endpoints have different values of $f$.
For a constant function this quantity is zero, while for the Boolean
parity function it is $n2^{n-1}$. Note that $I[f]$ equals the size of the cut
$I[f]=e(f^{-1}(0), f^{-1}(1))$ in the Hamming graph.

\begin{claim}[Nisan and Szegedy\cite{NiSz}]\label{cl:I}
$I[f]=\frac{q}{4}\sum\limits_{z\in Q_q^n}\wt(z)
|\widehat{(-1)^f}(z)|^2$.
\end{claim}
\begin{proof}
Let $M$ be the adjacency matrix of $Q^n_q$.  It is straightforward to verify that
 \begin{equation}\label{e:BF21}
2I[f]-(n(q-1)q^n-2I[f])=-(M(-1)^f,(-1)^f),
\end{equation}
since in the right-hand sum each edge is counted twice.

Recall that the characters are eigenfunctions of the matrix $M$, i.e.,
$M\phi_z=\lambda_z\phi_z$, where $\lambda_z=(q-1)n-q\wt(z)$
(see Section  \ref{10.3}).
We represent the function $(-1)^f$ as a linear combination of characters:
$(-1)^f=\frac{1}{q^{n/2}}\sum\limits_{z\in
Q_q^n}\widehat{(-1)^f}(z)\phi_z$. Substituting this expression into
(\ref{e:BF21}) and using the orthogonality of the characters, as well as
Parseval's equality $\sum\limits_{z\in
Q_q^n}|\widehat{(-1)^f}(z)|^2=q^n$, we obtain the equalities
$$4I[f]=n(q-1)q^n-\sum\limits_{z\in
Q_q^n}\lambda_z|\widehat{(-1)^f}(z)|^2=\sum\limits_{z\in
Q_q^n}q\wt(z) |\widehat{(-1)^f}(z)|^2.$$
\end{proof}

Define the density $\rho(f)=\frac{|C|}{q^n}$ of an indicator function $f={\bf 1}_C$.

\begin{theorem}\label{cor:X}   
  Let $f:Q_q^n\rightarrow \{0,1\}$, $f\neq \mathbf{0}$.\\
{\rm (a)} $\frac{I[f]}{q^{n+1}}\geq
(\cor(f)+1)(\rho(f)-\rho^2(f))$.\\
{\rm (b)} $\frac{I[f]}{q^{n+1}}=(\cor(f)+1)(\rho(f)-\rho^2(f))$ if and only  if $f$
is a perfect $2$-coloring.
\end{theorem}
\begin{proof}
(a) Since $(-1)^f={\bf 1}-2f$, the coefficients $\widehat{(-1)^f}(z)$ and
$\widehat{f}(z)$ vanish simultaneously for $\wt(z)>0$. By
Theorem \ref{t:cor-im} it follows that $\widehat{(-1)^f}(z)=0$ for
$0<\mathrm{wt}(z)\le \cor(f)$. Hence, from Proposition
\ref{cl:I} we have the inequality
\begin{equation}\label{eq:fei}
4I[f]\geq\sum\limits_{z\in
Q_q^n\setminus {\bar 0}}q(\cor(f)+1) |\widehat{(-1)^f}(z)|^2.
\end{equation}
Next we estimate the sum $\sum\limits_{z\in Q_2^n\setminus {\bar 0}}
|\widehat{(-1)^f}(z)|^2$. From Parseval's equality we have
$\sum\limits_{z\in Q_q^n} |\widehat{(-1)^f}(z)|^2 =q^n$, and moreover
$\widehat{(-1)^f}({\bar 0})=q^{n/2}(1-2\rho(f))$. Then
\begin{equation}\label{eq:fei1}
\sum\limits_{z\in Q_q^n\setminus {\bar 0}} |\widehat{(-1)^f}(z)|^2=
q^n-q^n(1-2\rho(f))^2=4q^{n}\rho(f)(1-\rho(f)).
\end{equation}
(b) Obviously, equality in   (\ref{eq:fei}) holds if and only if $(-1)^f-(\widehat{(-1)^f}(0))\phi_{\bar 0}$ is an eigenfunction with eigenvalue
$(q-1)n-q(\cor(f)+1)$. By Propositions \ref{c:colorlin}  and \ref{c:cor_immun_for_perfect_coloring}  it holds if and only if the function  $(-1)^f$ is a perfect $2$-coloring. 
\end{proof}

Finally, we consider Boolean functions. 
Let  $I_i[f]$ denote the probability  that flipping the $i$-th input bit changes the output value of  $f$.
The sum $\bar{I}[f]=\sum_i I_i[f]$  is called the {\sl total influence} (or {\sl average sensitivity}) of $f$.
It is easy to see that $\bar{I}[f]=I[f]/2^{n-1}$.
Consider also the normalized Walsh--Hadamard coefficients $\bar{W}_f(u)=\widehat{(-1)^f}(u)/2^{n/2}$. By Parseval's identity we have  $\sum_z\bar{W}_f^2(z)=1$  for every Boolean function $f:Q_2^n\rightarrow \{0,1\}$.
Hence the values  $\bar{W}_f^2(z)$ define a 
probability distribution on  $Q_2^n$. The {\sl Fourier entropy} of $f$ is defined by $H(f)=\sum\limits_{z\in
Q_2^n}\bar{W}_f^2(z)\log_2\frac{1}{\bar{W}_f^2(z)}$. By  Proposition \ref{cl:I} we have
$\bar{I}[f]=\sum\limits_{z\in Q_2^n}\bar{W}_f^2(z)wt(z)$.

The  celebrated  Fourier Entropy–Influence (FEI) conjecture was introduced by Friedgut and Kalai
 in 1996 \cite{FG}. The conjecture asserts the existence of a universal constant $C$ such that  
 $H[f]\leq C \cdot\bar{I}[f]$ for every
Boolean function $f$. \label{FEIprob}

\subsection{Bierbrauer--Friedman Inequality}

We derive Bierbrauer--Friedman inequality on the size of orthogonal arrays as a corollary of Theorem \ref{cor:X}.

\begin{theorem}[Bierbrauer--Friedman inequality \cite{Bier}, \cite{Friedman}]\quad\label{thBF} \quad \\
% --- proved for the binary case in '92, for the q-ary case in '98
    Let $f:Q_q^n\rightarrow \{0,1\}$, $f\neq \mathbf{0}$.\\
{\rm (a)}  Then  $\rho(f)\ge 1-\frac{n(q-1)}{q(\cor(f)+1)}$.\\
{\rm (b)} If $\rho(f)=1-\frac{n(q-1)}{q(\cor(f)+1)}$ 
    then $f$ is a perfect $2$-coloring.
\end{theorem}
\begin{proof}
Without loss of generality, we assume that $\rho(f) \leq 1/2$.\\
$\frac{I[f]}{q^{n}}\leq
(q-1)\rho(f) n$, since each vertex with value $1$ is adjacent to at most $n(q-1)$ vertices with value $0$.  By Theorem~\ref{cor:X}(a)  we obtain
$$\frac{ (q-1)\rho(f) n}{q}\geq  \frac{I[f]}{q^{n+1}}\geq (\cor(f)+1)(\rho(f)-\rho^2(f)).$$
This inequality is equivalent to part  (a). Part (b) follows directly from Theorem~\ref{cor:X}(b).
\end{proof}

Note that equality in the Bierbrauer--Friedman inequality cannot
be attained by a non-simple orthogonal array, since in that case
$$    (f,f)= \sum\limits_{x\in Q^n_q}|f(x)|^2>|\supp(f)|, $$ and, therefore,  equality (\ref{eq:fei1}) does not hold.

\begin{corollary}[Ostergard, Pottonen, Phelps, \cite{OPPh}]\quad
    Let $f:Q_q^n\rightarrow \{0,1\}$. The function $f$ is the
    characteristic function of a
    $1$-perfect code if and only if
    $\rho=\frac{1}{n(q-1)+1}$ and $\cor(f)=\frac{n(q-1)+1}{q}-1$.
\end{corollary}
\begin{proof}
    ($\Rightarrow$) The equality $\rho=\frac{1}{n(q-1)+1}$ follows immediately from the definition of a $1$-perfect code. 
    The value of $\cor(f)$ is obtained from Claim~\ref{c:cor_immun_for_perfect_coloring}
    by substituting the specific values $b=1$, $c=n(q-1)$.
    %//true, Delsarte proved it differently :)

    ($\Leftarrow$) Substituting these values into the Bierbrauer--Friedman inequality gives
       \[
    \frac{1}{n(q-1)+1}\ge 1-\frac{n(q-1)}{q\cdot
    \frac{n(q-1)+1}{q}}, \quad \mbox{i.e.,}   \quad  
    \frac{1}{n(q-1)+1}\ge \frac{1}{n(q-1)+1}.
    \]
    Hence, equality holds.
    Then by Theorem~\ref{thBF}(b)
    we obtain that $f$ is a perfect $2$-coloring.

    The value of the density implies that, on average, each radius-one ball contains exactly one vertex $x$ such that $f(x)=1$.
        Since $f$ is a perfect $2$-coloring, every radius-one ball contains vertices of both colors.
    Consequently, each ball contains exactly one $1$, and $f$ is the characteristic function of a $1$-perfect code. %//Ostergard-Pottonen-Phelps proved it differently, via the MacWilliams transform.
\end{proof}

It follows from the proof that equality in the Bierbrauer--Friedman inequality can occur only when
$(Mf,f)=0$, that is, when the perfect $2$-coloring $f$ has no adjacent vertices colored $1$.
However, one can show that every perfect $2$-coloring
of $Q^n_q$ is extremal in the sense of attaining
equality in some inequality similar to the
Bierbrauer--Friedman inequality.

\begin{claim}[\cite{Pot12}]
    Let $f:Q_q^n\rightarrow \{0,1\}$. The following inequality holds:
    \begin{equation}\label{e:BF1}
    \rho(f)\cdot q(\cor(f)+1)\le \bar s_{01}(f),
    \end{equation}
    where
$\bar s_{01}(f)$  denotes the average number of neighbors of color $1$ of a vertex of color $0$ (see Section \ref{3.3}). Equality
(\ref{e:BF1}) holds if and only if $f$
is a perfect coloring.
\end{claim}

The proof follows the same argument as that of the Bierbrauer--Friedman inequality.

\subsection{Fon-Der-Flaass Theorem}\label{11.6}

%\begin{definition}
    A Boolean correlation-immune function $f$ is called
    {\sl resilient} if $(-1)^f$ is balanced,
    i.e., $f$ takes the values  $0$ and $1$ equally often.
%\end{definition}
It is easy to see that a Boolean function $f$ is resilient only if $\widehat{(-1)^f}(\bar{0})=0$.

\begin{theorem}[Fon-Der-Flaass \cite{FdF2}]\label{thFdF}
    If a Boolean function $f: Q_2^n\rightarrow \{0,1\}$
    is neither constant nor resilient, then $\cor(f)\le \frac{2n}{3}-1$.
    Moreover, if equality is attained in this inequality,
    then $f$ is a perfect $2$-coloring.
\end{theorem}
\begin{proof}
    Let $g=\widehat{(-1)^f}$. Since $f$ is not constant,
    there exists $y\ne \bar{0}$ such that $g(y)=\widehat{(-1)^f}(y)\neq 0$. Since $f$
    is not resilient,
    $g(\bar{0})=\widehat{(-1)^f}(\bar{0})\neq 0$. By 
   Theorem \ref{th:Tits}, for the chosen $y$
    we have
    \[
    0=(g*g)(y)=\sum\limits_{x\in Q_2^n}g(x)g(y-x).
    \]
    Hence the sum must contain at least one additional nonzero term besides the two terms
     $g(\bar{0})g(y)$ and $g(y)g(\bar{0})$ since these two terms alone cannot sum to zero. Note that $\cor(f)=\cor((-1)^f)$.
    By Theorem \ref{t:cor-im} we have $g(x)=0$ for $0<\mathrm{wt}(x)\le
    \cor(f)$.
    Therefore,
    \[
    0=(g*g)(y)=2g(\bar{0})g(y)+\sum_{\substack{x\in Q_2^n\\ \mathrm{wt}(x)\geq\cor(f)+1\\
    \mathrm{wt}(y-x)\geq\cor(f)+1 }}g(x)g(y-x),
    \]
    and the sum  contains nonzero terms.
    Assume that $\cor(f)>\frac{2n}{3}-1$. Then
    $\mathrm{wt}(y)>\frac{2n}{3}$ and there exists $x$ for which $g(x)g(y-x)\neq 0$.
    From this it follows that $\mathrm{wt}(x)>\frac{2n}{3}$ and
    $\mathrm{wt}(y-x)>\frac{2n}{3}$. However, from the triangle
    inequality for the Hamming metric we have
    $\mathrm{wt}(y-x)=d(x,y)\leq
    d(\overline{1},x)+d(\overline{1},y)= n-\mathrm{wt}(x)+
    n-\mathrm{wt}(y)<\frac{2n}{3}$. Contradiction.

    Now let $\cor(f)=\frac{2n}{3}-1$. If $g(y)\neq 0$ at some vertex
    of weight $\mathrm{wt}(y)>\frac{2n}{3}$, then we obtain a contradiction in the same manner.
    Hence, $g$ is supported only at $\bar 0$ and at 
    vertices of weight
    $\frac{2n}{3}$. Then
    $f$ is a perfect coloring by
    Proposition \ref{c:colorlin}.
\end{proof}

The inequality in Theorem \ref{thFdF} is called the Fon-Der-Flaass
bound. The proof given above is due to
Khalyavin. As shown by Khalyavin \cite{Halyav}, the Fon-Der-Flaass
bound also remains valid  for binary orthogonal arrays.
In contrast to the Rao and Bierbrauer--Friedman inequalities,
no analogue of the Fon-Der-Flaass theorem  is currently  known for the $q$-ary hypercubes 
with $q\neq 2$.
The only known extension is to Q-polynomial distance-regular graphs (see \cite{CRCKP}).

\subsection{Problems}

\begin{exercise}
Let $f:Q_2^n\rightarrow \mathbb{C}$ and $\cor(f)=n-k$. Prove
that $f$ is uniquely determined by its values in the ball
of radius $k-1$ (i.e., from the values at vertices of weight at most $k-1$).
\end{exercise}

\begin{exercise}
Let $f:Q_q^n\rightarrow \mathbb{C}$ and $\cor(f)=n-k$. Prove
that $f$ is uniquely determined by its values outside the ball
of radius $n-k$ (i.e., from the values at vertices of weight at least
$n-k+1$).
\end{exercise}

\begin{exercise}
Let $q$ be prime. Regard $Q^n_q$ as the $n$-dimensional vector
space over the field $GF(q)$. Prove that if a function
$f:Q^n_q\rightarrow \mathbb{C}$ is balanced on all hyperplanes,
then $f$ is identically zero.
\end{exercise}

\begin{exercise}
Prove that the support size of a function $f:Q_2^n\rightarrow
\mathbb{C}$ that is balanced on faces of dimension $k$ is at least
$2^{n-k+1}$.
\end{exercise}

\begin{exercise}
Prove that the functions $f:Q_q^n\rightarrow \mathbb{C}$
balanced on faces of dimension $k$ form a linear
subspace, and that this subspace is the orthogonal
complement of the  span $L(n,k,q)$ of the set of indicators
of all  $k$- dimensional faces.
\end{exercise}

\begin{exercise}
Let $f\in L(n,k_1,q)$ and $g\in L(n,k_2,q)$. Prove that\\ $f\cdot g\in
L(n,k_2+k_1-n,q)$ if $k_2+k_1>n$.
\end{exercise}

\begin{exercise}
Let $f$ be a perfect coloring of $Q_2^n$, and let $\lambda$ be the largest eigenvalue of its quotient matrix other than  $n$. Prove that
$\cor(f)=\frac{n-\lambda}{2}-1$.
\end{exercise}

\begin{exercise}
Let $f$ be some coloring of $Q_2^n$, and let $\lambda$ be the largest eigenvalue of  the average quotient matrix $S_f$ (see Section \ref{3.3}) other than  $n$. Prove that $\cor(f)\leq\frac{n-\lambda}{2}-1$.
\end{exercise}

\begin{exercise}
A {\sl covering array of strength} $k$ is a table
whose rows are tuples from $P\subset Q^n_q$ such that the set
$P$ intersects all $k$-dimensional faces in $Q^n_q$. Prove
that the number of rows of a covering array of strength $2$ in $Q^n_2$ is at least
$\lceil\log_2 n\rceil$.
\end{exercise}

\begin{exercise}
Let a Boolean function $f: Q_2^n\rightarrow \{0,1\}$ be neither constant nor resilient and attain equality in the Fon-Der-Flaass  bound. Prove that in every face of dimension greater than $n/3$, the numbers of $x$ on which  $f(x)=1$ among even-weight and odd-weight vertices are  equal.
\end{exercise}

\section{Dual Code and Sarkar Identity}

\subsection{Dual Code}\label{12.1}

In this section, we regard $Q_q$ as the ring  $\mathbb{Z}_q$. So, for $u,v \in Q^n_q$
we write $v\perp u$ if $\langle v,u\rangle=\sum v_iu_i=0 \mod q$. Moreover,  for a code $C\subset Q_q^n$ we write   $v\perp C\subset Q_q^n$ if for all $u\in C$
the equality $\langle v,u\rangle=\sum v_iu_i=0 \mod q$ holds.
%\begin{definition}
 The set
$C^\perp=\{v\in Q^n_q \ |\ v\perp C\}$ is called the {\sl  dual code} of $C$.
%\end{definition}
It is easy to verify that the set $C^\perp$ is closed under
addition and multiplication in the ring $\mathbb{Z}_q$, and 
 $C\subseteq (C^\perp)^\perp$ by definition.
A code $C\subset Q_q^n$ is called {\sl additive} if it is closed under
addition and scalar multiplication in the ring $\mathbb{Z}_q$. In particular, an additive code in $Q^n_q$ is a subgroup of $\mathbb{Z}_q^n$. Moreover, for every additive code $C\subset \mathbb{Z}_q^n$ the dual code $C^\perp$ coincides with 
the subgroup $C^\perp\subset X(\mathbb{Z}_q^n)\simeq \mathbb{Z}_q^n $.

\begin{remark}
It is well-known that a field structure can be defined on a set of cardinality $q$ if and only if
$q=p^t$
is a prime power. For $t>1$ the addition and multiplication operations in the field $GF(p^t)$
will
not coincide with addition and multiplication modulo $q$.
Therefore, the definition of the dual code used here and in Section \ref{8.3}  coincides
only for prime $q$.
\end{remark}

\begin{claim}
Let $\Gamma$ be a $k$-dimensional face in $Q^n_q$, with $\bar{0}\in
\Gamma$, and let its free coordinates be indexed by  $I\subset \{1,\dots,n\}$. Then $\Gamma^{\perp}$ is an $(n-k)$-dimensional face in
$Q^n_q$, with $\bar{0}\in \Gamma^{\perp}$ and
$\{1,\dots,n\}\setminus I$
being the indices of its free coordinates.
\end{claim}
\begin{proof}
The sets of free coordinates of the two faces are disjoint. Therefore,
for any $i\in\{1,\dots,n\}$ and for any vertex $u$ of one face and  any vertex $v$
of the other face, $u_i=0$ or $v_i=0$. It remains to show that any tuple
$v$ having a nonzero coordinate $v_i\neq 0$, $i\in I$, is not
contained in $\Gamma^{\perp}$. Indeed, we have $e_i\in \Gamma$
and $\langle v,e_i\rangle=v_i\neq 0$.
\end{proof}

\begin{claim}\label{c:dual}
$v\perp C$ if and only if $\widehat{{\mathbf 1}_{C}}(v)=\frac{|C|}{q^{n/2}}$.
\end{claim}
\begin{proof}
We have the equality
$\widehat{{\mathbf 1}_{C}}(v)=\frac{1}{q^{n/2}}\sum\limits_{x\in
    C}\overline{\phi_v(x)}= \frac{1}{q^{n/2}}\sum\limits_{x\in
    C}\xi^{-\langle v, x\rangle}$. It is clear that $\sum\limits_{x\in
    C}\xi^{-\langle v, x\rangle}=|C|$
if and only if all terms are equal to $1$. And this is equivalent to
the condition $v\perp C$.
\end{proof}

\begin{claim}\label{cor:dualcode}
Let $C\subset Q_q^n$ be a cell of an equitable partition and let $v\in C^\perp$. Then $\lambda=n(q-1)-\wt(v)q$ is an eigenvalue of its quotient  matrix.
\end{claim}
\begin{proof}
From Proposition \ref{c:dual} we have $\widehat{{\mathbf 1}_{C}}(v)=\frac{|C|}{q^{n/2}}\neq 0$. By Proposition \ref{c:cor_immun_for_perfect_coloring0}  $\lambda=n(q-1)-\wt(v)q$ is an eigenvalue of the   quotient  matrix.
\end{proof}

\begin{corollary}\label{cor:face}
For any  $\Gamma\subset Q^n_q$, $\bar{0}\in \Gamma$, the equality
\[
    \widehat{{\mathbf 1}_{\Gamma}}=\frac{q^{\mathrm{dim}\,\Gamma}}{q^{n/2}}{\mathbf 1}_{{\Gamma^\perp}}
    \]
holds.
\end{corollary}
\begin{proof}
$\widehat{{\mathbf 1}_{\Gamma}}(v)=\frac{q^{\mathrm{dim}\,\Gamma}}{q^{n/2}}$
for $v\in \Gamma^\perp$ by Proposition \ref{c:dual}. From Parseval's
identity together with the fact $|\Gamma||\Gamma^\perp|=q^n$  it follows that 
$\widehat{{\mathbf 1}_{\Gamma}}(u)=0$ for $u\not\in
\Gamma^\perp$. \end{proof}

Every face in the hypercube can be represented as $\Gamma +a$,
where the face $\Gamma$ contains zero. We now determine the Fourier transform
of the characteristic function of an arbitrary face. By the properties
of the Fourier transform (Proposition \ref{cl:Four1-3}(c)) it follows that
$(\widehat{{\mathbf 1}_\Gamma(x+a)})(z)=
\frac{\phi_a(z)q^{\mathrm{dim}\,\Gamma}}{q^{n/2}}{\mathbf 1}_{\Gamma^\perp}(z)$.
Thus, for any $z\in \supp(\widehat{{\mathbf 1}_{\Gamma+a}})$
we have $\wt(z)\leq \dim(\Gamma^\perp)=n-\dim(\Gamma)$.
Note that Theorem \ref{t:cor-im}
    is an immediate consequence of Corollary  \ref{cor:face} and 
  Proposition
\ref{c:color_distr_proportion}.

If $q$
is prime, then $Q_q$ can be regarded as the finite field
$GF(q)$ with addition and multiplication modulo $q$. 
It is well known from linear algebra that the solution set of a homogeneous system of equations
\begin{equation}\label{eq:check_matrix}
C=\{x\in Q_q^n | Hx=\bar{0}\}
\end{equation}
 is a linear
subspace, and conversely, any linear subspace $C\subset
Q_q^n$ can be represented as the solution set of a system of equations
(\ref{eq:check_matrix}). Recall that the matrix $H=\{h_{ij}\}$ is called a
parity-check matrix of the code $C$.  The rows of $H$ belong to the
linear code $C^\perp$ and, when they are linearly independent,
form its basis. For the linear code $C^\perp$, the matrix $H$
is called a {\sl generator matrix}.

It is easy to see that the system of equations (\ref{eq:check_matrix}) determines an additive code over
the  ring $\mathbb{Z}_q$ for arbitrary $q$ and rows of $H$ generate the dual code. 

\begin{example}
Let $C\subset \mathbb{Z}_4^2$, $C=\{ (0,0), (2,2)\}$. By definition $C$ is an additive code.  $C^\perp=\{(0,0),(1,1), (2,2), (3,3), (1,3), (3,1)\}$ is also an additive code.
\end{example}

Next, we calculate the minimum distance of a binary code in terms of its generator matrix. 
 Let $H$ be a $t\times n$ binary generator matrix, let $C=\{uH : u\in Q_2^t\}$, and let $M(H)$ be the multiset of columns of $H$.  Define
\[
W_{M(H)}(u)=\sum_{a\in M(H)}(-1)^{\langle u, a\rangle}.
\]
\begin{claim}\label{p3_sidon}
The minimum distance of $C$ is
\[
 d(C)=\frac{n-\max_{u\ne0}W_{M(H)}(u)}{2}.
\]
\end{claim}
\begin{proof}
For a fixed $u\ne0$, a coordinate is one exactly when $\langle u, a\rangle=1$.  Hence
$\mathrm{wt}(uH)=\frac12\sum_{a\in M(H)}(1-(-1)^{u\cdot a})$.
Taking the minimum over nonzero $u$ gives the formula.
\end{proof}

\begin{claim}\label{cor:doubledual}
 If $C$ is an additive  code in $Q_q^n$, then
$C=(C^\perp)^\perp$.
\end{claim}
\begin{proof}
As mentioned above, $C\subseteq (C^\perp)^\perp\subset C$   for every additive code $C$. 
By Corollary \ref{claim:factor1} it follows $|(C^\perp)^\perp|= |C|$.
\end{proof}

\begin{corollary}\label{cor:face2}
Let $q$ be prime.  A code $C\subset Q_q^n$ is linear if and only if
$\widehat{{\mathbf 1}_{C}}=\frac{|C|}{q^{n/2}}{\mathbf 1}_{C^\perp}$.
\end{corollary}
\begin{proof}
From the definitions of linear and dual codes we have $\dim C+ \dim
C^\perp=n$. Hence $|C||C^\perp|=q^n$. Then necessity follows from Proposition
\ref{c:dual} and Parseval's equality.

From the equalities
$\widehat{{\mathbf 1}_{C}}=\frac{|C|}{q^{n/2}}{\mathbf 1}_{C^\perp}$  and $\widehat{\widehat{{\mathbf 1}_C}}={\mathbf 1}_{-C}$ we obtain
${\mathbf 1}_{-C}=\frac{|C|}{q^{n/2}}\widehat{{\mathbf 1}_{C^\perp}}$.
Since the code $(C^\perp)^\perp$ is linear, $-C$ is also a linear code.   Then $C=-C$.
\end{proof}

An affine subspace is a translate of a linear subspace
by some vector. The cardinality of any affine subspace in
$Q_q^n$ equals $q^s$ for some integer $s\geq0$. The intersection
of two affine subspaces is an affine subspace. In
particular, the following proposition holds.

\begin{claim}\label{pvs111} Let $C\subset Q_q^n$ be a linear code and
$\Gamma$ a face of the hypercube $Q_q^n$ with $\bar{0}\in
\Gamma$. Then the cardinality of the intersection $|C\cap (x+\Gamma)|$ equals $0$
    or $q^s$ for each $x\in Q_q^n$, where $s$ does not depend on $x$.
\end{claim}
\begin{proof}
A system of linear equations over  $GF(q)$ has  either no solutions  or  exactly $q^s$
solutions, depending on the right-hand side. 
The exponent $s$ depends only on the rank of the corresponding system of linear equations and is therefore independent of the right-hand side.
\end{proof}

\subsection{Hadamard Code and Plotkin Bound}

\begin{claim}\label{claim:Hadamard}
Let $q$ be prime and let $C\subset Q_q^n$ be  a Hamming code. Then the code $C^\perp$ has the minimum distance $n-\frac{n-1}{q}$ and cardinality $n(q-1)+1$.
\end{claim}
\begin{proof}
Since $C$ is a linear code, we have $|C^\perp|=q^n/|C|$. Since $C$ is a $1$-perfect code, the number $q^n/|C|$ equals the size of a ball of radius $1$, i.e.,  $n(q-1)+1$.  By  Proposition \ref{cor:dualcode} we obtain that  $n(q-1)-\wt(v)q=-1$ for nonzero $v\in C^\perp$  because $\lambda=-1$ is the only  nontrivial eigenvalue of a $1$-perfect code. Since  $C^\perp$ is a linear code, it contains all differences of pairs of codewords and its minimum distance equals the minimum weight of nonzero elements.
\end{proof}

The  dual code to the Hamming code is called the Hadamard code. Consider in more detail the binary case $q=2$. In the binary case the Hadamard code of length $n$ contains $n+1$ vectors, and all nonzero vectors of the Hadamard code of length $n$ have weight $\frac{n+1}{2}$. We now show that this is  optimal, i.e., there exists no code of larger cardinality for the same length with the same (or larger) minimum distance.

\begin{claim}[Plotkin bound]\label{c:plotkin}
Let the minimum distance of a code $C\subset Q^n_2$ be at least $\beta n$, where $\beta>\frac12$. Then $|C|\leq \frac{1}{2\beta-1}+1$.
\end{claim}
\begin{proof}
Consider the vectors $(-1)^v=((-1)^{v_1},\dots,(-1)^{v_n})$, where $v\in C$. The  squared norm $((-1)^v,(-1)^v)$ equals $n$ for any vector $v$. If $d(u,v)=r$, then $((-1)^v,(-1)^u)=n-2r$. We have the inequalities
$$0\leq(\sum\limits_{v\in C}(-1)^v,\sum\limits_{v\in C}(-1)^v)\leq
mn+2{m \choose 2}(n-2\beta n),$$
where $m=|C|$. Then $0\leq 1+(m-1)(1-2\beta)$.
\end{proof}

The inequality proved in Proposition \ref{c:plotkin} is known as the Plotkin bound. Substituting the minimum distance and cardinality of the Hadamard code into it, it is easy to see that the Hadamard code attains the Plotkin bound.

For the existence of a $1$-perfect code in $Q^n_2$ it is necessary and sufficient that $n=2^t-1$ (see Proposition \ref{c:codeHam}). Therefore, the minimum distance $\frac{n+1}{2}$ of the Hadamard code is even. Hence balls of radius $\frac{n-1}{4}$ with centers at the vertices of the Hadamard code do not intersect, and vertices at distance $\frac{n+1}{4}$ from a codeword are not contained in such balls. Consequently, the binary Hadamard code can correct $\frac{n-1}{4}$ errors and detect $\frac{n+1}{4}$ errors. 

A generalization of the Plotkin bound to $q$-ary codes will be proven in Section \ref{sec:LPB}.

\subsection{Extended 1-Perfect Code }

To each vector $v$ of a binary 1-perfect  code $C$, we add one more coordinate equal  to the parity of its weight, $|v|=v_1+\cdots+v_n\mod 2$. The  code $\overline{C}=\{(|v|,v) \ |\ v\in C\}$ is called the {\sl extended 1-perfect code}. Since all vectors of the code $\overline{C}$ have even weight, the minimum distance of $\overline{C}$ is even and at least the  minimum distance of  $C$. Therefore the code distance of $\overline{C}$ equals $4$.

The extended 1-perfect code can be regarded as a perfect coloring of the Boolean hypercube into three colors. Let $f={\mathbf 1}_C$ be a perfect coloring with quotient matrix $S=\begin{pmatrix}
0 & n\\
1 & n-1
\end{pmatrix}$. The coloring $f'(x_0,\overline{x})=
f(\overline{x})$, obtained by copying the coloring $f$ to the two parallel hyperfaces, has quotient matrix $S+I=\begin{pmatrix}
1 & n\\
1 & n
\end{pmatrix}$. The non-trivial eigenvalue $0$ of this matrix does not coincide with the non-trivial eigenvalue $-(n+1)$ of the  parity counter in $Q^{n+1}_2$.  Proposition \ref{c:color_distr_proportion} implies  that half of the vertices of each color have even weight, and the other half have odd weight. Since vertices of the same parity in the Boolean hypercube are not adjacent, we can split each of the colors into even and odd, and obtain a perfect coloring with quotient matrix $\begin{pmatrix}
0& 0 & 1 & n\\
0 & 0& 1 & n\\
1 & n & 0 & 0 \\
1 & n & 0 & 0 \\
\end{pmatrix}$. The first color class in the resulting perfect coloring is the extended 1-perfect code by construction. It follows from the matrix that  the first two  colors (as well as the two other ones) are equivalent. By Corollary \ref{cor:scleiv} after merging the  last two colors we obtain a perfect coloring with quotient matrix $\overline{S}=\begin{pmatrix}
0& 0 & n+1 \\
0 & 0& n+1 \\
1 & n & 0  \\
 \end{pmatrix}$.

The linearity of $C$ implies the linearity of  $\overline{C}$. It is easy to see that the parity-check matrix $\overline{H}$ of the code $\overline{C}$ is obtained from the parity-check matrix of the binary Hamming code by first adding a column of all zeros and then a row of all ones.

\begin{example}
The parity-check matrix of the extended Hamming code of length $8$ is
 $\overline{H}=\begin{pmatrix}
1 & 1& 1& 1& 1& 1& 1& 1\\
0 &0 & 0 & 0 & 1& 1 & 1 & 1\\
0& 0 & 1 & 1 & 0 & 0 & 1 & 1 \\
0 &1 & 0 & 1 & 0 & 1 & 0 & 1
\end{pmatrix}$.
\end{example}

Since the eigenvalues of the quotient matrix $\overline{S}$ are $0,\pm
(n+1)$,  Proposition \ref{cor:dualcode}  implies that every codeword of the dual code $\overline{C}^\perp$ has weight $0$, $m$, or $m/2$, where $m=n+1$ is the code length. Then for any vectors $u,v\in \overline{C}^\perp$ either $(-1)^u=-(-1)^v$, or $((-1)^v,(-1)^u)=m-2d(u,v)=0$. From the equalities $|\overline{C}||\overline{C}^\perp|=2^m$ and $|\overline{C}|=|C|=2^{m-1}/m$ we have $|\overline{C}^\perp|=2m$. Choosing one vector from each pair of opposite vectors of the code $\overline{C}^\perp$, we obtain a set of $m$ orthogonal vectors in $m$-dimensional space.   Thus, taking these vectors as  rows,  we obtain a Hadamard matrix  (see Section \ref{10.3}). It can be shown by induction that the Hadamard matrix obtained from the linear Hadamard code is the Hadamard--Sylvester matrix.

\subsection{Sarkar's Identity}

\begin{theorem}[Sarkar at al. \cite{CarSar}, \cite{VNTs}, and \cite{MacW})]\label{th:sarkar}
Let $f:Q^n_q\rightarrow \mathbb{C}$. Then\\
 (a)  For every face $\Gamma\subset Q^n_q$, $\bar{0}\in \Gamma$, the equality
$\widehat{f}*{\mathbf 1}_{\Gamma^\perp}=q^{n-\mathrm{dim}\,\Gamma}\widehat{f\cdot{\mathbf 1}_\Gamma}$ holds.\\
(b) Let $q$ be prime.  For every linear code $C \subset Q^n_q$,  the equality
$\widehat{f}*{\mathbf 1}_{C^\perp}=q^{n-\mathrm{dim}\,C}\widehat{f\cdot {\mathbf 1}_{C}}$ holds.
\end{theorem}
\begin{proof}
(a)   From Proposition \ref{c:convolution2} we have the equality
$\widehat{f}*\widehat{{\mathbf 1}_{\Gamma}}=q^{n/2}\widehat{f\cdot {\mathbf 1}_{\Gamma}}$.
Next we use the equality
$\widehat{{\mathbf 1}_{\Gamma}}=\frac{q^{\mathrm{dim}\,\Gamma}}{q^{n/2}}{\mathbf 1}_{\Gamma^\perp}$
(Corollary \ref{cor:face}).

(b)  Replacing Corollary \ref{cor:face} with Corollary \ref{cor:face2}, we obtain part (b) in the same way as part (a).
\end{proof}

From the Sarkar's identity or Theorem \ref{fsP} one can deduce

\begin{corollary}[\cite{Sar}]\label{cor:Sar}
Let $f:Q^n_q\rightarrow \mathbb{C}$ and let the face $\Gamma$ contain $\bar{0}$. Then $\sum\limits_{z\in\Gamma^\perp}\widehat{f}(z)=
q^{n/2-\mathrm{dim}\,\Gamma}\sum\limits_{x\in\Gamma}f(x)$.

\end{corollary}
\begin{proof}
By the Sarkar's identity we have the equalities
    \[
    \widehat{f}*{\mathbf 1}_{\Gamma^\perp}(0)=q^{n-\mathrm{dim}\,\Gamma}
    \widehat{f\cdot{\mathbf 1}_{\Gamma}}(0),\ \mbox{hence}
    \]
\[
\sum\limits_{z\in\Gamma^\perp}\widehat{f}(z)=q^{n-\mathrm{dim}\,\Gamma}\cdot
\frac{1}{q^{n/2}}\sum\limits_{x\in\Gamma}
f(x)\phi_0(x).
\]

Since $\phi_0(x)= 1$, we obtain the required equality.
\end{proof}

Note that this corollary can also be obtained by using Proposition \ref{c:plansh} (Plancherel's theorem).

\begin{theorem}[\cite{Carlet20}, \cite{Taran02}] \label{cor:Sar0}
    Let $f:Q^n_2\rightarrow \{0,1\}$ be a Boolean function, $\cor(f)=k$.
    Then the coefficients $\widehat{f}(z)$ have the form $2^{k-\frac{n}{2}}\cdot m$, where $m$ is an integer. If the function $f$ is resilient, then the coefficients $\widehat{f}(z)$  have the form $2^{1+k-\frac{n}{2}}\cdot m$, where $m$ is an integer.
\end{theorem}
\begin{proof}
By Theorem \ref{t:cor-im}, $\widehat{f}(z)=0$ for $0<\wt(z)\le k$.
We first prove the statement for coefficients satisfying $\widehat{f}(z)$ with $\wt(z)=k+1$.
    Consider an $(n-k-1)$-dimensional face $\Gamma$. The face $\Gamma^\perp$,
     $\bar{0}\in \Gamma$, has dimension
     $\mathrm{dim}\,\Gamma^\perp=k+1$ and
      contains exactly one $n$-tuple $z$ of weight $k+1$. Hence, by 
       Corollary \ref{cor:Sar} we have the equality
        \[
    \widehat{f}(z)+\widehat{f}(\bar
    0)=2^{\frac{n}{2}-(n-k-1)}\cdot\sum\limits_{x\in\Gamma}f(x).
    \]
    Let $m''=\sum\limits_{x\in\Gamma}f(x)$ be the number of values $1$ of  $f$
     in the face $\Gamma$.
         The number of ones  of $f$ in  every $(n-k)$-dimensional face is the same, since $\cor(f)=k$.    Let it be  $m'$. Then $\widehat{f}(\bar   0)=2^{-\frac{n}{2}+k}\cdot
    m'$.
      Thus, we obtain the equality
    \[
    \widehat{f}(z)=2^{k+1-\frac{n}{2}}\cdot m''-2^{-\frac{n}{2}+k}\cdot
    m'=2^{k-\frac{n}{2}}\cdot m .
    \]
    For a resilient function we have $\widehat{f}(\bar
0)=\frac{1}{2^{n/2}}\sum\limits_{x\in Q^n_2}f(x)=2^{\frac{n}{2}-1}$.
For $k\leq n-2$ it holds that $k+1-\frac{n}{2}\leq \frac{n}{2}-1$.
Boolean functions with $\cor(f)>n-2$ are only the parity counter, its negation and constants. In these cases we can verify the required
equality directly.

For the coefficients $\widehat{f}(z)$ with $\wt(z)=k+s$, $s>1$,
the statement can be proved by induction on $s$.  Indeed, consider a face
$\Gamma$ of dimension $n-k-s$. The face $\Gamma^\perp$
of dimension $k+s$ contains only one vertex of weight $k+s$,  for
vertices of smaller weight the statement has already been established by the induction
hypothesis.
\end{proof}

\begin{remark}[\cite{Taran02}]\label{cor:Sar1}
 If instead of the Boolean function $f$ we consider the function $(-1)^f$,
 then a similar argument shows that $\widehat{(-1)^f}(z)=2^{k-n/2+1}m$, where $m$ is an integer, 
 since the function $(-1)^f$ has even sums in all faces of nonzero dimension. If the function $f$ is
  resilient, then $\widehat{(-1)^f}(z)=2^{k-n/2+1}m$, where $m$ is even, when $k<n-1$, since in
  this case every $(n-k)$-dimensional face contains an even number of ones and  the second term $\widehat{(-1)^f}(\bar 0)$ equals zero.
\end{remark}

\begin{corollary}\label{corNS}
Let $f:Q^n_2\rightarrow \{0,1\}$. If $\widehat{(-1)^f}(z)=0$ for $\wt(z)\geq k$, then all coefficients $\widehat{(-1)^f}(z)$ have the form $m2^{\frac{n}{2}-k+1}$, where $m$ is an integer.
\end{corollary}
\begin{proof}
Let $g(x)=f(x)\oplus x_1\oplus \cdots\oplus x_n$.
$\widehat{(-1)^g}(z)=\widehat{(-1)^f}(z\oplus {\bar1})$ by Proposition \ref{cl:Four1-3} (c). Then
$\widehat{(-1)^g}(z)=0$ for $\wt(z)\leq n-k$. Consequently, by Theorem \ref{cor:Sar0},
$\widehat{(-1)^g}(z)=2^{n-k-n/2+1}m$, where $m$ is an integer.
\end{proof}

\subsection{Relevant Variables of Perfect Colorings}

The following theorem is a direct application of the Nisan--Szegedy theorem \cite{NiSz} to perfect colorings.

\begin{theorem}[Valyuzhenich \cite{Val26}]
 Let $f$ be a perfect coloring of  $Q^n_2$ with quotient matrix $S=\begin{pmatrix}
    n-b & b\\
    c & n-c
    \end{pmatrix}$. Then the number of relevant variables of $f$
    does not exceed $\frac{bc2^{(b+c)/2}}{b+c}$.
\end{theorem}
\begin{proof}
Without loss of generality, assume that $f$ takes values $\pm 1$.  Proposition \ref{c:cor_immun_for_perfect_coloring0} implies that $\widehat{f}(z)\neq 0$ only if $n-(b+c)=n-2\mathrm{wt}(z)$, i.e., $\mathrm{wt}(z)=(b+c)/2$ or $z=\bar 0$. Therefore, by Corollary \ref{corNS} every coefficient $\widehat{f}(z)$ has the form $m2^{\frac{n-(b+c)}{2}+1}$, where $m$ is an integer. Let $f_0$ and $f_1$ be the restrictions of $f$ to the hyperfaces defined by the equations $x_n=0$ and $x_n=1$, respectively. Denote the $(n-1)$-variable function $h=f_0 - f_1$. It is easy to verify directly that
$$\widehat{h}(z_1,\dots,z_{n-1})=\frac{1}{\sqrt{2}}\widehat{(f(x)  - f(x\oplus e_n))}(z_1,\dots,z_{n-1},1).$$
Therefore the coefficients $\widehat{h}$ are divisible by $2^{\frac{n-(b+c)+1}{2}}$.

 Since $h$ takes values  in $\{0, \pm 2\}$, from Parseval's identity we have
$$4|\supp(h)|=\|h\|^2=\|\widehat{h}\|^2\geq 2^{n-(b+c)+1}|\supp(\widehat{h})|.$$
By Proposition \ref{c:uncertainty_heisenberg} we have
$$|\supp(h)|\geq \frac{2^{n-1}}{|\supp(\widehat{h})|}.$$
Then 
$$|\supp(h)|\geq\frac{2^{2n-(b+c)-2}}{|\supp(h)|}.$$
Thus, $|\supp(h)|\geq 2^{n-1-(b+c)/2}$.

Consequently, if $f$ depends essentially on the $n$-th coordinate, then the values of $f$ differ in at least $2^{n-1-(b+c)/2}$ vertices adjacent along this coordinate. Hence if $f$ depends essentially on $t$ coordinates, then  $I[f]\geq t2^{n-1-(b+c)/2}$, where $I(f)$ is  the number of edges whose endpoints have different values. On the other hand, by Proposition \ref{cl:I} and Theorem \ref{cor:X} we have
$$I[f]= 2^{n-1}(b+c)(\rho(f)-\rho^2(f)),$$
where $\rho(f)=\frac{c}{b+c}$ or $\rho(f)=\frac{b}{b+c}$. Then $t\leq \frac{2^{(b+c)/2}bc}{b+c}$.
\end{proof}

\subsection{Testing Sets of Perfect Colorings}\label{secTest}

Consider the generalized Hadamard matrix $H_{q,n}$ whose rows are characters, i.e., $H_{q,n}(z,x)=\phi_z(x)$.  Let  $A_d$  be the square matrix obtained from $H_{q,n}$ by selecting the rows and columns corresponding to vertices of weight at least $d$.

\begin{claim}\label{c:doplinearball}
$(a)$ For any $d\geq 0$, the matrix $A_d$ is nondegenerate.

$(b)$ If $f(x)=0$ for all $x\in Q^n_q$ with $\wt(x)\leq d$ and $\widehat{f}(z)=0$ for all $z\in Q^n_q$ with $\wt(z)> d$, then $f={\bf 0}$.
\end{claim}
\begin{proof}
By  definition of $H_{q,n}$ we have  $\frac{1}{q^{n/2}}H_{q,n}f= \widehat{f}$.
The matrix $A_d$ is degenerate if and only if there exists a nonzero solution $A_du=0$. By definition of $A_d$ this solution is a function $f:Q_q^n\rightarrow\mathbb{C}$ such that $f(x)=0$ if $\wt(x)>d$ and $\widehat{f}(z)=0$ for $\wt(z)\leq d$.
So, the  conditions (a) and (b)  are equivalent. It therefore suffices to prove (b).

By the  condition of the proposition, $f\cdot{\mathbf 1}_{\Gamma}=0$ for any face $\Gamma$ of dimension $d$ containing $\bar{0}$. From  Sarkar's identity (Theorem \ref{th:sarkar} (a)) we have \\
$\widehat{f}*{\mathbf 1}_{\Gamma^\perp}=0$. Thus,  the sum of the values of $\widehat{f}$ over any face of dimension $n-d$ (not necessarily containing $\bar{0}$) equals zero. It is clear that $\widehat{f}$ is also balanced in faces of larger dimension. For each vertex $a$ of weight $d$ there exists a face $\Gamma_a=\{y\in Q^n_q : y_i=a_i,\ \text{if}\ a_i\neq 0\}$ of dimension $n-d$, which contains exactly one vertex of weight $d$ and the remaining vertices have larger weight. From the definition of $\widehat{f}$ we have $\widehat{f}(a)=0$. Thus, we have shown that $\widehat{f}$ equals zero on any vertex of weight $d$. Using faces of dimension $n-d+1$ we can similarly show that $\widehat{f}$ equals zero on vertices of weight $d-1$. Then, by induction we obtain $\widehat{f}={\bf 0}$ and, consequently,   ${f}={\bf 0}$. 
\end{proof}

\begin{claim}\label{c:doplinearball2}
An eigenfunction $f$ in $Q^n_q$ with eigenvalue $\lambda_k=n(q-1)-qk$ is uniquely determined by its values in the ball $B_k=\{x\in Q_q^n : \wt(x)\leq k\}$.
\end{claim}
\begin{proof}
Suppose that $f$ is not uniquely determined by its values in the ball $B_k=\{x\in Q_q^n : \wt(x)\leq k\}$, i.e., there exists an eigenfunction $f'$ with eigenvalue $\lambda_k$ such that $f|_{B_k}=f'|_{B_k}$. Consider the difference $g=f-f'$. By Proposition \ref{claim:ef}, $\widehat{g}(z)\neq 0$ only if $\wt(z)=k$. Consequently, $g|_{B_k}=0$ and $\widehat{g}|_{Q^n_q\setminus B_k}=0$. By Proposition \ref{c:doplinearball} (b), $g=0$.
\end{proof}

By Proposition \ref{c:dopeigenfun}   (see also Problem \ref{ex:Pnn})  and Proposition \ref{c:doplinearball2}, it follows

\begin{corollary}\label{c:doplinearball3}
If $k<n/2$, then an eigenfunction $f$ in $Q^n_q$ with eigenvalue $\lambda_k=n(q-1)-qk$ is uniquely determined by its values on the two spheres $A_{k}$ and $A_{k-1}$.
\end{corollary}

In the paper \cite{Vas15}, sufficient conditions are given for recovering an eigenfunction from a single sphere. Corollary \ref{c:doplinearball3} and Remark \ref{z:colorlin0} allow us to generalize Theorem \ref{th:2level} to the $q$-ary hypercube.

\begin{corollary}
Let $f$ be a perfect coloring of  $Q^n_q$ with quotient matrix $S$. Let $v$ be an eigenvector of $S$ corresponding to the eigenvalue $\lambda_k=n(q-1)-qk$, with $k<n/2$. If all coordinates of $v$ are distinct, then $f$ is uniquely determined by its values on the two spheres $A_{k}$ and $A_{k-1}$.
\end{corollary}

\subsection{Problems}

\begin{exercise}\label{exer2}
Let $f: Q_q^n \rightarrow\mathbb{C}$ and let $C\subset Q_q^n$ be a linear code. Prove that  $\frac{1}{q^{n/2}}\sum\limits_{x\in C}f(x)=\frac{1}{|C^\perp|}\sum\limits_{z\in C^\perp}\widehat{f}(z)$.
\end{exercise}

\begin{exercise}
Let $A\subset Q_q^n$, $q$ prime, and $C=\supp(\widehat{{\bf 1}_A})$. Prove that for any $v\in Q_q^n$ the equality $A+v=A$ holds if and only if $v\in C^\perp$.
\end{exercise}

\begin{exercise}
Let $C\subset Q_q^n$, $q$ prime, and $f=\widehat{{\bf 1}_C}$. Prove that $|f|=c{\bf 1}_A$ for some $c>0$ and $A\subset Q_q^n$ if and only if $C$ and $A$ are affine subspaces.
\end{exercise}

\begin{exercise}
Prove that  if every pair of vectors in an $n$-dimensional Euclidean space forms either an obtuse or a right angle, then their number does not exceed $2n$. Prove that the cardinality of a binary code of length $2n$ with code distance at least $n$ does not exceed $4n$.
\end{exercise}

\begin{exercise}
Prove that two distinct correlation-immune functions of order $k$ in the Boolean hypercube differ from each other on a set of cardinality at least $2^{k+1}$.
\end{exercise}

\begin{exercise}
Let $\Gamma\subset  Q_2^n$ be an arbitrary face containing $\bar{0}$. Suppose two Boolean functions $f$ and $g$ satisfy  $\widehat{f}|_\Gamma=\widehat{g}|_\Gamma$. Prove that the equality $\sum\limits_{x\in z\oplus\Gamma^\perp}(-1)^{f(x)}=\sum\limits_{x\in z\oplus\Gamma^\perp}(-1)^{g(x)}$ holds for any $z\in  Q_2^n$.
\end{exercise}

\begin{exercise}
Let $f:Q^n_2\rightarrow \{0,1\}$ be a perfect coloring with the largest eigenvalue other than  $n$ is  $\lambda=n-2k$. Prove that the function $f$ is completely determined by its values in the ball of radius $n-k$.
\end{exercise}

\section{MacWilliams  Identity and Weight Distributions}

\subsection{MacWilliams Transform}

Let $M$ be the adjacency matrix of the hypercube $Q^n_q$. Recall that $M_k=P_k(M)$, where $M_k$ is the distance-$k$ adjacency matrix and $P_k$ is the $q$-ary Krawtchouk polynomial of degree $k$. If $M\psi=\lambda\psi$, then from Proposition \ref{c:WDB_bound0} we have the equality
$\sum\limits_{x,d(x,a)=k}\psi(x)=(M_k\psi)(a)=P_k(\lambda)\psi(a)$.
Since the eigenvalue of a character $\phi_z$ depends only on the weight $\wt(z)$, we will further index the eigenvalues by the weights of their corresponding characters
$M\phi_z(x)=\lambda_{\wt(z)}\phi_z(x)$.

Consider an arbitrary function $f:Q^n_q\rightarrow \mathbb{C}$.
Introduce the notation $A_k[f]=\sum\limits_{x,d(x,\bar 0)=k}f(x)$ and
$\widehat{A}_k[f]=\sum\limits_{z,d(z,\bar 0)=k}\widehat{f}(z)$ for the {\sl weight distribution} of the function and its Fourier transform.

\begin{claim}[\cite{MacW}]\label{c:eqMW}
Let $f:Q^n_q\rightarrow \mathbb{C}$. Then $\widehat{A}_k[f]=
\frac{1}{q^{n/2}}\sum\limits_m P_k(\lambda_m)A_m[f]$.
\end{claim}
\begin{proof}
Taking into account that $\overline{\phi_z}=\phi_{-z}$, from the definitions we obtain the chain of equalities

$$\widehat{A}_k[f]=\frac{1}{q^{n/2}}\sum\limits_{z,d(z,\bar 0)=k}\sum\limits_x
f(x)\overline{\phi_z(x)} = \frac{1}{q^{n/2}}\sum\limits_x
f(x)\sum\limits_{z,d(z,\bar 0)=k}\phi_z(x)=$$
$$=\frac{1}{q^{n/2}}\sum\limits_x
f(x)\sum\limits_{z,d(z,\bar 0)=k}\phi_x(z)=\frac{1}{q^{n/2}}\sum\limits_x
f(x)P_k(\lambda_{\wt(x)})=\frac{1}{q^{n/2}}\sum_{m=0}^n
P_k(\lambda_m)A_m[f].$$
\end{proof}

Note that Proposition \ref{c:eqMW} can be regarded as a special case of Proposition \ref{c:plansh} (Plancherel's theorem), where the indicator of the sphere of radius $k$ is taken as the second function.

%\begin{definition}
The linear transform taking the coefficients $A_m[f]$ into the coefficients $\widehat{A}_k[f]$ is called the {\sl MacWilliams transform}.
%\end{definition}
From the formula $\widehat{\widehat{f(x)}}=f(-x)$ it follows that the MacWilliams transform is nondegenerate and idempotent for $q=2$.

Assume that  $q$ is a prime number and that $C$ is a linear code. Denote $A_k(C)= A_k[{\mathbf 1}_{C}]$. Then from Proposition \ref{c:eqMW} and Corollary \ref{cor:face2} we have
$$|C|A_k({C^\perp})=\sum\limits_m P_k(\lambda_m)A_m(C).$$
The latter equality is called the {\sl MacWilliams identity}. 
Note that the MacWilliams identity can be proved using Problems \ref{exer1} and \ref{exer2}.

Let $C$ be a code. Note that the number $f_{C}(a)=({\mathbf 1}_{C}\ast{\mathbf 1}_{-C})(a)$ is the number of ordered pairs $(x, x-a)$ of vertices of the code differing by the vector $a$. Define the weight spectrum $(A_0[C],\dots,A_n[C])$ of the code $C$ as $\frac{1}{|C|}(A_0[f_C],\dots,A_n[f_C])$. Note that $A_0[f_C]=\sum_x {\mathbf 1}^2_{C}(x)=|C|$ and
$$\sum\limits_mA_m[C]=\frac{1}{|C|}\sum_xf_{C}(x)=
\frac{1}{|C|}\sum_x\sum_a{\mathbf 1}_{C}(x){\mathbf 1}_{C}(x-a)=|C|.$$

 By the definition of a linear code, ${\mathbf 1}_{C}\ast{\mathbf 1}_{-C}=|C|{\mathbf 1}_{C}$.
 Therefore for a linear code $(A_0[C],\dots,A_n[C])= (A_0(C),\dots,A_n(C))$.

Then
$\widehat{A}_k[{\mathbf 1}_{C}\ast{\mathbf 1}_{-C}]=\widehat{A}_k[f_{C}]=
\frac{1}{q^{n/2}}\sum\limits_m P_k(\lambda_m)A_m[f_C]$. Moreover,
$\widehat{{\mathbf 1}_{C}\ast{\mathbf 1}_{-C}}=
q^{n/2}|\widehat{{\mathbf 1}_C}|^2$ by Propositions \ref{cl:Four1-3} and \ref{c:convolution}.
By Proposition \ref{c:eqMW} for any $k$, $0\leq k\leq n$, we have the equality
\begin{equation}\label{e:linprog}
q^n\cdot\sum\limits_{z,\wt(z)=k}|\widehat{{\mathbf 1}_C}(z)|^2=|C|\sum\limits_m
P_k(\lambda_m)A_m[C].
\end{equation}

\subsection{Linear Programming Bound}\label{sec:LPB}

The equalities (\ref{e:linprog}) provide the basis for applying linear programming methods to compute an upper bound on the size of a code with minimum distance $d$. Indeed, if the code $C$ does not contain two vertices at distance less than $d$, then $A_m[C]=0$ for $m=1,\dots,d-1$, since ${\mathbf 1}_{C}\ast{\mathbf 1}_{-C}(y)=0$ for $0<\wt(y)<d$. From this we obtain a {\sl linear programming problem} on the set of quantities $A_m[C]$ to find an upper bound on the size of a code with code distance $d$:
\begin{equation}\label{e:linprog1}
P_k(\lambda_0)+\sum\limits_{m\geq d} P_k(\lambda_m)A_m[C]\geq 0\
\text{for any} \ k
\end{equation}
$$|C|=1+\sum\limits_{m\geq d}A_m[C]\rightarrow \max.$$

It is known that
every solution of the dual minimization problem is an upper bound for the value of the linear programming problem, i.e., for the code size in this case. This is the basis for the linear programming bound on the size of a code with known minimum distance, which is proved in the following theorem.

\begin{theorem}[Delsarte \cite{Delsarte}]\label{thDelsart1}
Let $C$ be a code in  $Q^n_q$ with code distance $d$, and let 
 $(y_0,\dots,y_n)\in \mathbb{R}^{n+1}$ satisfy the conditions $y_0>0$, $y_i\geq 0$ 
 for $i=1,\dots,n$, and $\sum\limits_{k=0}^ny_kP_k(\lambda_m)\leq 0$ for $m\geq d$. 
 Then $|C|\leq \frac{1}{y_0}\sum\limits_{k=0}^n y_kP_k(\lambda_0)$.
\end{theorem}
\begin{proof}
As before, we formally assume that the polynomial $P_0(x)$ is identically equal to one. Then $|C|=P_0(\lambda_0)+\sum\limits_{m\geq d} P_0(\lambda_m)A_m[C]$. By the assumption of the theorem and  (\ref{e:linprog1}) we obtain the inequalities
$$ y_0|C|\leq
\sum\limits_{k=0}^ny_k\left(P_k(\lambda_0)+\sum\limits_{m\geq d}
P_k(\lambda_m)A_m[C]\right) $$
$$\leq\sum\limits_{k=0}^ny_kP_k(\lambda_0) +\sum\limits_{m\geq
d}A_m[C]\left(\sum\limits_{k=0}^ny_kP_k(\lambda_m)\right)\leq
\sum\limits_{k=0}^n y_kP_k(\lambda_0).$$
\end{proof}
Note that $P_i(\lambda_0)$ is the size of a ball of radius $i$ in $Q^n_q$. 

Let us discuss how to obtain upper bounds on code size using linear programming. Suppose that we wish to bound the size of a code whose pairwise distances between vertices take values only from a list $W=\{w_1,\dots,w_m\}$. Consider a polynomial $R$ of degree $m$ such that $R(\lambda_{w_i})\leq 0$, where $\lambda_{w_i}=n(q-1)-qw_i$ for $i=1,\dots,m$. Any polynomial of degree at most $m$ can be represented as a linear combination of the polynomials $P_i$ for $i=0,1,\dots,m$. Let $R=\sum\limits_{i=0}^my_iP_i$. If $y_0>0$ and $y_i\geq 0$ for $i=1,\dots,m$, then Theorem \ref{thDelsart1} implies  the inequality $|C|\leq \frac{1}{y_0}\sum\limits_{i=0}^m y_iP_i(\lambda_0)=R(\lambda_0)/y_0$. 

As an example, let us prove the Plotkin bound on the size of a code with large minimum distance. Let the minimum  distance of code be $d_C> (1-\frac1q)n$, i.e., if $w_i\in W$, then  $w_i\geq d_C$ and $\lambda_{w_i}=n(q-1)-qw_i\leq  n(q-1)-qd_C<0$. Consider the polynomial $R(x)=1-\frac{x}{n(q-1)-qd_C}$. 
$$R(\lambda_{w_i})= 1-\frac{n(q-1)-qw_i}{n(q-1)-qd_C}\leq 0$$ for all $w_i\in W$. Then
\begin{equation}\label{eqPlot}
|C|\leq R(\lambda_0)=
1-\frac{n(q-1)}{n(q-1)-qd_C}=\frac{qd_C}{qd_C-n(q-1)}.
\end{equation}
Inequality (\ref{eqPlot}) is  the Plotkin bound (Proposition \ref{c:plotkin}) for $q$-ary codes. In the case $q=2$ and $d_C=\frac{n+1}{2}$, $n$ odd, we obtain the Hadamard bound $|C|\leq n+1$.

Note that Theorem \ref{thDelsart1} remains true not only for codes but also for any functions $F:Q^n_q\rightarrow \mathbb{R}$, provided $A'_k[F]=\sum_{\wt(x)=k}F(x)\geq 0$ and $A'_0[F]=1$. In this case the role of the code size is played by the sum of the values of $F$, i.e., $\sum\limits_{k=0}^nA'_k[F]$. Consider a correlation-immune function $f={\bf 1}_C$, $\cor(f)+1>(1-\frac1q)n$, and define $F=\frac{q^{n/2}\widehat{f*f}}{|C|^2}$. Then $F(x)=\frac{q^n|\widehat{f}(x)|^2}{|C|^2}$, i.e., $A'_k[F]\geq 0$, $A'_0[F]=1$, $A'_k[F]= 0$ for $k\leq \cor(f)$. For this function the Plotkin bound (\ref{eqPlot}) also holds:
$$\sum\limits_{k=0}^nA'_k[F]\leq
\frac{q(\cor(f)+1)}{q(\cor(f)+1)-n(q-1)}.$$
Moreover, from Parseval's identity we have $\sum\limits_{k=0}^nA'_k[F]=q^n/|C|$. Then we obtain
$$\frac{|C|}{q^n}\geq 1-\frac{n(q-1)}{q(\cor(f)+1)}.$$
Thus, we recover the Bierbrauer–Friedman bound (Theorem \ref{thBF}), which is called the dual Plotkin bound in \cite{Bierbook}.

\subsection{Few-Weight Codes }

A {\sl few-weight code} is a code in which the set of nonzero codeword weights (or, for non-linear codes, the set of pairwise distances between codewords) contains only a small number of distinct values.
For example, the Hadamard code is an {\sl equidistant or constant-weight code}. In this section we estimate the cardinality of few-weight codes.

\begin{theorem}[Delsarte \cite{Delsarte73}]\label{th:Del1}
Let $f:Q^n_q\rightarrow\mathbb{C}$, and suppose there exist $k$ nonzero weights $W=\{w_1,\dots,w_k\}$ such that $\widehat{f}(x)=0$ for $\wt(x)\not \in W$ and $\widehat{f}(\bar 0)\neq 0$. Then $|\supp(f)|\geq q^n/b(q,k)$, where $b(q,k)$ is the cardinality of a ball of radius $k$.
\end{theorem}
\begin{proof}
Let $\lambda_i=n(q-1)-qw_i$ be the eigenvalue corresponding to weight $w_i$. By assumption, we have $f=\rho{\bf 1}+\sum\limits_{i=1}^k\varphi_i$, where $\varphi_i$ is an eigenfunction with eigenvalue $\lambda_i$ and $\rho=\frac{\widehat{f}(\bar 0)}{q^{n/2}}\neq 0$. Let $A=(M-\lambda_1I)\cdots(M-\lambda_kI)$ be a polynomial  in $M$ of degree $k$, where $M$ is the adjacency matrix of the hypercube $Q^n_q$. We have $Af= \rho A{\bf 1}$, since $A\varphi_i=\bar 0$. From the equality $M{\bf 1}=n(q-1)$ it follows that $A{\bf 1}= (n(q-1)-\lambda_1)\cdots(n(q-1)-\lambda_k){\bf 1}=\alpha{\bf 1}$, where $\alpha\neq 0$.

On the other hand,  since the $q$-ary hypercube is a distance-regular graph, every polynomial of $M$ is a linear combination of distance-$i$ adjacency matrices. Thus, $A=\sum\limits_{i=0}^k \beta_iM_i$, where $M_i$ is the distance-$i$ adjacency matrix. Hence $Af=g*f$, where $\supp(g)$ is contained in the ball of radius $k$. For any $x\in Q^n_q$ we have $0\neq Af(x)=g*f(x)$. Therefore every ball of radius $k$ contains at least one vertex
 $y\in Q^n_q$ such that $f(y)\neq 0$.
\end{proof}

We can consider Theorem \ref{th:Del1} as a generalization of the simple fact that every cell of an equitable 
$(k+1)$-partition intersects  every ball of radius $k$.

\begin{corollary}\label{corDel}
Let a linear code $C\subset Q^n_q$ contain words of only $k$ distinct nonzero weights. Then $|C|\leq b(q,k)$.
\end{corollary}
\begin{proof}
Consider the dual code $C^\perp$. By Delsarte's theorem we obtain $|C^\perp|\geq q^n/b(q,k)$. From the equality $|C||C^\perp|=q^n$ we obtain the desired inequality.
\end{proof}

However, the statement of Corollary \ref{corDel} remains true for nonlinear codes as well. We prove this in the case $q=2$.

\begin{theorem}[Delsarte \cite{Delsarte73}]
Let a code $C\subset Q^n_2$ consist of codewords whose pairwise distances take only $k$ distinct nonzero values. Then $|C|\leq b(2,k)$.
\end{theorem}
\begin{proof}
Let $D$ be the set of pairwise distances of the code $C$, $|D|=k$. Consider the set of functions $f:\{-1,1\}^n\rightarrow\mathbb{R}$ as a $2^n$-dimensional linear space $V_n$. It is easy to see that the monomials of the form $x_{i_1}\cdots x_{i_m}$, $m=0,\dots,n$, form an orthogonal basis of $V_n$, which corresponds to the basis of characters of the group $\mathbb{Z}_2^n$.

To each codeword $v\in C$ we associate the function $$f_v(x)=\prod\limits_{t\in D}(((-1)^v,x)-(n-2t)).$$ Since $((-1)^v,(-1)^u)=n-2d(u,v)$, we have $f_v((-1)^u)=0$ for any $u\in Q^n_2$ such that $d(u,v)\in D$, and $f_v((-1)^v)\neq 0$. Let us show that the set of functions $\{f_v : v\in C\}$ is linearly independent. Let $L=\sum\limits_{v\in C}\alpha_vf_v=0$. Then $L((-1)^v)=\alpha_vf_v((-1)^v)= 0$ and $\alpha_v=0$.

Each function $f_v $ is a polynomial  of degree at most $k$, in which the basis is formed by monomials of degree at most $k$.  The number of monomials of degree at most $k$ equals  $b(2,k)$.
\end{proof}

\subsection{Weight Distributions of Perfect Colorings}

The weight distribution of a $2$-coloring of the hypercube completely determines whether the coloring is perfect.

\begin{claim}
Let $C$ be a code in the hypercube $Q^n_q$, and let its weight distribution $(A_0[C],\dots,A_n[C])$ be the same as that of a cell $B$ of some equitable $2$-partition with quotient matrix $S$. Then ${\mathbf 1}_{C}$ is a perfect $2$-coloring with the same quotient matrix.
\end{claim}
\begin{proof}
The function ${\mathbf 1}_{C}$ is a perfect $2$-coloring if and only if $\supp( \widehat{{\mathbf 1}_{C}})$ contains only the zero vector and vectors of some fixed weight (see Propositions \ref{c:colorlin} and \ref{c:cor_immun_for_perfect_coloring0}). By the formula $q^n\cdot\sum\limits_{z,\wt(z)=k}|\widehat{{\mathbf 1}_{C}}(z)|^2=|C|\sum\limits_m P_k(\lambda_m)A_m[C]$ the Fourier weight distributions of $C$ and $B$ coincide. Thus, ${\mathbf 1}_{C}$ is a perfect $2$-coloring.  By the weight distribution  we can also determine its quotient matrix. Indeed, if $S=\begin{pmatrix}
    (q-1)n-b & b\\
    c & (q-1)n-c
    \end{pmatrix}$ then the nontrivial eigenvalue  equals $(q-1)n-b -c$.  It can be found from the Fourier weight distribution of $C$. Moreover, $\frac{b}{c}= \frac{|C|}{q^n-|C|}$ if $C$ is the first cell of the partition and 
    $\frac{b}{c}= \frac{q^n-|C|}{|C|}$ if $C$ is the second one.
\end{proof}

By Theorem \ref{th:Krot} the following statement follows directly.

\begin{claim}\label{cl:color_comp}
Let $f$ be a $k$-coloring of an amply regular graph. Suppose that for every $i$, $i=1,\dots,k$, the averaged color composition of the neighborhood of the $i$-th color $C_i$ of the coloring $f$ is the same as that of the $i$-th color $D_i$ of a perfect coloring $g$. If $A_2[C_i]=A_2[D_i]$ for every $i$, $i=1,\dots,k$, then $f$ is a perfect coloring with the same quotient matrix as $g$.
\end{claim}

For the case of $2$ colors Proposition \ref{cl:color_comp} implies that the first three weight
 coefficients  $A_0[C], A_1[C] ,A_2[C]$ of $C$ and   the first   three weight
 coefficients $(A_0[Q^n_q\setminus C],  A_1[Q^n_q\setminus C], A_2[Q^n_q\setminus C])$  of its complement are sufficient to determine whether the coloring is perfect.

By Proposition \ref{cl:Four1-3}, for any function taking only real values, the equalities $\widehat{f(-x)}=\overline{\widehat{f(x)}}$ and $\widehat{\widehat{f(-x)}}=f(x)$ hold. Therefore, the equality $\widehat{A}_k[f]= \frac{1}{q^{n/2}}\sum\limits_m P_k(\lambda_m)A_m[f]$ can be inverted:
$${A}_k[f]= \frac{1}{q^{n/2}}\sum\limits_m P_k(\lambda_m)\widehat{A}_m[f(-x)]=\frac{1}{q^{n/2}}\sum\limits_m P_k(\lambda_m)\overline{\widehat{A}_m[f]}.$$
In particular it means that the matrix $W=(P_k(\lambda_m))_{k,m}$ is nondegenerate.

Moreover, MacWilliams identity can be rewritten for the case of Cartesian products.
Suppose $K=Q^{n_1}_{q_1}\times Q^{n_2}_{q_2}$.  Proposition \ref{c:charact_cartesian_product} implies that the basis of characters on $K$ consists of products of characters of $Q^{n_1}_{q_1}$ and characters of $Q^{n_2}_{q_2}$. For a function $f:K\rightarrow \mathbb{C}$ we define $A_{k_1k_2}[f]=\sum\limits_{\wt(x_1)=k_1}\sum\limits_{\wt(x_2)=k_2}f(x_1,x_2)$ and $\widehat{A}_{k_1k_2}[f]=\sum\limits_{wt(z_1)=k_1}\sum\limits_{\wt(z_2)=k_2}\widehat{f}(z_1,z_2)$. Then, by  arguments 
 similar to the proof of Proposition \ref{c:eqMW}, we obtain the equality
\begin{equation}\label{e:eqMWproduct}
 \widehat{A}_{k_1k_2}[f]=
\frac{1}{q_1^{n_1/2}q_2^{n_2/2}}\sum\limits_{m_1,m_2}
P_{k_1}[n_1,q_1](\lambda_{m_1})P_{k_2}[n_2,q_2](\lambda_{m_2})A_{m_1m_2}[f].
\end{equation}
In equality $(\ref{e:eqMWproduct})$, the coefficients $\widehat{A}_{k_1k_2}[f]$ and ${A}_{k_1k_2}[f]$ can also be swapped (with the addition of complex conjugation if $q\neq 2$).

Let $f:Q^n_q\rightarrow \mathbb{C}$ and let $\Gamma\subset Q^n_q$ be some face. Define the local distribution of the function $f$ in the face $\Gamma$ as the weight distribution of the restriction $f|_{_\Gamma}$ of $f$ to the face $\Gamma$.

\begin{theorem}[Vasil'eva \cite{Vas09}]\label{th:Vas}
Let $g:Q^n_q\rightarrow \mathbb{C}$ be an eigenfunction of the hypercube $Q^n_q$. Then the local distribution of $g$ in the face $\Gamma^\perp$ linearly depends on the local distribution of $g$ in the face $\Gamma$.
\end{theorem}
\begin{proof}
Consider the hypercube $Q^n_q$ as the Cartesian product of faces: $Q^n_q=\Gamma\times \Gamma^\perp$. From equality (\ref{e:eqMWproduct}) we have
$$ A_{k_1k_2}[g]=
\frac{1}{q^{n/2}}\sum\limits_{m_1,m_2}
P_{k_1}[n_1,q](\lambda_{m_1})P_{k_2}[n_2,q](\lambda_{m_2})
\overline{\widehat{A}_{m_1m_2}[g]},$$ where $n_1=\dim \Gamma$ and $n_2=n-n_1=\dim\Gamma^\perp$.

By the definition of $\widehat{A}_{m_1m_2}[g]$, we have  $\widehat{A}_{m_1m_2}[g]=0$ for $m_1+m_2\neq m$, where $n(q-1)-mq=\lambda$ is the eigenvalue of $g$. Note that $P_{0}[n_2,q]=1$. Then for $k_1=0,\dots,\dim \Gamma$ we have the equalities
\begin{equation}\label{eq:thV}
A_{k_10}[g]=\frac{1}{q^{n/2}}\sum\limits_{m_1+m_2=m}
P_{k_1}[n_1,q](\lambda_{m_1}) \overline{\widehat{A}_{m_1m_2}[g]}.
\end{equation}

Since the MacWilliams transform is invertible,  the last equality expresses the coefficients $\overline{\widehat{A}_{m_1m_2}[g]}$ in the term of the numbers $A_{k_10}[g]$. Similarly, for $k_2=0,\dots,\dim \Gamma^\perp$ we have
$$
A_{0k_2}[g]=\frac{1}{q^{n/2}}\sum\limits_{m_1+m_2=m}
P_{k_2}[n_2,q](\lambda_{m_2}) \overline{\widehat{A}_{m_1m_2}[g]}.
$$
The set of quantities $A_{0k_2}[g]$ is the local distribution of $g$ in the face $\Gamma^\perp$. We have shown that it can be expressed linearly in terms of the set $A_{k_10}[g]$, which is the local distribution of $g$ in the face $\Gamma$.
\end{proof}

\begin{corollary}
Let $f:Q^n_q\rightarrow \{0,1\}$ be a perfect $2$-coloring with quotient matrix $S$. Then the weight distribution of the coloring $f$ in the face $\Gamma$ determines the weight distribution of $f$ in the face $\Gamma^\perp$.
\end{corollary}
\begin{proof}
Let $f_1$ and $f_2$ be two perfect $2$-colorings with quotient matrix $S$. Consider the function $g=f_1-f_2$. Since $f_1$ and $f_2$ are perfect $2$-colorings with the same parameters,  Proposition \ref{c:colorlin} implies that $g$ is an eigenfunction in $Q^n_q$. By Theorem \ref{th:Vas} it follows that if $g$ has zero local distribution in the face $\Gamma$, then it also has zero local distribution in the face $\Gamma^\perp$. Obviously, $A_{m_1m_2}[g]=A_{m_1m_2}[f_1]-A_{m_1m_2}[f_2]$. Therefore, if the local distributions of $f_1$ and $f_2$ coincide in the face $\Gamma$, then they also coincide in the face $\Gamma^\perp$.
\end{proof}

Vasil'eva's theorem states that the weight spectra of an eigenfunction $g$ in orthogonal faces are related by linear relations. Suppose that all coefficients in these relations are nonzero. In this case    the values of the eigenfunction $g:Q^n_2\rightarrow \mathbb{C}$ on vertices of weight less than 
$d$, $d\leq n/2$, are uniquely recovered from its values on vertices of weight $d$. Suppose eigenfunctions $g$ and $g'$ coincide on vertices of weight $d$. Consider the function $f=g-g'$. Let $u$ be a vertex of minimal weight $\wt(u)=d'<d$ for which $f(u)\neq 0$. Consider the face $\Gamma$ of dimension $d'$ containing the vertices $\bar 0$ and $u$. The local weight coefficients with respect  to this face satisfy $A_{00}[f]=0, \dots, A_{(d'-1)0}[f]=0$, $A_{d'0}[f]=f(u)\neq 0$. Consider the face $\Gamma^\perp$. On the one hand, $A_{0d}[f]=0$, since $f(v)=0$ for all vertices of weight $\wt(v)=d$. On the other hand, 
$A_{0d}[f]=\alpha A_{d'0}[f]$, where $\alpha\neq 0$ by assumption. It is a contradiction.

\subsection{Problems}

\begin{exercise}
Prove the equality \\ $2^{n-t}{t\choose k}=\sum\limits_{m=0}^{n-t}P_k[n,2](\lambda_m){n-t \choose m}$ for $0<k<t<n$.
\end{exercise}

\begin{exercise}
Using the generating function of Krawtchouk polynomials (\ref{eqKraw}), prove that for any function $f:Q^n_q\rightarrow \mathbb{C}$ the following polynomial equality holds:
$q^{n/2}\sum\limits_{k=0}^n\widehat{A}_k[f]z^k=\sum\limits_{m=0}^n A_m[f](1-z)^m(1+(q-1)z)^{n-m}.$ \\
In particular, for a linear code $C\subset Q^n_q$ the equality holds:\\
$|C|\sum\limits_{k=0}^nA_k({C^\perp})z^k=\sum\limits_{m=0}^n A_m(C)(1-z)^m(1+(q-1)z)^{n-m}.$
\end{exercise}

\begin{exercise}
For any function $f:Q^n_q\rightarrow \mathbb{C}$ and $i=0,\dots,n$ prove the polynomial equality\\
$\sum\limits_{k=0}^n\widehat{A}_k[f]{k\choose i}=q^{\frac{n}{2}-i}\sum\limits_{m=0}^n(-1)^m(q-1)^{i-m}A_m[f]{n-m\choose n-i},$ where binomial coefficients ${a \choose b}$ are taken to be zero when $a<b$.
\end{exercise}

\begin{exercise}
Prove that for any function $f:Q^n_q\rightarrow \mathbb{C}$ the following equality holds:
$\sum\limits_{k=0}^nk\widehat{A}_k[f]=q^{\frac{n}{2}-1}((q-1)nA_0[f]-A_1[f])$.
In particular, for a linear code $C\subset Q^n_q$ of dimension $m$ the equality holds:\\
$\sum\limits_{k=0}^nk{A}_k(C)=q^{{m}-1}((q-1)n-A_1({C^\perp}))$.
\end{exercise}

\begin{exercise}
For a linear code $C\subset Q^n_2$ prove the polynomial equality in two variables:\\
$|C|\sum\limits_{v\in {C^\perp}}x^{\wt{(v)}}y^{n-\wt{(v)}}=\sum\limits_{u\in {C}}(y-x)^{\wt{(u)}}(x+y)^{n-\wt{(u)}}$.
\end{exercise}

\begin{exercise}
Let $C\subset Q^n_2$ be a code of odd cardinality with distance $d$. Prove the inequality $P_k(\lambda_0)+\sum\limits_{m\geq d} P_k(\lambda_m)A_m[C]\geq \frac{{n\choose k}}{|C|}$ for every $k$, $0\leq k\leq n$.
\end{exercise}

\begin{exercise}
Let a code $C\subset Q^n_2$ consist of codewords whose pairwise distances take values only from a set $D$, with the numbers $t$ and $n-t$ both belonging  to $D$  and $0,n\not\in D$. Prove that\\
1) if $n/2\not \in D$, then $|C|\leq 1+{n\choose 2}+{n\choose 4}+\dots+{n\choose |D|}$;\\
2) if $n/2\in D$, then $|C|\leq {n\choose 1}+{n\choose 3}+\dots+{n\choose |D|}$.
\end{exercise}

\begin{exercise}
Let $A\subset Q_q^n$ be an affine  equidistant code and let $w$ be a nonzero weight of codewords. Prove that $|A|\leq q^n/L(w)$, where $L(w)$ is the minimum cardinality of the support of an eigenfunction with eigenvalue $\lambda=n(q-1)-qw$.
\end{exercise}

\begin{exercise}
Let $A\subset Q_2^n$ be an affine  equidistant code and let $w$ be a nonzero weight of codewords. Prove that $|A|\leq 2^{(n-|n-2w|)/2}$ (see Problem \ref{exer76}). Prove that this bound is attained for the affine codes $A=\{(x,x\oplus \bar 1,\bar 0) : x\in Q_2^{(n-i)/2}\}$, $i=n-2k$.
\end{exercise}

\begin{exercise}
Let ${\mathbf 1}_{C}$ be a linear perfect $2$-coloring in the hypercube. Prove that $C^\perp$ is an equidistant code.
\end{exercise}

\begin{exercise}
Let $C\subset Q^n_q$ be an equidistant code with code distance $d>n(q-1)/q$. Prove that $|C|\leq 1/(1-\frac{(q-1)n}{qd})$.
\end{exercise}

\section{Bent Functions}

\subsection{Definition and Basic Properties of Bent Functions}

We  regard the lines
$\{\alpha\phi_z \ |\ \alpha\in \mathbb{C}\}$, generated by the characters,
as coordinate axes. The distance between a vector $g$ and a coordinate axis
can be defined as the length of the vector $g'$, where $g=g'+\alpha_z\phi_z$ and
$(g',\phi_z)=0$. Then $\|g'\|^2_2=\|g\|^2_2-q^n|\alpha_z|^2$.
Since, by Parseval's identity, $q^n\sum_z|\alpha_z|^2=\|g\|^2_2$ 
the function that is farthest from all the coordinate axes is the one whose
Fourier coefficients have equal absolute values.

%\begin{definition}
A mapping $b:Q^n_q\rightarrow Q_q$ is called a {\sl bent function}
if $|\widehat{\xi^b}(z)|=\frac{1}{q^{n/2}}|(\xi^b,\phi_z)|= 1$ for
all $z\in Q^n_q$.
%\end{definition}
For every bent function $b$, we  have, by definition  
$q^{n/2}|\alpha_z|=|\widehat{\xi^b}(z)|=1$.
Consequently, bent functions are precisely those functions whose corresponding vectors are at equal distance from all coordinate axes generated by the characters.
The axes correspond, in terms of exponents, not only to linear functions $\xi^{\langle
x,a\rangle}$  but also to affine functions $\langle x,a\rangle+c$. Further,  we will prove that bent functions are 
at the maximum possible Hamming distance from the class of affine functions.

Recall that $\widehat{\xi^b}$ is called the Walsh--Hadamard transform of $b:Q^n_q\rightarrow Q_q$ and 
 $\widehat{\xi^b}(z)$ is the Walsh--Hadamard coefficients of $b$.
 
\begin{remark} For simplicity, throughout  this section we assume that $q$ is prime, so
 the hypercube $Q^n_q$ can be regarded
as the vector space $(GF(q))^n$.   In this case every affine function $f:Q^n_q\rightarrow Q_q$
can be represented in  the form $\langle x,a\rangle+c$. If
$q=p^k$ then it is necessary to use $Tr(x\cdot a)$ instead of $\langle x,a\rangle$, where $Tr:GF(q)\rightarrow  GF(p)$ is a linear function.
\end{remark}

Bent functions are used in block and stream ciphers to increase resistance to linear cryptanalysis, in which an attacker attempts to find the best affine approximation to an unknown function. A detailed survey of the theory of bent functions can be found in
\cite{Mesnager} and \cite{Tok12}.

\begin{claim}\label{kriteriibent}
 A function $b$ is a bent function
if and only if $\xi^b
*\xi^{\widetilde{b}}(z)=q^{n/2}\delta$, where $\tilde{b}(x)=-b(-x)$ and $\delta=q^{n/2}{\bf 1}_{\{\bar 0\}}$ is the Dirac delta function.
\end{claim}
\begin{proof}
 Apply the Fourier transform to the equality $\xi^b
 *\xi^{\widetilde{b}}(z)=q^{n/2}\delta$.
From the properties of the Fourier transform
$\widehat{f(-x)}(z)={\widehat{f(x)}}(-z)$ and
$\widehat{f*g}=q^{n/2}\widehat{f}\cdot\widehat{g}$ (see Propositions \ref{cl:Four1-3} and \ref{c:convolution}) we obtain the equivalent equality
$\widehat{\xi^b}\cdot\widehat{\xi^{\widetilde{b}}}=\widehat{\xi^b}\cdot\overline{\widehat{\xi^b}}=\widehat{\delta}=\mathbf{1}$.
Thus
$|\widehat{\xi^b}(z)|=1$
for every $z$, which is precisely the definition of a bent function.
\end{proof}

\begin{claim}\label{c:bent}
Let $b:Q^n_q\rightarrow Q_q$ be a bent function and suppose that
$\widehat{\xi^b}=\xi^f$ for some function $f:Q^n_q\rightarrow
Q_q$. Then $f$ is a bent function.
\end{claim}
\begin{proof}
We have the equalities
$\widehat{\xi^{f(x)}}=\widehat{\widehat{\xi^{b(x)}}}=\xi^{b(-x)}$.
Then $|\widehat{\xi^f}|=|\xi^{b(-x)}|=1$ and hence $f$ is a bent function.
\end{proof}

A bent function is called {\sl regular} if it satisfies the condition
of Proposition~\ref{c:bent}. For $q=2$ all bent functions are regular since  the Fourier transform maps real-valued functions to real-valued functions.
A bent function $b$ is called {\sl weakly regular} if
$\widehat{\xi^b}=\alpha\xi^f$ for some $q$-ary function $f$ and some
constant $\alpha$ with  $|\alpha|=1$.

\begin{example}\label{ex:quadbent}
 Define the quadratic form
$g:Q^{2n}_q\rightarrow Q_q$ by the equation $g(x,y)=\langle x,y\rangle=
\sum\limits_{i=1}^nx_iy_i \bmod q$. We show that $g$ is a
bent function. We have
$$\widehat{\xi^g}(u,v)=\frac{1}{q^n}\sum\limits_{x,y}\xi^{\langle x,y\rangle - \langle x, u\rangle -
\langle y, v\rangle}=\frac{1}{q^n}\sum\limits_{y}\xi^{-\langle y,
v\rangle}\sum\limits_{x}\xi^{\langle x,y-u\rangle }.$$ The sum
of a character over all vertices  vanishes unless  the
character is $\phi_0$ (see Proposition~\ref{c:char_sum}). Therefore,
$\sum\limits_{x}\xi^{\langle x,y-u\rangle}=0$ when $y\neq u$. Hence,
$\widehat{\xi^g}(u,v)=\xi^{-\langle u, v\rangle}$ and consequently $g$ is a
bent function by definition.
\end{example}

For an arbitrary function $f:Q^n_q\rightarrow Q_q$ we define the operator
$D_y[f](x)=f(x+y)-f(x)$, which can be called the {\sl discrete
derivative} of $f$ in the direction $y$.

\begin{claim}\label{cl:diff}
For any nonzero $y$, the function $D_y[f]$  takes each of the $q$
values equally often  if and only if $f$ is a
bent function.
\end{claim}
\begin{proof}
We have 
$$\xi^f *\xi^{\widetilde{f}}(z)=\sum\limits_{x\in
Q^n_q}\xi^{f(z-x)}\xi^{-f(-x)}=\sum\limits_{y\in
Q^n_q}\xi^{f(z+y)-f(y)}.$$

If $z\ne\bar 0$, the last sum is the sum of the $q^n$ values of
$\xi^{D_z[f]}$. If $D_z[f]$ takes each value of $Q_q$ equally often, then, since
$\sum_{a\in Q_q}\xi^a=0$,
we obtain
\[
\xi^f*\xi^{\widetilde f}(z)=0
\qquad\text{for }z\ne\bar 0.
\]
By Proposition~\ref{kriteriibent}, $f$ is therefore bent.

Conversely, if $f$ is bent, then
\[
\xi^f*\xi^{\widetilde f}(z)=0
\qquad\text{for }z\ne\bar 0.
\]
Hence
$\sum_{y\in Q_q^n}\xi^{D_zf}=0$
for every nonzero $z$. By Proposition~\ref{c:char_sum11}, this is equivalent to $D_z[f]$ taking each value of $Q_q$ equally often.
\end{proof}

\begin{claim}\label{splat12}
Let $f:Q_q^n\rightarrow Q_q$ be a  bent function, let
$L:Q_q^n\rightarrow Q_q^n$ be a non-degenerate affine transformation
and let $\ell: Q_q^n\rightarrow Q_q$ be an affine function. Then
$g=(f\circ L)+\ell$ is a bent function.
\end{claim}
\begin{proof}
A non-degenerate affine transformation of the argument of a function  implies a non-degenerate affine transformation of  the argument of its Walsh--Hadamard transform. Adding of a linear function permutes the Walsh--Hadamard coefficients. Adding of a constant function corresponds to multiplying the Walsh--Hadamard coefficients by a constant with absolute value  equal to one. Consequently, $g$ is a bent function by definition.
\end{proof}

The functions $f$ and $g$ satisfying the conditions of Proposition \ref{splat12} are called                          {\sl EA-equivalent}.

\subsection{Boolean Bent Functions and Plateaued Functions}

Let us consider the case $q=2$ in more detail. Recall that the Walsh--Hadamard transform of a Boolean function
$f:Q_2^n\to Q_2$ is defined as the Fourier transform
$\widehat{(-1)^f}$.  Consequently, a Boolean function is bent if and only if all its Walsh--Hadamard coefficients have absolute value $1$.

\begin{claim}\label{cl:bentsupp}
If $f$ is a Boolean bent function, then
\[
\widehat f(\bar 0)
=\frac{1}{2}\left(2^{n/2}-\widehat{(-1)^f}(\bar 0)\right)
=\frac{1}{2}\left(2^{n/2}\pm1\right).
\]
\end{claim}

\begin{proof}
It is easy to see that
$(-1)^f=\mathbf{1}-2f,$
and hence
$\widehat{(-1)^f}=\widehat{\mathbf{1}}-2\widehat f
=\delta-2\widehat f$.
In particular,
$\widehat{(-1)^f}(\bar 0)
=2^{n/2}-2\widehat f(\bar 0)$.
Since $f$ is bent, the definition of a bent function gives
$\left|\widehat{(-1)^f}(\bar 0)\right|=1$.
\end{proof}
So, the weight of every Boolean bent function equals either $2^{n-1}+2^{\frac{n}{2}-1}$ or $2^{n-1}-2^{\frac{n}{2}-1}$.

Let $\mathcal A$ be the set of affine Boolean functions, i.e., functions of the form
\[
f(x)=\langle x,z\rangle
\qquad\text{or}\qquad
f(x)=\langle x,z\rangle\oplus1.
\]

The {\sl nonlinearity} of a Boolean function $f$ is its Hamming distance from the set of affine functions:
\[
nl(f)=\min_{\ell\in\mathcal A}d_H(f,\ell),\
\mbox{where}\
d_H(f,\ell)=|\supp(f\oplus\ell)|.
\]

The Hamming distance between two Boolean functions $f$ and $g$ can be expressed in terms of the inner product of the corresponding functions $(-1)^f$ and $(-1)^g$. Namely,
$d_H(f,g)
=\frac12\left(2^n-((-1)^f,(-1)^g)\right).$
Indeed,
$2^n-2d_H(f,g)
=((-1)^f,(-1)^g),$
since the inner product is the difference between the number of vertices at which $f$ and $g$ coincide and the number at which they differ.

\begin{claim}\label{cl:nl}
The nonlinearity of a Boolean function $f$ is
\[
nl(f)
=\frac12\left(
2^n-2^{n/2}\max_z
\left|\widehat{(-1)^f}(z)\right|
\right).
\]
\end{claim}

\begin{proof}
A Boolean function $\ell$ is affine if and only if
\[
(-1)^\ell=\phi_z
\qquad\text{or}\qquad
(-1)^\ell=-\phi_z
\]
for some $z\in Q_2^n$. From the definition of the Fourier transform, we have
$((-1)^f,\pm\phi_z)
=\pm2^{n/2}\widehat{(-1)^f}(z)$.
Substituting this expression into the formula for the Hamming distance and minimizing over all affine functions gives the result.
\end{proof}

\begin{claim}[\cite{Rot}]
The inequality
$nl(f)\leq \frac12\left(2^n-2^{n/2}\right)$
holds for every Boolean function $f$. Equality is attained if and only if $f$ is a bent function.
\end{claim}

\begin{proof}
By Parseval's identity,
\[
\sum_{z\in Q_2^n}
\left|\widehat{(-1)^f}(z)\right|^2
=\sum\limits_{x\in Q_2^n}|(-1)^f(x)|^2
=2^n.
\]
Therefore,
$\max_z\left|\widehat{(-1)^f}(z)\right|^2\geq1$.
Equality can hold only if
$\left|\widehat{(-1)^f}(z)\right|=1$
for every $z\in Q_2^n$. Thus, by Proposition~\ref{cl:nl},
$nl(f)\leq\frac12\left(2^n-2^{n/2}\right)$,
with equality if and only if $f$ is bent.
\end{proof}

Since the nonlinearity of a Boolean function is an integer, Boolean bent functions cannot exist for odd $n$. On the other hand, Example~\ref{ex:quadbent} shows that bent functions exist for every even $n$.

A Boolean function $f$ is called {\sl plateaued} if its Walsh--Hadamard coefficients take only three distinct values,
$0,\quad \pm a$,
for some $a\in\mathbb R$.
It follows from Parseval's identity that $a^2$ divides $2^n$. Hence
$a^2=2^s$
for some nonnegative integer $s$. A plateaued function whose Walsh--Hadamard coefficients take the values
$0,\quad \pm2^{s/2}$
is called {\sl $s$-plateaued}. In particular, a $0$-plateaued function is a bent function.
It is obvious that EA-equivalent functions  either are or are not $s$-plateaued functions. 

\begin{claim}\label{plato}
A Boolean function $f$ is plateaued if and only if
$(-1)^f*(-1)^f*(-1)^f=
a^2 2^n(-1)^f$
for some $a\in\mathbb R$.
\end{claim}

\begin{proof}
By the convolution property of the Fourier transform (Proposition~\ref{c:convolution}), the equality
$(-1)^f*(-1)^f*(-1)^f
=
a^2 2^n(-1)^f$
is equivalent to
$2^n\left(\widehat{(-1)^f}\right)^3
=
a^2 2^n\widehat{(-1)^f}$.
Thus,
$\widehat{(-1)^f}
\left(
\left(\widehat{(-1)^f}\right)^2-a^2
\right)
=0.$
Therefore, every Walsh--Hadamard coefficient satisfies
$
\widehat{(-1)^f}(z)
\in\{0,\pm a\},$
which is precisely the definition of a plateaued function. The converse follows immediately by reversing the argument.
\end{proof}

\begin{theorem}[Tarannikov at al. \cite{Taran00}, \cite{SarMai}, \cite{Zheng}]\label{th:Tar}
Let $f$ be a Boolean function with $\cor(f)\leq n-2$. Then the following inequalities hold:\\
$nl(f)\leq 2^{n-1}-2^{\cor(f)}$ if $f$ is not resilient, and\\
$nl(f)\leq 2^{n-1}-2^{1+\cor(f)}$
if $f$ is resilient.

If equality is attained in either inequality, then $f$ is plateaued.
\end{theorem}

\begin{proof}
Suppose first that $f$ is not balanced. By Remark~\ref{cor:Sar1}, which follows from Sarkar's identity, the Walsh--Hadamard coefficients have the form
$m2^{1+\cor(f)-n/2},
$
where $m$ is an integer. Hence
\[
\max_z\left|\widehat{(-1)^f}(z)\right|
\geq 2^{1+\cor(f)-n/2}.
\]
By Proposition~\ref{cl:nl}, we obtain
$nl(f)
\leq
2^{n-1}-2^{\cor(f)}$.

If equality is attained, then all nonzero Walsh--Hadamard coefficients must have absolute value
$2^{\cor(f)+1-n/2}$.
Thus all Walsh--Hadamard coefficients belong to
$\{0,
\pm2^{\cor(f)+1-n/2}\}$,
and $f$ is plateaued.

The balanced case is analogous.
\end{proof}

We describe a method of constructing bent functions and plateaued functions
that generalizes Example \ref{ex:quadbent}.
Let $x\in Q^n_2$ and $s\leq n$. Denote by $x^s$ the vector of the
first $s$ coordinates of the vector $x$.

\begin{claim}[Maiorana--McFarland construction]\label{cplat133}
Let $x,y,u,v\in Q^n_2$, $z,w\in Q^s_2$, $s\leq n$. Consider
permutations $\tau:Q^n_2\rightarrow Q^n_2$ and
$\sigma:Q^s_2\rightarrow Q^s_2$ and an arbitrary function
$f:Q^n_2\rightarrow Q_2$.  Then $F(x,y,z)=\langle\tau(x),y\rangle +
\langle\sigma(x^s),z\rangle +f(x)$ is a $s$-plateaued function
of $2n+s$ variables.
\end{claim}
\begin{proof}
By the definition of the Fourier transform,
$$2^{n+\frac{s}{2}}\widehat{(-1)^F}(u,v,w)= \sum\limits_{x,y,z\in
Q^n_2}(-1)^{\langle\tau(x),y\rangle + \langle\sigma(x^s),z\rangle
+f(x)-\langle x,u\rangle-\langle v,y\rangle -\langle w,z\rangle}$$
$$=\sum\limits_{x\in Q^n_2}(-1)^{f(x)-\langle
x,u\rangle}\sum\limits_{y\in
Q^n_2}(-1)^{\langle\tau(x),y\rangle-\langle
v,y\rangle}\sum\limits_{z\in Q^s_2}(-1)^{\langle\sigma(x^s),z\rangle
-\langle w,z\rangle}.$$

The sums $\sum\limits_{z\in
Q^s_2}(-1)^{\langle\sigma(x^s),z\rangle -\langle w,z\rangle}$ are equal to
$2^s$ if $w=\sigma(x^s)$, and is equal to $0$ otherwise.
Similarly,
 the sums
$\sum\limits_{y\in Q^n_2}(-1)^{\langle\tau(x),y\rangle -\langle
v,y\rangle}$ is  equal to $2^n$ if $v=\tau(x)$ and is  equal to $0$ otherwise.

Consequently,
$2^{n+\frac{s}{2}}\widehat{(-1)^F}(u,v,w)=2^{n+s}(-1)^{f(x)-\langle
x,u\rangle}$ if $v=\tau(x)$ and $w=\sigma(x^s)$; otherwise
$\widehat{(-1)^F}(u,v,w)=0$.

Thus every nonzero Walsh--Hadamard coefficient of $F$ has absolute value
$\frac{2^{n+s}}{2^{n+s/2}}
=2^{s/2}$.
Hence $F$ is $s$-plateaued.
\end{proof}

\subsection{Nonlinearity of $q$-ary Functions}

Every affine function $g:Q_q^n\rightarrow Q_q$ can be represented as
$g(x)=\langle z,x\rangle+a$. For an arbitrary function
$f:Q_q^n\rightarrow Q_q$ we define the following quantities
$\delta^{z,a}_f=\frac{q-1}{q}- \frac{d_H(f,g)}{q^n}$, where
$g(x)=\langle z,x\rangle+a$.

The following properties of the quantities
$\delta^{z,a}_f$ follow directly from the definition.
\begin{claim}\label{rbent22}\ \quad\\
{\rm (a)} For any fixed  $a\in Q_q$ and $z\in Q^n_q$ 
$\sum\limits_{f}\delta^{z,a}_f=0$, where the sum is taken over all
functions $f: Q^n_q\rightarrow Q_q$.\\
{\rm (b)}  For any function $f$ and $z\in Q^n_q$ we have $\sum\limits_{a\in
Q_q}\delta^{z,a}_f=0$.\\
{\rm (c)}  The minimum value of $\delta^{z,a}_f$ is
$-\frac{1}{q}$, and the maximum value is $\frac{q-1}{q}$.
\end{claim}

As in the Boolean case, we define the nonlinearity of a $q$-ary
function $f:Q_q^n\rightarrow Q_q$ as the Hamming distance from $f$ to
the set of affine functions. From the definitions follows the equality
\begin{equation}\label{eqplat1}
nl(f)=q^{n-1}(q-1)-q^n\max\limits_{z,a}\delta^{z,a}_f.
\end{equation}

The following theorem can be viewed as another representation of Parseval's identity.

\begin{theorem}[Ryabov \cite{Ryabov20}]\label{rbent1}
$\sum\limits_{z\in Q^n_q ,a\in
Q_q}(\delta^{z,a}_f)^2=\frac{q-1}{q}$.
\end{theorem}
\begin{proof}
We prove the theorem by induction on $n$. For $n=0$ the equality is easily verified
directly. Let us prove the induction step. From the definitions we have
$$\delta^{z,a}_f=\frac1q\sum\limits_{c\in
Q_q}\delta_{f|_{x_n=c}}^{z',a+cz_n},\qquad \text{where}\quad
z=(z',z_n).$$ Therefore,
$$(\delta^{z,a}_f)^2=\frac{1}{q^2}\sum\limits_{c_1,c_2\in
Q_q}\delta_{f|_{x_n=c_1}}^{z',a+c_1z_n}\delta_{f|_{x_n=c_2}}^{z',a+c_2z_n}=
\frac{1}{q^2}\sum\limits_{c_1\neq c_2
}\delta_{f|_{x_n=c_1}}^{z',a+c_1z_n}\delta_{f|_{x_n=c_2}}^{z',a+c_2z_n}+
\frac{1}{q^2}\sum\limits_{c\in
Q_q}(\delta_{f|_{x_n=c}}^{z',a+cz_n})^2.$$

By the induction hypothesis we have
 $$\sum\limits_{z\in Q^n_q ,a\in
Q_q}\frac{1}{q^2}\sum\limits_{c\in
Q_q}(\delta_{f|_{x_n=c}}^{z',a+cz_n})^2=\frac{1}{q^2}\sum\limits_{c\in
Q_q}\sum\limits_{z_n\in Q_q}\sum\limits_{z'\in Q^{n-1}_q ,a\in
Q_q}(\delta_{f|_{x_n=c}}^{z',a+cz_n})^2=\frac{q^2}{q^2}\frac{q-1}{q}.$$

It remains to show that the contribution of the first term is zero. We have

$$\sum\limits_{z\in Q^n_q ,a\in
Q_q}\sum\limits_{c_1\neq c_2
}\delta_{f|_{x_n=c_1}}^{z',a+c_1z_n}\delta_{f|_{x_n=c_2}}^{z',a+c_2z_n}=
\sum\limits_{c_1\neq c_2 }\sum\limits_{z'\in Q^{n-1}_q
}\sum\limits_{z_n\in Q_q }\sum\limits_{a\in
Q_q}\delta_{f|_{x_n=c_1}}^{z',a+c_1z_n}\delta_{f|_{x_n=c_2}}^{z',a+c_2z_n}=$$
$$=\sum\limits_{c_1\neq c_2 }\sum\limits_{z'\in Q^{n-1}_q
}\left(\sum\limits_{a_1\in
Q_q}\delta_{f|_{x_n=c_1}}^{z',a_1}\right)\left(\sum\limits_{a_2\in
Q_q}\delta_{f|_{x_n=c_2}}^{z',a_2}\right)=0.$$

The penultimate equality follows from the fact that the system of equations
$
\left\{
\begin{array}{l}
a+c_1z_n=a_1;\\
a+c_2z_n=a_2.
\end{array}
\right.$ has a unique solution for the variables $a$ and $z_n$.
The last equality follows from Proposition \ref{rbent22}(b).
\end{proof}

\begin{claim}[]\label{rbent2}
For any function $f:Q_q^n\rightarrow Q_q$ , the following statements hold.\\
{\rm (a)}  $\widehat{\xi^f}(z)=q^{n/2}\sum\limits_{a\in Q_q}\delta^{z,a}_f\xi^a$.\\
{\rm (b)}  If $\delta^{z,a}_f=0$ for all $a\in Q_q$, then
$\widehat{\xi^f}(z)=0$. The converse holds if $q$ is prime.
\end{claim}
\begin{proof}
(a) Let $L_a=|\{x\in Q^n_q : f(x)-\langle x,z\rangle=a\}|$. Express
$\widehat{\xi^f}(z)$
 and $\delta^{z,a}_f$ in terms of $L_a$.
\begin{equation}\label{eqplat5}
\widehat{\xi^f}(z)=\frac{1}{q^{n/2}}\sum\limits_{x\in
Q^n_q}\xi^{f(x)-\langle
x,z\rangle}=\frac{1}{q^{n/2}}\sum\limits_{a\in Q_q}L_a\xi^a.
\end{equation}
$$\delta^{z,a}_f=\frac{q-1}{q}-\frac{q^n-L_a}{q^n}=\frac{L_a}{q^n}-\frac{1}{q}.$$
Consequently,
$L_a=q^n\delta_f^{z,a}+q^{n-1}$.
Using $\sum\limits_{a\in Q_q}\xi^a=0$ we obtain
$$\widehat{\xi^f}(z)=\frac{1}{q^{n/2}}\sum\limits_{a\in Q_q}(q^n\delta^{z,a}_f+q^{n-1})
\xi^a=q^{n/2}\sum\limits_{a\in Q_q}\delta^{z,a}_f\xi^a.$$ This proves the first statement.

(b) 
If $\delta_f^{z,a}=0$ for every $a\in Q_q$, then the first statement immediately gives
$\widehat{\xi^f}(z)=0$.

 Using equality (\ref{eqplat5}) and Proposition
\ref{c:char_sum11}, we conclude that if $q$
is prime and $\widehat{\xi^f}(z)=0$, then $L_a=q^{n-1}$ for every $a\in Q_q$.
Consequently, $\delta^{z,a}_f=0$ for every $a\in Q_q$.
\end{proof}

\begin{theorem}[Ryabov \cite{Ryabov21}]\label{thbent3}\quad \\ Let $q$ be prime. For any function $f:Q_q^n\rightarrow Q_q$ the following statement holds. \\
{\rm (a)} 
$nl(f)\leq (q-1)q^{n-1}-q^{\frac{n}{2}-1}$.\\
 {\rm (b)} If $nl(f)= (q-1)q^{n-1}-q^{\frac{n}{2}-1}$, then
 $f$ is a weakly regular bent function and $n$ is even.
\end{theorem}
\begin{proof}
(a)  To maximize
$nl(f)$, we need to solve the following optimization problem:

$\sum\limits_{z\in Q^n_q ,a\in
Q_q}(\delta^{z,a}_f)^2=\frac{q-1}{q},$

$\sum\limits_{a\in Q_q} \delta^{z,a}_f=0$ for every $z\in Q^n_q$,

$\max\limits_{z\in Q^n_q ,a\in Q_q}\delta^{z,a}_f\longrightarrow \min.$

Consider the following auxiliary problem: $\sum\limits_{a\in
Q_q}\sigma_{a}=0$, $\sum\limits_{a\in
Q_q}\sigma_{a}^2=A\frac{q-1}{q}$ and
$\max\limits_{a}\sigma_{a}\rightarrow \min.$ It is easy to see that, up to permutation, the unique solution  is $\sigma_0=-\frac{q-1}{q}\sqrt{A}$,
$\sigma_a=\frac{\sqrt{A}}{q}$ for $a\neq 0$.

The value $\max\limits_{a}\sigma_{a}$ increases with $A$. Then the optimal solution of the original problem
corresponds to the case when all the sums $\sum\limits_{a\in
Q_q}(\delta^{z,a}_f)^2$ are equal. Take $A=\frac{1}{q^n}$. Then
we obtain that  $\min\limits_f(\max\limits_{z\in Q^n_q,a\in
Q_q}\delta^{z,a}_f)\geq \frac{1}{q^{\frac{n}{2}+1}}$.
This proves  (a).

(b) Suppose now that equality is attained. Then equality must hold at every step of the above optimization argument. In particular, for every $z\in Q_q^n$ the values $\delta^{z,a}_f$, $a\in Q_q$,
 must consist of one value
$-\frac{q-1}{q^{\frac{n}{2}+1}}$  and $q-1$  values
$\frac{1}{q^{\frac{n}{2}+1}}$. Hence, by Propositions \ref{rbent2} (a),  $\widehat{\xi^f}(z)=-\xi^a$ for some $a\in Q_q$.
Then $f$ is a weakly regular bent function. Since the number
$nl(f)$ is an integer, $n$ must be even.
\end{proof}

\subsection{Problems}

\begin{exercise}
Let $f$ be a Boolean function and $F(y)=\sum\limits_{x\in Q^n_2}(-1)^{D_y[f](x)}$. Prove that each Fourier coefficient of $F$ is nonnegative. 
\end{exercise}

\begin{exercise}
Let $f$ be a Boolean function.  Prove that $$\sum\limits_{x,y,z\in Q^n_2}(-1)^{D_z[xD_y[f]](x)}=2^n\sum\limits_{v\in Q^n_2} 
(\widehat{(-1)^f}(v))^4.$$
\end{exercise}

\begin{exercise}
Let $f$ be a Boolean function.  Prove that  $$\sum\limits_{x,y,z\in Q^n_2}(-1)^{D_z[xD_y[f]](x)}\geq 2^{2n}$$ and this equality holds if and only if  $f$ is a bent function. 
\end{exercise}

\begin{exercise}
Prove that the Hamming distance between two distinct
bent functions in the  Boolean $n$-cube cannot be less than
$2^{n/2}$.
\end{exercise}

\begin{exercise}
For an arbitrary Boolean function $f:Q_2^n\rightarrow Q_2$ define
the entropy $H(f)$ by $H(f)=-\sum\limits_{x\in Q_2^n}p(x)\log_2
p(x)$, where $p(x)=(\widehat{(-1)^f}(x))^2/2^{n}$. Prove that for
even $n$ the maximum entropy equal to $n2^n$ is attained only by
bent functions.
\end{exercise}

\begin{exercise}
Prove that the restriction of a bent function to a hyperface is a
plateaued function.
\end{exercise}

\begin{exercise}
Prove that there are no Boolean functions of $n$ variables
that attain the Tarannikov bound when $\cor(f)\leq n/2$.
\end{exercise}

\section{Algebraic Normal Form of Boolean Functions}

\subsection{Möbius Transformation}

Two bases $\{u_1,\dots,u_n\}$ and $\{v_1,\dots,v_n\}$ in a vector
space over a field $F$ are called {\sl biorthogonal} if
$\langle u_i,v_j\rangle=0$ for $i\neq j$ and
$\langle u_i,v_i\rangle=1$ for all $i,j\in\{1,\dots,n\}$.
Let $U$ and $V$ be the square matrices whose rows are the
corresponding basis vectors. Then
$
U^{\mathsf T}V=I_n,
$
that is,
$
U^{\mathsf T}=V^{-1}.
$

Consider a finite set $T$ equipped with a partial order $\leq$.
For each $z\in T$, define the function
\[
\omega_z(x)=
\begin{cases}
1,&x\leq z,\\
0,&\text{otherwise}.
\end{cases}
\]
Consider the vector space $\mathbb{V}(T)$ of functions
$f:T\rightarrow F$. The number of functions $\omega_z$ is equal to
the dimension $|T|$ of the space $\mathbb{V}(T)$. Let $W$ be the
matrix whose rows are the vectors of values of the functions
$\omega_z$. It is easy to see that, if the elements of $T$ are
ordered compatibly with the partial order, then $W$ is a lower
triangular matrix with ones on the main diagonal. Consequently,
$W$ is nonsingular, and there exists a basis
$\{\nu_y:y\in T\}$ biorthogonal to the basis $\{\omega_z:z\in T\}$.

%\begin{definition}
Multiplication of a vector $f:T\rightarrow F$ by the matrix $W$
is called the {\sl Möbius transformation} of $f$. Multiplication
by the matrix $V^{\mathsf T}$, where $V$ is formed from the
functions $\nu_y$, is called the {\sl inverse Möbius transformation}.
%\end{definition}

The matrices $W$ and $V^{\mathsf T}$ are mutually inverse by definition, i.e.,
\[
V^{\mathsf T}(Wf)=(V^{\mathsf T}W)f=f.
\]
It follows that $V^{\mathsf T}$ is also lower triangular with ones
on the main diagonal, while $V$ is upper triangular. Thus,
\[
Wf(z)=\sum_{x\leq z}\omega_z(x)f(x)
\]
and
\[
V^{\mathsf T}g(y)=\sum_{z\leq y}\nu_z(y)g(z),
\]
since
\[
V^{\mathsf T}(y,z)=V(z,y)=\nu_z(y).
\]

Since the matrices $W$ and $V^{\mathsf T}$ are mutually inverse
lower triangular matrices, we also obtain the following.

\begin{claim}
For any distinct $y,z\in T$,
\[
\sum_{x:\,y\leq x\leq z}\nu_x(z)=0
\qquad\text{and}\qquad
\sum_{x:\,y\leq x\leq z}\nu_y(x)=0.
\]
\end{claim}

\begin{proof}
We have
\[
\begin{aligned}
\sum_{x:\,y\leq x\leq z}\nu_x(z)
&=\sum_{x:\,y\leq x\leq z}V^{\mathsf T}(z,x)\\
=\sum_{x:\,y\leq x\leq z}V^{\mathsf T}(z,x)W(x,y)&
=\sum_xV^{\mathsf T}(z,x)W(x,y)=0.
\end{aligned}
\]
The second equality follows in the same way by using that $WV^{\mathsf T}=I_n$.
\end{proof}

In the case of a linear order, the matrix $W$ has ones in all
positions on and below the main diagonal, while the matrix
$V^{\mathsf T}$ has only two nonzero diagonals: the main diagonal,
consisting of ones, and the first subdiagonal, consisting of $-1$'s.

As a partially ordered set, we can consider the Boolean hypercube
with the lexicographic order. In this case, the basis vectors
$\omega_z$ are the characteristic functions of the faces containing
$\bar 0$ and $z$. Let $F=\mathbb{C}$. It is easy to verify that the biorthogonal basis
$\{\nu_y\}$ consists of the characteristic functions of the faces
containing $\bar 1$ and $y$, multiplied by the parity character $(-1)^{x_1\oplus\dots\oplus x_n}$.
If the field $F$ has characteristic $2$, multiplication by the
parity character is not needed (see Proposition \ref{stsp0}).

\begin{claim}
Let $f:Q_2^n\rightarrow\mathbb{C}$ be an eigenfunction with
eigenvalue $\lambda=n-2k$. Then:

{\rm (a)} $f$ is determined by the sums of its values over the
$(n-k)$-dimensional faces containing the vertex $\bar 0$;

{\rm (b)}
the restriction of $f$ to the ball of radius $d\leq n-k$ centered
at the vertex $\bar 0$ is determined by the sums of its values over
the $d$-dimensional faces containing $\bar 0$.

\end{claim}

\begin{proof}

(a)  Consider Corollary \ref{cor:Sar} of Sarkar's identity
\[
\sum_{z\in\Gamma^\perp}\widehat{f}(z)
=
2^{n/2-\dim\Gamma}
\sum_{x\in\Gamma}f(x)
\]
for $\dim\Gamma=n-k$. By Proposition \ref{claim:ef},
$\widehat{f}(z)\neq0$ only if $\wt(z)=k$. Hence, the sum on the
left-hand side has exactly one nonzero term, and therefore
\[
\widehat{f}(z)
=
2^{k-n/2}\sum_{x\in\Gamma}f(x).
\]
The function $f$ can then be recovered from
$f=\widehat{\widehat{f}}$.

(b) We first show that, from the sums
$
\sum_{x\in\Gamma'}f(x)
$
over faces $\Gamma'$ of dimension $m\leq n-k$, one can compute the
sums
$
\sum_{x\in\Gamma}f(x)
$
over faces $\Gamma$ of dimension $m-1$.

There are $n-(m-1)$ faces $\Gamma'$ of dimension $m$ containing a
fixed face $\Gamma$ of dimension $m-1$. Any two of these faces
intersect exactly in $\Gamma$. Therefore,
\[
\sum_{\Gamma\subset\Gamma'}
  \sum_{x\in\Gamma'}f(x)
-(n-(m-1))\sum_{x\in\Gamma}f(x)
=
\sum_{x\in\Gamma}\sum_{y:\,d(x,y)=1}f(y)
-(m-1)\sum_{x\in\Gamma}f(x).
\]
Indeed, the left-hand side counts the values of $f$ on all faces
$\Gamma'$ of dimension $m$ containing $\Gamma$, after subtracting
the contribution of their common intersection $\Gamma$.
 The right-hand side counts the values of $f$ on all neighbors of vertices from  $\Gamma$.

Since $f$ is an eigenfunction,
$
\sum\limits_{y:\,d(x,y)=1}f(y)=\lambda f(x)$.
Consequently,
\[
\sum_{\Gamma\subset\Gamma'}
\sum_{x\in\Gamma'}f(x)
=
\bigl(n-2(m-1)+\lambda\bigr)
\sum_{x\in\Gamma}f(x).
\]

Let $m\leq d\leq n-k$. By the assumption on $\lambda$,
\[
n-2(m-1)+\lambda
\geq n-2(n-k-1)+n-2k
=2.
\]
Therefore, the sum of the values of $f$ over any face of dimension
less than $d$ can be computed recursively if the corresponding sums
over faces of dimension $d$ are known. 

Consider a natural partial order on the Boolean hypercube (see the next section). Then functions $\omega_z(x)$ are indicator functions of faces. The inverse Möbius
transformation then recovers the values of $f$ from  sums over faces.
\end{proof}

Note that, instead of the vertex $\bar 0$, any vertex of the
hypercube can be taken as the common intersection of all the faces.

\subsection{Algebraic Degree of  Boolean Functions}

We next consider a particular case of the Möbius transformation over
$GF(2)$ on the Boolean hypercube with the following natural order.  We define
\[
x\leq y\quad\Longleftrightarrow\quad
\supp(x)\subseteq\supp(y).
\]
Notice that $\{x\in Q_2^n:x\leq y\}$
is a face of dimension $\wt(y)$ containing the vertex $\bar 0$.

Every Boolean function $f:Q_2^n\rightarrow Q_2$ can be represented by 
in {\sl algebraic normal
form},
\begin{equation}\label{eqZhegal}
f(x_1,\dots,x_n)=\bigoplus\limits_{y\in Q_2^n}G[f](y)
x_1^{y_1}\cdots x_n^{y_n},
\end{equation}
where $x^0=1$, $x^1=x$, and
$G[f]:Q_2^n\rightarrow Q_2$.
The existence of such a representation follows from the fact that
addition and multiplication over $GF(2)$, together with the constants
$0$ and $1$, form a complete basis. In particular, multiplication
coincides with conjunction, while disjunction can be expressed as
\[
x\vee y=x\oplus y\oplus xy.
\]
A Boolean function has a unique representation in algebraic normal
form, since the operator $G$ maps Boolean functions to Boolean
functions in the same number of variables.

%\begin{definition}
The {\sl algebraic degree} of a Boolean function $f$ is the maximum
degree of a monomial in its algebraic normal form, i.e.,
\[
\deg(f)=\max_{G[f](y)=1}\wt(y).
\]
%\end{definition}

The algebraic degree of a set $C\subset Q_2^n$ is defined as the
algebraic degree of its characteristic function. 
The following statement holds.

\begin{claim}\label{stsp0}
For every Boolean function $f$,
\[
G[f](y)=\bigoplus\limits_{x\in Q_2^n,\;x\leq y}f(x),
\]
i.e., $G$ is the Möbius transformation of $f$ over $GF(2)$.
Moreover,
$G[G[f]]=f$.
\end{claim}

\begin{proof}
It is easy to see that if
$\supp(y)\not\subseteq\supp(x)$, then
$x_1^{y_1}\cdots x_n^{y_n}=0$,
whereas otherwise
$
x_1^{y_1}\cdots x_n^{y_n}=1$.
Fix $z\in Q_2^n$. We show that
\[
\bigoplus\limits_{x\in Q_2^n,\;x\leq z}
x_1^{y_1}\cdots x_n^{y_n}=1
\]
if and only if $z=y$. Indeed, if $y<z$, then the sum contains
$2^{\wt(z)-\wt(y)}$ ones and hence is equal to $0$. If $z$ and $y$
are incomparable or $z<y$, then all summands are equal to $0$.
Summing both sides of (\ref{eqZhegal}) over the face
\[
\{x\in Q_2^n:x\leq z\},
\]
we obtain
\[
\bigoplus\limits_{x\in Q_2^n,\;x\leq z}f(x)=G[f](z).
\]

Now fix $x\in Q_2^n$. On the right-hand side of (\ref{eqZhegal}),
the monomial
$x_1^{y_1}\cdots x_n^{y_n}$ is equal to $1$ precisely when
$y\leq x$. Therefore,
\[
f(x)=\bigoplus\limits_{y\in Q_2^n,\;y\leq x}G[f](y).
\]
Hence
$G[G[f]]=f $
for every Boolean function $f$.
\end{proof}

We say that a Boolean function is {\sl balanced modulo $2$} on a
face $\Gamma$ if the Boolean sum of its values on $\Gamma$ is equal
to $0$.

\begin{claim}\label{stsp1}
A Boolean function $f:Q_2^n\rightarrow Q_2$ is balanced modulo $2$
on all faces of dimension $m+1$ if and only if
$\deg(f)\leq m$.
\end{claim}
\begin{proof}
On the one hand, every monomial of degree less than the dimension of
a face $\Gamma$ takes the value $1$ an even number of times on
$\Gamma$. Therefore, $f$ is balanced modulo $2$ on every face whose
dimension is greater than $\deg(f)$.

On the other hand,
\[
\bigoplus\limits_{x\in Q_2^n,\;x\leq z}
x_1^{y_1}\cdots x_n^{y_n}=1
\]
if and only if $z=y$. Therefore, if the algebraic normal form
contains a monomial of degree $k$, $k>m$, then there exists a
$k$-dimensional face containing $\bar 0$ on which the sum of the
values of $f$ is equal to $1$.  Every $k$-dimensional face
can be partitioned into  disjoint $(m+1)$-dimensional faces.  Then  the sum of the
values of $f$ on one of these faces is equal to $1$. 
\end{proof}

\begin{claim}\label{stsp2}%[\cite{Taran11}]
Suppose that a Boolean function $f$ has an even number of ones on a
face of dimension
$m=n-\cor(f)$.
Then
$\deg(f)\leq m-1$.
If this number is odd, then
$\deg(f)=m$
and
$G[f](z)=1$
for every $z\in Q_2^n$ with $\wt(z)=m$.
\end{claim}

\begin{proof}
By the definition of correlation immunity, $f$ has the same number
of ones on all faces of dimension $m$. If this number is even, then
$\deg(f)\leq m-1$ by Proposition \ref{stsp1}. If it is odd, then every
face of dimension $m+1$ contains an even number of ones. Moreover,
since every face of dimension $m$ contains an odd number of ones,
\[
G[f](z)
=
\bigoplus\limits_{x\in Q_2^n,\;x\leq z}f(x)
=1
\]
for every $z\in Q_2^n$ with $\wt(z)=m$.
\end{proof}

\begin{corollary}[Siegenthaler inequality]\label{zi}
For every Boolean function $f:Q_2^n\rightarrow Q_2$,
\[
\deg(f)\leq n-\cor(f),
\]
and
\[
\deg(f)\leq n-\cor(f)-1
\]
for  resilient functions $f$, except for affine functions. 
\end{corollary}

By Proposition \ref{c:cor_immun_for_perfect_coloring}, a perfect
coloring of the $n$-dimensional Boolean hypercube with quotient matrix
\begin{equation}\label{matr}
\left(
\begin{array}{cc}
n-b & b\\
c & n-c
\end{array}
\right)
\end{equation}
is a correlation-immune function of order
\[
\frac{b+c}{2}-1.
\]
Therefore, Proposition \ref{stsp2} implies the following.

\begin{corollary}\label{stsp3}
Let $f:Q_2^n\rightarrow Q_2$ be a perfect coloring with quotient
matrix (\ref{matr}). Then
\[
\deg(f)\leq n-\frac{b+c}{2}+1.
\]
\end{corollary}

Boolean functions
$f:Q_2^n\rightarrow Q_2$ can be represented by their vectors of
values and hence regarded as elements of the Boolean cube of
dimension $2^n$. The set of value vectors of Boolean functions of
algebraic degree at most $m$ is called the {\sl Reed--Muller code} of type
$\mathcal{R}(m,n)$ in $Q_2^{2^n}$. The Reed--Muller code is linear,
since the degree of the sum of two polynomials does not exceed the
maximum of their degrees.

\begin{claim}[\cite{MacW}]\label{stki4}
For every nonzero Boolean function $f$,
\[
|\supp(f)|\geq 2^{n-\deg(f)}.
\]
If
$|\supp(f)|=2^{n-\deg(f)}$,
then $\supp(f)$ is an affine subspace.
\end{claim}

\begin{proof}
We use induction on $n$. If the function $f$ is equal to $1$ at a
single point, then its sum over the entire Boolean hypercube is odd, and
hence its algebraic normal form contains a term of degree $n$.

Suppose that $f$ takes the value $1$ at more than one point. Then,
in at least one direction, each of the two hyperfaces contains at
least one point at which $f$ takes the value $1$. The degree of the
restriction of $f$ to a hyperface does not exceed the degree of $f$.
By the induction hypothesis, each of the two hyperfaces contains at
least
$2^{(n-1)-\deg(f)}$
ones of $f$.

A Boolean function is an indicator function of an affine subspace  if and only if its intersection with
every face in every direction has cardinality $2^i$, $i=0,\dots,n$, or $0$. The
equality case can be proved by induction as well.
\end{proof}

\begin{corollary}
The minimum distance of the Reed--Muller code
$\mathcal{R}(m,n)$ is
$2^{n-m}$.
\end{corollary}

\subsection{Polynomial Representation of  Linear Codes}

By Theorem \ref{Markov} every isometry of a Boolean hypercube is the composition of a shift  and  permutation of variables.
A shift of a set $C\subset Q_2^n$ by any vector $y\in Q_2^n$
is performed by the change of variables $x'=x\oplus y$ of the indicator function ${\bf 1}_C$.
Therefore, isometric transformations do not change the degree of the
characteristic function of a set. An affine subspace $C$ can be transformed
into a linear subspace by adding to it any vector belonging to $C$.
Any linear subspace $C$ can be specified by a full-rank parity-check
matrix $H$, i.e.,
\[
C=\{x\in Q_2^n\mid Hx=\bar{0}\},
\]
where
$H=\{h_{ij}\}$ is a $k\times n$ matrix and
$\rank(C)=n-k$. 

\begin{claim}
The degree of the characteristic function of a linear code
$C\subset Q_2^n$ is equal to $n-\rank(C)$.
\end{claim}

\begin{proof}
Let
$
C=\{x\in Q_2^n\mid Hx=\bar{0}\}$,
where $H=\{h_{ij}\}$ is a $k\times n$ matrix and
$k=n-\rank(C)$. Consider the set
$C'=\{x\in Q_2^n\mid Hx=\bar{1}\}$.
Clearly, $C'$ is an affine subspace equivalent to the code $C$.
It is easy to see that
\[
{\mathbf 1}_{C'}(x)=\prod_i\left(\bigoplus_j h_{ij}x_j\right),
\]
i.e., the product of the sums of the rows of the matrix $H[x]$.
Clearly, $\deg({\mathbf 1}_{C'})\leq k$.

Consider a nonsingular $k\times k$ submatrix $A$ of $H$.
Since $A$ is a matrix over $GF(2)$, its nonsingularity is equivalent to
$\det(A)=1$. It follows from $\det(A)=1$ that the monomial corresponding
to the columns of $A$ occurs an odd number of times in the product
\[
\prod_i\left(\bigoplus_j h_{ij}x_j\right).
\]
Hence this product contains a monomial of degree $k$, and therefore
$\deg({\mathbf 1}_{C'})=k$.
Since $C'$ is equivalent to $C$, we obtain
$\deg({\mathbf 1}_C)=k=n-\rank(C)$.
\end{proof}

\begin{corollary}
Let $C\subset Q_2^n$ be a Hamming code, where $n=2^t-1$. Then
$\deg({\mathbf 1}_C)=t$.
\end{corollary}

\subsection{Algebraic Degrees of Bent and Plateaued Functions}

The following statement is  obtained from Sarkar's identity.

\begin{claim}\label{bitrpro8}
Let $f$ be a Boolean function and suppose that, for every
$v\in Q_2^n$,
$\widehat{(-1)^f}(v)=2^k m(v)$,
where $m(v)$ and $k$ are  integer. Then
$\deg(f)\leq\frac{n}{2}-k+1$.
\end{claim}

\begin{proof}
Consider a monomial
$x_1^{y_1}\cdots x_n^{y_n}$
of maximum degree in the algebraic normal form of $f$, i.e.,
\[
G[f](y)=1
\qquad\text{and}\qquad
\wt(y)=\deg(f).
\]
Then
\[
\sum_{x\leq y}f(x)\equiv1\pmod 2,
\]
and hence
\[
\sum_{x\leq y}(-1)^{f(x)}=2(m'+1),
\]
where $m'$ is even and may be negative.

By Sarkar's identity (Corollary \ref{cor:Sar}),
\begin{equation}\label{eq:MT1}
\sum_{v\in Q_2^n,\;v\leq y\oplus\bar 1}
\widehat{(-1)^f}(v)
=
2^{n/2-\wt(y)}
\sum_{v\in Q_2^n,\,x\leq y}(-1)^{f(x)}
=
(1+m')2^{n/2-\deg(f)+1}.
\end{equation}
By the assumption, the left-hand side is divisible by $2^k$.
Therefore,
$k\leq\frac n2-\deg(f)+1$,
which proves the proposition.
\end{proof}
Note that $k$ in this statement may also be negative.

\begin{theorem}\label{cl:platb} \quad \\
{\rm (a)} If $f$ is a bent function, then
$\deg(f)\leq n/2$
for $n\geq4$;\\
{\rm (b)}  If $f$ is resilient and
$nl(f)=2^{n-1}-2^{\cor(f)+1}$,
then
$\deg(f)=n-\cor(f)-1$;\\
{\rm (c)}  $f$ is not resilient and
$
nl(f)=2^{n-1}-2^{\cor(f)},
$
then
$
\deg(f)=n-\cor(f)$ and,
moreover,
$G[f](y)=1$
for every $y$ with $\wt(y)=n-\cor(f)$.
\end{theorem}

\begin{proof}
(a) For a bent function, we have
$\left|\widehat{(-1)^f}\right|=\mathbf{1}$.
It follows from Proposition \ref{bitrpro8} that
$\deg(f)\leq\frac{n}{2}+1$.
We now improve this bound. Let $\wt(y)=\deg(f)$. A face of dimension
$n-\deg(f)$ contains
$2^{n-\deg(f)}$ points, and
$n-\deg(f)>0$ for $n>2$. Therefore, the sum
\[
\sum_{v\in Q_2^n,\;v\leq y\oplus\bar 1}
\widehat{(-1)^f}(v)
\]
in (\ref{eq:MT1}) is even. It follows from (\ref{eq:MT1}) that
$2^{n/2-\deg(f)}$
is an integer. Hence
$\deg(f)\leq n/2$.

(b) By Corollary \ref{zi}, it is sufficient to prove that
\[
\deg(f)\geq n-\cor(f)-1.
\]
By Sarkar's identity (Remark \ref{cor:Sar1}) and Proposition
\ref{cl:nl}, all nonzero Walsh--Hadamard coefficients of $f$ have
the form
$\pm2^{2+\cor(f)-n/2}$.
By Theorem \ref{t:cor-im}, there exists
$y\in Q_2^n$ with
$\wt(y)=\cor(f)+1$
such that
\[
\left|\widehat{(-1)^f}(y)\right|
=2^{2+\cor(f)-n/2}.
\]
Similarly to (\ref{eq:MT1}), we have
\[
\pm2^{2+\cor(f)-n/2}
=
\sum_{v\in Q_2^n,\;v\leq y}
\widehat{(-1)^f}(v)
=
2^{\frac n2-\wt(y\oplus\bar 1)}
\sum_{x\leq y\oplus\bar 1}(-1)^{f(x)}\]
\[=
2^{\cor(f)-\frac n2+1}
\sum_{x\leq y\oplus\bar 1}(-1)^{f(x)}.
\]
Consequently,
\[
\sum_{x\leq y\oplus\bar 1}(-1)^{f(x)}=\pm2.
\]
This is possible only if
$\sum_{x\leq y\oplus\bar 1}f(x)\equiv1\pmod2$.
Therefore,
$$
\deg(f)\geq\wt(y\oplus\bar 1)
=n-\cor(f)-1.$$

(c)
By Corollary \ref{zi}, it is sufficient to prove that
\[
\deg(f)\geq n-\cor(f).
\]
 By Sarkar's identity (Remark \ref{cor:Sar1}) and Proposition
\ref{cl:nl}, all nonzero Walsh--Hadamard coefficients of $f$ have
the form
$\pm2^{\cor(f)+1-n/2}$.

Since $f$ is not resilient,
\[
\left|\widehat{(-1)^f}(\bar0)\right|
=2^{\cor(f)+1-n/2}.
\]
By Theorem \ref{t:cor-im},
$\widehat{(-1)^f}(z)=0$
for all $z$ satisfying
$0<\wt(z)\leq\cor(f)$.

Consider an arbitrary $y$ such that
\[
\wt(y)=n-\cor(f)
\qquad\text{and}\qquad
\wt(y\oplus\bar1)=\cor(f).
\]
Then, from the left-hand side of (\ref{eq:MT1}),
\[
2^{n/2-\wt(y)}
\sum_{x\leq y}(-1)^{f(x)}
=
\widehat{(-1)^f}(\bar0)
=
\pm2^{\cor(f)+1-n/2}.
\]
Therefore,
\[
\sum_{x\leq y}(-1)^{f(x)}=\pm2.
\]
This is possible only if
$\sum_{x\leq y}f(x)\equiv1\pmod2$.
Hence
$
G[f](y)=1$.

\end{proof}

Theorem \ref{cl:platb} means that a maximally nonlinear function for a
fixed correlation immunity is not only plateaued, but also attains
equality in the Siegenthaler inequality.

\subsection{Algebraic Immunity}

The following problem arises in the cryptanalysis of certain
ciphers. Suppose that a Boolean function
$f:Q_2^n\rightarrow Q_2$, an affine transformation
$L:Q_2^n\rightarrow Q_2^n$, and an output sequence
\[
b_0=f(x),\qquad b_1=f(L(x)),\qquad
\dots,\qquad b_i=f(L^i(x)),\qquad\dots
\]
are known. The task is to determine the key $x\in Q_2^n$.
The complexity of the resulting system of equations depends on the
degree of $f$. Therefore, if there exists a Boolean function $g$ such
that
\[
\deg(g)<\deg(f)
\qquad\text{and}\qquad
gf=\mathbf{0},
\]
then one can search for possible values of $x$ as solutions of the
simpler system
\[
g(x)b_0=0,\qquad\dots\qquad
g(L^i(x))b_i=0,\qquad\dots
\]
The algebraic immunity of $f$ is introduced to measure the resistance
of $f$ to this type of cryptanalysis.

%\begin{definition}
The {\sl algebraic immunity} of a Boolean function $f$ is defined by
\[
AI(f)=
\min_g
\left\{
\deg(g):
gf=\mathbf{0}
\ \text{or}\
g(f\oplus\mathbf{1})=\mathbf{0},
\quad g\neq\mathbf{0}
\right\}.
\]
%\end{definition}

\begin{theorem}[Courtois and Meier \cite{CMeier}]\label{th:Court}
For every Boolean function $f:Q_2^n\rightarrow Q_2$, there exists a
Boolean function $g$ such that
\[
\deg(g)\leq\left\lceil\frac n2\right\rceil \qquad
{\it and}\qquad
\deg(fg)\leq\left\lfloor\frac n2\right\rfloor.
\]
\end{theorem}

\begin{proof}
Let $A_d$ be the set of all monomials of degree at most $d$,
\[
|A_d|=\sum_{i=0}^d{n\choose i}.
\]
Take
$d_1=\left\lfloor\frac n2\right\rfloor$,
$d_2=\left\lceil\frac n2\right\rceil$.
Consider the set of functions
$
B=\{fw:w\in A_{d_2}\}$.
If $B$ is a linearly dependent set then  there exist $w_1,\dots, w_k\in A_{d_2} $ such that $fw_1\oplus\cdots\oplus fw_k=0$. Then we have found  a function $g=w_1\oplus \cdots \oplus w_k$ satisfying the conditions of the theorem. Otherwise,  $B$ is linearly independent   and $|B|=|A_{d_2}|$.

We have
\[
|B|+|A_{d_1}|=
\sum_{i=0}^{d_1}{n\choose i}
+
\sum_{i=0}^{d_2}{n\choose i}
=
\sum_{i=0}^{n}{n\choose i}
+
{n\choose\lfloor n/2\rfloor}
>2^n.
\]
Since the vector space of Boolean functions over $GF(2)$ has
dimension $2^n$, the set of functions
$B\cup A_{d_1}$ is linearly dependent. Hence there exists a linear
combination of the form
$h\oplus fg=\mathbf{0}$,
where $h$ is a linear combination of monomials from $A_{d_1}$ and
$g$ is a linear combination of monomials from $A_{d_2}$.
Therefore,
$fg=h$,
and so
\[
\deg(fg)\leq d_1,
\qquad
\deg(g)\leq d_2.
\]
\end{proof}

\begin{corollary}
For every Boolean function $f$ on the $n$-dimensional Boolean cube,
\[
AI(f)\leq\left\lceil\frac n2\right\rceil.
\]
\end{corollary}

\begin{proof}
By Theorem \ref{th:Court}, there exists a Boolean function $g$ such
that
\[
\deg(g\oplus fg)\leq\left\lceil\frac n2\right\rceil.
\]
Set
$h=g\oplus fg$.
Then
$fh=fg\oplus fg=\mathbf{0}$.
Hence
\[
AI(f)\leq\deg(h)\leq\left\lceil\frac n2\right\rceil.
\]
\end{proof}

\begin{claim}[Canteaut \cite{Canteaut}]
For every Boolean function $f$ on the $n$-dimensional Boolean cube,
\[
nl(f)\geq
\sum_{i=0}^{AI(f)-2}{n\choose i}.
\]
\end{claim}

\begin{proof}
We first estimate the weight
$\wt(f)=|\supp(f)|$.
Consider the system of $|\supp(f)|$ Boolean equations
\[
\sum_{w\in A_d}a_w w(x)=0
\]
in the unknowns $a_w$, where $x$ ranges over all elements of
$\supp(f)$. If this system has a nonzero solution $(a_w)$, then
$
f\cdot\left(\sum_{w\in A_d}a_w w\right)=\mathbf{0}$.
Choose $d$ so that
\[
AI(f)=d+1.
\]
Then the system has only the zero solution. Consequently, the number
of equations is at least the number of variables, i.e.,
\[
\sum_{i=0}^{AI(f)-1}{n\choose i}\leq\wt(f).
\]

Now consider the function
$f\oplus\ell$,
where $\ell$ is an arbitrary affine function. If
$(f\oplus\ell)g=\mathbf{0}$
for some function $g$, then
$
fg=g\ell$.
Therefore,
\[
fg(1\oplus\ell)
=
g\ell(1\oplus\ell)
=
g\ell\oplus g\ell
=
\mathbf{0}.
\]
Moreover,
\[
\deg(g(1\oplus\ell))\leq\deg(g)+1.
\]
Consequently,
$AI(f\oplus\ell)\leq AI(f)+1$.
Applying the same argument to $f\oplus\ell$ and $\ell$, we obtain
\[
AI(f)=AI(f\oplus\ell\oplus\ell)
\leq AI(f\oplus\ell)+1.
\]

It follows that
\[
\sum_{i=0}^{AI(f)-2}{n\choose i}
\leq
\sum_{i=0}^{AI(f\oplus\ell)-1}{n\choose i}
\leq
\wt(f\oplus\ell).
\]
By the definition of the nonlinearity,
$nl(f)=\min_\ell\wt(f\oplus\ell)$,
and hence the desired inequality follows.
\end{proof}

Exact relations between algebraic immunity and nonlinearity
(i.e., the distance from a function to affine functions or, more
generally, to functions of bounded degree) were obtained by
Lobanov \cite{Loban}.

\subsection{Problems}

\begin{exercise}
Prove that every
function
$f:Q_p^n\rightarrow Q_p$,
where $p$ is prime, can be represented uniquely as a polynomial over
$GF(p)$ of degree at most $n(p-1)$, with the degree of each variable
in every monomial not exceeding $p-1$.
\end{exercise}

\begin{exercise}
Prove that, for prime $q$, the degree of the characteristic function of
a linear code $C\subset Q_q^n$ is equal to $n-\rank(C)$.
\end{exercise}

\begin{exercise}
Prove that for a non-resilient plateaued function $f$ the following
property holds:
$
|\widehat{(-1)^f}(z)|=|\widehat{(-1)^f}(y)|
$
whenever $\wt(y)=\wt(z)$. In particular,
$
|\widehat{(-1)^f}(z)|
=2^{\cor(f)+1-\frac n2}
$
for $\wt(z)=\cor(f)+1$.
\end{exercise}

\begin{exercise}
Prove that a Boolean function
$f:Q_2^n\rightarrow Q_2$ depends linearly  on the variable $x_i$
(over $GF(2)$) if and only if
$\widehat{(-1)^f}(z)=0$
for all $z\in Q_2^n$ such that $z_i=0$.
\end{exercise}

\begin{exercise}
Prove that for every $k\geq0$ there exists $n(k)$ such that, for
$n\geq n(k)$, the following holds: if
$\cor(f)\geq n-k$,
then the Boolean function
$f:Q_2^n\rightarrow Q_2$ depends linearly on some of its variables.
\end{exercise}

\begin{exercise}
Prove that for every $k\geq0$ there exists $n(k)$ such that, for
$n\geq n(k)$, the following holds: if
$\widehat{f}(z)=0$
whenever $\wt(z)\geq k$, then the Boolean function
$f:Q_2^n\rightarrow Q_2$ does not depend on some of its variables.
\end{exercise}

\begin{exercise}
Prove that the size of the support of a Boolean function of degree $k$
in $n$ variables is divisible by
$2^{\lceil n/k\rceil-1}$.
(the McEliece theorem)
\end{exercise}

\begin{exercise}
Prove that the degree of the characteristic function of a binary
$1$-perfect code of length $n$ is at most
$\frac{n+1}{2}$.
Estimate the minimum Hamming distance between two binary $1$-perfect
codes.
\end{exercise}

\begin{exercise}
Estimate the minimum Hamming distance between two $s$-plateaued functions.
\end{exercise}

\begin{exercise}
Prove that the Reed--Muller codes $\mathcal{R}(m,n)$ and
$\mathcal{R}(n-m-1,n)$ are dual and
\[
\mathcal{R}(k,n)\subseteq\mathcal{R}(m,n)^\perp
\]
whenever $k+m<n$. Determine when
$\mathcal{R}(m,n)$ is self-dual.
\end{exercise}

\begin{exercise}
Prove that $\mathcal{R}(1,n)$ is the extended Hadamard code, while
$\mathcal{R}(n-2,n)$ is the extended Hamming code.
\end{exercise}

\begin{exercise}
Prove that if two Boolean functions $f$ and $g$, both of degree at most
$r$, agree on the ball $B_r$ of radius $r$, then $f=g$.
\end{exercise}

\begin{exercise}
Prove that if two Boolean functions $f$ and $g$ are EA-equivalent  and $\deg(f)\geq 2$ then $\deg(f)=\deg(g)$.
\end{exercise}

\begin{exercise}
 Prove that   every $n$-variable Boolean function of algebraic degree $2$ is EA-equivalent to quadratic  form $\bigoplus\limits_{0<i<k}x_{2i-1}x_{2i}$ for some $k\leq n/2$ (Dickson's theorem).
\end{exercise}

\section{Perfect Colorings of Multigraphs and Hypergraphs}

\subsection{Perfect Colorings of Hypergraphs}\label{hyp}

A {\sl multigraph} is a generalization  of a graph that admits multiple
edges and loops.
A graph or multigraph is called {\sl directed} if it consists of  ordered pairs of  vertices,  called {\sl arcs},   instead of edges.
 Any nonnegative integer square matrix  can be  the adjacency matrix of a directed multigraph.
The definition of a perfect coloring of a directed multigraph
is similar to the definition of a perfect coloring of a simple graph. However, in the case of directed multigraph, we consider as the neighborhood of a vertex  only those
vertices that are the endpoints of arcs originating at that vertex. The neighborhood or sphere of a vertex is a multiset, since some vertices can be the endpoints of several parallel arcs.  It is easy to see that the algebraic criterion
for  perfect coloring (Proposition \ref{s:perfcol_criteria}) also holds for multigraphs. Note that the quotient matrix of any perfect
coloring can be viewed as the adjacency matrix of some
directed multigraph. Then a perfect coloring of a graph with
this quotient matrix can be considered as a covering. Furthermore, a merging of colors  can be viewed as a
perfect coloring of this multigraph.

Let $G$ be a bipartite graph with parts $V$ and $U$. Consider an undirected multigraph $\mathcal{M}_{12}(G)$ whose vertices
are the vertices of the first part of $G$ and the vertices
$v_i,v_j\in V$ are connected by $m_{ij}$ edges if $G$ has
$m_{ij}$ distinct vertices from $U$ adjacent to both $v_i$ and $v_j$
simultaneously. For convenience, we assume that the multigraph $\mathcal{M}_{12}(G)$ contains loops at its vertices
of multiplicity equal to the vertex degree in $G$. If $G$ is biregular,
then the loops in $\mathcal{M}_{12}(G)$ have the same multiplicity. Let
the graph $G$ be amply biregular, i.e., pairs of vertices in both
parties that are at distance $2$ have the same number of adjacent vertices. Then all edges in the multigraph $\mathcal{M}_{12}(G)$ have the same multiplicity. 
Below, we assume that the graph $G$
is biregular (not necessarily regular) and consider only those perfect colorings of $G$ for which all colors in different
parts are distinct.

\begin{claim}[\cite{PA20}]\label{bipartity}
Let $f:V\cup U \rightarrow \{1,\dots,k\}$ be a perfect coloring of the bipartite graph
$G$ with parts $V$ and $U$. Then $f|_V$ is a perfect coloring of
$\mathcal{M}_{12}(G)$. \end{claim}
\begin{proof}
The adjacency matrix of $G$ can be represented as
$M=\begin{pmatrix}
0 & Y \\
Y^{\mathsf T} & 0
\end{pmatrix}$. Let $S=\begin{pmatrix}
0 & S_1 \\
S_2 & 0
\end{pmatrix}$ be the quotient matrix  of a perfect coloring. Then, from the algebraic criterion of perfect coloring (Proposition \ref{s:perfcol_criteria}), we have the equality $MF=FS$,
which is equivalent to the pair of equalities $YF_2=F_1S_1$ and $Y^{\mathsf T}F_1=F_2S_2$, where $F=\begin{pmatrix}
F_1 & 0 \\
0 & F_2
\end{pmatrix}$, and the matrices $F_1$ and $F_2$ correspond to the
colorings of the first and second parts of  $G$. We have $YY^{\mathsf T}F_1=YF_2S_2=F_1S_1S_2$. It remains to note that 
$YY^{\mathsf T}$ is the adjacency matrix of the multigraph
$\mathcal{M}_{12}(G)$ and, by Proposition \ref{s:perfcol_criteria},
 $S_1S_2$ is the quotient matrix of the perfect coloring
$f|_V$.
\end{proof}

The converse, generally speaking, is false. A perfect coloring
of  $\mathcal{M}_{12}(G)$ need not be the restriction of any
perfect coloring of  $G$. In particular, such examples are easy to
find when $\mathcal{M}_{12}(G)$ is a complete graph.

A {\sl hypergraph} is another generalization of the concept of a graph,  in which we consider  unordered sets of more than two vertices  as edges. A hypergraph is called {\sl $k$-uniform} if  its edge consists of $k$ vertices. The {\sl incidence matrix} of a hypergraph $G$ is a $(0,1)$-matrix $Y=(y_{ij})$ of size $n\times m$, where $n$ is the number of vertices and $m$ is the number of edges in $G$. Moreover, 
\[
y_{ij}= \left\{
\begin{array}{l}
1 \mbox{, if}\ v_i\in e_j \in E(G);\\
0\mbox{, else}.
\end{array}
\right.
\]

For a hypergraph $G$, we define a bipartite graph $\mathcal{D}(G)$,
whose first part is the vertices of $G$, and whose  second part is the edges
of $G$. A vertex $v\in V(G)$ of the first part is adjacent to a vertex $e\in E(G)$ of the second
part if $v\in e$. It is easy to see that the adjacency matrix
of $\mathcal{D}(G)$ can be represented as
$\begin{pmatrix}
0 & Y \\
Y^{\mathsf T} & 0
\end{pmatrix}$, where $Y$ is the incidence matrix of $G$.

%\begin{definition}
A vertex coloring of a hypergraph $G$ is called {\sl perfect} if
the induced (see Section \ref{Vizing2}) coloring of the bipartite graph $\mathcal{D}(G)$ is perfect.
%\end{definition}

In other words, a hypergraph coloring is
perfect if, for
any fixed edge color, any two vertices
of the same color are incident with an equal number of edges of that fixed 
color. 
This definition was first proposed by
Godsil \cite{Godsil95}. 

Similar to the algebraic criterion of perfect coloring (Proposition \ref{s:perfcol_criteria}), we obtain
that the definition of a perfect hypergraph coloring is equivalent to
the equation
$$
\begin{pmatrix}
0 & Y \\
Y^{\mathsf T} & 0
\end{pmatrix}
\begin{pmatrix}
0 & F_1\\
F_2 & 0
\end{pmatrix}=
\begin{pmatrix}
0 & F_1\\
F_2 & 0
\end{pmatrix}\begin{pmatrix}
0 & S_1\\
S_2 & 0
\end{pmatrix},
$$
where
the matrices $F_1$ and $F_2$ correspond to
colorings of the vertices and hyperedges of $G$, respectively. This
equality is equivalent to the pair of equalities $YF_2=F_1S_1$ and $Y^{\mbox{\rm \tiny
T}}F_1=F_2S_2$.

Consider the adjacency matrix $M$ of the multigraph
$\mathcal{M}_{12}(\mathcal{D}(G))$. It is easy to see that the entries
$m_{ij}$ of $M$ are equal to the number of edges of the hypergraph $G$ that
contain both the $i$-th and $j$-th vertices. In particular, $m_{ii}$
is the degree of the $i$-th vertex in $G$. Clearly,
$M=YY^{\mathsf T}$, where $Y$ is the incidence matrix of 
$G$.

By Proposition \ref{bipartity}, we have the following corollary.
\begin{corollary}
A perfect coloring of a hypergraph $G$ is a perfect coloring
of $\mathcal{M}_{12}(\mathcal{D}(G))$.
\end{corollary}

The definition of  perfect colorings of graphs
is equivalent to the definition of perfect colorings of
$2$-uniform hypergraphs. Indeed, a perfect coloring
of a graph induces a coloring of its edges into pairs of colors corresponding to their
endpoints. This is a perfect coloring of $\mathcal{D}(G)$ by
the definition of a perfect coloring. The converse is also true: the
induced colors of edges incident to a vertex uniquely
determine the colors of  vertices.

\subsection{Line Graphs and Transversals}

%\begin{definition}
Let $G$ be a $k$-uniform hypergraph. 
$\mathcal{E}(G)$ is a multigraph
whose vertices are the edges of $G$; two vertices of
$\mathcal{E}(G)$ are connected by $l$ edges if the intersection
of the corresponding edges of $G$ has cardinality $l$.
%\end{definition}

From the definition, it is clear that the adjacency matrix of the line graph
$\mathcal{E}(G)$ is $M=Y^{\mathsf T}Y-kI$, where $Y$
is the incidence matrix of $G$. If the multigraph $\mathcal{E}(G)$ does not contain multiple edges, then the multigraph $\mathcal{E}(G)$ is called a {\sl line graph} of $G$.
Of course, the line graph of a simple graph is always a simple graph.

A coloring $f'$ of $\mathcal{E}(G)$ is called {\sl induced} by a coloring $f$ of $G$  if the color of a vertex $e$ of $\mathcal{E}(G)$ is the multiset (taking  multiplicity into account) of colors of the vertices  of $G$ contained in  $e$, i.e., $f'(e) = \{f(v_1),f(v_2),\dots,f(v_k)\}$, where 
$e=\{v_1,v_2,\dots,v_k\}$.

\begin{claim}\label{cl:edge_coloring}
The coloring of the graph $\mathcal{E}(G)$ induced by a perfect
coloring of  $G$ is perfect.
\end{claim}
\begin{proof}
Let $f$ be a perfect $t$-coloring of $G$.
Consider an edge $e$ of color composition $\{i_1,\dots,i_t\}$, i.e., an edge
containing $i_j$ vertices of color $j$. By the definition of a perfect
coloring, a vertex of color $j$ is contained in a known
number of hyperedges of each color composition. We collect the vectors of these
colors for all vertices of the hyperedge to obtain the color composition
of its adjacent hyperedges. Perhaps some edge $e'$ intersects
edge $e$ at $s>1$ vertices; in which case its color composition will be added
$s$ times. However,  in the graph $\mathcal{E}(G)$, there are exactly $s$ edges between the vertices
corresponding to $e$ and $e'$.
\end{proof}

%\begin{definition}
A {\sl transversal}  in a hypergraph is a set of vertices that intersects (cover) every  hyperedge exactly once.
%\end{definition}

In particular, in a bipartite graph, either of its
parts is a transversal. It is easy to see that in a $k$-uniform $r$-regular
hypergraph, a transversal is a color class of a perfect $2$-coloring. Moreover, the following holds:

\begin{claim}[\cite{PA20}]\label{bipartity1}
A function $f:V\rightarrow \{0,1\}$ is the characteristic function of a transversal
of a $k$-uniform $r$-regular
hypergraph
$G$ if and only if it is a perfect coloring of the multigraph
$\mathcal{M}_{12}(\mathcal{D}(G))$ with quotient matrix
$S=\begin{pmatrix}
r & (k-1)r \\
r & (k-1)r
\end{pmatrix}$.
\end{claim}
\begin{proof}
Since $G$ is $r$-regular, the vertices of the first part
of $\mathcal{D}(G)$ have degree $r$. Then each vertex
of  $\mathcal{M}_{12}(\mathcal{D}(G))$ is incident to $r$
loops. Vertices of the transversal are not adjacent to other vertices from the transversal in $\mathcal{M}_{12}(\mathcal{D}(G))$. Moreover,
each vertex not from the transversal is adjacent to $r$ vertices of the transversal in 
$\mathcal{M}_{12}(\mathcal{D}(G))$,
since each vertex contained in exactly $r$ edges, each of which
contains one vertex from the transversal. Therefore, $f$ is a
perfect coloring of the multigraph
$\mathcal{M}_{12}(\mathcal{D}(G))$. Moreover, this fact also follows from Proposition \ref{bipartity}, since the transversal
induces a $3$-coloring of the graph $\mathcal{D}(G)$.

If $f$ is a perfect coloring with quotient matrix
$S$, then any vertex of the first color is adjacent only to  itself  and is not
adjacent to other vertices of the first color in 
$\mathcal{M}_{12}(\mathcal{D}(G))$. This means that each edge of 
$G$ contains at most one vertex of the first color. Furthermore, any vertex
of the second color is  adjacent to $r$ vertices of the first color, which means that every
hyperedge incident to it contains a vertex of the first color.
\end{proof}

\begin{remark}
For $k$-uniform hypergraphs, it is natural to define the concept of a $k$-dimensional adjacency matrix. 
Methods for determining of ordinary eigenvalues and eigenvectors of multidimensional matrices are known,
 preserving some. A detailed exposition of the theory of perfect
hypergraph colorings can be found in \cite{Taranenko22}.
  In particular, analogous of 
 Theorems \ref{th:Vizing}, \ref{Angl}, Proposition \ref{s:perfcol_criteria} (algebraic criterion of perfect coloring), and Corollary \ref{p:perf_color_properties} (Lloyd's theorem)
 were proved to be valid for $k$-uniform hypergraphs as well.
\end{remark}

\subsection{Problems}

\begin{exercise}
Let $f$ be an eigenfunction of the bipartite graph $G$. Prove
that the restriction of  $f$ to the first part, is an eigenfunction
of $\mathcal{M}_{12}(G)$.
\end{exercise}

\begin{exercise}\label{exerDH0}
Prove that the eigenvalues of the multigraph
$\mathcal{M}_{12}(\mathcal{D}(G))$ are nonnegative.
\end{exercise}

\begin{exercise}
Let $F_1$ and $F_2$ be perfect colorings of the vertices and hyperedges 
of a hypergraph with quotient matrices $S_1$ and $S_2$, respectively. Prove the equality
$F_1^{\mathsf T}F_1S_1=S_2^{\mathsf T}F_2^{\mathsf T}F_2$.
\end{exercise}

\begin{exercise}\label{exerDH1}
Let $T$ be a transversal of a uniform regular hypergraph $G$.
Prove that $T$ is an independent set in the multigraph
$\mathcal{M}_{12}(\mathcal{D}(G))$. Prove that if
$\mathcal{D}(G)$ is amply regular, then  $T$ 
attains the Delsarte--Hoffman bound in the graph obtained from
$\mathcal{M}_{12}(\mathcal{D}(G))$ by removing loops and edge multiplicities.
\end{exercise}

\begin{exercise}
Let $G$  be a regular bipartite graph.
Prove that the restriction of to the first part  of $G$ an eigenfunction $f$ with eigenvalue
$\lambda$ 
is an eigenfunction of the multigraph $\mathcal{M}_{12}(G)$ with
eigenvalue $\lambda^2$. Conversely,
an eigenfunction of the multigraph $\mathcal{M}_{12}(G)$ with a positive eigenvalue
extends
to an eigenfunction of $G$.
\end{exercise}

\begin{exercise}
Let $G$ be a biregular bipartite graph and let $G'$ be the graph obtained from
$\mathcal{M}_{12}(\mathcal{D}(G))$ by removing loops and edge multiplicities. 
Prove that  any perfect coloring of $G$ is 
  a perfect coloring $G'$.
\end{exercise}

\begin{exercise}
Let $G_1$ be a regular graph.  Let  $\Gamma(G_1)$ be a multigraph  in which
edges connect vertices at distances $1$ and $2$ in $G_1$, with taking multiplicities into account. Prove that if $\varphi$ is a covering
from $G_2$ to $G_1$, then $\varphi$ is a covering
from $\Gamma(G_2)$ to $\Gamma(G_1)$.
\end{exercise}

\begin{exercise}
Formulate and prove an analogue of Theorem 
\label{thAF} for uniform regular hypergraphs.
\end{exercise}
\section{Latin Squares and $1$-Factorizations of Graphs}

\subsection{Latin Squares}

%\begin{definition}
A {\sl  Latin rectangle} is an  $m_1\times m_2$ table
filled with symbols (for example, elements of $Q_q$) such that each symbol occurs at most once in
every line (row or column). 
%\end{definition}

If the number of rows, columns, and distinct symbols are  all equal to $q$, then such
a Latin rectangle is called a {\sl Latin square of order} $q$.
The Latin coloring of the graph $Q^2_q$ introduced above can be
regarded as an equivalent definition of a Latin
square. The theory of Latin squares can be found in
\cite{DenesK}.

\begin{example}
Latin squares of orders $3$, $4$, and $5$.\\
 $\begin{array}{|c|c|c|}
\hline 0&1&2 \\
\hline
1&2&0 \\
 \hline
 2&0&1\\
 \hline
\end{array}$\qquad
 $\begin{array}{|c|c|c|c|}
\hline 0&1&2&3 \\
\hline
1&0&3&2 \\
 \hline
 2&3&1&0\\
\hline 3&2&0&1\\ \hline
\end{array}$\qquad
$\begin{array}{|c|c|c|c|c|}
\hline 0&1&2&3&4 \\
\hline
1&0&3&4&2 \\
 \hline
 2&3&4&0&1\\
 \hline
 3&4&1&2&0\\
\hline 4&2&0&1&3\\ \hline
\end{array}$

\end{example}

%\begin{definition}
 A binary operation
 $\circ:A^2\rightarrow A$ is called a {\sl 
    quasigroup} if   for every  $a,b\in A$ the equations
    $x\circ a=b$ and $a\circ x=b$ have   unique solutions.
    The cardinality
of the set $A$ is called the {\sl order} of the quasigroup.
%\end{definition}
More precisely, a quasigroup is a pair consisting of a set and a binary operation
defined on it, in this case
$(A, \circ)$.
Throughout this book we 
consider only quasigroups of finite order.  Since the nature of
 the underlying set is irrelevant,
 we will always assume that $A=Q_q$ for some $q\in
\mathbb{N}$. The following proposition is an immediate consequence of  the definitions.

\begin{claim}
A table of size $q\times q$, filled with $q$ distinct symbols,
is a Latin square if and only if it is the Cayley table of a quasigroup.
\end{claim}

\begin{claim}
If a quasigroup $(Q_q,\circ)$ is a semigroup (i.e., the operation
$\circ$ is associative), then it is a group.
\end{claim}
\begin{proof}
Fix $a\in Q_q$. By the definition of a quasigroup, there exists $e\in
Q_q$ such that $a\circ e= a$. Then, for every $b\in Q_q$, associativity yields  $a\circ (e\circ
b)= (a\circ e)\circ b= a\circ b$. Since left division is uniquely defined in a quasigroup, it follows that 
$e\circ b=b$ and, in particular, $e\circ a=a$.
Applying the same argument on the other side, we obtain
$b\circ e=b$ for all $b\in Q_q$. Thus, $e$ is an
identity element, and therefore $(Q_q,\circ)$ is a group.
\end{proof}

Let $\tau_i$, $i=0,1,2$, be permutations of the elements of the set
$Q_q$.
It is straightforward to verify that the operation
  $x*y=\tau_0(\tau_1(x) \circ
\tau_2 (y))$ is a quasigroup  operation  whenever $\circ$ is.  The permutations $\tau_1$ and $\tau_2$, acting on the arguments
of the quasigroup, correspond to permutations of rows and columns of  the associated
Latin square, and  $\tau_0$ renames its symbols.

Furthermore,  the inverse operation of a quasigroup with respect to any argument is also a quasigroup. More precisely, the operations $*$ and $*'$,
defined by the equality $x*z=y \Leftrightarrow x\circ y=z$ or
the equality $z*'y=x \Leftrightarrow x\circ y=z$, are also
quasigroup operations.

\subsection{Theorems of Vizing and Ryser}\label{7.2}

The adjacency matrix of a bipartite graph has a block
structure $ \left(\begin{array}{ll}
0 & M \\
M^{\mathsf T} & 0   \\
\end{array} \right)$, where $M$ is a rectangular matrix of size $n_1\times n_2$, where $n_1$ and $n_2$ are the sizes of two parts of the graph.
The matrix $M$ is called the {\sl reduced adjacency matrix} of the bipartite
graph.

Let $A=\{a_{ij}\}$ be an  $n_1\times n_2$ matrix.
Define the graph $g(A)$ as follows. Its vertex set
$V(g(A))=\{(i,j) \ |\ a_{ij}\neq 0\}$ consists of the nonzero entries
of the matrix.
 Two vertices $(i,j)$ and $(i',j')$
 are adjacent by an edge of the first
direction if $i=i'$, and by an edge of the second direction if $j=j'$.
If $A$ is the reduced adjacency matrix of some bipartite
graph $G$, then  $g(A)$ is  the line graph of
$G$. The following observation is immediate.

\begin{claim}\label{mper031} If two vertices of $g(A)$ belongs
 to the same connected component, then there exists a path between them whose edges alternate between the first and second directions.
\end{claim}

Next we prove a weak version of Vizing's theorem. 

\begin{theorem}[Vizing]\label{mper041}
Suppose that every row and every column of a rectangular $(0,1)$-matrix $A$ contains at most $q$ ones. Then the ones of the matrix can be properly colored with $q$ colors, that is, adjacent ones receive distinct colors.
\end{theorem}
\begin{proof}
Assume, to the contrary, that such a coloring does not exist.
Let $f:W\rightarrow
Q_q$
be a proper coloring of a subset $W\subset V(g(A))$ of vertices of the graph $g(A)$, where the subset $W$
has maximal possible cardinality.
 Then there exists a
vertex $x\in V(g(A))\setminus W$,  whose neighborhood contains
vertices of all $q$ colors.

Since $x$ has at most $q-1$ neighbors in each direction, there exist two distinct colors that are absent among the neighbors in one of the two directions.
Assume that color $0$ occurs only among the neighbors of $x$ in the first direction and color $1$ only among its neighbors in the second direction. Let $y$ and $z$ denote the neighbors of $x$ colored $0$ and $1$, respectively. Consider the subgraph of $g(A)$ induced by all vertices of colors $0$ and $1$. Since the coloring is proper, this subgraph is bipartite. If $y$ and $z$ belonged to the same connected component of this subgraph, then adjoining the vertex $x$ to a path connecting $y$ and $z$ would produce a cycle whose edge directions alternate by Proposition~\ref{mper031}.   Then the length of the  path between $y$ and $z$ is even.  It is impossible because $y$ and $z$ have different colors. Hence $y$ and $z$ belong to different connected components. Now interchange colors $0$ and $1$ throughout the connected component containing $z$. The resulting coloring remains proper. Consequently, color $1$ is now available for $x$, and we may extend the coloring to $W\cup\{x\}$, contradicting the maximality of $W$. \end{proof}

Let $A$ be the  reduced adjacency matrix of a bipartite
graph $G$.  By  Theorem \ref{mper041}, if the maximum  degree of the vertices of $G$  is at most $r$, then there exists a proper edge coloring of $G$
 with $r$ colors. This is a special case of the well-known theorem
of Vizing \cite{Viz64}, which states that
every graph of maximum degree $r$ admits a proper edge coloring with at most $r+1$ colors.

\begin{theorem}[Ryser]
A Latin rectangle $P$ of size $m_1\times m_2$ can be extended to a
Latin square of order $q$ if and only if
each of the $q$  symbols occurs at least $m_1+m_2-q$
times.
\end{theorem}
\begin{proof}
($\Rightarrow$)  Each 
letter occurs exactly $q-m_1$ times in the last $q-m_1$ rows, and exactly  $q-m_2$
times in the last $q-m_2$ columns  of the Latin square. Hence it  occurs at most $2q-m_1-m_2$ times outside $P$. Since every letter  occurs $q$ times in total, it occurs in $P$ at  least $m_1+m_2-q$
times.

($\Leftarrow$) First suppose that $m_2=q$. In this case, the
condition 
is automatically satisfied, since every symbol
occurs exactly once in each of the $m_1$ rows, and therefore exactly $m_1$ times.

 Construct
from the Latin rectangle $P$ a table $A=Inv(P)$ of size $q\times
q$ according to the rule
\[ A({i,j})= \left\{
\begin{array}{l}
k\mbox{, if there exists}\ k\in \{1,\dots,m_1\} \mbox{ such that}\ P(k,j)=i;\\
0\mbox{, otherwise}.
\end{array}
\right.
\]

Observe that if  $P$ can be extended to a Latin
square, then the operation $Inv$ corresponds to the inversion with respect to the first
argument of the quasigroup defining this Latin square.

Here it is convenient to regard the alphabet of the Latin rectangle as $\{1,2,\ldots,q\}$, while the symbol $0$ denotes an empty cell. By $\supp(A)$ we denote the $(0,1)$-matrix obtained from $A$ by replacing every nonzero entry with $1$.

For example,

 $P=\begin{array}{|c|c|c|c|}
\hline 4&1&2&3 \\
\hline
1&4&3&2 \\
 \hline
 0&0&0&0\\
\hline 0&0&0&0\\ \hline
\end{array}$\qquad
$A=\begin{array}{|c|c|c|c|}
\hline 2&1&0&0 \\
\hline
0&0&1&2 \\
 \hline
 0&0&2&1\\
\hline 1&2&0&0\\ \hline
\end{array}$
\qquad $\supp(A)=\begin{array}{|c|c|c|c|}
\hline 1&1&0&0 \\
\hline
0&0&1&1 \\
 \hline
 0&0&1&1\\
\hline 1&1&0&0\\ \hline
\end{array}$.

In every column and every row of  $\supp(A)$ there are exactly $m_1$ ones.
Let $B=J-\supp(A)$, where $J$ is the all-ones matrix. 
Then every row and every column of $B$ contains exactly $q-m_1$ ones. By Vizing's theorem, the ones of $B$ can be properly colored with the colors $m_1+1,\ldots,q$. It is immediate that the union of the matrices $A$ and $B$ forms a Latin square $L$.

Now compute $Inv(L)$. By the definition, the operation $Inv$ maps
 Latin squares to Latin squares, and satisfies
$Inv(Inv(L))=L$. Therefore, the first $m_1$ rows of the Latin square
$Inv(L)$ coincide with the original Latin rectangle $P$.

Now assume that $m_2<q$. Applying $Inv$ to $P$ yields a $q\times m_2$ rectangle $A=Inv(P)$. Each column of $\supp(A)$ contains exactly $m_1$ ones, while each row contains at least $m_1+m_2-q$ ones by assumption.
 Let $B=J-\supp(A)$, where $J$
is the $q\times m_2$  of all-ones matrix. Each
column of $B$ contains $q-m_1$ ones, and  each row contains at most  $m_2-(m_1+m_2-q)=q-m_1$ ones.
 Therefore, by Vizing's theorem, the ones of $B$ can be properly colored with $q-m_1$ colors.
 Together with $A$, this yields a Latin rectangle of size $q\times m_2$. By the first part of the proof, this rectangle can be extended to a Latin square $L$, and the Latin square $Inv(L)$ contains the original rectangle $P$.
 \end{proof}

\begin{corollary}\label{cor:latangle}
Every Latin rectangle of size $m_1\times q$ containing $q$
distinct symbols can be extended to a Latin square.
\end{corollary}

Note that a direct generalization of Corollary \ref{cor:latangle} to
Latin parallelepipeds is impossible, since there exists a Latin
parallelepiped of size $2\times 5\times 5$ which cannot be
extended even to a Latin parallelepiped of size $3\times
5\times 5$ over an alphabet of $5$ elements \cite{MKW}.

\subsection{Theorems of Hall, Kőnig and Egerváry}

%\begin{definition}
A {\sl matching} in a graph $G$ is a set of edges such that every
vertex of $G$ is incident to at most one edge of the set.
%\end{definition}

It follows from Vizing's theorem that if the maximum degree of a
bipartite graph is at most $r$, then its edge set can be partitioned
into $r$ matchings.
We now consider the problem of determining the maximum cardinality of
a matching.

%\begin{definition}
A {\sl vertex cover} of a graph $G$ is a set of vertices such that
every edge of $G$ is incident to at least one vertex of the set.
%\end{definition}

\begin{theorem}[Kőnig--Egerváry]
In every bipartite graph, the cardinality of a maximum matching equals
the cardinality of a minimum vertex cover.
\end{theorem}

We now formulate the Kőnig--Egerváry theorem in terms of matrices.

\begin{theorem}[Kőnig--Egerváry]
Let $A$ be a matrix.
Then the minimum total number of rows and columns covering all nonzero
entries of $A$ equals the cardinality of a maximum independent set in
the graph $g(A)$.
\end{theorem}

The equivalence of these formulations is immediate if we regard
$\supp(A)$ as the reduced adjacency matrix of a bipartite graph.
Rows and columns covering all nonzero entries correspond to a vertex
cover of the graph, while an independent set in $g(A)$ corresponds to
a matching.
We now prove the theorem.

\begin{proof}
It is clear that each row or column of a cover of $\supp(A)$ contains
at most one vertex of an independent set.
Therefore, the cardinality of an independent set cannot exceed the
number of rows and columns in the cover.

Let $C$ be a maximum independent set in $g(A)$.
Then every vertex of $V(g(A))\setminus C$ is adjacent to at least one
vertex of $C$.

Suppose that some vertex $c\in C$ has two neighbors
$a_1,a_2\in V(g(A))\setminus C$ in different directions, neither of
which is adjacent to any other vertex of $C$.
Then
$C'=(C\setminus\{c\})\cup\{a_1,a_2\}$
is an independent set of larger cardinality, contradicting the
maximality of $C$.

Hence every vertex $c\in C$ may have neighbors that are adjacent to no
other vertex of $C$ only in one of the two directions, namely, either
along its row or along its column.
Let $C_1$ denote the set of vertices of the first type and $C_2$ the
set of vertices of the second type.

Suppose that a vertex $c_1\in C_1$ is connected to a vertex
$c_2\in C_2$ by an alternating path whose vertices alternately belong
to $C$ and to $V(g(A))\setminus C$.
There exist vertices
$a_i\in V(g(A))\setminus C$, $i=1,2$, adjacent only to $c_i$.
Replacing the vertices of $C$ on the path by the remaining vertices of
the path, including $a_1$ and $a_2$, yields an independent set whose
cardinality exceeds that of $C$, a contradiction.

Thus the vertices of $C$ can be partitioned into three classes:
$C'_1$, consisting of vertices connected by alternating paths to
vertices of $C_1$, $C_1\subseteq C'_1$;
$C'_2$, consisting of vertices connected by alternating paths to
vertices of $C_2$, $C_2\subseteq C'_2$;
and the remaining vertices.

For each vertex of $C'_1$, select its row, and for each vertex of
$C'_2$, select its column.
Let $X\subseteq V(g(A))$ denote the set of vertices covered by these
lines.

Now consider a vertex
$a\in V(g(A))\setminus X$.
 There exist vertices from $C$ adjacent to $a$ both by
row and by column by definition of $C_i$. Both
these vertices belong to $C\setminus(C'_1\cup C'_2)$,
since if one of them belongs to $C'_i$, then the other also
belongs to $C'_i$  by the construction. For each vertex of $C\setminus(C'_1\cup
C'_2)$, take the row containing it. 
The selected rows and columns together cover all nonzero entries of
$\supp(A)$, and their total number is exactly $|C|$.
\end{proof}

Consider a bipartite graph with parts $U$ and $W$.
For every subset $U'\subseteq U$, let
$S(U')\subseteq W$
denote the set of vertices adjacent to at least one vertex of $U'$.

\begin{theorem}[Hall]
A bipartite graph with parts $U$ and $W$ has a matching of cardinality
$|U|$ if and only if
$|U'|\le |S(U')|$
for every subset $U'\subseteq U$.
\end{theorem}

We now formulate Hall's theorem in terms of matrices.

\begin{theorem}[Hall]
Let $A$ be an $n_1\times n_2$ matrix.
The graph $g(A)$ contains an independent set of cardinality $n_1$ if
and only if $\supp(A)$ contains no all-zero submatrix of size
$k_1\times k_2$ satisfying
$k_1+k_2>n_2$.
\end{theorem}

The equivalence of these formulations is immediate if we regard
$\supp(A)$ as the reduced adjacency matrix of a bipartite graph.
 By definition,  sizes of largest all-zero submatrices are $|U'|\times (n_2- |S(U')|)$, where  $U'\subset U$. Then the condition
$|U'|\leq |S(U')|$ is equivalent to $k_2= n_2-|S(U')|\leq n_2-k_1$ for every all-zero submatrix.  Now let us prove the theorem.

\begin{proof}
Suppose that $A$ contains an all-zero submatrix of size
$k_1\times k_2$ with
$k_1+k_2>n_2$.
Then among these $k_1$ rows one can choose at most
$n_2-k_2<k_1$
independent vertices of the graph $g(A)$.

Conversely, let $k_1\times k_2$ be a largest zero submatrix of $A$.
Then a minimum cover of $\supp(A)$ consists of
\[
n_1-k_1+n_2-k_2
\]
rows and columns.
Since
\[
k_1+k_2\le n_2,
\]
we have
\[
n_1-k_1+n_2-k_2\ge n_1.
\]
Therefore, by the Kőnig--Egerváry theorem, the graph $g(A)$ contains an
independent set of cardinality at least $n_1$.
\end{proof}

\subsection{Perfect Matchings}\label{7.3}

%\begin{definition}
A {\sl perfect matching ($1$-factor)} in a graph (or hypergraph) $G$
is a set of edges (hyperedges) such that every vertex of $G$ is
incident to exactly one edge of the set.
%\end{definition}

Clearly, a necessary condition for the existence of a perfect matching
in a graph is that the number of vertices be even. In the case of a
bipartite graph, it is additionally necessary that the two parts have
equal cardinality.

%\begin{definition}
A {\sl $1$-factorization} of a graph $G$ is a partition of its edge set
into perfect matchings.
%\end{definition}

A necessary condition for the existence of a $1$-factorization is that
the graph be regular.

We first consider an $r$-regular bipartite graph $G$.
Regularity implies that the two parts of $G$ have equal cardinality,
and its reduced adjacency matrix $M$ contains exactly $r$ ones in every
row and every column.
Diagonals consisting entirely of ones in $M$ correspond to perfect
matchings of $G$, while a proper coloring of $g(M)$ with $r$
colors is equivalent to a $1$-factorization of $G$.
Therefore, Vizing's theorem immediately yields the following result.

\begin{corollary}\label{konig}
Every regular bipartite graph admits a $1$-factorization.
\end{corollary}

We now turn to arbitrary (not necessarily bipartite) graphs.
A perfect matching in a $r$-regular graph $G$ corresponds to a symmetric
(with respect to the main diagonal)  diagonal consisting of ones in its adjacency
matrix $M$.
Consequently, a $1$-factorization of $G$ is equivalent to a symmetric
proper (and perfect at the same time) coloring of the ones of $M$, i.e., a perfect coloring of $g(M)$ with quotient matrix 
\begin{equation}\label{latmatrix}
2(J-I)=
\left(
\begin{array}{cccc}
0&2&\dots&2\\
2&0&2&\dots\\
\vdots&\vdots&\ddots&\vdots\\
2&\dots&2&0
\end{array}
\right).
\end{equation}
 A line graph $\mathcal{E}(G)$ differs from the graph  $g(M)$ in that we consider opposite points to be the same, i.e.,  $g(M)$ covers $\mathcal{E}(G)$. 
Since $G$ is $r$-regular, $\mathcal{E}(G)$ is
$2(r-1)$-regular.
By definitions, a perfect matching of $G$ is a cell of an equitable $2$-partition of  $\mathcal{E}(G)$ with quotient matrix 
\[
\left(
\begin{array}{cc}
0&2(r-1)\\
2&2(r-2)
\end{array}
\right).
\]
Moreover, every $1$-factorization of $G$ induces a perfect coloring of $\mathcal{E}(G)$ with quotient matrix (\ref{latmatrix}).

The concepts of a perfect matching and $1$-factorization  are naturally generalized to uniform hypergraphs. 
It is easy to see that a perfect matching and  $1$-factorization of a regular hypergraph  $G$ also correspond to perfect colorings of $\mathcal{E}(G)$. Consider a hypergraph $G'$ obtained from $G$ as follows: $V(G')=E(G)$, $E(G')=V(G)$, and $v\in e$ in $G'$ if and only if  $e\in v$ in $G$.
Then perfect matchings of $G$ one-to-one correspond to transversals of $G'$.

\subsection{Permanents and Doubly Stochastic Matrices}

%\begin{definition}
The {\sl permanent} $\per(A)$ of a $k\times m$ matrix
$A$, where $k\le m$, is defined by
\[
\per(A)=
\sum_{\sigma}
a_{1\sigma(1)}a_{2\sigma(2)}\cdots a_{k\sigma(k)},
\]
where the sum is taken over all injective mappings
\[
\sigma:\{1,\ldots,k\}\rightarrow\{1,\ldots,m\}.
\]
%\end{definition}

Observe that, for every nonnegative matrix $A$,
\[
\per(A)=0
\quad\Longleftrightarrow\quad
\per(\supp(A))=0.
\]

If $M$ is the reduced adjacency matrix of a bipartite graph $G$, then
$\per(M)$ equals the number of matchings of maximum possible
cardinality.
In particular, when $k=m$, $\per(M)$ equals the number of perfect
matchings.

\begin{theorem}[\cite{AlonF}]
Let $M$ be the adjacency matrix of a graph $G$.
Then the number of perfect matchings in $G$ does not exceed
$\sqrt{\per(M)}$.
\end{theorem}

\begin{proof}
Consider an arbitrary pair of perfect matchings in $G$.
Their union forms a $2$-regular spanning (containing all graph's vertices) subgraph (allowing multiple
edges), called a {\sl $2$-factor}.
If every cycle of the $2$-factor $f$ is even (counting a pair of parallel
edges as a $2$-cycle), then its edges can be partitioned into two
perfect matchings in exactly $2^{s(f)}$ ways, where $s(f)$ denotes the
number of cycles of length at least four.
If the $2$-factor contains an odd cycle, such a decomposition is
impossible.
Hence the square of the number of perfect matchings equals
$
\sum_f2^{s(f)},
$
where the sum ranges over all $2$-factors whose cycles are all even.

Now consider a diagonal consisting entirely of ones in the adjacency
matrix $M$. Each vertex of the graph is adjacent to two edges along the diagonal:
one $1$ in the row and another $1$ in the column corresponding to the
vertex. Therefore, a diagonal of $1$s corresponds to some
$2$-factor in the graph $G$. 
Such a diagonal determines a  $2$-factor of $G$.
Moreover, on every cycle of length at least three we may replace each
selected entry by its symmetric counterpart without changing the
corresponding $2$-factor.
Consequently,
$
\per(M)\ge
\sum_f2^{s(f)},
$
where the sum ranges over all $2$-factors.
The theorem follows immediately.
\end{proof}

For a bipartite graph $G$ with equally sized parts, the square root of
the permanent of its adjacency matrix equals the permanent of its
reduced adjacency matrix and therefore equals the number of perfect
matchings.

\begin{claim}
Let $G$ be a bipartite graph with reduced adjacency matrix $M$ of size
$k\times m$.
Then the total number of matchings in $G$ equals
$\per(A)$,
where
$A=(I_k\mid M)$
is the $k\times(k+m)$ matrix obtained by adjoining the identity matrix
to $M$.
\end{claim}

\begin{proof}
We establish a bijection between diagonals of $k$ ones in $A$ and
matchings in $G$.
Each such diagonal consists of a partial diagonal in $M$,
representing the edges of the matching, together with diagonal entries
of $I_k$, representing those vertices of the first part not incident
to the matching.
Since every edge of a bipartite graph is incident to a unique vertex of
the first part, this correspondence is bijective.
\end{proof}

A nonnegative matrix is called {\sl doubly stochastic} if the sum of
the entries in every row and every column equals $1$.
In particular, every permutation matrix
\[
P_{ij}(\tau)=
\left\{
\begin{array}{ll}
1,&j=\tau(i),\\
0,&\text{otherwise},
\end{array}
\right.
\]
is doubly stochastic.

\begin{theorem}[Kőnig--Frobenius]\label{mper043}
For every doubly stochastic matrix $A$, we have
\[
\per(A)>0.
\]
\end{theorem}

\begin{proof}
Suppose that a doubly stochastic matrix $A$ of size $m\times m$
contains exactly $m$ nonzero entries
($|\supp(A)|=|V(g(A))|=m$).
Then $A$ is a permutation matrix, and therefore
$\per(A)=1$.

Assume, to the contrary, that there exists a doubly stochastic matrix
of size $m\times m$ with zero permanent.
Among all such matrices, choose one, denoted by $B$, having the minimum
number of entries satisfying
$0<b_{ij}<1.
$
As observed above,
$|V(g(B))|>m$,
so such an entry necessarily exists.

For every entry satisfying $0<b_{ij}<1$, there exist neighboring
entries
\[
0<b_{i'j}<1,
\qquad
0<b_{ij'}<1,
\]
where $i'\ne i$ and $j'\ne j$.
Hence every vertex of the subgraph of $g(B)$ induced by the entries
$0<b_{ij}<1$ has degree at least two.
Consequently, this subgraph contains an even cycle $c$ whose edge
directions alternate.

Let
\[
\alpha=b_{i_0j_0}
=\min\{b_{ij}\}
\le
\beta=\min\{1-b_{ij}\},
\]
where both minimums are taken over the vertices of the cycle $c$.

Construct a matrix $C$ whose entries are zero except on the cycle $c$,
where they are equal to $-\alpha$ at vertices lying at even distance
(including distance zero) from $(i_0,j_0)$ and to $\alpha$ at vertices
lying at odd distance.
The sum of the entries of $C$ in every row and every column is zero.
Therefore,
$B+C$
is again doubly stochastic.

Moreover, $B+C$ contains fewer entries satisfying
$0<b_{ij}<1$
than $B$.
Since $\supp(B+C)\subseteq \supp(B)$ and
$\per(B)=0$,
we also have
$\per(B+C)=0$.
This contradicts the minimality of the number of fractional entries in
$B$.
The case
$\beta<\alpha$
is treated analogously.
\end{proof}

It follows from the Kőnig--Frobenius theorem (or from Vizing's theorem) that every $(0,1)$-matrix
having the same number of ones in each row and in each column has
positive permanent.  Van der Waerden conjectured that the minimum permanent among all $n \times n$ doubly stochastic matrices is  achieved by the matrix for which all entries are equal to $1/n$. Proofs of this conjecture were obtained  by Egorychev \cite{Egor} and Falikman \cite{Falik}.

\begin{theorem}[Birkhoff]
Every doubly stochastic matrix is a convex combination of permutation
matrices.
\end{theorem}

\begin{proof}
If a doubly stochastic matrix $A$ of size $m\times m$ contains exactly
$m$ nonzero entries, then $A$ is a permutation matrix.

Assume, to the contrary, that there exists a doubly stochastic matrix
that cannot be expressed as a convex combination of permutation
matrices.
Choose such a matrix $B$ with the minimum number of nonzero entries.

By the Kőnig--Frobenius theorem,
$\per(B)>0$.
Hence there exists a diagonal $d$ such that
$\alpha=\min\{b_{ij}\mid(i,j)\in d\}>0$.]
Let $P$ be the permutation matrix corresponding to the diagonal $d$.
Since $\frac{1}{1-\alpha}-\alpha=1$,  the matrix
\[
\frac{1}{1-\alpha}(B-\alpha P)
\]
is doubly stochastic and contains fewer nonzero entries than $B$.
By the minimality of $B$, there exist permutation matrices
$P_i$ and positive coefficients $\beta_i$ satisfying
\[
\sum_i\beta_i=1,
\qquad
\frac{1}{1-\alpha}(B-\alpha P)
=
\sum_i\beta_iP_i.
\]
Therefore,
$B=\alpha P+(1-\alpha)\sum_i\beta_iP_i$,
which expresses $B$ as a convex combination of permutation matrices, a
contradiction.
\end{proof}

A {\sl convex polytope} is a bounded subset of $\mathbb{R}^n$
defined by a system of linear inequalities and equalities,
\[
M_1x\le b_1,\qquad M_2x=b_2,
\]
where $b_1$ and $b_2$ are real vectors and $M_i$ are
$m_i\times n$ real matrices.
For example, the set of doubly stochastic matrices forms a convex
polytope.

A {\sl vertex} (or {\sl extreme point}) of a convex polytope is a point
that cannot be expressed as a convex combination of other points of the
polytope.
Every convex polytope is the convex hull of its vertices, and,
conversely, the convex hull of any finite set of points is a convex
polytope.

By Birkhoff's theorem, the vertices of the polytope of doubly
stochastic matrices are precisely the permutation matrices.

Let $F$ be an $n\times k$ $(0,1)$-matrix corresponding to a
$k$-coloring of a graph.
By Proposition~\ref{claimprojector},
$P=F(F^{\mathsf T}F)^{-1}F^{\mathsf T}$
is the orthogonal projector onto the subspace spanned by the indicator
functions of the color classes.
By~(\ref{eq:projector}), the matrix $P$ is doubly stochastic.

By Proposition~\ref{claimprojector1}, the coloring defined by $F$ is a
perfect coloring of a graph $G$ if and only if the corresponding
projector $P$ commutes with the adjacency matrix of $G$.
Since the set of matrices commuting with a fixed matrix forms a linear
subspace, it follows that the set of doubly stochastic matrices
commuting with the adjacency matrix of $G$ is a convex polytope.
We denote this polytope by $DSP(G)$.

\begin{claim}[Godsil \cite{Godsil97}]
If every vertex of $DSP(G)$ is a permutation matrix, then every perfect
coloring of $G$ is an orbit coloring.
\end{claim}

\begin{proof}
Let $F$ be the $n\times k$ matrix corresponding to a perfect coloring
of $G$.
By Proposition~\ref{claimprojector1}, the corresponding projector
$P$ belongs to $DSP(G)$.

Hence
\[
P=\sum_i\alpha_iP_i,\
\mbox{where}\qquad
\alpha_i\ge0,\quad
\sum_i\alpha_i=1,
\]
and the matrices $P_i$ are permutation matrices.
Since $P_i\in DSP(G)$, each $P_i$ commutes with the adjacency matrix of
$G$ and therefore represents an automorphism of $G$ (see Section \ref{sec_cov}).

By~(\ref{eq:projector}),
$P(x,y)>0$
if and only if the vertices $x$ and $y$ have the same color.
Since
$P(x,y)=\sum_i\alpha_iP_i(x,y)$,
there exists an index $i$ such that
$P_i(x,y)=1$.
Thus, for every pair of vertices of the same color, there exists an
automorphism of $G$ mapping $x$ to $y$.
Hence the coloring defined by $F$ is an orbit coloring.
\end{proof}

A trivial $1$-partition is equitable in every regular graph.  Consequently, we obtain the  following statement. 
\begin{corollary}
If $G$ is a
regular graph,  and  every vertex of $DSP(G)$ is a permutation matrix, then  $G$ is vertex-transitive.
\end{corollary}

\subsection{Symmetric Latin Squares}

%Let $Q_q$ denote the complete graph on $q$ vertices.

\begin{claim}
If $q$ is even, then a $1$-factorization of the complete graph $Q_q$
is equivalent to a symmetric Latin square of order $q$ with zeros on
the main diagonal.
\end{claim}

\begin{proof}
Consider a $1$-factorization of the graph $Q_q$.
Each edge $e$ corresponds to a pair of symmetric ones in the adjacency
matrix of $Q_q$.
Replace this pair by the symbol $i$ whenever $e$ belongs to the
$i$-th $1$-factor.
The resulting matrix is symmetric of order $q$, and each symbol
$0,1,\ldots,q-1$ appears exactly once in every row and every column,
since every vertex is incident to exactly one edge of each
$1$-factor.
Conversely, reversing this construction yields a $1$-factorization of
the complete graph from any symmetric Latin square with zeros on the
main diagonal.
\end{proof}

A $1$-factorization of the complete graph $Q_q$ may be viewed as a
schedule for a round-robin tournament.
The participants correspond to the vertices of the graph, and if the
edge $\{u,v\}$ belongs to the $i$-th matching, then participants $u$
and $v$ play each other in round $i$.

%\begin{definition}
A quasigroup $f:Q_q^2\rightarrow Q_q$ and its Cayley table (Latin
square) are called {\sl idempotent} if
$f(a,a)=a$
for every $a\in Q_q$.

A quasigroup
$f:Q_{2q}^2\rightarrow Q_{2q}$
and its Cayley table are called {\sl semi-idempotent} if
\[
f(a,a)=a,\quad a\in Q_q,\quad
\mbox{and}\quad
f(a,a)=a-q,\quad
a\in Q_{2q}\setminus Q_q.
\]
%\end{definition}

We now prove the existence of a $1$-factorization of every complete
graph of even order, as well as the existence of a symmetric
idempotent Latin square of every odd order.

\begin{theorem}[\cite{Rosa}]\label{th:LSsym}
$1.$ For every even $q\ge2$, there exists a $1$-factorization of the graph
$Q_q$ and a semi-idempotent symmetric Latin square of order $q$. \\
$2.$
For every odd $q\ge1$, there exists an idempotent symmetric Latin
square of order $q$.
\end{theorem}

\begin{proof}
Let $A_{2q}[Q]$ denote a symmetric Latin square of order $2q$ over the
alphabet $Q$ with constant main diagonal,
let $B_{2q}[Q]$ denote a symmetric semi-idempotent Latin square of order $2q$ over
$Q$, and let
$C_{2q+1}[Q]$ denote a symmetric idempotent Latin square of order $2q+1$ over $Q$.

We construct these squares by induction.

The initial squares are
\[
A_2[Q_2]=B_2[Q_2]=
\left(
\begin{array}{cc}
0&1\\
1&0
\end{array}
\right),
\]
and
\[
C_3[Q_3]=
\left(
\begin{array}{ccc}
0&2&1\\
2&1&0\\
1&0&2
\end{array}
\right).
\]

Assume that
$A_{2q}[Q_{2q}]$,
$B_{2q}[Q_{2q}]$,
and
$C_{2q+1}[Q_{2q+1}]$
have already been constructed.

Then
\[
A_{4q}[Q_{4q}]
=
\left(
\begin{array}{cc}
A_{2q}[Q_{2q}]&
A_{2q}[Q_{4q}\setminus Q_{2q}]\\
A_{2q}[Q_{4q}\setminus Q_{2q}]&
A_{2q}[Q_{2q}]
\end{array}
\right),
\]
and
\[
B_{4q+2}[Q_{4q+2}]
=
\left(
\begin{array}{cc}
C_{2q+1}[Q_{2q+1}]&
C_{2q+1}[Q_{4q+2}\setminus Q_{2q+1}]\\
C_{2q+1}[Q_{4q+2}\setminus Q_{2q+1}]&
C_{2q+1}[Q_{2q+1}]
\end{array}
\right).
\]

Next, let $E$ be an arbitrary Latin square of order $2q+1$ over the
alphabet $Q_{2q+1}$ with zeros on the main diagonal.
Starting from the Latin square
\[
\left(
\begin{array}{cc}
C_{2q+1}[Q_{4q+2}\setminus Q_{2q+1}]&E\\
E^{\mathsf T}&
C_{2q+1}[Q_{4q+2}\setminus Q_{2q+1}]
\end{array}
\right),
\]
interchange the main diagonals of the two copies of
$C_{2q+1}[Q_{4q+2}\setminus Q_{2q+1}]$
with those of the blocks
$E$
and
$E^{\mathsf T}$.
The resulting Latin square is
$A_{4q+2}[Q_{4q+2}]$.

Similarly, starting from
\[
\left(
\begin{array}{cc}
B_{2q}[Q_{2q}]&
B_{2q}[Q_{4q}\setminus Q_{2q}]\\
B_{2q}[Q_{4q}\setminus Q_{2q}]&
B_{2q}[Q_{2q}]
\end{array}
\right),
\]
interchange the second halves of the main diagonals of the two diagonal
blocks.
The resulting Latin square is
$B_{4q}[Q_{4q}]$.

It remains to construct
$C_{4q-1}[Q_{4q-1}]$
and
$C_{4q+1}[Q_{4q+1}]$.

Renaming the symbols if necessary, we may assume that the symbols in
the first row and first column of
$A_{4q}[Q_{4q}]$
appear in increasing order.
Replace the diagonal entry in position $(i,i)$ by the symbol $i$ for
every
$i\in Q_{4q}$.
After deleting the first row and the first column, we obtain the Latin
square
$C_{4q-1}[Q_{4q}\setminus\{0\}]$.
Finally, renaming the symbols yields
$C_{4q-1}[Q_{4q-1}]$.

The construction of
$C_{4q+1}[Q_{4q+1}]$
from
$A_{4q+2}[Q_{4q+2}]$
is completely analogous.
This completes the induction.
\end{proof}

\subsection{Problems}

\begin{exercise}
Prove that every bipartite graph contains a matching whose edges are
incident to all vertices of maximum degree.
\end{exercise}

\begin{exercise}
Prove that a doubly stochastic matrix of size $n\times n$ cannot
contain a zero submatrix of size $m_1\times m_2$, where
$m_1+m_2>n$.
\end{exercise}

\begin{exercise}
Prove that a nonnegative matrix of size $n\times n$ has zero permanent
if and only if it contains a zero submatrix of size
$m_1\times m_2$, where $m_1+m_2>n$.
\end{exercise}

\begin{exercise}\cite{Lovasz}
Let $A$ be a $(0,1)$-matrix of size $m\times n$ containing $ka$ ones
in each row and $kb$ ones in each column, where $k>0$ is an integer.
Prove that one can select $a$ ones from each row of $A$ so that exactly
$b$ ones are selected from each column.
\end{exercise}

\begin{exercise}\cite{Lovasz}
Prove that every bipartite graph can be extended by adding vertices and
edges to one of its parts so that the maximum degree does not increase
and the resulting bipartite graph is regular.
\end{exercise}

\begin{exercise}
Prove that every regular graph has a covering graph admitting a
$1$-factorization.
\end{exercise}

\begin{exercise}
Prove that every partial Latin square of order $q$, that is, a
$q\times q$ array in which  cells are filled with at most $q$
distinct symbols and no symbol occurs more than once in any row or
column, can be completed to a Latin square of order $2q$.
\end{exercise}

\begin{exercise}
Prove that the maximum of the permanents of doubly stochastic matrices is equal to $1$ and it is attained only on permutation matrices.
\end{exercise}

\begin{exercise}
Prove (by using van der Waerden's conjecture on permanents)  that the asymptotics  of the number of Latin squares of order $n$ is equal to $n^{n^2(1-o(1))}$ as $n\rightarrow \infty$. 
\end{exercise}

\section{MDS Codes }

\subsection{Basic Properties of MDS Codes}\label{MDS}

Let $C\subseteq Q_q^n$. Recall that the minimum distance (code distance) of the code $C$ is
defined by
\[
d_C=\min_{\substack{x,y\in C\\x\neq y}} d(x,y).
\]

A set $C\subseteq Q_q^n$ is called a {\sl maximum distance separable
(MDS) code} with distance $d$ if
$|C\cap\Gamma|=1$
for every $(d-1)$-dimensional face $\Gamma\subset Q_q^n$.

Consider the hypergraph whose vertices are the vertices of the
hypercube and whose hyperedges are its $(d-1)$-dimensional faces.
It follows directly from the definitions that an MDS code with
distance $d$ is precisely a transversal of this hypergraph.

As mentioned above (see Example \ref{MDSexam}), when $d=2$, the
characteristic function of an MDS code in $Q_q^n$ is a perfect
$2$-coloring with quotient matrix
\[
\left(
\begin{array}{cc}
n & n(q-2)\\
0 & n(q-1)
\end{array}
\right).
\]
Its nontrivial eigenvalue $\lambda=-n$ is the smallest eigenvalue of
$Q_q^n$.

The cardinality of an MDS code is determined directly from the
definition.

\begin{claim}[Singleton bound]\label{proSingl}\quad

{\rm (a)} For every $C\subseteq Q_q^n$,
$ |C|\le q^{\,n-d_C+1}$.

{\rm (b)} A code $C$ with minimum distance $d_C$ is an {\rm MDS} code
if and only if
\[
|C|=q^{\,n-d_C+1}.
\]
\end{claim}

\begin{proof}
Part (a) follows from the fact that every face of dimension
$d_C-1$ contains at most one codeword, while
$Q_q^n$ is partitioned into $q^{\,n-d_C+1}$ disjoint faces of
dimension $d_C-1$.

(b) If $C$ is an MDS code, then each of these
$q^{\,n-d_C+1}$ faces contains exactly one codeword.
Conversely, if $|C|=q^{\,n-d_C+1}$, then every face of
dimension $d_C-1$, in every direction, must contain a codeword.
\end{proof}

%\begin{definition}
Let $\Gamma$ be an $m$-dimensional face of $Q_q^n$.
The {\sl retract} of a set $C\subseteq Q_q^n$ onto $\Gamma$
is the subset $C[\Gamma]\subseteq Q_q^m$ consisting of all
tuples in $\Gamma\cap C$ after deleting the coordinates fixed in
$\Gamma$.
%\end{definition}
Thus, the characteristic function ${\mathbf 1}_{C[\Gamma]}$ is simply the
restriction of ${\mathbf 1}_C$ to the face $\Gamma$.

%\begin{definition}
Let $\Gamma$ be an $m$-dimensional face of $Q_q^n$.
The {\sl projection} of a set $C\subseteq Q_q^n$ onto $\Gamma$
is the subset of $Q_q^m$ obtained by deleting the coordinates fixed
in $\Gamma$ from every tuple of $C$.
%\end{definition}

Thus, the projection onto a face is the union of the retracts over all
parallel faces of the same direction. A retract of a code is called a
{\sl punctured code}, while a projection is called a {\sl shortened
code}. The class of MDS codes is closed under both retracts and
projections.

\begin{claim}\label{mds99}
A subset of $Q_q^n$ is an {\rm MDS} code with distance $d$
($2\le d\le n-1$) if and only if

{\rm (a)} all of its retracts of dimension at least $d$ are
{\rm MDS} codes with distance $d$; or

{\rm (b)} all of its projections onto faces of a fixed dimension
$m\ge n-d+1$ are {\rm MDS} codes with distance $d-n+m$.
\end{claim}

\begin{proof}
(a) The definition of an MDS code is equivalent to the statement that
every face of dimension at least $d$ has the property that each of its
$(d-1)$-dimensional subfaces contains exactly one codeword.

(b) This follows from Proposition~\ref{proSingl}. Under projection
onto a face of dimension $m$, all vertices of each face of dimension
$n-m$ are mapped to a single point. If $n-m\le d-1$, then the
projection has the same cardinality as the original code, while its
minimum distance can decrease by at most $n-m$. Proposition
\ref{proSingl}(a) therefore implies that the minimum distance
decreases by exactly $n-m$, and Proposition
\ref{proSingl}(b) shows that the projection is an MDS code.

Conversely, suppose that the projection onto some
$m$-dimensional face is an MDS code with distance $d-n+m$.
Then every $(d-1)$-dimensional face containing all directions along
which the projection is taken contains a codeword. Since the
projection direction is arbitrary, every $(d-1)$-dimensional face of
the hypercube contains a codeword. Finally, if the original code
contained two codewords at distance $d'<d$, then some projection would
have minimum distance smaller than $d-n+m$, contradicting the
assumption that the projection is an MDS code.
\end{proof}

\begin{claim}\label{proAO2}
Let $C\subseteq Q_q^n$ be an {\rm MDS} code with minimum distance
$d_C\ge3$. Then
\[
|C|\le q^{\,q-1}
\qquad\text{and}\qquad
n\le q+d_C-2.
\]
\end{claim}

\begin{proof}
By Proposition~\ref{mds99}, there exists a projection $C'$ of $C$
onto a face of dimension
$m=n-d_C+3$,
whose minimum distance is $3$. We have
\[
|C'|=|C|=q^{\,n-d_C+1}=q^{\,m-2}.
\]

By the Hamming bound (see Section \ref{sec:it}),
$|C'|\bigl(m(q-1)+1\bigr)\le q^m$.
Hence,
$q^2\ge m(q-1)+1$,
which implies $m\le q+1$. Therefore,
$
n\le q+d_C-2$.
\end{proof}

\subsection{$n$-ary Quasigroups}

An {\sl $n$-ary quasigroup} of order $q$ is a pair consisting of a set
of elements (which, without loss of generality, we take to be $Q_q$)
and an $n$-ary operation
$f:Q_q^n\rightarrow Q_q $
such that
$f(x)\neq f(y)$
whenever the Hamming  distance between $x$ and $y$ is equal to $1$.
In other words, any choice of $n$ coordinates in the equation $f(x_1,\dots,x_n)=x_0$
uniquely determines the remaining coordinate. 
The function $f$ itself is also called an
{\sl $n$-ary quasigroup}.

The notion of an $n$-ary quasigroup extends the notion of a binary
($2$-ary) quasigroup introduced above.

\begin{example}
Let $\circ$ be a group operation on the set $Q_q$. Then the function
\[
f(x_1,\dots,x_n)=x_1\circ\cdots\circ x_n
\]
is an $n$-ary quasigroup, called an {\sl iterated group}.
\end{example}

It follows immediately from the definition that an $n$-ary quasigroup
defines a perfect coloring of $Q^n_q$ with quotient matrix
\[
n(J_q-I_q)=
\left(\begin{array}{ccccc}
0 & n & n & \dots & n\\
n & 0 & n & \dots & n\\
\vdots & \vdots & \vdots & \ddots & \vdots\\
n & n & \dots & n & 0
\end{array}\right),
\]
whose unique nontrivial eigenvalue is $\lambda=-n$.

Each color class of the coloring defined by an $n$-ary quasigroup is
an MDS code with minimum distance $2$. Moreover, the graph of an
$n$-ary quasigroup,
\[
\{\bar{x}\in Q_q^{\,n+1}\mid
x_{n+1}=f(x_1,\dots,x_n)\},
\]
is itself an MDS code in $Q_q^{\,n+1}$.

The table of values (Cayley table) of an $n$-ary quasigroup is called
a {\sl Latin $n$-cube}; for $n=2$, it is simply a {\sl Latin square}.

\begin{remark}
It is natural to generalize the notion of a doubly stochastic matrix to
multidimensional arrays. A function
$f:Q_q^n\rightarrow[0,1]$
is called a {\sl polystochastic matrix} if the sum of its values over
every line (that is, every $1$-dimensional face) equals $1$. The set
of all polystochastic matrices is a multidimensional analogue of the
Birkhoff polytope.

It is easy to see that the indicator function of an MDS code with
minimum distance $2$ is a polystochastic matrix. Such an MDS code is
sometimes called a {\sl multidimensional permutation}. Every
multidimensional permutation is a vertex of the generalized Birkhoff
polytope. However, an analogue of Birkhoff's theorem does not hold:
there exist vertices of the generalized Birkhoff polytope that are not
multidimensional permutations (see \cite{LL14,PotTar}). The
description of the vertices of the generalized Birkhoff polytope and
the estimation of their number remain open problems.\label{Berprob}
\end{remark}

\subsection{Orthogonal Systems of Quasigroups}

Any {\rm MDS} code with minimum distance greater than $2$ can be
represented as the set of solutions of a system of equations involving
$m$-ary quasigroups \cite{EMull}. Indeed, consider an {\rm MDS} code of
length $n$ and minimum distance $d$. By the definition of an {\rm MDS}
code, every face of the hypercube of dimension
$n-m=d-1$
contains exactly one codeword. In particular, the first $m$
coordinates uniquely determine $(n-m)$-dimensional face, and, therefore,  the remaining $n-m$ coordinates.
Hence, the {\rm MDS} code is the solution set of the system
\begin{equation}\label{eMDS11}
x_{m+1}=f_1(x_1,\dots,x_m),\ \dots,\,
x_n=f_{n-m}(x_1,\dots,x_m).
\end{equation}
It follows immediately from the definition of an {\rm MDS} code that
the functions $f_1,\dots,f_{n-m}$ are $m$-ary quasigroups. Indeed, since all
coordinates of an MDS code play symmetric roles, any choice of $m$ coordinates
uniquely determines the remaining coordinates by means of a system of
the form (\ref{eMDS11}), although in general with different
quasigroups.

If the solution set of a system  (\ref{eMDS11}) is an MDS
code, then the corresponding collection of $m$-ary quasigroups
$f_1,\dots,f_{n-m}$ is called {\sl orthogonal}. Proposition
\ref{mds99}(a) implies that every subset of an orthogonal system is
again orthogonal.

A pair of Cayley tables (Latin squares) corresponding to binary
($2$-ary) quasigroups $f$ and $g$ is called {\sl orthogonal} if all
pairs
\[
\bigl(f(x_1,x_2),\,g(x_1,x_2)\bigr),
\qquad (x_1,x_2)\in Q_q^2,
\]
are distinct. Thus, a system of mutually orthogonal Latin
squares (MOLS)  is a special case of an orthogonal system of
quasigroups. According to the arity $m$, one speaks of orthogonal
Latin squares, orthogonal Latin cubes, or orthogonal Latin
hypercubes.

\begin{example}
Mutually orthogonal Latin squares of order $4$.

\smallskip

\begin{tabular}{|c|c|c|c|}
 \hline
  0&{\color[rgb]{0.1,1,0.1} \bf 2}&3&1\\
  \hline
  3&1&0&{\color[rgb]{0.1,1,0.1} \bf 2}\\
  \hline
  1&3&{\color[rgb]{0.1,1,0.1} \bf 2}&0\\
  \hline
  {\color[rgb]{0.1,1,0.1} \bf 2}&0&1&3\\
 \hline
 \end{tabular}
\begin{tabular}{|c|c|c|c|}
   \hline
  0&{\color[rgb]{0.1,1,0.1} \bf 2}&3&1\\
   \hline
  2&0&1&{\color[rgb]{0.1,1,0.1} \bf 3}\\
   \hline
  3&1&{\color[rgb]{0.1,1,0.1} \bf 0}&2\\
   \hline
  {\color[rgb]{0.1,1,0.1} \bf 1}&3&2&0\\
 \hline
 \end{tabular}
 \begin{tabular}{|c|c|c|c|}
 \hline
  0&{\color[rgb]{0.1,1,0.1} \bf 2}&3&1\\
   \hline
  1&3&2&{\color[rgb]{0.1,1,0.1} \bf 0}\\
   \hline
  2&0&{\color[rgb]{0.1,1,0.1} \bf 1}&3\\
   \hline
  {\color[rgb]{0.1,1,0.1} \bf 3}&1&0&2\\
 \hline
 \end{tabular}

\end{example}

\begin{example}
Orthogonal Latin cubes of order $4$.

\begin{minipage}{1\textwidth}
%\vskip15mm
\includegraphics[width=1.0\textwidth]{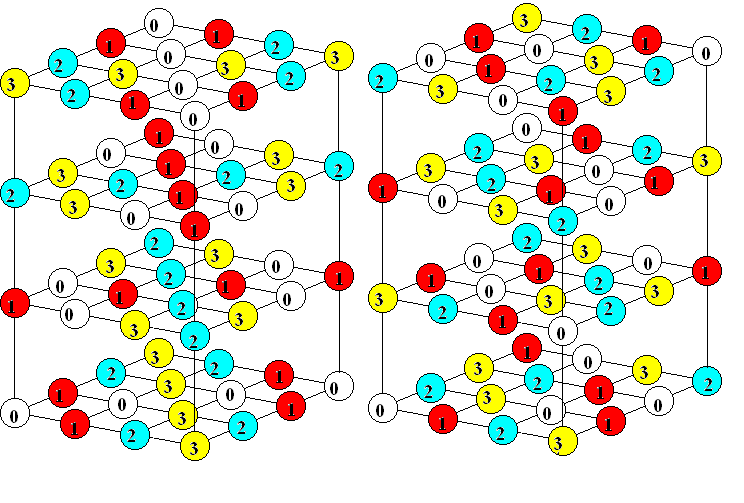}
% \center{\includegraphics[width=0.07\textwidth]{E3.png}}
\end{minipage}

\end{example}

\begin{claim}\label{proAO5}
If $f_1,\dots,f_m$ is a system of mutually orthogonal Latin squares
(binary quasigroups) of order $q$, then
$m\le q-1$.
\end{claim}

\begin{proof}
It is clear that renaming entries of a Latin square preserves the property of  orthogonality.
Then, without loss of generality, we may assume that the first rows
($x_1=0$) of all Latin squares are identical and that
$f_i(0,0)=0$ for every $i\in\{1,\dots,m\}$.
Consequently, every pair of Latin squares already contains all pairs
of the form $(a,a)$ in the first row. Therefore, the values
$f_i(1,0)$ must be nonzero and pairwise distinct. Since there are only
$q-1$ nonzero elements of $Q_q$, it follows that $m\le q-1$.
\end{proof}

The preceding proof immediately implies the following result.

\begin{corollary}\label{AP01}
Let $C\subset Q_q^{\,q+1}$ be an {\rm MDS} code with minimum distance
$q$. Then $C$ contains no pair of codewords at distance $q+1$, i.e., no  pair of codewords differs in every coordinate.
\end{corollary}

\begin{proof}
Such an MDS code is equivalent to a system of mutually orthogonal
Latin squares (see (\ref{eMDS11})). Since the isometry group of the Hamming graph acts
transitively on the vertices, we may assume that $\bar{0}\in C$ and
that the first rows of the corresponding Latin squares are identical.
The proof of the previous proposition shows that the vertices
$(i,\dots,i)$ do not belong to $C$, because
$(0,i,\dots,i)\in C$ for every $i\in Q_q$. Hence no codeword of $C$ is
at distance $q+1$ from $\bar{0}$.
\end{proof}

\begin{claim}\label{proAO9}
Let $C\subset Q_q^n$ be an {\rm MDS} code with minimum distance
$d_C$, where $n>d_C\ge3$. Then
$d_C\le q$ and $n\le2q-2$.
\end{claim}

\begin{proof}
Let $d_C=m+1$, where $m\ge2$. Then the code $C$ corresponds to an
orthogonal system of $m$ $(n-m)$-ary quasigroups (see (\ref{eMDS11})).

By Proposition~\ref{mds99}, every retract of $C$ of dimension $m+2$
also determines an orthogonal system, namely a system of $m$
mutually orthogonal Latin squares. Proposition~\ref{proAO5}
therefore gives
$m\le q-1$,
which is equivalent to $d_C\le q$.
Finally, Proposition~\ref{proAO2} yields
$n\le q+d_C-2\le2q-2$.
\end{proof}

For $q=2$, Proposition~\ref{proAO9} immediately implies the following.

\begin{corollary}\label{proAO7}
The Boolean hypercube $Q_2^n$ contains no {\rm MDS} code with minimum
distance $d$ satisfying
$2<d<n$.
\end{corollary}

Thus, for a fixed alphabet size $q$, MDS codes of arbitrarily large
length exist only for minimum distances $2$ and $n$. MDS codes with
minimum distance $2$ correspond to Latin hypercubes, whereas an MDS
code in $Q_q^n$ with minimum distance $n$ is simply a diagonal of the
hypercube.

Since $|Q_q\times Q_p|=qp$, we identify the Cartesian product $Q_q\times Q_p$ with the set $Q_{qp}$.
Let   $x\in Q_q$ and $y\in Q_p$; we will consider 
the pair $(x,y)$  as an element of the set $Q_{pq}$.
The following theorem states that the direct product of MDS codes over
alphabets $Q_p$ and $Q_q$ is an MDS code over the alphabet $Q_{pq}$.

\begin{theorem}[MacNeish]\label{thMN}
If the sets $M_1\subset Q_q^n$ and $M_2\subset Q_p^n$ are {\rm
MDS} codes with distance $d$, then the set
$$M=\{(x,y)=((x_1,y_1),\dots,(x_n,y_n)) \ |\ x\in M_1, y\in M_2\}$$ is also an
 {\rm MDS} code with distance $d$ in $Q^n_{pq}$.
\end{theorem}
\begin{proof}
From the condition of the theorem and Proposition \ref{proSingl} it follows that
$|M|=q^{n-d+1}p^{n-d+1}$. Then $M$ attains the Singleton bound for codes with minimal distance $d$. It remains to prove that the set $M$ contains no
elements at distance less than $d$. It follows from the inequality
$$d((x,y),(x',y'))=\max\{d(x,x'),d(y,y')\}.$$
\end{proof}

MacNeish originally formulated this construction for systems of mutually orthogonal Latin squares.

\begin{remark}
Keevash \cite{Keevash18} proved that for any fixed $n\in
\mathbb{N}$ and code distance $d\leq n$, there exists $q_0\in
\mathbb{N}$ such that for $q\geq q_0$, the hypercube $Q^n_q$ contains an {\rm
MDS} code with distance $d$.
\end{remark}

As mentioned above, the case  $d=n-1$ corresponds to the system of $n-2$ MOLS.
In this particular case, Keevash's result means that for any $k$
there exists $q_0(k)\in \mathbb{N}$ such that for $q\geq q_0(k)$
there is a system of $k$ orthogonal Latin squares of order
$q$.  Bose,   Shrikhande and Parker \cite{BShP} proved that $q_0(2)=7$, more precisely,
for any $q\neq 2,6$ there exists a pair of orthogonal Latin
squares. We will prove this theorem in the next section.

\subsection{Systems of mutually orthogonal Latin squares}

Let $K$ be an abelian group of order $q$. Consider permutations
$\pi_1,\dots,\pi_k$ of $K$ satisfying the condition that, for every
$i\neq j$, the pointwise difference
\[
(\pi_i-\pi_j)(x)=\pi_i(x)-\pi_j(x)
\]
is also a permutation of $K$.

\begin{claim}
The functions
$f_i(x,y)=x+\pi_i(y)$, $i=1,\dots,k,$
form an orthogonal system.
\end{claim}

\begin{proof}
It is sufficient to show that, for every $i\neq j$, all $q^2$ pairs
\[
\bigl(f_i(x,y),f_j(x,y)\bigr)
\]
are distinct. Suppose that
\[
f_i(x_1,y_1)=f_i(x_2,y_2), \qquad
f_j(x_1,y_1)=f_j(x_2,y_2).
\]
Subtracting the two equalities gives
\[
(\pi_i-\pi_j)(y_1)=(\pi_i-\pi_j)(y_2).
\]
Since $\pi_i-\pi_j$ is a permutation by assumption, we obtain
$y_1=y_2$. Substituting this into
$f_i(x_1,y_1)=f_i(x_2,y_1)$ yields $x_1=x_2$.
\end{proof}

%From claim \ref{MDS-33} we have

\begin{corollary}\label{MOLS} If $q\neq 2$, $q=p^s$ is a prime
power, then there exists a system $f_1,\dots,f_{q-1}$ of mutually orthogonal
$2$-quasigroups (Latin squares) of order $q$.
\end{corollary}
\begin{proof}
We can identify  $Q_q$ with the Galois field  $GF(q)$.
Consider a set of permutations $\{\pi_i : i\in  GF(q)\backslash \{0\} \}$, where $\pi_i(x)=ix$.  Obviously, the difference $\pi_i-\pi_j=\pi_{i-j}$ is also a permutation. 
\end{proof}

\begin{claim}\label{clLS2}
Let the number $k$ have the following decomposition into prime factors:
$k=2^{\delta_1}3^{\delta_2}\cdots p_i^{\delta_i}\cdots$. Then
there exists a system of\\ $\min\limits_{i, \delta_i\neq 0}
(p_i^{\delta_i}-1)$ mutually orthogonal Latin squares
of order $k$.
\end{claim}
\begin{proof}
From Corollary \ref{MOLS} we obtain that there exists a system of
$p_i^{\delta_i}-1$ mutually orthogonal Latin squares
of order $p_i^{\delta_i}$. By Theorem \ref{thMN} it follows that if
there exist systems of $m$ MOLS
of orders $k_1$ and $k_2$, then there exists a system of $m$ MOLS of order $k_1k_2$.
\end{proof}

\begin{corollary}\label{corLS}
If $k\not\equiv 2\mod 4$, then there exists a pair of orthogonal Latin
squares of order $k$.
\end{corollary}

\begin{lemma}[\cite{DenesK}]\label{lemLS}
Suppose there exist pairs of orthogonal Latin squares of orders $k_1$,
$k_1+1$, $k_3$ and three mutually orthogonal Latin squares
of order $k_2$, $k_2\geq k_3$. Then there exists a pair of orthogonal
Latin squares of order $k_1k_2+k_3$.
\end{lemma}
\begin{proof}
Let $f_i:Q_{k_2}^2\rightarrow Q_{k_2}$, $i=0,1,2$, be three mutually
orthogonal Latin squares of order $k_2$;
$g_i:Q_{k_1}^2\rightarrow Q_{k_1}$, $i=1,2$, be a pair of orthogonal
Latin squares of order $k_1$;
$h_i:(Q_{k_1}\cup\{q\})^2\rightarrow Q_{k_1}\cup\{q\}$, $i=1,2$, be a
pair of orthogonal Latin squares of order $k_1+1$ such that
$h_1(q,q)=h_2(q,q)=q$; $t_i:Q_{k_3}^2\rightarrow Q_{k_3}$, $i=1,2$,
be a pair of orthogonal Latin squares of order $k_3$.

The alphabet $R$ of cardinality $k_1k_2+k_3$ will be represented in the form $(Q_{k_1}\times
Q_{k_2})\cup Q_{k_3}$, where $Q_{k_3}=\{q_0,\dots,q_{k_3-1}\}$.
Define a pair of orthogonal Latin squares $F_1$, $F_2$ in the
alphabet $R$ according to the following rule.

If $x_1,x_2\in Q_{k_1}$, $y_1,y_2\in Q_{k_2}$ and $f_0(y_1,y_2)\geq
k_3$, then $$F_i((x_1,y_1),(x_2,y_2))=(g_i(x_1,x_2),f_i(y_1,y_2)), \qquad
i=1,2.$$

If $x_1,x_2\in Q_{k_1}$, $y_1,y_2\in Q_{k_2}$ and $f_0(y_1,y_2)<
k_3$, then when $h_i(x_1,x_2)\neq q$ we define
$$F_i((x_1,y_1),(x_2,y_2))=(h_i(x_1,x_2),f_i(y_1,y_2)), \qquad i=1,2,$$ and
when $h_i(x_1,x_2)=q$ we define $$F_i((x_1,y_1),(x_2,y_2))= q_j, \qquad \mbox{\rm where}\quad
f_0(y_1,y_2)=j. $$ Note that by definition, $h_1(x_1,x_2)$ and
$h_2(x_1,x_2)$ cannot equal $q$ simultaneously.

 If $x_1\in Q_{k_1}$, $y_1\in Q_{k_2}$ and $q_j\in
\{q_0,\dots,q_{k_3-1}\}$, then we define
$$F_i((x_1,y_1),q_j)=(h_i(x_1,q),f_i(y_1,y)), \qquad i=1,2,\qquad \mbox{\rm where}\quad
f_0(y_1,y)=j.$$

Symmetrically, if $x_2\in Q_{k_1}$, $y_2\in Q_{k_2}$ and $q_j\in
\{q_0,\dots,q_{k_3-1}\}$, then we define
$$F_i(q_j,(x_2,y_2))=(h_i(q,x_2),f_i(y,y_2)), \qquad i=1,2, \qquad \mbox{\rm where}\quad
f_0(y,y_2)=j.$$

For $q_{j_1},q_{j_2}\in \{q_0,\dots,q_{k_3-1}\}$, we define
$F_i(q_{j_1},q_{j_2})=t_i(q_{j_1},q_{j_2})$, $i=1,2$.

It suffices to check that $F_1$ and $F_2$ are Latin
squares in the alphabet $R$ and for any pair $(a,b)\in R$ there exists
a pair $u,v\in R$ such that $a=F_1(u,v)$, $b=F_2(u,v)$. 

Let us verify that the  row of $F_i$ does not contain two equal elements.  Since $f_0$ is a quasigroup, there is at most  one $q_j\in Q_{k_3}$ in a row.  Equalities $(g_i(x'_1,x_2),f_i(y'_1,y_2))=(g_i(x_1,x_2),f_i(y_1,y_2))$  or $(h_i(x'_1,x_2),f_i(y'_1,y_2))=(h_i(x_1,x_2),f_i(y_1,y_2))$ imply $x'_1=x_1$ and $y'_1=y_1$  since $g_i$, $h_i$ and $f_i$ are quasigroups.  If  $(g_i(x'_1,x_2),f_i(y'_1,y_2))=(h_i(x_1,x_2),f_i(y_1,y_2))$ then $y'_1=y_1$. Therefore, the row contains only one type element, either  $(g_i(x'_1,x_2),f_i(y_1,y_2))$ or $(h_i(x_1,x_2),f_i(y_1,y_2))$ depending on whether $f_0(y_1,y_2)\geq k_3$.

 Let $a=(x_3,y_3)$ and
$b=(x_4,y_4)$, where $x_3,x_4\in Q_{k_1}$, $y_3,y_4\in Q_{k_2}$.
There exist $y_1,y_2 \in Q_{k_2}$ such that $f_1(y_1,y_2)=y_3$ and
$f_2(y_1,y_2)=y_4$. Let $f_0(y_1,y_2)\geq k_3$, there exist
$x_1,x_2\in Q_{k_1}$ such that $g_1(x_1,x_2)=x_3$,
$g_2(x_1,x_2)=x_4$. Then $u=(x_1,y_1)$, $v=(x_2,y_2)$. Let
$f_0(y_1,y_2)< k_3$, either there exist $x_1,x_2\in Q_{k_1}$ such that
$h_1(x_1,x_2)=x_3$, $h_2(x_1,x_2)=x_4$ and then $u=(x_1,y_1)$,
$v=(x_2,y_2)$; or there exists $x_1\in Q_{k_1}$ such that $h_1(x_1,q)=x_3$,
$h_2(x_1,q)=x_4$ and then $u=(x_1,y_1)$, $v=q_j$, where
$j=f_0(y_1,y_2)$; or there exists $x_2\in Q_{k_1}$ such that
$h_1(q,x_2)=x_3$, $h_2(q,x_2)=x_4$ and then $u=q_j$, where
$j=f_0(y_1,y_2)$, $v=(x_2,y_2)$.

Let $a=(x_3,y_3)$ and $b=q_j$, where $x_3\in Q_{k_1}$, $y_3\in
Q_{k_2}$, $q_j\in Q_{k_3}$. There exist $y_1,y_2 \in Q_{k_2}$ such
that $f_0(y_1,y_2)=j$ and $f_1(y_1,y_2)=y_3$. Let
$f_2(y_1,y_2)=y_4$. There exist $x_1,x_2\in Q_{k_1}$ such that
$h_1(x_1,x_2)=x_3$ and $h_2(x_1,x_2)=q$. Then $u=(x_1,y_1)$,
$v=(x_2,y_2)$. The case $a=q_j$ and
$b=(x_4,y_4)$, where $x_4\in Q_{k_1}$, $y_4\in Q_{k_2}$, $q_j\in
Q_{k_3}$, can be considered analogously.

If $a,b\in Q_{k_3}$, then the required follows from the orthogonality of the
Latin squares $t_1$ and $t_2$.
\end{proof}

Note that in Lemma \ref{lemLS} one can take $k_3=1$.

\begin{lemma}[\cite{DenesK}]\label{lemLS1}
For $k\equiv 1\mod 4$, $k\geq 5$, there exists a triple of mutually
orthogonal Latin squares of order $k$.
\end{lemma}
\begin{proof}
If $k$ is not divisible by $3$, then a triple of mutually orthogonal
Latin squares of order $k$ exists by Proposition \ref{clLS2}.

Let $k=4t+1=3p$. Then $t$ is not divisible by $3$. If $t=2$ or
$t=4$, then the required follows from Proposition  \ref{clLS2}. If
$t=2^\delta$, $\delta\geq 3$, or all divisors of $t$ are at least $5$, then
there exists a system of $4$ mutually orthogonal Latin squares
$f_0,f_1,f_2,f_3$ of order $t$ by Proposition \ref{clLS2}. Let
$g_1,g_2,g_3$ and $h_1,h_2,h_3$ be triples of mutually orthogonal Latin
squares of orders $4$ and $5$, respectively. From the proof of Lemma
\ref{lemLS} ($k_1=4$, $k_2=t$, $k_3=1$) it is clear that one can obtain three
pairs $F_1, F_2$, $F_2, F_3$ and $F_1, F_3$ of orthogonal Latin squares of order $k$, using the triple
$f_0,f_i,f_j$, and two pairs $g_i,g_j$ and $h_i,h_j$, where
$i,j\in\{1,2,3\}$.
  Thus, a triple $F_1, F_2, F_3$ of mutually
orthogonal Latin squares of order $k$ is obtained.
\end{proof}

\begin{theorem}[Bose, Shrikhande and Parker, \cite{BShP}]\label{thLS1}
For $k\neq 2,6,10,14$, there exists a pair of orthogonal Latin
squares of order $k$.
\end{theorem}
\begin{proof}
By Corollary \ref{corLS}, it remains to investigate the case $k\equiv 2\mod
4$. Three cases $k\equiv 2,6,10\mod 12$  are possible. Then 
$k=3(4t+1)+11$, $k=3(4t+1)+3$, or $k=3(4t+1)+7$. From Lemmas
\ref{lemLS} and \ref{lemLS1} we obtain the required statement for $k_1=3$,
$k_2=4t+1$, and $k_3=11,3,7$, with the exception of the cases
$k=2,6,10,14,22,26,38$. For $k=22$ and $k=26$, one can apply Lemma
\ref{lemLS} with $k_1=3$, $k_2=7$, $k_3=1$ and $k_3=5$. For $k=38$,
one can apply Lemma \ref{lemLS} with $k_1=7$, $k_2=5$, and $k_3=3$.
\end{proof}

Pairs of orthogonal Latin squares of orders $2$ and $6$ do not
exist. Pairs of orthogonal Latin squares of orders $10$ and
$14$ have been constructed by   computer search. In \cite{Ji}, pairs of
orthogonal cubes of all orders $k$, with the exception of $k=3$ and $k\equiv
2\mod 4$, have been constructed.
 The  existence of a triple of
mutually orthogonal Latin squares of order $10$ and of a pair of
orthogonal Latin cubes of order $10$ remain
open problems. \label{LSprob}

An {\sl affine plane} is a system of points and lines
satisfying the following properties.

1) Any two distinct points lie on a unique line.

2) Given any line and any point not on that line there is a unique line which contains the point and does not meet the given line. 

3) There exist four points such that no three are collinear (points not on a single line).

It can be shown that in a finite affine plane, each line
contains the same number of points $k$, called the {\sl order} of the
plane. Moreover, all lines are partitioned into $k+1$ classes of $k$
parallel lines each. Thus, an affine plane
contains $k^2$ points.

\begin{claim} Affine planes of order $k$ correspond bijectively
to {\rm MDS} codes in $Q^{k+1}_k$ with distance $k$
(systems of $k-1$ mutually orthogonal Latin squares
of order $k$).
\end{claim}
\begin{proof} Suppose an {\rm MDS} code $C\subset Q^{k+1}_k$ with minimum distance
$k$ is given. We call the elements of the code points, and we take as lines the
intersections of the {\rm MDS} code $C$ with hyperplanes. The 
property 1) is satisfied by Corollary \ref{AP01}. Properties 2) and 3)
hold by the construction. Let us prove that an affine plane generates an
{\rm MDS} code. Enumerate  $k$ lines of each parallel class
with the numbers $0,\dots,k-1$. To each point we associate a tuple from $Q^{k+1}_k$,
where in the $i$-th place is the number of the line from the $i$-th class to which
this point belongs. By construction, the resulting set has
cardinality $k^2$, and from property 1) it follows that the distances between
points are $k$. Indeed, if two tuples coincide in two positions then there exist two lines containing both points.  The resulting subset of the hypercube is an {\rm
MDS} code by Proposition \ref{proSingl} (b).
\end{proof}

No system of $k-1$ mutually orthogonal Latin
squares of order $k$ is known when $k$ is not a prime
power. The conjecture that such systems, and consequently,
affine planes of order not equal to a prime power, do not
exist is one of the classical problems of combinatorics.
\label{PGprob}

Recall that a strongly regular graph is a distance-regular graph of
diameter $2$ (if $k\neq v-1$) with parameters $(b_0, b_1; c_1, c_2)=
(k, k-\lambda-1; 1, \mu)$.

Let $C\subset Q^n_q$ be some code with code distance $d$.
Consider the minimum distance graph $G(C)$, whose vertices
are the vertices of the code $C$, and two vertices are joined by an edge if
the Hamming distance between them equals $d$.

\begin{claim}
If $C\subset Q^n_q$ is an {\rm MDS} code with distance $n-1$, then
$G(C)$ is a strongly regular graph.
\end{claim}
\begin{proof}
Since an {\rm MDS} code with distance $n-1$ in $Q^n_q$ is defined by
a system of $n-2$ mutually orthogonal Latin squares, by
Proposition \ref{proAO5} we have $n\leq q+1$. From Corollary \ref{AP01}
it follows that for $n=q+1$, the graph $G(C)$ is complete. We
further assume that $n\leq q$. 
$G(C)$ is a $k$-regular graph, where $k=n(q-1)$. Indeed, every codeword $u$ belongs to $n$ hyperplanes and every hyperplane contains $q$ codewords. The intersection of two hyperplanes with different directions  is a
face of dimension $n-2$. Then this intersection contains only the  codeword $u$. 

If two vertices of the {\rm MDS} code are at
distance $n-1$, then their common neighbors are all the vertices
lying in the $(n-1)$-dimensional face common to them. Then
$\lambda=\frac{|C|}{q}-2=q-2$. Suppose two vertices $u,v\in C$
are at distance $n$. Let $u=(u_1,\dots,u_n)$ and $v=(v_1,\dots,v_n)$ with $d(u,v)=n$. 
A common neighbor must agree with $u$ in exactly one coordinate and with $v$ in exactly one different coordinate. There are $n(n-1)$ ordered choices of these coordinates.
The intersection of any two $(n-1)$-dimensional
faces of different directions  is a
face of dimension $n-2$. By the definition of an {\rm MDS} code, any
$(n-2)$-dimensional face of the hypercube contains exactly one code vertex.
Then $\mu$ equals the number of pairs of such $(n-1)$-dimensional faces, i.e.,
$\mu=n(n-1)$.
\end{proof}

Every Latin square of order $q$ corresponds to an MDS code with
distance $2$ in $Q^3_q$. Moreover, every system of $n-2$
orthogonal Latin squares of order $q$
 defines an {\rm MDS} code with distance $n-1$ in the hypercube $Q^n_q$.
 Consequently, every Latin square and every system of MOLS correspond to a strongly regular graph via the minimum-distance graph of the corresponding MDS code.

\subsection{Transversals of Latin Hypercubes}

%\begin{definition}
A {\sl transversal} of a Latin hypercube (and, in particular, of a
Latin square) is a diagonal of the hypercube whose entries are pairwise
distinct.
%\end{definition}

It follows immediately that a Latin square admits an orthogonal Latin
square if and only if it can be partitioned into transversals.
For Latin hypercubes, however, the situation is different: a necessary
condition for the existence of an orthogonal Latin hypercube is that
every two-dimensional face (that is, every Latin square contained in
the hypercube) can be partitioned into transversals. A partition of the
entire hypercube into transversals is not sufficient.

\begin{claim}\label{transver}
If $q$ and $n$ are even, then the Latin hypercube defined by
the equality $f(x_1,\dots,x_n)=x_1+\dots+x_n \mod q$ has no
transversals.
\end{claim}
\begin{proof}
We have $\sum\limits_{j=0}^{q-1}j= \frac{q(q-1)}{2}\mod q$. Suppose the
Latin hypercube contains a transversal $y^0,y^1,\dots,y^{q-1}\in
Q^n_q$. Then $$\sum_{i=1}^{n}\sum_{j=0}^{q-1} y^j_i =
\sum_{i=1}^{n} \frac{q(q-1)}{2} = q\frac{n(q-1)}{2}=0\mod q.$$ While
$$\sum_{j=0}^{q-1} f(y^j) = \sum_{j=0}^{q-1}j= \frac{q(q-1)}{2}=\frac{q}{2}\mod
q,$$ which is nonzero because $q$ is even. We obtain a contradiction with the equality that defines the
Latin hypercube.
\end{proof}

\begin{remark} Ryser's well-known conjecture states that any Latin square of
odd order has a transversal. By
Proposition \ref{transver} it is clear that not all Latin squares of
even order contain a transversal. However, one can ask the question
of the existence in Latin squares of a partial transversal of the greatest
possible length. Montgomery \cite{Mont} proved that Latin
squares of sufficiently large even order $q$ contain a partial
transversal of length $q-1$. For Cayley  tables of groups, a criterion  of the existence  of transversals is known.  Hall and Paige conjectured that the  Cayley table of a finite group  has transversals if and only if  its Sylow $2$-subgroups are trivial or noncyclic. This conjecture has now been  proved (see \cite{Bray}).

Wanless \cite{Wanless}
generalized  Ryser's  conjecture to Latin hypercubes of odd orders. It is true for orders $1$,  $3$ and $5$ (see \cite{PPV}). Taranenko \cite{Taranenko18} proved that for order $4$ and every even $n$  there is only one  Latin $n$-cube without transversal,   up to equivalence. This is described in Proposition \ref{transver}. \label{Tprob}|
\end{remark}

\subsection{Linear MDS Codes}\label{8.3}

Further, we 
construct linear MDS codes over the Galois field $GF(q)$, where $q$ is
a prime power. The hypercube $Q^n_q$ will now be
considered as an $n$-dimensional vector space over 
$GF(q)$. In this sense, the concepts of linearity, affine hull, dimension, etc.,
apply to subsets of the hypercube.
%\begin{definition}
 Recall that   the { rank} of a code is the dimension of its
affine hull.
%\end{definition}

By Proposition \ref{proSingl}  the rank  of a linear {\rm MDS} code $C\subset
Q_q^n$ with minimum distance $d$ equals $r=n-d+1$. From Proposition
\ref{proAO2} we note that the rank of a linear MDS code with distance
greater than 2 does not exceed $q-1$.

We will use the representation of a linear {\rm MDS} code in the form of the
solution of a system of linear equations $C=\{x\in Q_q^n\ |\ Ax=\bar{0}\}$. Recall that the matrix $A$ 
is called a parity-check matrix for the code $C$. Without loss of generality we admit that the parity-check 
$r\times n$ matrix $A$ has rank $r$.

\begin{claim}\label{MDS-33} A linear code $C$ with parity-check $r\times n$ matrix $A$
is an {\rm
    MDS} code of rank $n-r$ with minimum distance $d=r+1$
if and only if any $r$ columns of the matrix
$A=\{a_{ij}\}$ form a nonsingular matrix.\end{claim}
\begin{proof} Any linear code is the solution set of a system of linear equations.
 It suffices to note that any face of dimension $r$ in the hypercube
 contains exactly one codeword of the code $C$ if and only
if the corresponding minor in the matrix $A$ is non-zero, i.e., any
set of $r$ columns of the matrix $A$ is a linearly independent
basis. Indeed, fix $n-r$ coordinates.
The remaining $r$ coordinates satisfy an $r\times r$ linear system.
There is a unique solution if and only if the corresponding submatrix is nonsingular.
\end{proof}

 The code $C^\perp=\{x\in Q_q^n : \sum_i x_iy_i=0\
\mbox{for all}\ y\in C\}$ is called {\sl dual} to $C$. Obviously, code $C^\perp$  is linear. 
A linear code and its dual code are MDS codes at the same time.
In Section \ref{12.1} we considered  dual codes over  rings $\mathbb{Z}_q$. These notions of duality coincide in the case of prime $q$.

\begin{claim} If $C\subset Q_q^n$ be a linear {\rm MDS} code of rank $n-r$,
then $C^\perp$ is a linear {\rm MDS} code of rank $r$.
\end{claim}
\begin{proof} Let $C=\{x\in Q_q^n\ |\ Ax=\bar{0}\}$, where
$A$ is an $r\times n$ matrix and the linear combinations of the rows
of the matrix $A$ form the code $C^\perp$. Consider some set consisting of
$r$ coordinates. The submatrix of $A$ consisting of these $r$ columns forms a nonsingular square matrix by Proposition \ref{MDS-33}.  Without loss of generality we can  take the first $r$ coordinates. Consider an arbitrary  $u\in Q_q^r$.  There is a unique linear combination $w$ of rows of $A$ such that the first $r$ coordinates of $w$ coincide with $u$.   Then  there is a unique codeword from  $C^\perp$ in every $(n-r)$-dimensional face.  
Thus, the set $C^\perp$ is a linear MDS code of rank
$r$.
\end{proof}

\begin{corollary}  Linear {\rm MDS} codes  in $Q_q^n$ with distances $d$
and $n-d+2$ exist simultaneously. \end{corollary}

\begin{claim}\label{RS} If $q=p^s$ is a prime power, then in
$Q_q^{q+1}$ there exist linear {\rm MDS} codes with distance $d$,
for any $d$, $3\leq d\leq q$. \end{claim} \begin{proof} The matrix
$V_r=\left(\begin{array}{cccccc}
0&1&1& 1&\dots&1  \\
0&0&1& 2&\dots& q-1 \\
0&0&1& 2^2 &\dots & (q-1)^2\\
\dots &\dots &\dots & \dots &\dots & \dots \\
1&0&1& 2^{r-1} &\dots & (q-1)^{r-1}\\
\end{array} \right)$
satisfies the conditions of Proposition \ref{MDS-33} for $r< q$,
since all its minors of order $r$ are non-zero. More precisely, any $r$ of
the last $q-1$ columns of the matrix constitute a Vandermonde
determinant, and other sets of $r$ columns can be expressed by Vandermonde determinant of smaller order.
\end{proof}

\begin{claim} \label{MDS-333} Let $q=2^t$, $t\geq 2$, $n=q+2$. There exists a
 linear {\rm MDS} code $C\subset Q_q^n$ of rank $n-3$.
\end{claim}
\begin{proof} Note that $a^2-b^2=(a-b)^2$ in $GF(2^t)$, i.e., $a^2\neq b^2$ for $a\neq b$. It is easy to check that
all minors of dimension $3$ of the matrix $\left(\begin{array}{ccccccc}
1&0&0& 1&\dots&1  \\
0&1&0& 2&\dots& q-1 \\
0&0&1& 2^2 &\dots & (q-1)^2
\end{array} \right)$
are non-zero. Then the proposition follows from Proposition \ref{MDS-33}.
\end{proof}

The {\rm MDS} codes proposed in Proposition \ref{MDS-333} are called {\sl
    Asturian codes}.

\begin{remark}
Ball \cite{Ball} proved that if $q$ is prime, then the length of a
linear MDS code does not exceed $q+1$ (when $q=2^t$, $t>1$, the length of a
linear MDS code does not exceed $q+2$), except for the cases
when the code distance equals $2$ or the length of the code. \label{MDSprob}
\end{remark}

Ball's result strengthens Proposition \ref{proAO9} for linear
codes. A well-known MDS conjecture states that the same statement holds
for arbitrary MDS codes when $q\neq 2^t$. In particular, Bespalov
\cite{Bespalov} proved that for odd $q$ there do not exist
MDS codes of length $q+2$ or more with minimum distance $4$.

Let $\alpha$ be a primitive element in the field $GF(q)$. Consider the
matrix

$RS= \left(\begin{array}{ccccc}
1& 1&     1&       \dots&1  \\
1& \alpha &\alpha^2&\dots& \alpha^{q-2} \\
1& \alpha^2 &\alpha^4&\dots& \alpha^{2(q-2)} \\
\dots& \dots& \dots& \dots & \dots\\
1& \alpha^r &\alpha^{2r}& \dots & \alpha^{r(q-2)}
\end{array} \right).$

Every set of $r+1$ ($r<q-2$) columns of the matrix $RS$ constitutes
a Vandermonde matrix. By Proposition \ref{MDS-33}, the matrix $RS$
is a parity-check matrix of an {\rm MDS} code. It is easy to see that
this code is {\sl cyclic}, i.e., words of the form
$(y_1,y_2,\dots,y_n)$ and $(y_2,\dots,y_n,y_1)$ are either both contained or
both not contained in the code. The {\rm MDS} code with parity-check
matrix $RS$ is called the {\sl  Reed--Solomon code}. 

The main results of this section were generalized to Galois rings in \cite{Dong}.

\subsection{Applications of MDS Codes}

MDS codes have numerous applications in information processing.

Consider first the transmission of information over a communication
channel with erasures. During transmission, a message
$(a_1,a_2,\dots,a_n)$ may be transformed into a word of the form
$(a_1,*,a_3,\dots,*,\dots,a_n)$, where the symbol $*$ denotes an erased
letter. If the transmitted words are codewords of an MDS code with
minimum distance $d$, then up to $d-1$ erased symbols can be recovered.
In fact, every code with minimum distance $d$ allows the recovery of
$d-1$ erasures, since every face of dimension less than $d$ contains at
most one codeword. What distinguishes MDS codes is that every face of
dimension $d-1$ contains exactly one codeword.

MDS codes are also widely used in distributed storage systems. Suppose
that information is stored across $n$ servers and that it must be
possible to recover the entire file after the failure of fewer than
$d$ servers. This can be achieved by representing the information as
codewords of an MDS code with minimum distance $d$, storing the
$i$-th coordinate of each codeword on the $i$-th server. The failure
of $d-1$ servers is then equivalent to the erasure of $d-1$
coordinates, from which the original information can be recovered
uniquely.

Another important application of MDS codes is secret sharing. Suppose
there are $n-1$ participants. The goal is to distribute shares of a
secret key so that any $m$ participants together can reconstruct the
key, whereas any group of at most $m-1$ participants has no
information  about it.

Assume there exists an MDS code of length $n$ with minimum distance
$n-m+1$. The code itself is assumed to be public. Let
$(a_1,a_2,\dots,a_n)$ be an arbitrary codeword. The value
$a_i\in Q_q$ is assigned as the private share of the $i$-th
participant, while $a_n\in Q_q$ is taken as the secret key. Since any
$m$ coordinates uniquely determine the remaining ones, every group of
$m$ participants can reconstruct the entire codeword, and hence the
secret key $a_n$, by solving a system of the form
(\ref{eMDS11}). On the other hand, any group of only $m-1$
participants knows the values of only $m-1$ variables. Consequently,
the $m$-ary quasigroup
\[
f_{n-m}(x_1,\dots,x_m)=x_n
\]
takes each value of $Q_q$ exactly once as the remaining unknown
variable ranges over $Q_q$. Therefore, every possible value of the
secret key is equally likely, and the participants obtain no
information about it.

Reed--Solomon codes over the
field $GF(2^t)$ are especially useful for applications. These codes are widely used as error-correcting codes.  In particular, for $t=8$, the elements of the alphabet are
bytes.  The Reed--Solomon code  with length $2^t-1$ and minimum distance $r+1$ can correct up to
$r/2$ errors,  i.e., corrupted  binary blocks of length $t$.
 The Reed--Solomon code is used in RAM controllers, when writing information to
hard drives, and on CD/DVD disks.

\subsection{Problems}

\begin{exercise}
Prove that if a Latin parallelepiped of size $q\times
q\times (q-1)$ is filled with elements of $Q_q$ such that in each line
no element occurs twice, then this parallelepiped can be
completed to a Latin cube.\end{exercise}

\begin{exercise}
Prove that for every $n\geq 2$, there is,
up to isometry of  $Q^n_3$,  exactly one {\rm MDS} code with distance $2$.
\end{exercise}

\begin{exercise}
Prove that  every  Latin $n$-cube of order $3$ has transversals.
\end{exercise}

\begin{exercise}
 Prove that if there exists a pair of
orthogonal Latin hypercubes in a $q$-ary alphabet, then their
dimensions do not exceed $q-1$.
\end{exercise}

\begin{exercise}
Prove that a set of $q-2$ mutually orthogonal Latin
squares of order $q$ can always be completed
 to a set of $q-1$ mutually orthogonal Latin squares of order $q$.
\end{exercise}

\begin{exercise}[\cite{Pot22}]
Prove that every partial Latin $n$-cube of order $q$, 
 can be completed to a Latin $n$-cube of order not greater than $q^n$.
\end{exercise}

\begin{exercise}
Prove that an MDS code with distance $3$ is a completely
regular code in $Q^n_q$, more precisely the first cell of a perfect
$3$-coloring with parameters $\begin{pmatrix}
0 & n(q-1) & 0 \\
1 & (n-2)(n-1)+q-2 & (n-1)(q-n+1) \\
0 & n(n-1) & n(q-n)
\end{pmatrix}$.
\end{exercise}

\begin{exercise}
Let an $r\times n-r$ matrix $A$ have only non-zero
minors of all sizes over the field $GF(q)$. Prove that the matrix
$(I_r|A)$ of size $r\times n$ is a parity-check matrix of an
MDS code.
\end{exercise}

\begin{exercise}
Represent the matrix $V_r$ from Proposition \ref{RS} in the form $V_r=(B|A)$,
where $B$ is of size $r\times r$. Prove that the matrix $B^{-1}A$ does not
contain any zero minors of any orders.
\end{exercise}

\begin{exercise}
Let $a_1,\dots,a_n,b_1\dots,b_n$ be pairwise distinct elements of the field
$GF(q)$. Prove that the Cauchy matrix $C$, where
$c_{ij}=\frac{1}{a_i-b_j}$, has no zero minors of any
orders.
\end{exercise}

\begin{exercise}
Let $C\subset Q^n_q$ and suppose each ball of radius $\rho$ contains no
more than $L$ elements of $C$. Prove that $|C|\leq Lq^{n-\lfloor
\frac{(L+1)\rho}{L}\rfloor}$ (the generalized Singleton bound).
\end{exercise}

\section{Combinatorial $t$-Designs}

\subsection{Designs as Perfect Colorings of Graphs and Hypergraphs}

Recall that the Johnson graph $J(n,k)$ is the graph whose vertices
are all binary tuples of length $n$ and weight $k$. Two
vertices of the graph are joined by an edge if the Hamming distance between
the corresponding binary tuples is two. Thus, the
adjacency matrix of the Johnson graph is a submatrix of the distance-$2$
matrix in the Boolean hypercube.
% $M_2=\frac12(M^2-nI)$, where $M$ is the adjacency matrix of the Boolean cube.
The vertices of the Johnson graph $J(n,k)$
can be viewed as points of the integer simplex
$\{(x_1,\dots,x_k) \ |\ x_i\in \mathbb{N},
\sum\limits_{i=1}^kx_i\leq n\}$, where $\sum\limits_{i=1}^jx_i$
is the index of the coordinate containing the $j$-th one of the binary tuple.

%\begin{definition}
A {\sl combinatorial design} (more precisely,
a $t$-design) with parameters $t$-$(n,k,\lambda)$
 is a collection  $k$-element
subsets (blocks) of an $n$-element set, such that any $t$-element
subset is contained in exactly $\lambda$ blocks.
%\end{definition}

The blocks of a design with parameters $t$-$(n,k,\lambda)$ can be considered
as vertices of the graph $J(n,k)$. A combinatorial $t$-design is conveniently
represented as a list of blocks. A block can be specified by the
indices of the elements of the block or by the characteristic function of the block, i.e.,
as a Boolean vector of length $n$ and weight $k$. A detailed survey of the theory
of combinatorial designs can be found in \cite{CDinitz}.

\begin{example} Table of blocks of a design with parameters 2-(13,4,1).

$\begin{array}{cccc}
  2 & 3 & 5 & 11\\
  3 & 4 & 6 & 12\\
  4 & 5 & 7 & 13\\
  1 & 5 & 6 & 8\\
  2 & 6 & 7 & 9\\
  3 & 7 & 8 & 10\\
  4 & 8 & 9 & 11\\
  5 & 9 & 10 & 12\\
  6 & 10 & 11 & 13\\
  1 & 7 & 11 & 12\\
  2 & 8 & 12 & 13\\
  1 & 3 & 9 & 13\\
  1 & 2 & 4 & 10
\end{array}$
\end{example}

\begin{claim}\label{c:design1}
The characteristic function of a design $D$ with parameters
$(k-1)$-$(n,k,\lambda)$
 is a perfect coloring of the graph
$J(n,k)$ with quotient matrix $S=\left(\begin{array}{cc}
k(\lambda-1)& k(n-k-\lambda+1) \\
k\lambda&  k(n-k-\lambda)\\
\end{array} \right)$.
\end{claim}
\begin{proof}
 Two blocks $u,v\in V(J(n,k))$ are adjacent if and only
if their intersection contains $(k-1)$ elements. Any
$(k-1)$-element set is contained in $\lambda$ blocks of the design
$D$. Consequently,
 a block $u\in D$ is adjacent to $\lambda-1$ vertices $v\in D$ that intersect
   it in some $(k-1)$-element
set. Meanwhile, a $k$-element block $u\not \in D$ has
$\lambda$ adjacent vertices $v\in D$ that intersect
   it in a fixed $(k-1)$-element
set. Since  $u\in D$ includes exactly $k$ subsets
of cardinality $k-1$, we have $s_{11}=k(\lambda-1)$ and $s_{21}=k\lambda$.
The degree of the graph $J(n,k)$ is $k(n-k)$. Hence, we have
$s_{12}=k(n-k-\lambda+1)$ and $s_{22}=k(n-k-\lambda)$.
\end{proof}

For $\lambda=1$, combinatorial designs with parameters
$t$-$(n,k,\lambda)$ are usually denoted $S(t,k,n)$, and for $t=k-1$
they are additionally called {\sl Steiner systems} (for $k=3$, {\sl Steiner triple
systems}) of order $n$.

For $\lambda=1$, a combinatorial design with parameters $t$-$(n,k,1)$
can be regarded  as a transversal in the hypergraph $G_{n,k,t}$,
whose vertices are all possible $k$-element blocks, and
hyperedges
are the sets of blocks containing a fixed $t$-element
set.

In the case $t=k-1$, the multigraph
$\mathcal{M}_{12}(\mathcal{D}(G_{n,k,t}))$ of vertex adjacency corresponding to
hypergraph $G_{n,k,t}$ coincides with
the graph $J(n,k)$  after removing all loops.

 As a special case of Proposition
\ref{c:design1} with $\lambda=1$ (see also Problem \ref{exerDH1} and
Proposition \ref{PCclaim11} on the Delsarte--Hoffman bound) we have the following corollary.

\begin{corollary}
Every Steiner system $S(k-1,k,n)$ is a maximum independent
sets  of $J(n,k)$.
\end{corollary}

Define the multigraph $(J(n,k))_{k-t}$ as follows. We take the vertices of the Johnson graph 
$J(n,k)$; two $k$-blocks  $u$ and $v$ are adjacent if the distance between $u$ and $v$ does not exceed  $k-t$ in $J(n,k)$ and the multiplicity of the edge equals the number of common  $t$-subsets of $u$ and $v$.
$(J(n,k))_{k-t}$ is the same as 
$\mathcal{M}_{12}(\mathcal{D}(G_{n,k,t}))$ (see Section \ref{hyp}) without loops.

\begin{claim}[\cite{PA20}]\label{c:design11}
The characteristic function of a design $D$ with parameters
$t$-$(n,k,\lambda)$ is a perfect coloring of the multigraph
$(J(n,k))_{k-t}$ with quotient matrix
 $S=\left(\begin{array}{cc}
{k\choose t}(\lambda-1)& {k\choose t}({{n-k}\choose{k-t}}-\lambda+1) \\
{k\choose t}\lambda&  {k\choose t}({{n-k}\choose{k-t}}-\lambda)\\
\end{array} \right)$
and the design $D$ is a maximum independent
set in $(J(n,k))_{k-t}$ if $\lambda=1$.
\end{claim}
\begin{proof}
We show that if $\lambda=1$ then the design $D$ is a maximum independent set.
Consider the hypergraph $G_{n,k,t}$. The vertices of the graph $(J(n,k))_{k-t}$
are joined by an edge if and only if they are contained in
some edge of the hypergraph $G_{n,k,t}$. An independent set
in $(J(n,k))_{k-t}$ intersects an edge of the hypergraph in at most
one element, therefore its cardinality is not greater than the cardinality of a
transversal of the hypergraph. Since the design $D$ is a transversal
in $G_{n,k,t}$, it is a maximum independent set in
$(J(n,k))_{k-t}$. The proof that the characteristic
function of the design $D$ with parameters $t$-$(n,k,\lambda)$ is a
perfect $2$-coloring of the multigraph $(J(n,k))_{k-t}$ is identical to Proposition \ref{c:design1}, replacing $k$ $(k-1)$-subsets of a block by its ${k \choose t}$ $t$-subsets.
\end{proof}

The converse is also true: the first cell of a equitable partition of the multigraph
$(J(n,k))_{k-t}$ with quotient matrix $S$ is a
$t$-$(n,k,\lambda)$ design. Let us prove this in the particular case
$\lambda=1$.

\begin{claim}[\cite{PA20}]\label{c:design19}
A set $D$ is a design with parameters $t$-$(n,k,1)$ if and
only if its characteristic function is a
perfect coloring of the multigraph $(J(n,k))_{k-t}$ with quotient
matrix $S=\left(\begin{array}{cc}
0 & {k\choose t}{{n-k}\choose{k-t}} \\
{k\choose t}&  {k\choose t}({{n-k}\choose{k-t}}-1)\\
\end{array} \right)$.
\end{claim}
\begin{proof}
($\Rightarrow$) Follows from Propositions 
\ref{c:design11}.

($\Leftarrow$) Let  $D$ be the first cell in an equitable partition
 with quotient matrix $S$. Any $t$-element set
is contained in at most one $k$-block from $D$, since
$s_{00}=0$. Consider a $k$-block $u\not \in D$.   Every adjacency between $u$ and a block of $D$ arises from one of the  ${k\choose t}$ $t$-element subsets of $u$. 
Since $s_{10}={k\choose t}$,  every one of these $t$-element subsets belongs to a $k$-block from $D$.  We conclude that any $t$-element subset
is contained in a $k$-block from $D$. \end{proof}

Note that Proposition \ref{c:design19} is essentially a special
case of Proposition \ref{bipartity1}.

\begin{remark}
The problem of finding perfect colorings of the Johnson graphs is very complicated.
For example, only trivial $e$-perfect codes in the Johnson graph are currently known. Trivial examples include the entire vertex set ($e=0$), a single vertex ($e$ is the diameter of the graph), or specific single-element shifts when parameters align trivially. Unlike perfect codes, the existence of 
 $t$-$(n,k,\lambda)$ designs has been proven if $n$ is sufficiently large   (see Remark \ref{rm:Keevash} below).
\end{remark}

\subsection{Steiner Systems}

Recall that a perfect code with minimum distance $2e+1$ or
an $e$-perfect code in a graph is a subset $C$
of the vertices of the graph such that the balls of radius $e$ centered at the codewords are pairwise disjoint and cover all vertices  of the graph. 
$e$-Perfect codes in hypercubes will be studied  in Section \ref{15.3}.

\begin{claim}
Let $C$ be an $e$-perfect code in $Q^n_2$ and $\bar{0}\in C$. Then
the set $\{u\in C \ |\ \wt(u)=2e+1\}$ is an $S(e+1,2e+1,n)$
design.
\end{claim}
\begin{proof}
Consider an arbitrary vertex of weight $e+1$. By the definition of an
$e$-perfect code, it belongs to exactly one ball centered at a
point $u\in C$, and $\wt(u)\leq e+1 +e=2e+1$ by the triangle inequality. There are no vertices of weight less than $2e+1$ in the code $C$, since
the balls of radius $e$ centered at such a vertex and at the vertex $\bar{0}\in
C$ intersect.  Then $\wt(u)=2e+1$.
\end{proof}

In particular, the set of vertices of weight $3$ in a $1$-perfect binary
code forms a Steiner triple system. The Steiner triple system
contained in the Hamming code can be defined as follows:
$\{\{x,y,x\oplus y\} \ |\ x,y\in Q^t_2\setminus \{\bar{0}\}, x\neq
y\}$, where the length of the Hamming code is $2^t-1$. The indices $x$ and $y$
of the coordinates of vertices in the  Boolean $n$-cube 
can be viewed as elements of the set $Q^t_2\setminus
\{\bar{0}\}$, i.e.,
 the columns of the parity-check matrix of the Hamming code.

\begin{claim}\label{claim:DesRet}
Let $D$ be a combinatorial design with parameters $t-(n,k,\lambda)$.
Then its retract $D'=\{u'\in Q^{n-1}_2 \ | \ (u',1)\in D\}$
is a combinatorial design with parameters
$(t-1)$-$(n-1,k-1,\lambda)$.
\end{claim}
\begin{proof}
Consider the blocks of the design $D$ that contain the $n$-th element.
Any $(t-1)$-element subset not containing the $n$-th element
is contained in $\lambda$ such blocks by the definition of a design. Thus,
removing the $n$-th element from these blocks, we obtain a design with
the required parameters.
\end{proof}

The design $D'$ is called a {\sl shortening} of the original design $D$.

\begin{claim}\label{c:conddesign}
If a combinatorial design with parameters $t$-$(n,k,\lambda)$
exists, then the numbers $\lambda{{n-i} \choose {t-i}}/{{k-i} \choose
{t-i}}$ for $i=0,\dots,t-1$ are integers.
\end{claim}
\begin{proof}
Let $N$ be the number of blocks in a design with parameters
$t$-$(n,k,\lambda)$. Consider a bipartite graph, one part of
which is the blocks of the design, and the other all possible
$t$-element sets. Join each $t$-element set by edges
to the blocks in which it is contained. From the definition of a design it follows that
the number of edges in the graph equals $\lambda{{n} \choose {t}}=N{{k}
\choose {t}}$. 
Applying the same counting argument to the shortening obtained by fixing one point yields the corresponding formula for $i=1$. Repeating this process gives the formulas for all $i=0,\dots,t-1$.
\end{proof}
 The proof of the following theorem is taking from \cite{Rosa}.

\begin{theorem}[Kirkman]\label{th:STS}
A Steiner triple system $S(2,3,n)$ exists if and only
if $n\equiv 1\mod 6$ or $n\equiv 3\mod 6$.
\end{theorem}
\begin{proof}
From Proposition \ref{c:conddesign} it follows that a necessary condition
for the existence of a Steiner triple system of order $n$ is follows. 
The number $n(n-1)$ is divisible by $6$ and the number $n-1$ is divisible by $2$. This condition holds
only if $n\equiv 1\mod 6$ or $n\equiv 3\mod 6$.

Let us prove the existence of Steiner triple systems for $n\equiv 3\mod 6$.
Note that the number $m=n/3$ is odd. By Theorem \ref{th:LSsym},
there exists an idempotent symmetric quasigroup
$f:Q_m^2\rightarrow Q_m$ of order $m$. Consider the $n$-element
set $Q_m\times Q_3$, whose elements we will write in the form
$xi$, where $x\in Q_m$, and $i\in\{0,1,2\}$. Define the set $S$
consisting of triples of elements of $Q_m\times
Q_3$ of two types\\ $\{\{x0,x1,x2\} \ |\ x\in Q_m\}$ and\\
$\{\{xi,yi,f(x,y)(i+1\mod 3)\} \ |\ i\in Q_3, x,y\in Q_m, x\neq
y\}$.\\ Let us check that $S$
is a Steiner triple system. Consider a pair of elements $xi,zj\in Q_m\times
Q_3$. If $x=z$, then the pair is contained in a triple of the first type, and if
$i=j$ then pair is contained in a triple  of the second type. Suppose $x\neq z$ and $i\neq j$. Then $j=i+1\mod 3$
or conversely $i=j+1\mod 3$ (which does not matter, since the pair
$\{xi,zj\}$ is unordered). Without loss of generality, set $j=i+1\mod
3$. From the definition of a quasigroup, there exists a unique $y\in Q_m$
such that $z=f(x,y)$; moreover, $y\neq x$ by the idempotence
of the quasigroup. Then there is a triple of the second type $\{xi,yi,zj\}$.
The triple containing an arbitrary pair of elements is found
uniquely; therefore, $S$ is a Steiner system.

Let us prove the existence of Steiner triple systems for $n\equiv 1\mod 6$.
Note that the number $m=(n-1)/3$ is even. By Theorem \ref{th:LSsym},
there exists a semi-idempotent symmetric quasigroup
$f:Q_m^2\rightarrow Q_m$ of order $m$. Consider the $n$-element
set $(Q_m\times Q_3)\cup\{ \delta\}$, where $\delta$ is a
special additional symbol. Define three types  of triples from the set  $(Q_m\times Q_3)\cup\{ \delta\}$  which constitute the set $S$:\\
$\{\{x0,x1,x2\} \ |\ x\in Q_{m/2}\}$,\\ $\{\{xi,yi,f(x,y)(i+1\mod
3)\} \ |\ i\in Q_3, x,y\in Q_m, x\neq y\}$ and \\
$\{\{(x+m/2)i,x(i+1\mod 3),\delta\} \ |\ i\in Q_3, x\in Q_{m/2}\}$, 
Pairs of elements $\{xi,\delta\}$ are contained in triples of the third type.
The proof that the remaining pairs of elements, with the exception of
pairs\\ $\{(x+m/2)i,(x+m/2)(i+1\mod 3)\}$ and \\ $\{(x+m/2)i,x(i+1\mod 3)\}$, $x\in Q_{m/2}$, are covered by triples of the second type is analogous to
the previous case $n\equiv 3\mod 6$. To prove the
uniqueness of the covering of pairs, it suffices to show that the pairs
$\{(x+m/2)i,x(i+1\mod 3)\}$, $x\in Q_{m/2}$, are not contained in triples of
the second type. While the pairs $\{(x+m/2)i,(x+m/2)(i+1\mod 3)\}$,
$x\in Q_{m/2}$, on the contrary, are contained in triples of the second type.
Indeed, $f(x+m/2,y)\neq x$ for $y\neq x+m/2$, since
$f(x+m/2,x+m/2)= x$ by the semi-idempotence of the quasigroup $f$. However, for
any $x\in Q_{m/2}$ there exists $y\neq x+m/2$ such that
$f(x+m/2,y)=x+m/2$.
\end{proof}

The constructions of Steiner triple systems used in Theorem \ref{th:STS} were devised by
Bose \cite{Bose} and Skolem \cite{Skolem}.

\begin{remark}\label{rm:Keevash}
 Hanani \cite{Hanani} proved the existence of Steiner quadruple
systems for any $n$ satisfying the necessary conditions
$n\equiv 2\mod 6$ or $n\equiv 4\mod 6$ (Proposition
\ref{c:conddesign}). Wilson \cite{Wilson} (for $t=2$) and Keevash
\cite{Keevash14} (for arbitrary $t$)
 proved that for  fixed
$t,k,\lambda$ there exists $n_0\in \mathbb{N}$ such that a combinatorial
design with parameters $t$-$(n,k,\lambda)$ exists for $n\geq n_0$
provided the necessary conditions from Proposition \ref{c:conddesign} hold. There is an alternative proof of the last theorem \cite{Glock}.
\end{remark}

Let $D$ be a Steiner system $S(2,k,n)$ and $M\subset Q_k^k$ be an
{\rm MDS} code with distance $k-1$, containing the subcode
$\{\underbrace{(i,\dots,i)}_{k} \ |\ i=1,\dots,k\}$. The code
distance $k-1$ of the MDS code $M$ guarantees that its vectors, with the
exception of vectors of the form $(i,\dots,i)$, do not
contain identical symbols. 
Indeed, if a codeword contains the same symbol twice, then it is at distance at most $k-2$ from the corresponding constant word.
Consequently, every nonconstant codewords is a permutation of $Q_k$. 
For every block 
 $\{a_1,\dots,a_k\}\in D$ and every nonconstant codeword $(\pi_1,\dots,\pi_k)\in M$, include the ordered $k$-tuple $(a_{\pi_1},\dots,a_{\pi_k})$ in $E[D,M]$.

\begin{claim}\label{e:SteMDS}
The set $C=E[D,M]\cup \{\underbrace{(i,\dots,i)}_{k} \ |\
i=1,\dots,n\}$ is an {\rm MDS} code with minimum distance $k-1$ in the
hypercube $Q^k_n$.\end{claim}
\begin{proof}
The code distance of the set $E[D,M]$ equals $k-1$.
Indeed, since $D$
is a Steiner system, each pair is contained in only one
block, consequently, codewords corresponding to distinct blocks are at distance at least $k-1$.
Since $M$ is an {\rm MDS} code with distance $k-1$, 
permutations of the elements of one block differ in at least $k-1$
positions. Since $E[D,M]$ consists of permutations, the vertices
of the form $(i,\dots,i)$ are at distance $k-1$ from the vertices of
$E[D,M]$.
 The number
of elements in the code equals $n+k(k-1){n \choose 2}/{k \choose 2}=n^2$. By Singleton bound (Proposition \ref{proSingl}) it follows that $C$
is an {\rm MDS} code.
\end{proof}

If $k=3$, one can take the idempotent
Latin square of order $3$ as $M$. The MDS codes with
distance $2$ obtained from a Steiner triple system are graphs of quasigroups, which
are called {\sl totally symmetric}. For $k>3$, the
construction from Proposition \ref{e:SteMDS} yields an MDS code
equivalent to a system of $k-2$ mutually orthogonal Latin squares.

\subsection{Symmetric Designs and Projective Planes}

A {\sl projective plane} is a system of points and lines
satisfying the following properties.

1) Given any two distinct points, there is exactly one line incident with both of them.

2) Given any two distinct lines, there is exactly one point incident with both of them.

3) There are four points such that no line is incident with more than two of them.

It can be shown that in a finite projective plane, each line
contains the same number  $k+1$ of points, where $k$ is the order
of the plane. From the symmetry of points and lines in the definition of a projective plane it follows that
 through each point there pass $k+1$ lines. In total, a
projective plane of order $k$ contains $k^2+k+1$ points,
since each of the $k+1$ lines passing through a fixed
point contains $k$ points, besides the fixed point. From properties
1)-2) it follows that all these points are distinct and there are no other points on the
plane. Projective planes of order $k$ are Steiner
 designs $S(2,k+1,k^2+k+1)$, where the points correspond to elements, and
 the lines to blocks.

\begin{example}
The projective plane of order $2$ (Fano plane) contains $7$
lines $\{1,2,3\},\{1,4,5\},\{1,6,7\},\{2,4,6\},\{2,5,7\},\{3,4,7\},
\{3,5,6\}$.

\begin{minipage}{1\textwidth}
%\vskip15mm
\includegraphics[width=0.5\textwidth]{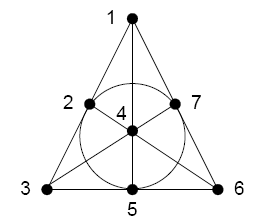}
% \center{\includegraphics[width=0.07\textwidth]{E3.png}}
\end{minipage}

\end{example}

An affine plane of order $k$ can be obtained from a projective
plane of the same order by removing one block and all points
of that block from the projective plane. The reverse operation: adding
to each line of the affine plane one point (the same point for
parallel lines) and an additional line consisting of all the
added points allows one to construct a projective plane   from every affine plane 
of the same order. An affine plane of order
$k$ can be viewed as a design $S(2,k,k^2)$.

From the equality ${k^2+k+1 \choose 2}=(k^2+k+1){k+1 \choose 2}$ it follows that
a projective plane contains as many lines, namely
$k^2+k+1$, as points. Any two lines in a projective plane
intersect in exactly one point, therefore, by reversing the incidence
relation and considering points as blocks and lines as elements
of the design, we obtain a design with the same parameters. Projective
planes that map to each other under this transformation are called
{\sl dual}.

\begin{claim} From any design $D$ of cardinality $v$ with parameters
$2$-$(v,k,\lambda)$, after replacing vertices with blocks, i.e., after
reversing the incidence relation, one obtains a design with the same
parameters.
\end{claim}
\begin{proof} Let $M$ be the incidence matrix of the design $D=\{b_j : j=1,\dots,v\}$, i.e.
\[
M_{v\times v}=(m_{ij})= \left\{
\begin{array}{l}
1 \mbox{{, if}}\ x_i \in b_j\in D;\\
0\mbox{, otherwise}.
\end{array}
\right.
\]
From the definition of a design it follows that each column of the matrix $M$
contains $k$ ones. By Proposition \ref{claim:DesRet} it follows that
each row of the matrix $D$ contains the same number of ones, because all retracts of the design are designs having the same parameters.
Since $D$ is a square matrix, its rows also contain $k$ ones. Consider an
arbitrary pair of rows; by the definition of a design there are exactly
$\lambda$ columns in which both these rows have ones. Let
the numbers $a$ and $b$ satisfy the equation $\lambda a^2 +
2(k-\lambda)ab+(v-2k-\lambda)b^2=0$. Consider the matrix $Q=aM+b(J-M)$, which
is obtained from $M$ by replacing $1$ with $a$ and $0$ with $b$. From the choice of the numbers
$a$ and $b$ it follows that the matrix $Q$ is orthogonal. On the other hand,
an arbitrary orthogonal matrix $Q$ whose elements take
only two values corresponds to the incidence matrix of some
design. It is clear that the matrix $Q^{\mathsf T}$ corresponds to the design obtained
by reversing the incidence relation, and the matrices $Q$ and $Q^{\mathsf T}$
are orthogonal simultaneously.
\end{proof}

$2$-Designs in which the number of blocks equals the number of points are called {\sl
symmetric}.

In \cite{Chowla}, \cite{Bruck} the following conditions for the
    existence of symmetric
    $2$-designs were proved.
\begin{theorem}[Bruck, Ryser, Chowla]
Suppose there exists a symmetric $2$-$(v,k,\lambda)$ design.

If $v$ is even, then $k-\lambda=x^2$ for some integer number $x$.

If $v$ is odd, then there exist integer numbers $x$, $y$, and
$z$ (not all zero simultaneously) such that the equality holds
$$z^2=(k-\lambda)x^2+(-1)^{(v-1)/2}\lambda y^2.$$
\end{theorem}

This theorem yields an infinite series of arithmetically
    admissible design parameters for which designs nevertheless do not exist.
 Note that the statement on the existence of $2$-$(v,k,\lambda)$ with arithmetically admissible parameters \cite{Wilson} becomes
        true only for sufficiently large $v$, for fixed $k$ and $\lambda$.

\begin{remark}
From the Bruck--Ryser--Chowla theorem one can prove that there do not
exist finite projective planes of orders 14, 21, and 22.
The non-existence of a projective plane of order 6 follows from the
absence of a pair of orthogonal Latin squares of order 6.
 The non-existence of a projective plane of order 10 was proved by using computer \cite{Lam}. 
Exhaustive computer classification \cite{LamKol} showed
 that for each of the orders 3, 4, 5, 7, and 8 there exists a
unique projective plane, corresponding to the complete system of
linear orthogonal Latin squares; and for order 9 there are
four projective planes, two of which are dual.
\end{remark}

\subsection{Hadamard Matrices}\label{10.3}

Consider the matrix $H_{q,n}$ of size $q^n\times q^n$ whose rows are the characters $\phi_z(x)$ of ${\mathbb Z}^n_q$, i.e.\ $H_{q,n}(z,x)=\phi_z(x)$. From the orthogonality properties of characters (Proposition \ref{c:charact_orthonorm_basis}) it follows that $H_{q,n}H_{q,n}^*=mI_m$, where $m=q^n$. 

\begin{claim}
Consider a function $f:Q^n_q\rightarrow {\mathbb C}$ as a vector of $ {\mathbb V}({\mathbb Z}^n_q)$.  Then 
$H_{q,n}f=\widehat{f}$. 
\end{claim}
Proof follows directly from definitions.

%\begin{definition}
  An $m\times m$ matrix $H$ is called a {\sl Hadamard matrix} if it satisfies
  $HH^{\mathrm{T}}=mI_m$ and consists only of $\pm 1$.
  Matrices consisting of entries with absolute value $1$ and satisfying
  $HH^{*}=mI_m$ are called {\sl generalized Hadamard matrices}. If all entries  are
    powers of $\xi=e^{2\pi i/q}$, then the matrix is called {\sl Butson-type}  Hadamard matrix. 
%\end{definition}

Obviously, a Hadamard matrix becomes orthogonal (unitary) after multiplying by $\frac{1}{\sqrt m}$. From the definitions of the matrices and characters it is clear that the matrices $H_{q,n}$ are symmetric.
The $2^n\times 2^n$ matrix $H_{2,n}$ is called the {\sl Hadamard--Sylvester matrix}.

It is clear that permuting the rows or columns of a generalized Hadamard matrix does not take it out of the class of such matrices. Moreover, multiplying a column or a row of a Butson-type Hadamard matrix consisting of powers of $\xi=e^{2\pi i/q}$ by $\xi^j$ does not take the matrix out of this class.
Generalized Hadamard matrices obtained from one another by such transformations are called {\sl equivalent}.

\begin{claim}
If $H_{n_1}$ and $H_{n_2}$ are generalized Hadamard matrices, then $H_{n_1}\otimes H_{n_2}$ is a generalized Hadamard matrix.
\end{claim}
\begin{proof}
$(H_{n_1}\otimes H_{n_2})(H_{n_1}\otimes H_{n_2})^*=(H_{n_1}\otimes
H_{n_2})(H^*_{n_1}\otimes H^*_{n_2})= (H_{n_1}H_{n_1}^*)\otimes
(H_{n_2}H_{n_2}^*)=n_1I_{n_1}\otimes n_2I_{n_2}=n_1n_2I_{n_1n_2}$.
\end{proof}

\begin{example}
The Hadamard-–Sylvester matrix
 $$H_{2,3}=\begin{pmatrix}
1& 1& 1 &1 &1 &1 &1 & 1\\
1 &-1 &1 &-1 &1 &-1 &1 &-1 \\
1 & 1& -1 &-1 &1 & 1& -1& -1 \\
1 &-1 &-1 &1 &1 &-1 &-1 &1\\
1& 1& 1 &1 &-1 &-1 &-1 & -1\\
1 &-1 &1 &-1 &-1 &1 &-1 &1 \\
1 & 1& -1 &-1 &-1 & -1& 1& 1 \\
1 &-1 &-1 &1    &-1 &1 &1 &-1
\end{pmatrix}.$$
\end{example}

\begin{lemma}\label{l:matrix}
Let $f: Q^n_q\rightarrow \mathbb{C}$ be a function.
Define the $q^n\times q^n$ matrix $A$ by 
$A_f(x,y)=f(x-y)$. Then $A_fA_g=A_{f*g}$.
\end{lemma}
\begin{proof}
We have  $$A_fA_g(x,z)=\sum\limits_{y\in
Q^n_q}f(x-y){g(y-z)}= \sum\limits_{u\in
Q^n_q}f(u)g(x-z-u)=(f*g)(x-z)=A_{f*{g}}(x,z).$$
\end{proof}

\begin{claim}\label{cl:Hadam}
A function $b$ is a bent function if and only if
$H(x,y)=\xi^{b(x-y)}$ is a generalized Hadamard matrix.
\end{claim}
\begin{proof}
We have $H^*(x,y)=\xi^{-b(y-x)}=\xi^{\tilde{b}(x-y)}$. Then from Lemma~\ref{l:matrix} we obtain the equality
$HH^*(x,y)=[\xi^b*\xi^{\tilde{b}}](x-y)$. By Proposition~\ref{kriteriibent}, the equation $HH^*=NI_N$, where $N=q^n$, is equivalent to $b$ being bent.
\end{proof}

For $q=2$, the matrix $H(x,y)=(-1)^{b(x-y)}$ is a Hadamard matrix.

By multiplying rows and columns by $-1$, any Hadamard matrix can be transformed so that its first column and first row consist entirely of ones. Such matrix is called  normalized Hadamard matrix. Then in the other rows and columns there will be exactly half of the entries equal to $1$ and half equal to $-1$. Consider the $(4m+3)\times (4m+3)$ matrix $A_m$ obtained from the normalized Hadamard matrix by deleting the first row and the first column and replacing $-1$ by $0$. The rows of the matrix $A_m$ can be regarded as binary vectors of length $4m+3$ and weight $2m+1$.

\begin{claim}
The set of rows of the matrix $A_m$ forms a combinatorial design with parameters $2$-$(4m+3,2m+1,m)$.
\end{claim}
\begin{proof}
It suffices to verify that for every pair of columns of the matrix $A_m$ there are exactly $m$ rows in which both columns contain a~$1$. Indeed, in the original Hadamard matrix each column (except the first) contained exactly $2m$ entries equal to $1$ and $2m$ entries equal to $-1$. Then $s(1,1)+s(1,-1)=s(-1,1)+s(-1,-1)=2m+2$, where $s(\delta_1,\delta_2)$ is the number of pairs $(\delta_1,\delta_2)$ in the two columns. The orthogonality of the Hadamard matrix implies $s(1,1)+s(-1,-1)=s(1,-1)+s(-1,1)$. Hence $s(1,1)=m+1$. Since the matrix is normalized, every non-first column begins with $1$, then by deleting the first row from consideration yields the required result.
\end{proof}

As can be seen from the construction, the number of blocks in a design with parameters $2-(4m+3,2m+1,m)$ equals the number of points, i.e.\ $4m+3$. Therefore such a design is symmetric. It is easy to see that if we form a matrix from the blocks of a design with parameters $2$-$(4m+3,2m+1,m)$, replace $0$ by $-1$ in it, and add a row and a column of ones, we obtain a Hadamard matrix. The orthogonality of the columns of the resulting matrix follows from the parameters of the design.

It is not difficult to show that $m\times m$ Hadamard matrices can exist only when $m$ is divisible by $4$
(except for orders $m=1, 2$, which are trivial cases). 
Indeed, let $H$ be a normalized Hadamard  matrix of size $m\times m$. Consider the  second and third rows of $H$ as a collection of $m$ pairs. Let  $x_1$ be the number of pairs  $(1,1)$, $x_2$ be the number of pairs  $(-1,1)$, $x_3$ be the number of pairs  $(1,-1)$, $x_4$ be the number of pairs  $(-1,-1)$. The orthogonality condition of the first three rows yields the following system of linear equations:
\begin{equation}\label{eq:H4}
  \left\{
    \begin{aligned}
      x_1 + x_2 + x_3 + x_4 &= m,\\
      x_1 + x_2 - x_3 - x_4 &= 0,\\
       x_1 - x_2 + x_3- x_4 &= 0,\\
        x_1 - x_2 - x_3 + x_4 &= 0.
    \end{aligned}
  \right.
\end{equation}
The matrix of the system (\ref{eq:H4}) is a  Hadamard--Sylvester matrix of order $4$. The unique solution of  (\ref{eq:H4}) is $x_1=x_2=x_3=x_4=m/4$. Therefore, $m/4$ is an integer.
The famous conjecture states that this condition is sufficient for the existence of Hadamard matrices.
A detailed survey of the theory of Hadamard matrices is available in \cite{Horadam}. \label{Hprob}

\subsection{Problems}

\begin{exercise} Prove that  designs with parameters
$t$-$(n,k,\lambda)$ and designs with parameters
$i$-$(n,k,\lambda_i)$ for $i=1,\dots,t-1$ and $\lambda_i=\lambda{{n-i}
\choose {t-i}}/{{k-i} \choose {t-i}}$ determine each other uniquely.
\end{exercise}

\begin{exercise} Prove that every design with parameters
$2$-$(n,k,\lambda)$ generate a design with
parameters $2$-$(n,n-k,\lambda')$, where $\lambda'=\lambda{{n}
\choose {2}}/{{k} \choose {2}}- 2\lambda\frac{n-1}{k-1}+ \lambda$.
\end{exercise}

\begin{exercise}
Prove that a maximum clique in the graph $J(n,k)$ is
Delsarte.
\end{exercise}

\begin{exercise}
Prove that a design with parameters $S(k-2,k,n)$ is a completely
regular code in the graph $J(n,k)$.
\end{exercise}

\begin{exercise}
In a construction analogous to the construction of the set $E(D,M)$ (see
Proposition \ref{e:SteMDS}), take as
 $D$ a Steiner quadruple system of order $n$, and as $M$ the set
  of even permutations. Prove that the
set $E[D,M]\cup \{(i,i,j,j), (i,j,i,j), (i,j,j,i) \ |\ \
i,j=1,\dots,n\}$ is an {\rm MDS} code with minimum  distance $3$ in the
hypercube $Q^4_n$.
\end{exercise}

\begin{exercise}
Prove that an MDS code in $Q^{k+1}_k$ with distance $k$ does not contain
vertices at distance $k+1$. Prove that any such code generates
a design with parameters $2$-$(k^2,k,1)$.
\end{exercise}

\begin{exercise}
Let $M$ be incidence matrix of a design $D$ with parameters
$2$-$(n,k,\lambda)$. Prove the equality $MM^{\rm
T}=\lambda(J_n+\frac{n-k}{k-1}I_n)$.
\end{exercise}

\begin{exercise}
Consider a Hadamard matrix $A$ of order $n$ whose first row
and first column consist of $1$'s. Replace all occurrences of $-1$'s in the matrices $A$ and
$-A$ by $0$. Prove that the orthogonal array consisting of $2n$
rows of these matrices has strength $2$. Prove that if $n=2^k$, then
an orthogonal array of strength $2$ of smaller cardinality does not exist.
\end{exercise}

\section{Perfect Codes and Tilings}

\subsection{Tilings of Groups}

%\begin{definition}
A  {\sl factorization} or {\sl tiling} of a group $(Gr,\circ)$ is a
pair $(C,B)$ of subsets $C\subset Gr$ and $B\subset Gr$ such that,
for every $x\in Gr$, there exists a unique pair of elements
$c\in C$ and $b\in B$ such that
$x=c\circ b$,
i.e.,
$Gr=C\circ B$.
%\end{definition}

If $C\leq Gr$ is a subgroup, then $B$ can be chosen to contain one
representative from each coset of $C$ in $Gr$. It is easy to see that,
for a finite group, uniqueness of the representation implies
$|Gr|=|C||B|$.

If $(C,B)$ is a tiling, then $(C,B\circ a)$ is also a tiling for every
$a\in Gr$. Therefore, we may always assume that the identity element
$\varepsilon$ of the group belongs to $C\cap B$. In this case, it is
clear that
$C\cap B=\{\varepsilon\}$.

\begin{claim}
A pair $(C,B)$ is a tiling of the group $(Gr,\circ)$ if and only if
$(C^{-1}\circ C)\cap(B\circ B^{-1})=\{\varepsilon\}$ 
and
$|Gr|=|C||B|$.
\end{claim}

\begin{proof}
The equality
$c_1\circ b_1=c_2\circ b_2
$
is equivalent to
$
c_2^{-1}\circ c_1=b_2\circ b_1^{-1}$.
Thus, uniqueness of the representation is equivalent to
\[
(C^{-1}\circ C)\cap(B\circ B^{-1})=\{\varepsilon\}.
\]
Under this condition, all elements of $C\circ B$ are distinct, and
their number is $|C||B|=|Gr|$. Hence every element of $Gr$ has a
representation of the form $c\circ b$.
\end{proof}

In the case of an abelian group, this condition takes the more
symmetric form
\[
(C-C)\cap(B-B)=\{0\}.
\]

Let $A$ be the connection set of the Cayley graph
$G=\operatorname{Cay}(Gr,A)$. Consider the unit ball
$B=A\cup\{\varepsilon\}$.
It is easy to see that a subset $C\subset Gr$ is a $1$-perfect code in
$G$ if and only if $(C,B)$ is a tiling of the group $Gr$.

\begin{example}
Let $C\subset Q_q^n$ be a Hamming code and let
\[
B=\{x\in Q_q^n\mid \operatorname{wt}(x)\leq1\}
\]
be the unit ball. Then $(C,B)$ is a tiling when
$n=\frac{q^t-1}{q-1}$.
\end{example}

\begin{claim}\label{cl:tiling}
A pair $(C,B)$ is a tiling of an abelian group $K$ if and only if
$
{\mathbf 1}_C*{\mathbf 1}_B={\bf 1}$,
or, equivalently,
$\widehat{{\mathbf 1}_C}\cdot\widehat{{\mathbf 1}_B}=\frac{\delta}{|K|^{1/2}}={\bf 1}_{\{0\}}$.
\end{claim}

\begin{proof}
The first equality follows directly from the definition of a tiling.
The second follows from Proposition \ref{c:convolution} and the
equality
$\widehat{\bf 1}=\delta$.
\end{proof}

As a generalization of tilings of an abelian group $K$, we can consider
pairs of functions $(f,g)$, $f,g:K\rightarrow\mathbb C$, satisfying
\[
f*g=\alpha{\bf 1}
\]
for some $\alpha\in\mathbb C$. Similarly to Proposition \ref{cl:tiling}, a
criterion for a pair $(f,g)$ to be a {\sl multiple tiling} is
\[
\supp(\widehat f)\cap\supp(\widehat g)=\{0\}.
\]

A tiling of a group can be used for compact storage of the group.
To perform the group operation, its Cayley table is required, which
has size $|Gr|^2$. If a tiling $(C,B)$ of $Gr$ is available, it is
sufficient to store the multiplication tables of the elements of $Gr$
by the elements of $C$ and by the elements of $B$, i.e., two tables
of sizes $|Gr||C|$ and $|Gr||B|$.

A tiling $(C,B)$ of the group $\mathbb Z_2^k$ can be used to construct
$1$-perfect codes in $Q_2^n$. We regard the elements of $\mathbb Z_2^k$
as binary vectors of length $k$. Let $\bar 0\in C$ and $|C|=n+1$.
Taking the nonzero elements of $C$ as columns, construct a matrix $H$
of size $k\times n$.

\begin{claim}[\cite{Cohen}]
The union, over all $b\in B$, of the solution sets of the systems
$
Hx=b
$
over $GF(2)$ is a $1$-perfect code in $Q_2^n$.
\end{claim}

\begin{proof}
Let $x\in Q_2^n$. We first prove that there is a codeword at distance
at most $1$ from $x$. Consider
\[
z=x_1c^1+\cdots+x_nc^n\in\mathbb Z_2^k,
\]
where $c^1,\ldots,c^n$ are the columns of $H$. Since $(C,B)$ is a
tiling, there exist $c\in C$ and $b\in B$ such that
$z=c+b$.
If $c=\bar0$, then $x$ is a codeword. Otherwise, $c$ is one of the
columns of $H$, say $c^i$, and hence
$
z+c^i=b$.

Therefore, the vector
$
x+(0\ldots0\underset{i}{1}0\ldots0)
$
is a codeword.

It remains to prove that the distance between distinct codewords is
at least $3$. Suppose that $x$ and $y$ are two codewords that differ
only in positions $i$ and $j$. Then
\[
x_1c^1+\cdots+x_nc^n=b^1,\qquad
y_1c^1+\cdots+y_nc^n=b^2
\]
for some $b^1,b^2\in B$, and hence
$
c^i+c^j=b^1+b^2$.
Since $c^i,c^j\in C$, this contradicts the definition of a tiling.
If $x$ and $y$ differ in only one position, we obtain an analogous
contradiction using $\bar0\in C$.
Thus, the minimum distance of the code is at least $3$, and every
vector is at distance at most $1$ from a codeword. Hence the code is
$1$-perfect.
\end{proof}

\subsection{Tilings of Graphs and the Code-Anticode bound}

The notion of a tiling can be extended to vertex-transitive graphs.

%\begin{definition}
A {\sl tiling} of a vertex-transitive graph $G$ is a pair $(A,\Pi)$,
where $A\subset V(G)$ and $\Pi\subset  \operatorname{Aut}(G)$, such that
\[
V(G)=\bigsqcup_{\varphi\in \Pi}\varphi(A).
\]
%\end{definition}

Clearly, the conditions
$|\Pi||A|=|V(G)|$
and
\[
\varphi_1(A)\cap\varphi_2(A)=\varnothing
\qquad\text{for all }\varphi_1\neq\varphi_2,\quad
\varphi_1,\varphi_2\in \Pi,
\]
are necessary and sufficient for $(A,\Pi)$ to be a tiling.

Consider a hypergraph $\Gamma$ whose set of vertices is $V(G)$ and whose hyperedges  are $\pi(A)$, where $\pi\in \operatorname{Aut}(G)$. By definition,  the set $\{\pi(A) : \pi\in \Pi\}$ is a perfect matching of  $\Gamma$.

A subset $A\subset V(G)$ whose internal distance does not exceed a given maximum diameter $d$ is called an {\sl anticode}. If $G$ is the Cayley graph of an abelian group $K$,  then we have $|A||C|\leq |K|$ for every code $C$ with minimal distance $d+1$.  This inequality is called a {\sl code-anticode bound}.
It is easy to see that if the  code-anticode bound is attained for $C$ and  $A$ then $(C,A)$ is a
tiling of $K$.  

The code-anticode bound is generalized to  vertex-transitive graphs. 

\begin{claim}
Let $G$ be a  vertex-transitive graph and the internal distance of  $A\subset V(G)$ is at most  $d$. Then, for every code $C\subset V(G)$ with minimum distance $d+1$, the bound $|A||C|\leq |V(G)|$ holds.
\end{claim}
\begin{proof}
Let $Aut(G)$ cover $V(G)$ $\alpha=| \operatorname{Aut}(G)|/|V(G)|$ times. 
The proof is based on double-counting the number of pairs
$(\pi, c)$ such that $\pi \in  \operatorname{Aut}(G)$, $c \in C$, and $c \in \pi(A)$. On one hand, this number is not greater than $| \operatorname{Aut}(G)|$, but on the other hand, it is equal to $\alpha|C||A|$, because for every pair $(c,a)$ there are $\alpha$ distinct $\pi\in \operatorname{Aut}(G)$ such that $c=\pi(a)$.  Consequently, $|\operatorname{Aut}(G)|\geq \alpha|C||A|$.
\end{proof}

It follows from the proof  that if the code-anticode bound is reached, then $(A, \operatorname{Aut}(G))$ is a multiple tiling.

\subsection{Perfect Codes with Minimum Distance Greater Than $3$}\label{PC}

%\begin{definition}

A {\sl perfect code with minimum distance $2e+1$} (an
{\sl $e$-perfect code}) in a graph $G$ is a set $C\subset V(G)$ such
that the balls $B(x,e)$ centered at vertices $x\in C$ and of radius
$e$ form a partition of the vertices of the graph, i.e.,
\[
V(G)=\bigsqcup_{x\in C}B(x,e).
\]

%\end{definition}

A perfect code with minimum distance $2e+1$ in a Cayley graph
$G=Cay(Gr,A)$ is equivalent to a tiling $(C,A^e\cup\{\varepsilon\})$,
where
\[
A^e=\{a_1\circ\dots\circ a_k\mid a_i\in A,\ k=1,\dots,e\},
\]
as well as to a $1$-perfect code in the graph
$G=Cay(Gr,A^e)$.

In a distance-regular graph with parameters
$(b_0,\dots,b_{d-1},c_1,\dots,c_d)$, an $e$-perfect code $C$
induces a perfect coloring with quotient matrix
\[
S_C=
\begin{pmatrix}
a_0 & b_0 & 0 & 0& \dots & 0\\
c_1 & a_1 & b_1 & 0 & \dots & 0\\
0  & c_2 & a_2 & b_2 &  \dots & 0 \\
\dots & \dots & \dots & \dots & \dots & \dots \\
0 & \dots & 0 & 0 & c_e & a_e+b_e
\end{pmatrix},
\]
where the vertices of the code $C$ are assigned the first color,
and the colors of the remaining vertices are determined by their
minimum distance from  the code vertices.

It is easy to see that a coloring with the quotient matrix $S_C$ is a distance coloring with respect to the code $C$. Then every $e$-perfect code is a completely regular code by definition.

\begin{claim}[Lloyd's condition]\label{cl:LC}

Let $C$ be an $e$-perfect code  in a
distance-regular graph $G$. Then, for every eigenvalue $\lambda$ of
the matrix $S_C$, except for $\lambda_0=\deg G$, we have
\[
1+P_1(\lambda)+\dots+P_e(\lambda)=0,
\]
where $P_t(x)$ is the  Krawtchouk-type polynomial associated with the graph $G$
(see Section \ref{4.1}).

\end{claim}

\begin{proof}

Let $f$ be the characteristic function of the code $C$. By the
definition of an $e$-perfect code, we have
\[
If+Mf+M_2f+\dots+M_ef=\mathbf{1},
\]
where $M_i=P_i(M)$ is the distance-$i$ adjacency matrix of the
graph $G$. Let $\lambda_i$, $i=0,\dots,e$, be the eigenvalues of
the quotient matrix $S_C$.

By Corollary \ref{p:perf_color_properties1}, there exist suitable
coefficients $\alpha_i$ such that
\[
f=\sum_{i=0}^e\alpha_i\varphi_i,
\]
where $\varphi_i$ is an eigenfunction of the graph $G$ with
eigenvalue $\lambda_i$.

We show that $\alpha_i\neq0$ for every $i=0,\dots,e$. By
Corollary \ref{p:perf_color_properties}, for every eigenvalue of
the matrix $S_C$, an eigenfunction of the graph with the same
eigenvalue can be obtained as a linear combination of the
characteristic functions of the colors. The characteristic
function of the color at distance $i$ from the code is $M_if$.
Therefore, all characteristic functions of the colors are linear
combinations of the same eigenfunctions of the graph as the first
color.

By Proposition \ref{c:color_distr_proportion}, if $\alpha_i=0$ for some $i=0,\dots,e$, then all
characteristic functions of the colors would be orthogonal to the
corresponding eigenspace of the matrix $M$, which is a
contradiction.

Thus, by Proposition  \ref{c:WDB_bound0},
\[
\mathbf{1}
=If+Mf+M_2f+\dots+M_ef
=\sum_{i=0}^e
\alpha_i\bigl(1+P_1(\lambda_i)+\dots+P_e(\lambda_i)\bigr)\varphi_i.
\]
Since $\mathbf{1}$ is an eigenfunction corresponding to the
eigenvalue $\lambda_0$, we obtain
\[
1+P_1(\lambda_i)+\dots+P_e(\lambda_i)=0,
\qquad i=1,\dots,e.
\]

\end{proof}

If a distance-regular graph is a Cayley graph of an abelian group,
then Lloyd's condition is equivalent to Proposition  \ref{cl:tiling}.
Indeed,
\[
\mathbf{1}=If+Mf+M_2f+\dots+M_ef={\mathbf 1}_B*f,
\]
where $f={\mathbf 1}_C$ and $B$ is the ball of radius $e$ centered at the
zero element.

\begin{corollary}\label{cor:Lloyd}

The matrix $S_C$ has no repeated eigenvalues, and
\[
\det(S_C-xI)
=
(-1)^{e+1}
\left(\prod_{t=1}^e c_t\right)
(x-\lambda_0)
\sum_{t=0}^eP_t(x).
\]

\end{corollary}

\begin{proof}

Since the matrix $S_C-xI$ is tridiagonal with positive
off-diagonal entries, its rank is at least $e$. Consequently, the
matrix $S_C$ has no repeated eigenvalues.

By Lloyd's condition, the polynomials
\[
\det(S_C-xI)
\quad\text{and}\quad
(x-\lambda_0)\sum_{t=0}^eP_t(x)
\]
have the same roots. By (\ref{eqST0}), the leading coefficient of
the polynomial $P_e(x)$ is
\[
\frac{1}{\prod_{t=1}^e c_t}.
\]
The leading coefficient of the polynomial
$\det(S_C-xI)$ is $(-1)^{e+1}$.
This proves the stated equality.

\end{proof}

\subsection{Perfect Codes in Hamming Graphs}\label{15.3}

\begin{claim}\label{cor:Lloydcub}

Let $C\subset Q^n_q$ be an $e$-perfect code.
Then
\[
(n-\lambda_1)\cdots(n-\lambda_e)=\frac{e!q^n}{|C|},
\]
where $\lambda_i$, $i=1,\dots,e$, are all the eigenvalues of the
parameter matrix $S_C$ except for $\lambda_0=(q-1)n$.

\end{claim}

\begin{proof}

The degrees of the vertices of the $n$-dimensional Hamming graph are
$(q-1)n$, so $(q-1)n$ is the largest eigenvalue of the matrix $S_C$.
For Hamming graphs, we have $c_i=i$. As was noted in Corollary
\ref{cor:Lloyd}, the matrix $S_C$ has no repeated eigenvalues.
Therefore,
\[
\det(S_C-xI)=(-1)^{e+1}(x-\lambda_0)\cdots(x-\lambda_e).
\]
On both sides of the equality
\[
\det(S_C-xI)
=
(-1)^{e+1}e!(x-(q-1)n)
\sum_{t=0}^eP_t(x),
\]
which follows from Corollary \ref{cor:Lloyd}, we cancel the factor
$(x-(q-1)n)$. Here $P_t=P_t[n,q]$ are the $q$-ary Krawtchouk
polynomials. By Proposition \ref{c:WDB_bound0},
$
\sum_{t=0}^eP_t(n)$
is the size of a ball of radius $e$.

Since $C\subset Q^n_q$ is an $e$-perfect code, the balls of radius $e$ centered at its codewords partition
$Q^n_q$. Hence
$
|C|\sum_{t=0}^eP_t(n)=q^n$.
The desired equality follows.

\end{proof}

Let us consider the Boolean hypercube  in more detail.

\begin{claim}\label{cl:spher}

Let $A_r\subset Q^n_2$ be the sphere of radius $r$ centered at
$\bar 0$. For any fixed $r\in\mathbb{N}$, and sufficiently large $n$,
we have
$
\widehat{{\mathbf 1}_{A_r}}(z)>0
$
for all $z$ such that $\wt(z)<n/2$.

\end{claim}

\begin{proof}

Let $\wt(z)=m$, where $m$ may depend on $n$. We have
\[
2^{n/2}\widehat{{\mathbf 1}_{A_r}}(z)
=
\sum_{\wt(x)=r}(-1)^{\langle x,z\rangle}
=
\sum_{i=0}^r
(-1)^i
{m\choose i}
{{n-m}\choose {r-i}}.
\]

If
$2mr+m+r\leq n$,
then
\[
{m\choose 1}{{n-m}\choose {r-1}}
\leq
\frac12{{n-m}\choose r}.
\]
Therefore,
\[
\sum_{i=0}^r
(-1)^i
{m\choose i}
{{n-m}\choose {r-i}}
>
\frac12{{n-m}\choose r}.
\]

If
$
2mr+m+r>n$,
then $m\to\infty$ as $n\to\infty$. Applying
\[
m(m-1)\cdots(m-i+1)=m^i(1+o(1))
\]
as $m\to\infty$, together with the binomial theorem, we obtain
\[
\begin{aligned}
\sum_{i=0}^r
(-1)^i
{m\choose i}
{{n-m}\choose {r-i}}
&=
\sum_{i=0}^r
(-1)^i
\frac{m^i(n-m)^{r-i}}{i!(r-i)!}(1+o(1))\\
&=
\frac{(n-2m)^r(1+o(1))}{r!}.
\end{aligned}
\]
Since $m<n/2$, the last expression is positive for sufficiently
large $n$. This proves the proposition.

\end{proof}

By Proposition \ref{c:WDB_bound0},
$
2^{n/2}\widehat{{\mathbf 1}_{A_r}}(z)
=
P_r[n,2](\lambda_{\wt(z)})$. Recall that $\lambda_{\wt(z)}= n-2\wt(z)$.
Consequently, for sufficiently large $n$, the Krawtchouk polynomials
$P_r[n,2](x)$ are nonnegative at
$x=n,n-2,\dots,2$ 
or at
$x=n,n-2,\dots,1$.
Actually,
it is possible to prove that $P_r[n,2](x)$ is nonnegative throughout the interval
$x\in[1,n]$ for $n>n_0(r)$.

\begin{remark}\label{lastremark}

If $r$ grows sufficiently slowly (in particular, if
$2^{r/2}\leq n$), Proposition \ref{cl:spher} remains valid. Since
$
\widehat{{\mathbf 1}_{B_r}}(z)>0$
for $\wt(z)<n/2$ and sufficiently large $n$, Proposition \ref{cl:tiling}
implies that
$
\widehat{{\mathbf 1}_C}(z)=0$
for $0<\wt(z)<n/2$, where $C$ is an arbitrary $r$-perfect code.
Thus, the characteristic function ${\mathbf 1}_C$ is orthogonal to the
eigenspaces of the Hamming graph corresponding to positive
eigenvalues (except for $n$).
As follows from the proof of Proposition \ref{cl:LC} (Lloyd's condition), the decomposition
of the characteristic function ${\mathbf 1}_C$ as a linear combination of
eigenfunctions contains nonzero components corresponding to every
eigenvalue of the matrix $S_C$. Therefore, the existence of an  $r$-perfect code $C$ implies that  $S_C$ has no positive eigenvalues other than $n$,   for sufficiently large
$n$.

\end{remark}

$1$-Perfect linear Hamming codes in $Q^n_q$ were
constructed for prime power $q$ and  $n=\frac{q^t-1}{q-1}$  (see Proposition \ref{c:codeHam}). In addition, in hypercubes $Q_q^n$ of odd dimension
$n=2s+1$, there are perfect codes with maximal
distance $s$, consisting of all vertices of the form
$(a,\dots,a)$, $a\in Q_q$. A
one-point set in any graph is formally a perfect code with covering
radius equal to (or greater than) the diameter of the graph, but this
case is usually not taken into account when listing perfect codes. Next we prove that the number of other parameters of perfect codes in $Q^n_2$  is finite. Our proof follows the book
\cite{Lint}.

\begin{theorem}[Tiet\"av\"ainen \cite{Tiet}, Zinoviev and Leontiev
\cite{ZL72}]\label{th:percod}

There are only finitely many dimensions $n$ for which the binary
hypercube $Q^n_2$ contains an $e$-perfect code with
$1<e<\frac{n-1}{2}$.

\end{theorem}

\begin{proof}

From Proposition \ref{c:WDB_bound0} and the formula
$\lambda_z=n-2\wt(z)$ for the eigenvalues of  the
Boolean  hypercube, we straightforwardly obtain
$P_t[n,2](-n)=P_t(-n)=(-1)^t{n\choose t}$ and
$P_t(2-n)=(-1)^t\left({n-1\choose t}-
{n-1\choose t-1}\right)$.

We use the well-known combinatorial identity (see Problem
\ref{exer})
\[
2n\sum\limits_{t=0}^s(-1)^{t}{n \choose t}
=(-1)^{s}\left((s+1){n \choose {s+1}}
+(n-s){n \choose {s}}\right)
\]
to calculate $e!\sum\limits_{t=0}^eP_t(-n)$ and
$e!\sum\limits_{t=0}^eP_t(2-n)$. For $e\geq 2$, we obtain
\[
e!\sum\limits_{t=0}^eP_t(-n)
=(-1)^e(n-1)\cdots(n-e),
\]
\[
e!\sum\limits_{t=0}^eP_t(2-n)
=(-1)^e(n-2)\cdots(n-e-1)
-(-1)^ee(n-2)\cdots(n-e)
\]
\[
=(-1)^e(n-2)\cdots(n-e)(n-2e-1).
\]

Suppose that there exists an $e$-perfect code $C\subset Q^n_2$, and
let $\lambda_i$, $i=0,\dots,e$, be the eigenvalues of the quotient
matrix $S_C$.

Then, by Corollary \ref{cor:Lloyd}, taking into account that
$\lambda_0=n$, we have
\[
(\lambda_1-x)\cdots(\lambda_e-x)
=(-1)^e e!\sum\limits_{t=0}^eP_t(x).
\]

Substituting $x=-n$ and $x=2-n$, respectively, gives
\begin{equation}\label{eq:percod0}
(n+\lambda_1)\cdots(n+\lambda_e)
=(n-1)\cdots(n-e),
\end{equation}
and
\begin{equation}\label{eq:percod}
(n+\lambda_1-2)\cdots(n+\lambda_e-2)
=(n-2)\cdots(n-e)(n-2e-1).
\end{equation}

Both sides  of  (\ref{eq:percod0}) are nonzero. Both sides of   (\ref{eq:percod})  can be zero only
when $n=2e+1$. Define $\Delta(n)$ by $n=2^{\Delta(n)}m$, where $m$ is
odd. Among any $k$ consecutive integers, at most every second integer
is divisible by $2$, every fourth integer by $4$, and so on. Therefore,
\begin{equation}\label{eq:percod1}
\Delta((m+1)\cdots(m+k))
\leq k+\max\limits_{i=1,\dots,k}\Delta(m+i).
\end{equation}

If $C$ is an $e$-perfect code, then the eigenvalues of the quotient
matrix $S_C$ are eigenvalues of the adjacency matrix of the Boolean
hypercube (see Corollary \ref{p:perf_color_properties}). Since
$\lambda_i=n-2t$ for some $t=1,\dots,n$, we have
$\Delta(n+\lambda_i)\geq 1$ for $i=1,\dots,e$. Moreover,
\[
\Delta((n+\lambda_i)(n+\lambda_i-2))\geq 3,
\]
since one of the factors is divisible by $4$.

Multiplying the left- and right-hand sides of (\ref{eq:percod0}) and
(\ref{eq:percod}), and applying inequality (\ref{eq:percod1}), we
obtain
\[
\begin{aligned}
3e
&\leq
\Delta((n+\lambda_1)\cdots(n+\lambda_e)
(n+\lambda_1-2)\cdots(n+\lambda_e-2))\\
&=
\Delta((n-2e-1)(n-1)(n-2)^2\cdots(n-e)^2)\\
&\leq
\Delta(n-2e-1)+2e+
\max\limits_{i=1,\dots,e}\Delta(n-i).
\end{aligned}
\]
Hence
$e\leq
\max\limits_{i=1,\dots,e}\Delta(n-i)+\Delta(n-2e-1)$,
and therefore $n\geq 2^{e/2}$.

For $q=2$ the matrix $S_C$ contains coefficients $a_i=0$,
$i=0,\dots,e$, and $b_e=n-e$. Hence the sum of the eigenvalues of
$S_C$ is
$\tr(S_C)=n-e$.
Thus
\[
\sum\limits_{i=1}^e\lambda_i=-e,
\]
and all the numbers $\lambda_i$ are distinct (see Corollary
\ref{cor:Lloyd}).

Therefore, for $e>1$, at least one of them is positive. The 
inequality   $n\geq 2^{e/2}$ and Remark \ref{lastremark} imply that, for sufficiently
large $n$, there can be no $e$-perfect codes in the Boolean
$n$-cube, except for the cases $e=1$ and
$e=\frac{n-1}{2}$.

\end{proof}

Thus, in order to determine all possible parameters of
$e$-perfect codes in Boolean $n$-cubes, it suffices to consider
finitely many values of $n$ and $e$. The following conditions are
necessary for the existence of an $e$-perfect code: the size of a
ball,
\[
\sum\limits_{t=0}^e{n\choose t},
\]
divides the number of vertices of the graph (i.e., the latter is equal
to $2^n$); the eigenvalues of the quotient matrix $S_C$ are
eigenvalues of the adjacency matrix of the Boolean $n$-cube (integers
of the same parity as $n$); and the equality in Proposition
\ref{cor:Lloydcub} holds.

The listed conditions, for
$1<e<\frac{n-1}{2}$, are satisfied only for $n=23$, with
\[
S_C=
\begin{pmatrix}
0 & 23 & 0 & 0\\
1 & 0 & 22 & 0\\
0 & 2 & 0 & 21\\
0 & 0 & 3 & 20
\end{pmatrix},
\]
whose eigenvalues are
$\lambda=23,\ 7,\ -1,\ -9$.
Golay \cite{Golay} discovered a linear code with these parameters.
The Golay code is unique up to isometry.

Zinoviev, Leontiev, and Tiet\"av\"ainen proved
Theorem \ref{th:percod}, and also proved the nonexistence of perfect codes with
nontrivial parameters (except for the ternary Golay code
$n=11$, $e=2$, $q=3$ \cite{Golay}) when $q=p^s$, where $p$ is prime.
It is now known (\cite{ZL73},\cite{BZLF}) that for arbitrary $q$ and
$2<e<\frac{n-1}{2}$, there are no $e$-perfect codes in $Q^n_q$.

\subsection{Problems}

\begin{exercise}

Let $A$ be a generating set of the Cayley graph
$G=Cay(K,A)$ of an abelian group $K$, and let
$B=A\cup\{0\}$. Prove that if $(C,B)$ is a tiling of the group $K$,
then the coloring $f:V(G)\rightarrow B$ in which each vertex
$x\in K$ is colored by $b$ whenever $x=c+b$ for some
$c\in C$ and $b\in B$, is a perfect coloring.
\end{exercise}

\begin{exercise}

Prove that if a distance-regular graph contains an $e$-perfect code, then for every eigenvalue $\lambda$ of the
quotient matrix $S_C$, it holds
$
|A_{e+1}|P_e(\lambda)=|A_e|P_{e+1}(\lambda)$,
where $A_e$ is the sphere of radius $e$.

\end{exercise}

\begin{exercise}

Prove the equality
\[
\det(S_C-xI)
=(-1)^{e+1}
(c_{e+1}P_{e+1}(x)-b_eP_e(x))
\prod\limits_{t=1}^e c_t.
\]

\end{exercise}

\section{Difference Sets and Other Combinatorial Configurations}

\subsection{Difference Sets}

Let $K$ be a finite abelian group.
%\begin{definition}
    A set $D\subseteq K$ is called a {\sl difference set}
     with parameters $(v,k,\lambda)$,
    if $|K|=v$, $|D|=k$ and for every nonzero $a\in K$
    there are exactly $\lambda$  ordered pairs $d_1, d_2\in D$ such that $d_1-d_2=a$.
%\end{definition}

\begin{claim}\label{c:diffset}
A set $D\subseteq K$ is a difference set
 with parameters $(v,k,\lambda)$ if and only if
  ${\mathbf 1}_{_D} *{\mathbf 1}_{-D}=\lambda\cdot\mathbf{1}
+(k-\lambda)\delta/v^{1/2}$, where $\delta=v^{1/2}{\bf 1}_{\{0\}}$.
\end{claim}
\begin{proof}
Compute the convolution $\sum\limits_{x\in
K}{\mathbf 1}_{D}(x)\cdot{\mathbf 1}_{D}(-(y-x))$. For $y\ne 0$, the summand is equal to $1$ precisely when
$x\in D$  and $x-y\in D$. Thus, the number of such $x$ is exactly the number of pairs
$x,x-y\in D$ satisfying
$x-(x-y)=y$.
 By the
definition of a difference set, this number is $\lambda$.
The value of the convolution as $y=0$
equals $\sum\limits_{x\in K}{\mathbf 1}_{D}(x)\cdot{\mathbf 1}_{D}(x)=k$.
The converse follows immediately from the same calculation.
\end{proof}

\begin{corollary}[\cite{Mesnager}]\label{cor:bent_diffset} Let $b$ be a Boolean bent function such that $\widehat{(-1)^b}(\bar0)=1$.
Consider the set $B=\{x\in\{0,1\}^n\mid b(x)=1 \}$. Then
$B$ is a difference set with parameters
$(2^n,2^{n-1}-2^{n/2-1},2^{n-2}-2^{n/2-1})$. Conversely, if
$B$ is a difference set with these parameters, then ${\mathbf 1}_{B}$
is a bent function.
\end{corollary}
\begin{proof}
We have $|K|=2^n$ by assumption and $|B|=\frac12(2^n-2^{n/2})$ by
Proposition \ref{cl:bentsupp}. Similarly to the proof of Proposition
\ref{cl:bentsupp}, we have the equalities
$b*b=\frac14(\mathbf{1}-(-1)^b)*(\mathbf{1}-(-1)^b)=\frac14(2^n\mathbf{1}+2^{n/2}\delta-
2\cdot2^{n/2}\mathbf{1})=(2^{n-2}-2^{n/2-1})\mathbf{1}+2^{n/2-2}\delta$.
Here we used the linearity of convolution and the fundamental property of a
bent function (Proposition \ref{kriteriibent}).  Proposition
\ref{c:diffset} 
shows that $B$ is a difference set with the stated parameters.
Relying on Proposition
\ref{kriteriibent}, it is easy to prove the converse.
\end{proof}

\begin{remark}
If $b$ is a Boolean bent function and $\widehat{(-1)^b}(\bar0)=-1$, then,
 similarly,  its support is a difference set
with parameters $(2^n,2^{n-1}+2^{n/2-1},2^{n-2}+2^{n/2-1})$. These  parameters, together with those in Corollary \ref{cor:bent_diffset} are called {\sl McFarland parameters}.
\end{remark}
It will be shown below that any difference set with McFarland parameters can be transformed into a Hadamard matrix.

We next show that any difference set gives rise to an orthogonal matrix whose entries take only two values.
 Let there be a difference set $D$
with parameters $(v,k,\lambda)$. For arbitrary $d_1,d_2\in K$, consider the characteristic functions
$f_1(x)={\mathbf 1}_{D}(x-d_1)$ and $f_2(x)={\mathbf 1}_{D}(x-d_2)$. By the definition of a difference set, if $d_1\neq d_2$,
the functions $f_1$ and $f_2$ are simultaneously equal to $1$ in exactly
$\lambda$ points. By symmetry, they differ at $2(k-\lambda)$
points, and they are simultaneously equal to $0$ in the remaining
$v-\lambda-2(k-\lambda)$ points. Let $a$ and $b$ satisfy
the equation
\[
\lambda a^2+2(k-\lambda)ab+(v-\lambda-2(k-\lambda))b^2=0.
\]
The preceding counting shows that the functions $g_{d_1}=af_1+b({\bf
1}-f_1)$ and $g_{d_2}=af_2+b({\bf 1}-f_2)$  are orthogonal if and only if the above equation holds. Thus,
for any distinct $d_1,d_2\in K$, the functions $g_{d_1}$ and
$g_{d_2}$ are orthogonal.  Let us show that $a=1$ and $b=-1$ for difference
sets with McFarland parameters. Indeed,
\[
\lambda -2(k-\lambda)+(v-\lambda-2(k-\lambda))=v-4(k-\lambda)=
2^n-4\cdot2^{n-2}=0.
\]

Write the value vectors of the functions $g_d$, $d\in K$, as rows of a matrix
$T$. The resulting $v\times v$ matrix $T$ is orthogonal up to a scalar factor, since all rows have equal norms. 
For difference sets with McFarland parameters, the entries of $T$ are $\pm1$. Hence $T$ is a Hadamard matrix. This also shows directly that
 $(-1)^{{\mathbf 1}_{D}}*(-1)^{{\mathbf 1}_{D}}(x)=0$ for $x\neq
\bar 0$. Therefore, by Proposition \ref{cl:diff} (or equivalently by
Proposition \ref{cl:Hadam}) we obtain the following result.

\begin{corollary}[\cite{Mesnager}]\label{cor:bent_diffset1} Let $B$ be a difference set with parameters
$(2^n,2^{n-1}\pm2^{n/2-1},2^{n-2}\pm2^{n/2-1})$. Then $b={\mathbf 1}_{B}$
is a Boolean bent function.
\end{corollary}

From the orthogonality of $T$, it follows that the pairs $(a,a)$,  $(b,a)$, $(a,b)$, $(b,b)$ occur in any pair of rows
(or columns) of $T$ the same number of times.
 By
construction $T(x,y)=a{\mathbf 1}_{D}(x-y)+b(\mathbf{1} -{\mathbf 1}_{D})(x-y)$.
Replace  the entries $a$ of $T$ by $1$ and the entries $b$ by $0$. Thus, 
consider the matrix $Q_{_D}(x,y)={\mathbf 1}_{D}(x-y)$. Then each of
the pairs $(1, 1)$, $(0, 1)$, $(1, 0)$, $(0, 0)$ occurs equally often  in any pair of rows or
columns.
   If we regard the rows of $Q_{_D}$ as the characteristic
   vectors of $k$-subsets of a $v$-element set, we obtain
   a $2$-design with parameters $2$-$(v,k,\lambda)$.
 Recall that 
 $2$-designs with the same number of points and blocks  are called symmetric.

%\begin{definition}
  A {\sl partial difference set}
     with parameters $(v,k,\lambda,\mu)$ is a set $D\subseteq K$
    such that $|K|=v$, $|D|=k$ and for any nonzero $a\in D$
    there are exactly $\lambda$ ordered pairs of elements $d_1, d_2\in D$ satisfying
    $d_1-d_2=a$, whereas for every nonzero $a\in K\setminus D$
    there are exactly $\mu$ such pairs.
%\end{definition}

Similarly to Proposition \ref{c:diffset} it is easy to prove that the
definition of a partial difference set implies  the equality

\begin{equation}\label{eq:partdiff} {\mathbf 1}_{D}
*{\mathbf 1}_{-D}=\lambda\cdot{\mathbf 1}_{D}+\mu\cdot(\mathbf{1}-{\mathbf 1}_{D})+
\frac{k-\mu}{v^{1/2}}\delta.\end{equation}
Indeed, at a nonzero element of $D$, the left-hand side has value $\lambda$; at a nonzero element outside $D$, it has value $\mu$; and at $0$, it has value $k$.

In other words, if the convolution  of the characteristic function
of a set with  the characteristic function of its negative is a linear combination of the characteristic function
of the set, the constant function, and the Dirac $\delta$ function, then the set
is a partial difference set. Its parameters $\mu,\lambda,k$
can be determined from the coefficients of this linear combination. 

Throughout 
we additionally assume that $D=-D$ and $0\notin D$.
Applying the Fourier transform to
(\ref{eq:partdiff}), we obtain the following quadratic equation for $f=\widehat{{\mathbf 1}_{D}}$
$$
f \cdot f=\lambda\cdot f+\mu\cdot(\delta-f)+
\frac{k-\mu}{v^{1/2}}\mathbf{1}.$$

Consequently,  $f$ can take only the values
$$\frac{\lambda-\mu\pm\Delta}{2},\quad \mbox{\rm where}\quad
\Delta=\sqrt{(\lambda-\mu)^2+4(k-\mu)/v^{1/2}}$$ on the nonzero
elements of $K$, and $f(0)=k/v^{1/2}$.

Let $D^+=\{x\in K \ | \ f(x)=\frac{\lambda-\mu +\Delta}{2}\}$.

\begin{claim}[Delsarte \cite{Ma}]
$D^+$ is a partial difference set.
\end{claim}
\begin{proof}
It suffices to verify that   the characteristic function of $D^+$  satisfies an identity of the form  (\ref{eq:partdiff})  for
some $k'$, $\lambda'$, and $\mu'$. Since $D=-D$,  we have $f(x)=f(-x)$. Hence $D^+=-D^+$. 
By the inversion formula for the Fourier transform,
 $\widehat{\widehat{f(x)}}=f(-x)$ and therefore
$\widehat{f}={\mathbf 1}_{D}$. 
Consequently,
$\widehat{f}\cdot\widehat{f}=\widehat{f}$. Applying the Fourier
transform once more and using the convolution property   (Proposition \ref{c:convolution}) gives $f*
f=v^{1/2}f$. By the definition of $D^+$ the equality
$f=\alpha_1{\mathbf 1}_{{D^+}} + \alpha_2\mathbf{1}+ \alpha_3\delta$ holds for
some coefficients $\alpha_1, \alpha_2, \alpha_3$.

Substituting
this representation of $f$ into $f* f=v^{1/2}f$ and taking into account
${\mathbf 1}_{D^+}*\delta=v^{1/2}{\mathbf 1}_{D^+}$ and
${\mathbf 1}_{D^+}*\mathbf{1}=|D^+|\mathbf{1}$, we obtain that the convolution
${\mathbf 1}_{D^+} *{\mathbf 1}_{D^+}$ can be expressed as the required linear
combination. Thus, there exist some coefficients
$\mu',\lambda',k'$ for which
the equality
$${\mathbf 1}_{D^+}
*{\mathbf 1}_{D^+}=\lambda'\cdot{\mathbf 1}_{D^+}+\mu'\cdot(\mathbf{1}-{\mathbf 1}_{D^+})+
\frac{k'-\mu'}{v^{1/2}}\delta$$
holds. By the characterization of partial difference sets, $D^+$ is therefore a partial difference set.
\end{proof}

\subsection{Difference Sets and Strongly Regular Graphs as Perfect Colorings}\label{17.5}

Consider the Cayley graph $\mathrm{Cay}(K,D)$ of a group $K$ with connection set $D$, where $D=-D$ and $0\notin D$. It follows from the definition of a partial difference set with parameters $(v,k,\lambda,\mu)$ that $\mathrm{Cay}(K,D)$ is a strongly regular graph with parameters $(v,k,\lambda,\mu)$.

Indeed, let $x,y\in K$ be adjacent vertices of $\mathrm{Cay}(K,D)$; thus, $y=x+a$ for some $a\in D$. By the definition of a partial difference set, there are exactly $\lambda$ ordered pairs $(d_1,d_2)\in D\times D$ such that
$a=d_1-d_2$.
For each such pair,
$x+d_1=y+d_2$
is a common neighbor of $x$ and $y$ in $\mathrm{Cay}(K,D)$. Hence, $x$ and $y$ have exactly $\lambda$ common neighbors. Similarly, if $x$ and $y$ are nonadjacent, then they have exactly $\mu$ common neighbors. Therefore, $\mathrm{Cay}(K,D)$ is strongly regular with parameters $(v,k,\lambda,\mu)$.

The converse also holds: if $\mathrm{Cay}(K,D)$ is strongly regular, then $D$ is a partial difference set in $K$ with the corresponding parameters.

As shown in the previous section, any difference set $D$ with parameters $(v,k,\lambda)$ gives rise to a symmetric $2$-$(v,k,\lambda)$ design. If $D=-D$ and $0\notin D$, then the matrix
$Q_{D}(x,y)={\mathbf 1}_D(x-y)
$
is the adjacency matrix of $\mathrm{Cay}(K,D)$. Moreover, it follows immediately from the definitions that the rows of the adjacency matrix of a strongly regular graph with parameters $(v,k,\lambda,\lambda)$ form the incidence vectors of the blocks of a symmetric $2$-$(v,k,\lambda)$ design. Indeed, any two distinct rows have scalar product $\lambda$, since two distinct vertices have exactly $\lambda$ common neighbors.

The converse is also true. Suppose that the rows of a matrix are the characteristic vectors of the $k$-subsets (blocks) of a symmetric $2$-$(v,k,\lambda)$ design, and suppose that the matrix is symmetric and has zeros on its main diagonal. Then this matrix is the adjacency matrix of a strongly regular graph with parameters $(v,k,\lambda,\lambda)$.

We next show that every strongly regular graph gives rise to a perfect coloring of a certain hypergraph.

Consider the complete graph $K_n$. Let $\Gamma_n$ be the hypergraph whose vertices are the edges of $K_n$ and whose hyperedges are the triples of edges that form triangles in $K_n$. Let $G$ be a strongly regular graph on the same vertex set as $K_n$, with parameters $(n,k,\lambda,\mu)$. Define a $2$-coloring of the vertices of $\Gamma_n$ by
$$
f(e)=
\begin{cases}
1,& e\in E(G),\\
0,& e\notin E(G).
\end{cases}
$$
The color of a hyperedge is the unordered composition of the colors of its three vertices.

We verify that this coloring is perfect. Let $e=\{x,y\}$ be a vertex of $\Gamma_n$. Every hyperedge of $\Gamma_n$ containing $e$ corresponds to a vertex $z\in V(K_n)\setminus \{x,y\}$ and consists of the three edges
$\{x,y\}, \{x,z\}, \{y,z\}$.

First suppose that $e=\{x,y\}$ has color $1$, so that $x$ and $y$ are adjacent in $G$. Among the $n-2$ possible vertices $z$, the numbers for which $z$ is adjacent to both $x$ and $y$, to exactly one of $x,y$, or to neither of them are, respectively,
$\lambda, 2(k-\lambda-1), n-2k+\lambda$.
Hence, the numbers of hyperedges of colors $111$, $110$, and $100$ containing $e$ are, respectively,
$\lambda, 2(k-\lambda-1), n-2k+\lambda$.

These numbers depend only on the parameters of $G$ and not on the particular color-$1$ vertex $e$.

Now suppose that $e=\{x,y\}$ has color $0$, so that $x$ and $y$ are nonadjacent in $G$. In this case, the corresponding numbers are
$\mu, 2(k-\mu), n-2k+\mu-2$.

Consequently, the numbers of hyperedges of colors $011$, $010$, and $000$ containing $e$ are, respectively,
$\mu, 2(k-\mu), n-2k+\mu-2$.
Again, these numbers depend only on the parameters of $G$ and not on the particular color-$0$ vertex $e$.

Thus, for each vertex color and each hyperedge color, every vertex of that color is incident with the same number of hyperedges of that color. Hence, the coloring is perfect.
This construction shows that every strongly regular graph determines a perfect $2$-coloring of $\Gamma_n$.

The converse requires some care. A perfect $2$-coloring of $\Gamma_n$ imposes regularity conditions on the graph formed by the edges of color $1$. Apart from certain degenerate cases (disconnected graphs or 
star graphs), these conditions imply that the resulting graph is strongly regular. In particular, the numbers of common neighbors of two vertices are determined by whether the corresponding edge of $\Gamma_n$ has color $1$ or color $0$.
The correspondence between perfect $2$-colorings of $\Gamma_n$ and strongly regular graphs is not one-to-one. In particular, the definition of a perfect coloring of a hypergraph depends only on the unordered composition of the colors of the vertices in a hyperedge and therefore does not distinguish between different orderings of these vertices.

\subsection{Association Schemes}

It is well known (see \cite{CRC}) that the above definition of a distance-regular graph is equivalent to the following one, which is also frequently taken as the definition. A connected graph $G$ is distance-regular if and only if, for every pair of vertices $x,y$ with $d(x,y)=i$ and every pair of integers $j,k$, the number of vertices
$\{z\in V(G)\mid d(x,z)=j,\ d(y,z)=k\}$
depends only on $i,j,k$ and not on the particular choice of the vertices $x$ and $y$.

The latter definition is closely related to the notion of an association scheme.

An {\sl association scheme} $\mathcal{A}$ on a set $X$ is a set of binary relations
$R_0,R_1,\dots,R_m\subseteq X\times X$
satisfying the following conditions:\\
$
1)\quad R_0=\{(x,x)\mid x\in X\}.$\\
$2)\quad \text{If }R_i\in\mathcal{A},\text{ then }
R_i^*=\{(x,y)\mid (y,x)\in R_i\}\in\mathcal{A}$.\\
$3)\quad \text{For every }(x,y)\in R_k,\text{ the number of vertices }z\in X$

$\text{such that }(x,z)\in R_i\text{ and }(z,y)\in R_j
\text{ depends only on }i,j,k.$

A distance-regular graph gives rise to an association scheme in the following way. Let $X=V(G)$ be the vertex set of a distance-regular graph $G$, and let
$R_i=\{(x,y)\in X\times X\mid d(x,y)=i\}$.
The relations $R_0,R_1,\dots,R_m$ form an association scheme. A detailed exposition of the theory of association schemes can be found in \cite{BIto}.

Consider an association scheme
$\mathcal{A}=\{R_0,R_1,\dots,R_m\}$
on a set $X$ of $n$ elements, where $R_0$ is the relation consisting of the pairs $(x,x)$. Suppose that the association scheme is symmetric, i.e.,
$(x,y)\in R_i\quad\Longrightarrow\quad (y,x)\in R_i $
for every $i\in{1,\dots,m}$. Then the association scheme determines a perfect coloring of the hypergraph $\Gamma_n$ with colors
$\{R_1,\dots,R_m\}$.

Indeed, the vertices of $\Gamma_n$ are the $2$-element subsets of $X$, that is, the edges of the complete graph $K_n$, and its hyperedges are the triples of pairs that form a triangle. We color a vertex $\{x,y\}$ of $\Gamma_n$ by the relation $R_i$ containing $(x,y)$. This is well defined because the scheme is symmetric.

Consider a vertex ${x,y}$ of color $R_k$. A hyperedge containing this vertex is determined by a third element $z\in X\setminus \{x,y\}$. The colors of the other two vertices of this hyperedge are determined by the relations containing $(x,z)$ and $(z,y)$. Therefore, the number of hyperedges of each color composition containing ${x,y}$ is determined by the numbers
$$
p_{ij}^{k}=
\left|\{z\in X\mid (x,z)\in R_i,\ (z,y)\in R_j\}\right|.
$$
By condition (3) of the definition of an association scheme, these numbers depend only on $i,j,k$ and not on the particular choice of $(x,y)\in R_k$. Hence, the resulting coloring of $\Gamma_n$ is perfect.

Despite the close relationship between the definitions of a symmetric association scheme and a perfect coloring of $\Gamma_n$, the two notions do not coincide.

\begin{example}
Consider the hypergraph $\Gamma_8$ and the following coloring of its vertices (the edges of $K_8$):

%\vskip50mm

\begin{minipage}{0.24\textwidth}
\center{\includegraphics[width=2.5\textwidth]{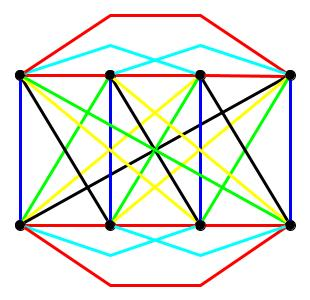}}
\end{minipage}

It is easy to see that the coloring shown in the figure is perfect for the hypergraph $\Gamma_8$, in which two vertices are considered adjacent if the corresponding edges of $K_8$ belong to the same triangle. At the same time, this perfect coloring does not define an association scheme. Indeed, for the edges of $K_8$ of the red color, the two endpoints are not equivalent. Namely, from one endpoint of a red edge there is a blue edge that completes a triangle with a black edge, whereas from the other endpoint there is a green edge. In a symmetric association scheme, the two endpoints of an edge belonging to the same relation must be equivalent.
\end{example}

\subsection{Bent Functions as Perfect Colorings}

Next, we consider in more detail the abelian group $\mathbb{Z}_2^n$. Denote by $\Delta_n$ the $3$-uniform hypergraph whose vertices are all nonzero elements of the group $\mathbb{Z}_2^n$. A triple of vertices ${a_1,a_2,a_3}$ forms a hyperedge if
$a_1+a_2+a_3=0$.
Note that this equality implies that the elements $a_i$, $i=1,2,3$, are nonzero and pairwise distinct.

Let $D$ be a partial difference set in $\mathbb{Z}_2^n$ with parameters $(v,k,\lambda,\mu)$. Then the coloring ${\mathbf 1}_D$ of the hypergraph $\Delta_n$ is perfect. Indeed, in $\mathbb{Z}_2^n$ the equality
$a_1+a_2+a_3=0$
is equivalent to
$a_1=a_2-a_3$.

For each $a_1\in D$, there are $\lambda$ ordered pairs $(d_1,d_2)\in D\times D$ such that
$a_1=d_1-d_2$.
Since the pairs $(d_1,d_2)$ and $(d_2,d_1)$ determine the same hyperedge, the vertex $a_1$ is contained in $\lambda/2$ hyperedges consisting entirely of elements of $D$. Similarly, if $a_1\notin D$, then $a_1$ is contained in $\mu/2$ hyperedges containing two elements of $D$.

Every pair of distinct nonzero elements $a,b\in\mathbb{Z}_2^n$ belongs to exactly one hyperedge of $\Delta_n$, namely
$\{a,b,a+b\}$.
Therefore, the numbers of hyperedges of each possible color composition containing a vertex $a$ can be determined from $\lambda$ and $\mu$. These numbers depend only on the color of $a$. Consequently, the coloring ${\mathbf 1}_D$ is perfect.

The converse is also true. Let $D$ be the set of vertices of color $1$. For a vertex $a\in D$, the number of hyperedges of color composition $(1,1,1)$ containing $a$ is $\lambda/2$, where $\lambda$ is the corresponding partial-difference-set parameter. Similarly, for a vertex $a\notin D$, the number of hyperedges of color composition $(0,1,1)$ containing $a$ is $\mu/2$. Hence, the corresponding partial-difference-set parameters are determined by the parameters of the perfect coloring.

Then, from Proposition \ref{cor:bent_diffset}, we obtain the following claim.

\begin{claim}\label{avgdop}
Boolean bent functions correspond one-to-one to perfect $2$-colorings of the hypergraph $\Delta_n$ with the corresponding quotient matrix.
\end{claim}

By Proposition \ref{cl:edge_coloring}, a perfect coloring of the hypergraph $\Delta_n$ induces a perfect coloring of its line graph $\mathcal{E}(\Delta_n)$. Since the intersection of two hyperedges of $\Delta_n$ contains at most one vertex, $\mathcal{E}(\Delta_n)$ is a simple graph. Note that a $2$-coloring ${\mathbf 1}_D$ of the hypergraph $\Delta_n$ induces a $4$-coloring of the graph $\mathcal{E}(\Delta_n)$, because each hyperedge of $\Delta_n$ consists of three vertices and can contain $0$, $1$, $2$, or $3$ vertices of color $1$.

Consider the elements of the group $\mathbb{Z}_2^n$ as elements of the vector space $(GF(2))^n$. A hyperedge of $\Delta_n$ corresponds to a $2$-dimensional subspace of $(GF(2))^n$. Moreover, two hyperedges of $\Delta_n$ have a common vertex if and only if the corresponding $2$-dimensional subspaces intersect in a $1$-dimensional subspace. Thus, the graph $\mathcal{E}(\Delta_n)$ is the Grassmann graph $J_2(n,2)$.

The {\sl Grassmann graph} $J_q(n,k)$ is the graph whose vertices are the $k$-dimensional subspaces of the vector space $(GF(q))^n$. Two vertices are adjacent if the corresponding subspaces intersect in a $(k-1)$-dimensional subspace. Grassmann graphs are distance-regular (see \cite{BCN}).

Consider the graph $J_2(n,2)$ whose vertices are the $2$-dimensional subspaces of $(GF(2))^n$. Let
$f:Q_2^n\to Q_2$
be a Boolean function, and assume that $f(\bar 0)=1$. Each $2$-dimensional subspace of $(GF(2))^n$ consists of four elements. The function $f$ induces a coloring of the vertices of $J_2(n,2)$ into four colors according to the number of ones of $f$ on the corresponding subspace. Thus, every Boolean function gives rise to a $4$-coloring of the graph $J_2(n,2)$.

\begin{theorem}[\cite{PA20}]\label{avg}
A bent function $b:Q_2^n\to Q_2$ satisfying the normalization conditions $\widehat{(-1)^b}(\bar 0)=1$, %[\emph{to be checked against the Walsh-transform convention}]
$b(\bar 0)=1$
corresponds to a perfect $4$-coloring of the graph $J_2(n,2)$ with quotient matrix
{\tiny
\[
\begin{pmatrix}
3(2^{n-3}-2^{(n/2)-2}-1)& 3\cdot2^{n-2}-3 & 3(2^{n-3}+2^{(n/2)-2}) & 0\\
2^{n-2}-2^{(n/2)-1} & 5\cdot2^{n-3}-2^{(n/2)-2}-5& 2^{n-1}+2^{(n/2)-1} & 2^{n-3}+2^{(n/2)-2}-1\\
2^{n-3}-2^{(n/2)-2} & 2^{n-1}-2^{(n/2)-1}-1 & 5\cdot2^{n-3}+2^{(n/2)-2}-3 & 2^{n-2}+2^{(n/2)-1}-2 \\
0 & 3(2^{n-3}-2^{(n/2)-2}) & 3\cdot2^{n-2} & 3(2^{n-3}+2^{(n/2)-2}-2)\\
\end{pmatrix}.
\]
}
\end{theorem}

\begin{proof}
By Proposition \ref{cl:bentsupp},
$
|\supp(b)|
=\frac12\left(2^n+2^{n/2}\right)$.

%The sign is determined by the value of the Walsh transform at $\bar 0$; this sign should be chosen %consistently with the normalization in the statement of the theorem.

For $y\ne\bar 0$, Proposition \ref{cl:diff}  implies
\[
|\{x\in Q_2^n:b(x)\oplus b(x+y)=0\}|
=
|{x\in Q_2^n:b(x)\oplus b(x+y)=1}|.
\]
Consider a $2$-dimensional subspace containing ${\bar 0,y}$. Denote by $A_{\alpha\beta}$ the number of pairs of values
$(\alpha,\beta)=(b(x),b(x+y))$,
$x\in (GF(2))^n$.

Then
\[
A_{01}+A_{10}=A_{00}+A_{11},
\]
and
\[
A_{01}+A_{10}+A_{00}+A_{11}=2^n.
\]
Moreover, by the change of variable $x\mapsto x+y$,
\[
A_{01}=A_{10}.
\]
Finally,
\[
A_{10}+A_{11}=|\operatorname{supp}(b)|.
\]
Using these equalities and the chosen sign of the Walsh coefficient at $\bar 0$, we obtain
\[
A_{01}=A_{10}=2^{n-2},
\]
\[
A_{11}=2^{n-2}+2^{(n/2)-1},
\qquad
A_{00}=2^{n-2}-2^{(n/2)-1}.
\]

Consequently, if the color of a $2$-dimensional subspace is known, then the numbers of adjacent subspaces of each color are determined.

For example, a subspace of color $4$ is not adjacent to a subspace of color $1$. It is adjacent to
$\frac{3A_{00}}2$
subspaces of color $2$, since each pair of vertices $(x,x+y)$ with
$(b(x),b(x+y))=(0,0)$
is counted twice, namely as $(x,x+y)$ and $(x+y,x)$. It is adjacent to
$\frac{3A_{11}-12}{2}
$
subspaces of color $4$, because it contains four elements, including $\bar 0$, at which the function $b$ takes the value $1$. Fixing each of the three nonzero vertices gives
$\frac{A_{11}-4}{2}$
adjacent subspaces of color $4$.

Thus, the coloring is perfect, with quotient matrix
\[
\frac12
\begin{pmatrix}
3A_{00}-6&6A_{01}-6&3A_{11}&0\\
2A_{00}&4A_{01}+A_{00}-10&2A_{01}+2A_{11}&A_{11}-2\\
A_{00}&2A_{01}+2A_{00}-2&4A_{10}+A_{11}-6&2A_{11}-4\\
0&3A_{00}&6A_{01}&3A_{11}-12
\end{pmatrix}.
\]
Substituting the values of $A_{00}$, $A_{01}$, and $A_{11}$ gives the matrix in the statement.
\end{proof}

The converse is also true: every perfect $4$-coloring with the quotient matrix given in Theorem \ref{avg} corresponds to a Boolean bent function (see \cite{PA20}).

It is easy to see that a difference set with the corresponding parameters defines a perfect coloring of the graph $J_2(n,2)$ with the parameters from Theorem \ref{avg}. Hence, we obtain the following corollary.

\begin{corollary}
The characteristic function of a difference set with parameters
\[
\left(2^n,
2^{n-1}-2^{n/2-1},
2^{n-2}-2^{n/2-1}\right)
\]
is a bent function.
\end{corollary}

It is easy to verify that, by merging the colors pairwise according to their parity in the perfect coloring of $J_2(n,2)$ constructed above, we obtain a perfect coloring with quotient matrix
\[
\begin{pmatrix}
3\cdot2^{n-2}-3&3\cdot2^{n-2}-3\\
3\cdot2^{n-2}&3\cdot2^{n-2}-6
\end{pmatrix}.
\]
However, one cannot recover a bent function from a perfect $2$-coloring with these parameters, because adding an affine function to a bent function does not change the coloring.

\subsection{Spreads}

Let $q=p^r$ be a prime power. Consider $Q_q^n$ as a vector space over the finite field $GF(q)$.
A {\sl spread}, more precisely a $k$-spread, in the vector space $Q_q^n$ is a collection of $k$-dimensional linear subspaces that covers every nonzero vector exactly once.

The definition of a $k$-spread is equivalent to the definition of a $1$-factor (perfect matching) in the hypergraph whose vertices are the elements of $Q_q^n\setminus{\bar 0}$ and whose hyperedges are $k$-dimensional subspaces of $Q_q^n$. Thus, $k$-spreads are perfect colorings of hypergraphs (see Section \ref{7.3}). But sometimes they also generate perfect colorings of Grassmann graphs. 

In the Grassmann graph $J_q(n,k)$, the $k$-dimensional subspaces correspond to vertices. The distance between two $k$-dimensional subspaces $U$ and $V$ in $J_q(n,k)$ is
$
d(U,V)=k-\dim(U\cap V)$.
Therefore, a $k$-spread is equivalent to a code in $J_q(n,k)$ with the maximum possible minimum distance $k$ and the maximum possible number of vertices,
$c=\frac{q^n-1}{q^k-1}$.

Indeed, subspaces at distance $k$ intersect only in the zero vector, and their nonzero vectors are therefore pairwise disjoint. Since each $k$-dimensional subspace contains $q^k-1$ nonzero vectors, the number of such subspaces cannot exceed $(q^n-1)/(q^k-1)$.

\begin{claim}
$2$-spreads in $J_2(n,2)$ correspond one-to-one to the first cells of equitable $2$-partitions with quotient matrix
\[
\begin{pmatrix}
0&6\cdot2^{n-2}-6\\
3&6\cdot2^{n-2}-9
\end{pmatrix}.
\]
\end{claim}

\begin{proof}
Every vertex of $J_2(n,2)$ corresponds to a $2$-dimensional subspace containing three nonzero vectors. Two elements of a $2$-spread are not adjacent in $J_2(n,2)$, since they intersect only in the zero vector.

Now consider a vertex of $J_2(n,2)$ that does not belong to the spread. The corresponding $2$-dimensional subspace contains three nonzero vectors. Each of these vectors belongs to exactly one element of the spread. These three elements of the spread are distinct, and each intersects the given $2$-dimensional subspace in a $1$-dimensional subspace. Hence, the vertex is adjacent to exactly three vertices of the spread.

Thus, the partition into the spread and its complement is equitable, with quotient matrix
\[
\begin{pmatrix}
0&6\cdot2^{n-2}-6\\
3&6\cdot2^{n-2}-9
\end{pmatrix}.
\]

Conversely, suppose that an equitable $2$-partition with this quotient matrix is given, and let $c$ and $c'$ denote the sizes of its first and second cells, respectively. The total number of vertices of $J_2(n,2)$ is
$
\frac{(2^n-1)(2^n-2)}{6}$.
Hence
\[
c+c'=\frac{(2^n-1)(2^n-2)}{6}.
\]
Counting edges between the two cells in two ways gives
\[
c(6\cdot2^{n-2}-6)=3c'.
\]
Solving this system yields
$c=\frac{2^n-1}{3}$.

Each vertex in the first cell is nonadjacent to every other vertex in the first cell. Thus, the corresponding $2$-dimensional subspaces intersect only in the zero vector. Since there are $(2^n-1)/3$ such subspaces and each contains three nonzero vectors, they cover all nonzero vectors of $Q_2^n$. Hence, they form a $2$-spread.
\end{proof}

It is well known that every $2$-spread and every $3$-spread in $Q_q^n$ gives a completely regular code in $J_q(n,2)$ and $J_q(n,3)$, respectively (see \cite{Martin98}). Other perfect colorings of Grassmann graphs are considered in \cite{Mog22}.

Clearly, a $k$-spread in $Q_q^n$ can exist only if
$\frac{q^n-1}{q^k-1}$
is an integer.

\begin{claim}
If $k$ divides $n$, then there exists a $k$-spread in $Q_q^n$.
\end{claim}

\begin{proof}
Suppose that $n=kt$. Regard $Q_q^n$ as a $t$-dimensional vector space over the field $GF(q^k)$. The $1$-dimensional linear subspaces of this vector space form a $1$-spread. Indeed, every nonzero vector belongs to a unique $1$-dimensional subspace over $GF(q^k)$.

Each $1$-dimensional subspace over $GF(q^k)$ is a $k$-dimensional subspace over $GF(q)$. To see this, let
$
A_v=\{\alpha v:\alpha\in GF(q^k)\}
$
be the $1$-dimensional subspace containing a nonzero vector $v$. Since $GF(q^k)$ is a $k$-dimensional vector space over $GF(q)$, the set $A_v$ is a $k$-dimensional vector space over $GF(q)$. Thus, the $1$-dimensional subspaces over $GF(q^k)$ become $k$-dimensional subspaces over $GF(q)$.

Moreover, two distinct $1$-dimensional subspaces over $GF(q^k)$ intersect only in the zero vector. Therefore, these subspaces form a $k$-spread in $Q_q^n$.
\end{proof}

A {\sl partial spread} is a collection of linear subspaces of a fixed dimension that intersect pairwise only in the zero vector.

\begin{claim}[Dillon, see \cite{Mesnager}]\label{Dillon}
Let
$
\{V_i:i=1,\dots,2^{n-1}\}
$
be a partial $n$-spread in $Q_2^{2n}$. Then the characteristic function $b$ of the set
$B=\bigcup_{i=1}^{2^{n-1}}V_i\setminus{\bar 0}$
is a bent function.
\end{claim}

\begin{proof}
Since the subspaces $V_i$ intersect pairwise only in $\bar 0$, the sets $V_i\setminus{\bar 0}$ are pairwise disjoint. Hence,
$b=\sum_{i=1}^{2^{n-1}}{\mathbf 1}_{V_i}-\frac12\delta$. Consequently,
$$(-1)^b= 2b-{\mathbf 1}=
2\sum_{i=1}^{2^{n-1}}{\mathbf 1}_{V_i}-\delta-\mathbf 1.
$$
By using 
$
\widehat{{\mathbf 1}_{V_i}}={\mathbf 1}_{V_i^\perp}$ (Corollary \ref {cor:face2}),
and  
$
\widehat{\delta}=\mathbf 1$, $\widehat{\mathbf 1}=\delta$ ( Section \ref{10.4}),
we obtain
$$\widehat{(-1)^b}=
2\sum_{i=1}^{2^{n-1}}{\mathbf 1}_{V_i^\perp}
-\mathbf 1-\delta.
$$

For $i\ne j$,
$
V_i^\perp\cap V_j^\perp=(V_i+V_j)^\perp.
$
Since
$
V_i\cap V_j=\{\bar 0\}
$
and both subspaces have dimension $n$ in the $2n$-dimensional space $Q_2^{2n}$, we have
$V_i+V_j=Q_2^{2n}$.
Therefore,
$
V_i^\perp\cap V_j^\perp={\bar 0}$.
Thus, the sets $V_i^\perp\setminus\{\bar 0\}$ are pairwise disjoint. Consequently, the sum
$\sum_{i=1}^{2^{n-1}}{\mathbf 1}_{V_i^\perp}$
takes only the values $0$ and $1$ away from $\bar 0$, and hence $\widehat{(-1)^b}$ takes only the values $\pm 1$. Therefore, $b$ is a bent function.
\end{proof}

\subsection{Problems}

\begin{exercise}
Prove that every partition of a set of $2n$ elements  into pairs  is a completely regular code in $J(2n,2)$.
Prove that every  partition of a set of $3n$ elements into triples is a completely regular code in $J(3n,3)$.
\end{exercise}

\begin{exercise}
Prove that the complement of a difference set is also a
difference set.
\end{exercise}

\begin{exercise}
Let $M$ be the adjacency matrix of a strongly regular graph with
parameters $(v,k,\lambda, \lambda+2)$. Prove that the rows of the matrix
$M+I$ are the characteristic
vectors of $(k+1)$-blocks of a symmetric $2$-$(v,k+1,\lambda+2)$
design.
\end{exercise}

\begin{exercise}
Let $D$ be a  $2$-$(v,k,\lambda)$ design
of size $v$. Let the rows of the matrix $Q_{_D}$ be the
characteristic vectors of the blocks of the design $D$. Prove that
the columns of $Q_{_D}$ are the characteristic vectors of the blocks of a
design with the same parameters.
\end{exercise}

\begin{exercise}
Let $\{V_i : i=1,\dots,2^{n-1}+1\}$ be a partial $n$-spread in
$Q_2^{2n}$. Prove that the characteristic function of the set
$B=\cup_{i=1}^{2^{n-1}+1}V_i$ is a bent function.
\end{exercise}

\begin{exercise}
Prove that, using the Dillon construction (Proposition
\ref{Dillon}), one can construct at least ${2^n+1 \choose 2^{n-1}}$  distinct bent functions in $Q^{2n}_2$.
\end{exercise}

\begin{exercise}
Prove that, using the Maiorana--McFarland construction
(Proposition \ref{cplat133}), one can construct at least
$(2^n)!2^{2^n}$ distinct bent functions in $Q^{2n}_2$.
\end{exercise}

\begin{exercise}
Let there be a difference set $D$ with parameters $(v,k,1)$.
Prove that the $v\times v$ matrix $Z$ whose rows
are the characteristic vectors of the sets $a-D$, $a\in K$, does
not contain a $2\times 2$ submatrix consisting of four ones.
\end{exercise}

\begin{exercise}
Prove that for any prime $p$ there exists a square matrix
of size $p(p+1)+1$ containing $(p+1)(p(p+1)+1)$ ones and having
no $2\times 2$ submatrix consisting entirely of  ones  (the
{\sl Zarankiewicz problem}).
\end{exercise}

\begin{exercise}
Prove that, for any $p\ge 1$, a square matrix of size
$p(p+1)+1$ containing
more than $(p+1)(p(p+1)+1)$ ones must contain a $2\times 2$
submatrix consisting  entirely of  ones.
\end{exercise}

\begin{exercise}
Prove that for any $p\ge 1$ a bipartite graph with
$p(p+1)+1$ vertices in each part and $R$ edges contains a cycle of length $4$
whenever $R>(p+1)(p(p+1)+1)$. If $R=(p+1)(p(p+1)+1)$ and $p$ is
prime, prove that there exists such a graph without cycles of length less than $6$.
\end{exercise}

\begin{exercise}
Let $f$ be a perfect $2$-coloring of the hypergraph
$\Gamma_n$ (see Section \ref{17.5}).
Prove that if the edges of either color of $f$ induce a
connected regular graph $G$, then $G$ is strongly regular.
\end{exercise}

\pagebreak


\begin{thebibliography}{999}
\addcontentsline{toc}{section}{References}

\bibitem{AlonCh}
N. Alon and F. R. K. Chung,
``Explicit construction of linear sized tolerant networks,''
in \textit{Proceedings of the First Japan Conference on Graph Theory and Applications},
Hakone, 1986, pp. 15--19, 1988.

\bibitem{AlonM}
N. Alon and V. D. Milman,
``$\lambda_1$, isoperimetric inequalities for graphs, and superconcentrators,''
\textit{J. Comb. Theory, Ser. B}, vol. 38, no. 1, pp. 73--88, 1985.

\bibitem{AlonF}
N. Alon and S. Friedland,
``The maximum number of perfect matchings in graphs with a given degree sequence,''
\textit{Electron. J. Comb.}, vol. 15, \#13, 2008.

\bibitem{Angluin}
D. Angluin and A. Gardiner,
``Finite common coverings of pairs of regular graphs,''
\textit{J. Comb. Theory, Ser. B}, vol. 30, pp. 184--187, 1981.

\bibitem{Avgust}
S. V. Avgustinovich,
``On a property of perfect binary codes,''
\textit{Diskretn. Anal. Issled. Oper.}, vol. 2, no. 1, pp. 4--6, 1995 (Russian).

\bibitem{AM10}
S. V. Avgustinovich and I. Yu. Mogilnykh,
``Perfect 2-colorings of Johnson graphs $J(8,3)$ and $J(8,4)$,''
\textit{J. Appl. Industr. Math.}, vol. 5, no. 1, pp. 19--30, 2011.

\bibitem{AM11}
S. V. Avgustinovich and I. Yu. Mogil'nykh,
``Induced perfect colorings,''
\textit{SEMR}, vol. 8, pp. 310--316, 2011.

\bibitem{Bab}
L. Babai,
``Spectra of Cayley graphs,''
\textit{J. Comb. Theory, Ser. B}, vol. 27, no. 2, pp. 180--189, 1979.

\bibitem{Ball}
S. Ball,
``On sets of vectors of a finite vector space in which every subset of basis size is a basis,''
\textit{J. Eur. Math. Soc. (JEMS)}, vol. 14, no. 3, pp. 733--748, 2012.

\bibitem{BIto}
E. Bannai and T. Ito,
\textit{Algebraic Combinatorics I: Association Schemes},
Mathematics Lecture Note Series, Menlo Park, CA: Benjamin/Cummings, 1984.

\bibitem{BZLF}
L. A. Bassalygo, V. A. Zinov'ev, V. K. Leont'ev, and N. I. Fel'dman,
``Nonexistence of perfect codes over some composite alphabets,''
\textit{Probl. Peredaci Inform.}, vol. 11, no. 3, pp. 3--13, 1975 (Russian).

\bibitem{Besp}
E. A. Bespalov, D. S. Krotov, A. A. Matiushev, A. A. Taranenko, and K. V. Vorob'ev,
``Perfect 2-colorings of Hamming graphs,''
\textit{J. Combin. Des.}, vol. 29, no. 6, pp. 367--396, 2021.

\bibitem{Bespalov}
E. Bespalov,
``On the non-existence of extended 1-perfect codes and MDS codes,''
\textit{J. Comb. Theory, Ser. A}, vol. 189, Article ID 105607, 11 pp., 2022.

\bibitem{Bier}
J. Bierbrauer,
``Bounds on orthogonal arrays and resilient functions,''
\textit{J. Combin. Des.}, vol. 3, no. 3, pp. 179--183, 1995.

\bibitem{Bierbook}
J. Bierbrauer,
\textit{Introduction to Coding Theory}, 2nd ed.,
Boca Raton, FL: Taylor \& Francis, 2017.

\bibitem{Bose}
R. C. Bose,
``On the construction of balanced incomplete block designs,''
\textit{Ann. Eugenics}, vol. 9, pp. 353--399, 1939.

\bibitem{BShP}
R. C. Bose, S. S. Shrikhande, and E. T. Parker,
``Further results on the construction of mutually orthogonal Latin squares and the falsity of Euler's conjecture,''
\textit{Canad. J. Math.}, vol. 12, pp. 189--203, 1960.

\bibitem{BCN}
A. E. Brouwer, A. M. Cohen, and A. Neumaier,
\textit{Distance-Regular Graphs},
Ergebnisse der Mathematik und ihrer Grenzgebiete (3), vol. 18,
Berlin: Springer-Verlag, 1989.

\bibitem{BGKM}
A. E. Brouwer, C. D. Godsil, J. H. Koolen, and W. J. Martin,
``Width and dual width of subsets in polynomial association schemes,''
\textit{J. Comb. Theory, Ser. A}, vol. 102, no. 2, pp. 255--271, 2003.

\bibitem{BHaemers}
A. E. Brouwer and W. H. Haemers,
\textit{Spectra of Graphs},
Universitext, Berlin: Springer, 2012.

\bibitem{Bray}
J. N. Bray, Q. Cai, P. J. Cameron, P. Spiga, H. Zhang,
``The Hall-Paige conjecture, and synchronization for affine and
diagonal groups'',
 \textit{J. Algebra}, vol. 545. pp. 27--42, 2020.

\bibitem{Bruck}
R. H. Bruck and H. J. Ryser,
``The nonexistence of certain finite projective planes,''
\textit{Canad. J. Math.}, vol. 1, pp. 88--93, 1949.

\bibitem{CCD}
P. Camion, B. Courteau, and P. Delsarte,
``On $r$-partition designs in Hamming spaces,''
\textit{Appl. Algebra Engrg. Comm. Comput.}, vol. 2, no. 3, pp. 147--162, 1992.

\bibitem{Canteaut}
A. Canteaut,
``Open problems related to algebraic attacks on stream ciphers,''
in \textit{International Workshop on Coding and Cryptography (WCC 2005)},
Bergen, March 2005, Lecture Notes in Computer Science, vol. 3969,
Berlin: Springer, 2006, pp. 1--11.

\bibitem{CCZ}
C. Carlet, P. Charpin, and V. Zinoviev,
``Codes, bent functions and permutations suitable for DES-like cryptosystems,''
\textit{Des. Codes Cryptography}, vol. 15, no. 2, pp. 125--156, 1998.



\bibitem{CarSar}
 C. Carlet and P. Sarkar,  ``Spectral domain analysis of correlation immune and resilient Boolean functions", \textit{Finite fields and Applications}, vol.  8, pp. 120--130, 2002.


\bibitem{Carlet15}
C. Carlet,
``Boolean and vectorial plateaued functions, and APN functions,''
\textit{IEEE Trans. Inf. Theory}, vol. 61, no. 11, pp. 6272--6289, 2015.

\bibitem{CMes}
C. Carlet and S. Mesnager,
``Four decades of research on bent functions,''
\textit{Des. Codes Cryptography}, vol. 78, pp. 5--50, 2016.

\bibitem{Carlet20}
C. Carlet,
\textit{Boolean Functions for Cryptography and Coding Theory},
Cambridge: Cambridge University Press, 2020.

\bibitem{Chowla}
S. Chowla and H. J. Ryser,
``Combinatorial problems,''
\textit{Canad. J. Math.}, vol. 2, pp. 93--99, 1950.

\bibitem{Cohen}
G. Cohen, S. Litsyn, A. Vardy, and G. Zémor,
``Tilings of binary spaces,''
\textit{SIAM J. Discrete Math.}, vol. 9, no. 3, pp. 393--412, 1996.
DOI: 10.1137/S0895480195280137.

\bibitem{CMeier}
N. Courtois and W. Meier,
``Algebraic attacks on stream ciphers with linear feedback,''
in \textit{Proceedings of Eurocrypt 2003},
Lecture Notes in Computer Science, vol. 2656, pp. 345--359, 2003.

\bibitem{25}
M. B. Crawford, D. J. Marchette, W. Maxwell, and S. S. Mendelson,
``Spectral properties of random graphs with fixed equitable partition,''
E-print 2311.07675, arXiv.org, 2023.
DOI: 10.48550/arXiv.2311.07675.

\bibitem{CDSachs}
D. M. Cvetković, M. Doob, and H. Sachs,
\textit{Spectra of Graphs: Theory and Application},
New York: Academic Press, 1980.

\bibitem{Delsarte}
P. Delsarte,
``An algebraic approach to the association schemes of coding theory,''
\textit{Philips Research Reports Supplements}, vol. 10, 1973.

\bibitem{Delsarte73}
P. Delsarte,
``Four fundamental parameters of a code and their combinatorial significance,''
\textit{Information and Control}, vol. 23, pp. 407--438, 1973.

\bibitem{DelLev}
P. Delsarte and V. I. Levenshtein,
``Association schemes and coding theory,''
\textit{IEEE Trans. Inf. Theory}, vol. 44, no. 6, pp. 2477--2504, 1998.

\bibitem{DenesK}
J. Dénes and A. D. Keedwell,
\textit{Latin Squares: New Developments in the Theory and Applications},
Annals of Discrete Mathematics, vol. 46,
Amsterdam: North-Holland, 1991.

\bibitem{Dev}
K. Devriendt and P. Van Mieghem,
``Tighter spectral bounds for the cut size, based on Laplacian eigenvectors,''
\textit{Linear Algebra Appl.}, vol. 572, pp. 68--91, 2019.

\bibitem{Dong}
X. D. Dong, C. B. Son, and E. Gunawan, ``Matrix characterization of MDS linear codes over
modules,'' \textit{Linear Algebra Appl.},  vol. 277, no. 1, pp. 57--61, 1998.


\bibitem{Egor}
G. P. Egorychev
``The solution of van der Waerden’s problem for permanents.,''
\textit{ Adv. Math}. vol. 42, pp. 299--305, 1981.

\bibitem{EMull}
J. T. Ethier and G. L. Mullen,
``Strong forms of orthogonality for sets of hypercubes,''
\textit{Discrete Math.}, vol. 312, no. 12--13, pp. 2050--2061, 2012.

\bibitem{EVardy}
T. Etzion and A. Vardy,
``Perfect binary codes and tilings: problems and solutions,''
\textit{SIAM J. Discrete Math.}, vol. 11, no. 2, pp. 205--223, 1998.

\bibitem{Etzion}
T. Etzion,
\textit{Perfect Codes and Related Structures},
Singapore: World Scientific, 2022.

\bibitem{Falik}
D. I. Falikman, 
``Proof of the van der Waerden conjecture regarding the permanent of a doubly stochastic matrix,''
\textit{Math. Notes}, vol. 29, pp. 475--479, 1981;  translation from \textit{Mat. Zametki}, vol. 29, pp. 931--938, 1981.

\bibitem{FdF1}
D. G. Fon-Der-Flaass,
``Perfect 2-colorings of a hypercube,''
\textit{Sib. Math. J.}, vol. 48, no. 4, pp. 740--745, 2007.
Translated from \textit{Sib. Mat. Zh.}, vol. 48, no. 4, pp. 923--930, 2007.

\bibitem{FdF2}
D. G. Fon-Der-Flaass,
``A bound on correlation immunity,''
\textit{SEMR}, vol. 4, pp. 133--135, 2007.

\bibitem {FG}
 E. Friedgut, G. Kalai,  ``Every monotone graph property has a
sharp threshold'', \textit{Proceedings of the American Mathematical Society},
vol.124, no. 10, pp. 2993--3002, 1996.


\bibitem{Friedman}
J. Friedman,
``On the bit extraction problem,''
in \textit{Proc. 33rd IEEE Symposium on Foundations of Computer Science},
pp. 314--319, 1992.

\bibitem{Glock}
S. Glock, D. Kühn, A. Lo, and D. Osthus,
``The existence of designs via iterative absorption: hypergraph $F$-designs for arbitrary $F$,''
\textit{Mem. Amer. Math. Soc.}, vol. 284, no. 1406, v+131 pp., 2023.

\bibitem{Godsil0}
C. D. Godsil,
``Equitable partitions,''
in \textit{Combinatorics, Paul Erdős is Eighty}, vol. 1,
Bolyai Soc. Math. Stud., Budapest: János Bolyai Math. Soc., pp. 173--192, 1993.






\bibitem{Godsil}
C. D. Godsil,
\textit{Algebraic Combinatorics},
Chapman and Hall Mathematics Series, New York: Chapman \& Hall, 1993.

\bibitem{Godsil95}
C. D. Godsil,
``Tools from linear algebra,''
in \textit{Handbook of Combinatorics}, vol. 2,
R. L. Graham, M. Grötschel, and L. Lovász, Eds.,
Amsterdam: Elsevier; Cambridge, MA: MIT Press, pp. 1705--1748, 1995.

\bibitem{Godsil97}
 C. Godsil, ``Compact graphs and equitable partitions,'' \textit{ Linear Algebra and Its Applications},
vol. 255, pp. 259-266,  1997.



\bibitem{GodsilR}
C. Godsil and G. Royle,
\textit{Algebraic Graph Theory},
Graduate Texts in Mathematics, vol. 207,
New York: Springer-Verlag, 2001.

\bibitem{GM}
C. Godsil and K. Meagher,
\textit{Erdős--Ko--Rado Theorems: Algebraic Approaches},
Cambridge: Cambridge University Press, 2016.

\bibitem{Golay}
M. J. E. Golay,
``Notes on digital coding,''
\textit{Proc. IRE}, vol. 37, no. 6, p. 657, 1949.

\bibitem{Golubev}
K. Golubev,
``Graphical designs and extremal combinatorics,''
\textit{Linear Alg. Appl.}, vol. 604, pp. 490--506, 2020.

\bibitem{Haemers}
W. H. Haemers,
\textit{Eigenvalue Techniques in Design and Graph Theory},
Dissertation, Technische Hogeschool Eindhoven, 1979.
Mathematical Centre Tracts, vol. 121,
Amsterdam: Mathematisch Centrum, 1980.

\bibitem{Haemers1}
W. H. Haemers,
``Interlacing eigenvalues and graphs,''
\textit{Linear Algebra Appl.}, vols. 226--228, pp. 593--616, 1995.

\bibitem{Hatami}
H. Hatami, P. Hatami, and S. Lovett,
``Higher-order Fourier analysis and applications,''
\textit{Foundations and Trends in Theoretical Computer Science},
vol. 13, no. 4, pp. 247--448, 2019.

\bibitem{Hamming50}
R. W. Hamming,
``Error detecting and error correcting codes,''
\textit{Bell System Tech. J.}, vol. 29, pp. 147--160, 1950.

\bibitem{Hamming}
R. W. Hamming,
\textit{Coding and Information Theory},
Englewood Cliffs, NJ: Prentice-Hall, 1980.

\bibitem{Hanani}
H. Hanani,
``On quadruple systems,''
\textit{Can. J. Math.}, vol. 12, pp. 145--157, 1960.

\bibitem{CDinitz}
C. J. Colbourn and J. H. Dinitz, Eds.,
\textit{Handbook of Combinatorial Designs}, 2nd ed.,
Discrete Mathematics and Its Applications,
Boca Raton, FL: Chapman \& Hall/Taylor \& Francis, 2007.

\bibitem{Horadam}
K. J. Horadam,
\textit{Hadamard matrices and their applications}, 
Princeton, NJ: Princeton University Press, 263 p., 2007.


\bibitem{Hoff}
A. J. Hoffman,
``On eigenvalues and colorings of graphs,''
in \textit{Graph Theory and Its Applications}
(Proc. Advanced Sem., Math. Research Center, Univ. of Wisconsin, Madison, WI, 1969),
pp. 79--91, 1970.

\bibitem{Huang}
H. Huang, B. Xia, and S. Zhou,
``Perfect codes in Cayley graphs,''
\textit{SIAM J. Discrete Math.}, vol. 32, no. 1, pp. 548--559, 2018.

\bibitem{Isaev}
M. Isaev, T. Makai,  B. McKay,  P. Pralat, J. Tan and M. Zhukovskii
``Canonical labelling of random regular graphs,''
arXiv preprint arXiv:2602.17567,  2026.
     

\bibitem{Ji}
L. Ji and J. Yin,
``Constructions of new orthogonal arrays and covering arrays of strength three,''
\textit{J. Combin. Theory Ser. A}, vol. 117, no. 3, pp. 236--247, 2010.


\bibitem{Keevash18}
P. Keevash,
``The existence of designs II,''
arXiv preprint arXiv:1802.05900, 2018.

\bibitem{Keevash14}
P. Keevash,
``A short proof of the existence of designs,''
arXiv preprint arXiv:2411.18291 , 2024.


\bibitem{Kelly}
S. Kelly,
``Constructions of intriguing sets of polar spaces from field reduction and derivation,''
\textit{Des. Codes Cryptogr.}, vol. 43, no. 1, pp. 1--8, 2007.

\bibitem{Halyav}
A. V. Khalyavin,
``Estimates of the capacity of orthogonal arrays of large strength,''
\textit{Mosc. Univ. Math. Bull.}, vol. 65, no. 3, pp. 130--131, 2010 ; translation from Vest. Mosk. Univ. Mat. Mekh. vol. 65, no. 3, pp. 49-51, 2010.

\bibitem{Horosh}
D. B. Khoroshilova,
``On the parameters of perfect 2-colorings of circulant graphs,''
\textit{Diskretn. Anal. Issled. Oper.}, vol. 18, no. 6, pp. 82--89, 2011 (Russian).

\bibitem{Kir}
T. E. Kireeva,
``Perfect orientation colorings of cubic graphs,''
\textit{Sib. Èlektron. Mat. Izv.}, vol. 15, pp. 1353--1360, 2018.

\bibitem{Koolen}
J. H. Koolen, W. S. Lee, and W. J. Martin,
``Characterizing completely regular codes from an algebraic viewpoint,''
arXiv preprint arXiv:0911.1828 [math.CO].

\bibitem{Krotov11}
D. S. Krotov,
``On weight distributions of perfect colorings and completely regular codes,''
\textit{Des. Codes Cryptogr.}, vol. 61, no. 3, pp. 315--329, 2011.

\bibitem{Krotov12}
D. S. Krotov,
``On the binary codes with parameters of triply-shortened 1-perfect codes,''
\textit{Des. Codes Cryptography}, vol. 64, no. 3, pp. 275--283, 2012.

\bibitem{Krotov14}
D. S. Krotov,
``On calculation of the interweight distribution of an equitable partition,''
\textit{J. Algebraic Combin.}, vol. 40, no. 2, pp. 373--386, 2014.

\bibitem{KMP}
D. S. Krotov, I. Yu. Mogilnykh, and V. N. Potapov,
``To the theory of $q$-ary Steiner and other-type trades,''
\textit{Discrete Mathematics}, vol. 339, no. 3, pp. 1150--1157, 2016.

\bibitem{CRCKP}
D. S. Krotov and V. N. Potapov,
``Completely regular codes and equitable partitions,''
in \textit{Completely Regular Codes in Distance-Regular Graphs},
M. Shi and P. Solé, Eds.,
Boca Raton, FL: Chapman \& Hall/CRC, 2025, pp. 1--84.



\bibitem{CRC}
M. Shi and P. Solé, Eds.,
\textit{Completely Regular Codes in Distance-Regular Graphs},
Boca Raton, FL: Chapman \& Hall/CRC, 2025.

\bibitem{LamKol}
C. W. H. Lam, G. Kolesova, and L. Thiel,
``A computer search for finite projective planes of order 9,''
\textit{Discrete Math.}, vol. 92, no. 1--3, pp. 187--195, 1991.

\bibitem{Lam}
C. W. H. Lam, L. Thiel, and S. Swiercz,
``The nonexistence of finite projective planes of order 10,''
\textit{Canad. J. Math.}, vol. 41, no. 6, pp. 1117--1123, 1989.

\bibitem{Lei}
F. T. Leighton,
``Finite common coverings of graphs,''
\textit{J. Combin. Theory Ser. B}, vol. 33, no. 3, pp. 231--238, 1982.

\bibitem{Lang}
S. Lang,
\textit{Algebra},
Graduate Texts in Mathematics, vol. 211, rev. 3rd ed.,
New York: Springer-Verlag, 2002.

\bibitem{Lev}
V. I. Levenshtein,
``Universal bounds for codes and designs,''
in \textit{Handbook of Coding Theory},
V. S. Pless and W. C. Huffman, Eds.,
Amsterdam: Elsevier, 1998, ch. 6, pp. 499--648.

\bibitem{LL14}
N. Linial and Z. Luria, 
``On the vertices of the {{\(d\)}}-dimensional {Birkhoff} polytope.''
\textit{Discrete Comput. Geom.}, vol. 51, no. 1, pp. 161--170, 2014.


\bibitem{Lint75}
J. H. van Lint,
``A survey of perfect codes,''
\textit{Rocky Mt. J. Math.}, vol. 5, pp. 199--224, 1975.

\bibitem{CameronLint}
J. H. van Lint and R. M. Wilson,
\textit{A Course in Combinatorics}, 2nd ed.,
Cambridge: Cambridge University Press, 2001.

\bibitem{Lint}
J. H. van Lint,
\textit{Introduction to Coding Theory}, 3rd ed.,
Graduate Texts in Mathematics, vol. 86,
Berlin: Springer-Verlag, 1999.

\bibitem{Lis}
M. A. Lisitsyna and S. V. Avgustinovich,
``Perfect colorings of the prism graph,''
\textit{Sib. Èlektron. Mat. Izv.}, vol. 13, pp. 1116--1128, 2016.

\bibitem{Lloyd}
S. P. Lloyd,
``Binary block coding,''
\textit{Bell Syst. Tech. J.}, vol. 36, no. 2, pp. 517--535, 1957.

\bibitem{Loban}
M. S. Lobanov,
``Exact relation between nonlinearity and algebraic immunity,''
\textit{Discrete Math. Appl.}, vol. 16, no. 5, pp. 453--460, 2006.
Translated from \textit{Diskretn. Mat.}, vol. 18, no. 3, pp. 152--159, 2006.

\bibitem{LSSYa}
O. A. Logachev, A. A. Salnikov, and V. V. Yashchenko,
\textit{Boolean Functions in Coding Theory and Cryptography},
Translations of Mathematical Monographs, vol. 241,
Providence, RI: American Mathematical Society, 2012.
Translated from Russian by S. Nikova.

\bibitem{Lovasz}
L. Lovász,
\textit{Combinatorial Problems and Exercises},
corrected reprint of the 1993 2nd ed.,
Providence, RI: AMS Chelsea Publishing, 2007.

\bibitem{Ma}
S. L. Ma,
``A survey of partial difference sets,''
\textit{Des. Codes Cryptography}, vol. 4, pp. 221--261, 1994.

\bibitem{MacW}
F. J. MacWilliams and N. J. A. Sloane,
\textit{The Theory of Error-Correcting Codes}, Parts I and II,
North-Holland Mathematical Library, vol. 16,
Amsterdam: North-Holland, 1977.




\bibitem{Markov}
A. A. Markov, Jr.,
\textit{Selected Works}, Vol. II:
\textit{Theory of Algorithms and Constructive Mathematics. Mathematical Logic. Information Science and Related Questions},
2003 (Russian).

\bibitem{Martin}
W. J. Martin,
\textit{Completely Regular Subsets},
Ph.D. thesis, University of Waterloo, 1992.

\bibitem{Martin98}
W. J. Martin, ``Completely regular designs,'' \textit{J. Combin. Des.}, vol. 4, pp. 261–273, 1998.


\bibitem{Mesnager}
S. Mesnager,
\textit{Bent Functions: Fundamentals and Results},
Cham: Springer, 2016.

\bibitem{McKay:80}
B. D. McKay,
``Practical graph isomorphism,''
in \textit{Proceedings of the 10th Manitoba Conference},
Winnipeg/Manitoba, 1980,
Congressus Numerantium, vol. 30, pp. 45--87,
Winnipeg, MB: Combinatorial Press, 1981.

\bibitem{MKW}
B. D. McKay and I. M. Wanless,
``A census of small Latin hypercubes,''
\textit{SIAM J. Discrete Math.}, vol. 22, pp. 719--736, 2008.

\bibitem{Miraf}
S. M. Mirafzal, M. Ziaee,
``A note on the automorphism group of the Hamming graph,'' 
\textit{Trans. Comb.}, vol. 10,  no. 2, pp. 129--136, 2021.





\bibitem{Mog20}
I. Yu. Mogilnykh,
``Perfect codes from PGL(2,5) in star graphs,''
\textit{Sib. Èlektron. Mat. Izv.}, vol. 17, pp. 534--539, 2020.

\bibitem{Mog22}
I. Yu. Mogilnykh, ``Completely regular codes in Johnson and Grassmann graphs with small covering radii,``
\textit{Electron. J. Comb.}, vol. 29, no. 2, Paper No. P2.57, 15 p. 2022.

\bibitem{Mont}
R. Montgomery,
``A proof of the Ryser--Brualdi--Stein conjecture for large even $n$,''
arXiv:2310.19779.

\bibitem{Morgan:65}
H. L. Morgan,
``The generation of a unique machine description for chemical structures: A technique developed at Chemical Abstracts Service,''
\textit{J. Chem. Doc.}, vol. 5, no. 2, pp. 107--113, 1965.

\bibitem{Morris}
S. A. Morris,
\textit{Pontryagin Duality and the Structure of Locally Compact Abelian Groups},
London Mathematical Society Lecture Note Series, vol. 29,
Cambridge: Cambridge University Press, 1977.

\bibitem{NiSz}
N. Nisan and M. Szegedy,
``On the degree of Boolean functions as real polynomials,''
\textit{Comput. Complexity}, vol. 4, no. 4, pp. 301--313, 1994.
Special issue on circuit complexity (Barbados, 1992).

\bibitem{OPPh}
P. R. J. Östergård, O. Pottonen, and K. T. Phelps,
``The perfect binary one-error-correcting codes of length 15:
Part II -- Properties,''
\textit{IEEE Trans. Inf. Theory}, vol. 56, no. 6, pp. 2571--2582, 2010.

\bibitem{ParLis}
O. G. Parshina and M. A. Lisitsyna,
``The perfect 2-colorings of infinite circulant graphs with a continuous set of odd distances,''
\textit{Sib. Electron. Math. Reports}, vol. 17, pp. 590--603, 2020.

\bibitem{PPV}
A. L. Perezhogin, V. N. Potapov, S. Yu. Vladimirov, ``Every latin hypercube of order 5 has transversals,'' \textit{Journal of Combinatorial designs},  vol. 32, no. 11, 679–699,  2024.


\bibitem{Preparata}
F. P. Preparata,
``A class of optimum nonlinear double-error-correcting codes,''
\textit{Inform. and Control}, vol. 13, pp. 378--400, 1968.

\bibitem{Pot12}
V. N. Potapov,
``On perfect 2-colorings of the $q$-ary $n$-cube,''
\textit{Discrete Math.}, vol. 312, no. 6, pp. 1269--1272, 2012.

\bibitem{PA20}
V. N. Potapov and S. V. Avgustinovich,
``Combinatorial designs, difference sets, and bent functions as perfect colorings of graphs and multigraphs,''
E-print 2403.02904, arXiv.org, 2024.

\bibitem{Pot22}
V. N. Potapov,
``Embedding in MDS codes and Latin cubes,''
\textit{J. Comb. Des.}, vol.  30, no. 9, pp. 626-633, 2022.

\bibitem{PotTar}
V. N. Potapov and A. A. Taranenko, 
``Asymptotic bounds on the numbers of vertices of polytopes of polystochastic matrices,''
\textit{Discrete Math.},  349, No. 1, Article ID 114653, 7 p. 2026.

\bibitem{Puz}
S. A. Puzynina,
``Periodicity of perfect colourings of an infinite rectangular grid,''
\textit{Diskretn. Anal. Issled. Oper., Ser. 1}, vol. 11, no. 1, pp. 79--92, 2004 (Russian).

\bibitem{Rao}
R. C. Rao,
``Factorial experiments derivable from combinatorial arrangements of arrays,''
\textit{Suppl. J. Roy. Statist. Soc.}, vol. 9, pp. 128--139, 1947.

\bibitem{RS}
H. Reiter and J. D. Stegeman,
\textit{Classical Harmonic Analysis and Locally Compact Groups}, 2nd ed.,
London Mathematical Society Monographs, New Series, vol. 22,
Oxford: Clarendon Press, 2000.

\bibitem{Rosa}
A. Rosa,
\textit{Combinatorial Designs with Applications},
Notes, Banská Bystrica: Belianum, 2015.

\bibitem{Roth}
R. M. Roth,
``Higher-order MDS codes,''
\textit{IEEE Trans. Inf. Theory}, vol. 68, no. 12, pp. 7798--7816, 2022.

\bibitem{Rot}
O. S. Rothaus,
``On 'bent' functions,''
\textit{J. Combin. Theory, Ser. A}, vol. 20, no. 3, pp. 300--305, 1976.

\bibitem{Ryabov20}
V. G. Ryabov,
``Approximation of restrictions of $q$-valued logic functions to linear manifolds by affine analogues,''
\textit{Diskr. Mat.}, vol. 32, no. 4, pp. 89--102, 2020.
English translation: \textit{Discrete Math. Appl.}, vol. 31, no. 6, pp. 409--419, 2021.

\bibitem{Ryabov21}
V. G. Ryabov,
``Criteria for maximal nonlinearity of a function over a finite field,''
\textit{Diskr. Mat.}, vol. 33, no. 3, pp. 79--91, 2021.
English translation: \textit{Discrete Math. Appl.}, vol. 33, no. 2, pp. 117--126, 2023.

\bibitem{Sar}
 P. Sarkar,
``Spectral domain analysis of correlation immune and resilient Boolean functions,''
\textit{Cryptology ePrint Archive}, Report 2000/049, September 2000.



\bibitem{SarMai}
P. Sarkar and S. Maitra, ``Nonlinearity Bounds and Constructions of Resilient
Boolean Functions,''
\textit{Proceedings of CRYPTO 2000, Lecture Notes in Computer
Science,} vol. 1880, pp. 515--532, 2000.

\bibitem{Shapiro}
H. S. Shapiro and D. L. Slotnick,
``On the mathematical theory of error-correcting codes,''
\textit{IBM J. Res. Develop.}, vol. 3, pp. 25--34, 1959.

\bibitem{Skolem}
Th. Skolem,
``Some remarks on the triple systems of Steiner,''
\textit{Math. Scand.}, vol. 6, pp. 273--280, 1958.

\bibitem{Sol}
F. I. Solov'eva,
\textit{On Perfect Codes and Related Topics},
Lecture Notes, Pohang University of Science and Technology, 2004.

\bibitem{Sot}
E. Sotnikova and A. Valyuzhenich,
``Minimum supports of eigenfunctions of graphs: a survey,''
\textit{The Art of Discrete and Applied Mathematics}.
DOI: 10.26493/2590-9770.1404.61e.

\bibitem{Tanner}
R. M. Tanner,
``Explicit concentrators from generalized $n$-gons,''
\textit{SIAM J. Algebraic Discrete Methods}, vol. 5, no. 3, pp. 287--293, 1984.

\bibitem{Tao}
T. Tao,
``An uncertainty principle for cyclic groups of prime order,''
\textit{Math. Res. Lett.}, vol. 12, no. 1, pp. 121--127, 2005.

\bibitem{Tao12}
T. Tao,
\textit{Higher Order Fourier Analysis},
vol. 142, Providence, RI: American Mathematical Society, 2012.

\bibitem{Taranenko18}
 A. A. Taranenko, 
 ``Transversals in completely reducible multiary quasigroups
and in multiary quasigroups of order 4,'' 
\textit{Discrete Math.},  vol. 341,  pp. 405--420,  2018. 


\bibitem{Taranenko}
A. A. Taranenko,
``Algebraic properties of perfect structures,''
\textit{Linear Algebra Appl.}, vol. 607, pp. 286--306, 2020.

\bibitem{Taranenko22}
A. A. Taranenko,
``Perfect colourings of hypergraphs,''
\textit{Linear Multilinear Algebra}, vol. 73, no. 8, pp. 1566--1590, 2025.

\bibitem{Taran00}
Y. Tarannikov,
``On resilient Boolean functions with maximal possible nonlinearity,''
in \textit{Progress in Cryptology -- INDOCRYPT 2000} (Calcutta),
Lecture Notes in Computer Science, vol. 1977,
Berlin: Springer, pp. 19--30, 2000.

\bibitem{Taran02}
 Yu. V. Tarannikov,
``On correlation-immune and stable Boolean functions,''
\textit{Mat. Vopr. Kibern. } vol.11,  pp. 91-148, 2002  (Russian). 

\bibitem{Tiet}
A. Tietäväinen,
``On the nonexistence of perfect codes over finite fields,''
\textit{SIAM J. Appl. Math.}, vol. 24, pp. 88--96, 1973.

\bibitem{Tits}
R. Titsworth. \textit{Correlation properties of cyclic sequences}, PhD thesis, California Institute of Technology,
1963.


\bibitem{Tok12}
N. Tokareva,
\textit{Bent Functions: Results and Applications to Cryptography},
Amsterdam: Elsevier/Academic Press, 2015.

\bibitem{VNTs}
M. Tsfasman, S. Vlăduţ, and D. Nogin,
\textit{Algebraic Geometric Codes: Basic Notions},
Mathematical Surveys and Monographs, vol. 139,
Providence, RI: American Mathematical Society, 2007.

\bibitem{ValVor}
A.A. Valyuzhenich, K.V. Vorob’ev, ``Minimum supports of functions on the Hamming graphs with
spectral constraints,'' \textit{Discrete Math.}, vol. 342, pp. 1351--1360, 2019. 



\bibitem{Val26}
A. Valyuzhenich,
``An upper bound on the number of relevant variables for Boolean functions on the Hamming graph,' 
\textit{Discrete Math. } vol. 349, no. 2, Article ID 114745, 6 p.,  2026.

\bibitem{Vas09}
A. Yu. Vasil'eva,
``Local and interweight spectra of completely regular codes and of perfect colorings,''
\textit{Probl. Inf. Transm.}, vol. 45, no. 2, pp. 151--157, 2009.
Translated from \textit{Probl. Peredachi Inf.}, vol. 45, no. 2, pp. 84--90, 2009.

\bibitem{Vas12}
A. Yu. Vasil'eva,
``On reconstructive sets of vertices in the Boolean cube,''
\textit{J. Appl. Ind. Math.}, vol. 6, no. 3, pp. 393--402, 2012.

\bibitem{Vas15}
A. Yu. Vasil'eva,
``Reconstruction of eigenfunctions of a $q$-ary $n$-dimensional hypercube,''
\textit{Probl. Inf. Transm.}, vol. 51, no. 3, pp. 231--239, 2015.
Translated from \textit{Probl. Peredachi Inf.}, vol. 51, no. 3, pp. 31--40, 2015.

\bibitem{Vas62}
Yu. L. Vasil'ev,
``On closely-packed nongroup codes,''
\textit{Probl. Kibern.}, vol. 8, pp. 337--339, 1962.
German translation: \textit{Probl. Kybernetik}, vol. 8, pp. 375--378, 1965.

\bibitem{Viz64}
V. G. Vizing,
``On an estimate of the chromatic class of a $p$-graph,''
\textit{Diskret. Analiz}, Novosibirsk, vol. 3, pp. 25--30, 1964 (Russian).

\bibitem{Vizing}
V. G. Vizing,
``Distributive coloring of graph vertices,''
\textit{Diskretn. Anal. Issled. Oper.}, vol. 2, no. 4, pp. 3--12, 1995 (Russian).





\bibitem{Wanless}
I. M. Wanless, Transversals in Latin squares: a survey. Surveys in
combinatorics 2011, 403-437, London Math. Soc. Lecture Note Ser.,
392, Cambridge Univ. Press, Cambridge, 2011.



\bibitem{Wegener}
I. Wegener, \textit{The complexity of Boolean functions}. Wiley-Teubner
Series in Computer Science. John Wiley \& Sons, Ltd., Chichester; B.
G. Teubner, Stuttgart, 1987. xii+457.

\bibitem{WL:68}
B.~{Yu}. Weisfeiler and A.~A. Leman.
``The reduction of a graph to canonical form and the algebra which
  appears therein.''
\textit{Nauchno-Technicheskaya Informatsia}, vol. 2,  no 9.  pp. 12--16, 1968. (Russian)
%\newblock In Russian, translated in
%  \url{https://www.iti.zcu.cz/wl2018/pdf/wl_paper_translation.pdf.}



\bibitem{Wel}
J. Wellens, ``Relationships between the number of inputs and other
complexity measures of Boolean functions,'' 
\textit{Discrete Analysis},  19--21 pp.  2022,  arXiv:2005.00566v2


\bibitem{Wilson}
 R. M. Wilson,  ``An existence theory for pairwise balanced designs.
III. Proof of the existence conjectures.,'' 
\textit{J. Combinatorial Theory
Ser. A}, vol. 18.  pp. 71-79. 1975.


\bibitem{Zheng}
Y. Zheng and X.-M. Zhang, ``Improving upper bound on the nonlinearity of high
order correlation immune functions,''
\textit{Proceedings of Selected Areas in Cryptography
2000, Lecture Notes in Computer Science,} vol. 2012, pp. 262--274, 2001.


\bibitem{ZL72}
 V. A. Zinov’ev,  V. K. Leont’ev,
``On perfect codes,''
\textit{Probl. Peredaci Inform}. vol 8, mo. 1, pp. 26-35, 1972.  (Russian)

\bibitem{ZL73}
 V. A. Zinov’ev, V. K.  Leont’ev, 
``On non-existence of perfect codes over Galois fields,''
\textit{Probl. Control Inf. Theory,} vol. 2, no. 2,  pp.16--24.  1973; translation from Probl. Upravl. Teor. Inform. vol. 2, no. 2, pp. 123--132, 1973.





\end{thebibliography}
\end{document}